\documentclass[reqno]{amsart}

\usepackage[hmargin=2.5cm,top=1.5cm,bottom=1.5cm,includehead]{geometry}

\usepackage{showkeys}
\usepackage{amsmath,mathrsfs}
\usepackage{amssymb, amsmath}
\usepackage{graphicx}
\usepackage{enumerate,color,xcolor}
\usepackage{upgreek}
\usepackage{pdfsync}

\usepackage[normalem]{ulem}

\IfFileExists{ulem.sty}{\usepackage[normalem]{ulem}}{\let\uwave\underline}
\usepackage{hyperref}
\usepackage{colonequals}

\newtheorem{thm}{Theorem}[section]
\newtheorem{lem}{Lemma}[section]
\newtheorem{rmk}{Remark}[section]
\newtheorem{cor}{Corollary}[section]
\newtheorem{prop}{Proposition}[section]

\newcommand{\vv}[1]{\boldsymbol{#1}}

\numberwithin{equation}{section}

\newcommand{\beq}{\begin{equation}}
	\newcommand{\eeq}{\end{equation}}
\newcommand{\ben}{\begin{eqnarray}}
	\newcommand{\een}{\end{eqnarray}}
\newcommand{\beno}{\begin{eqnarray*}}
	\newcommand{\eeno}{\end{eqnarray*}}
\let\f=\frac

\newcommand{\bit}{\begin{itemize}}
	\newcommand{\eit}{\end{itemize}}

\newcommand{\N}{{\mathbb N}}
\newcommand{\Z}{{\mathbb Z}}
\newcommand{\R}{{\mathbb R}}
\newcommand{\T}{{\mathbb T}}
\newcommand{\pa}{\partial}
\newcommand{\vf}{\varphi}
\newcommand{\eps}{\varepsilon}
\newcommand{\ov}{\overline}

\newcommand{\ba}{\begin{aligned}}
	\newcommand{\ea}{\end{aligned}}
\def\na{\nabla}

\let\wt=\widetilde
\let\wh=\widehat
\let\udl=\underline
\def\fA{{\mathfrak A}}
\def\fB{{\mathfrak B}}

\def\fM{{\mathfrak M}}

\def\sE{{\mathscr E}}
\def\sD{{\mathscr D}}
\def\sC{{\mathscr C}}
\def\sF{{\mathscr F}}
\def\cL{{\mathcal L}}
\def\cN{{\mathcal N}}
\def\sD{{\mathscr D}}

\def\sR{{\mathscr R}}
\def\sM{{\mathscr M}}
\def\sP{{\mathscr P}}
\def\tf{{\mathsf f}}
\def\tg{{\mathsf g}}
\def\th{{\mathsf h}}
\def\tQ{{\mathsf Q}}
\def\ueps{{\udl{\eps}}}
\def\uh{{\udl{h}}}
\def\uth{{\udl{\mathsf h}}}
\def\rm{\mathrm{r}}
\def\kAo{\mathrm{k}_1}
\def\kAt{\mathrm{k}_2}
\def\virgp{\raise 2pt\hbox{,}}
\def\cdotpv{\raise 2pt\hbox{;}}

\def\<{\langle}
\def\>{\rangle}

\def\gs{\gtrsim}
\def\ls{\lesssim}
\newcommand{\rr}[1]{\left( #1 \right)}
\newcommand{\nr}[1]{\left| #1 \right|}
\newcommand{\nnr}[1]{\left\| #1 \right\|}
\def\lta{\ell^2_\alpha}

\def\vect{\pa_v^t}
\def\ent{e^{\nu t}}
\def\entau{e^{\nu\tau}}
\def\emnt{e^{-\nu t}}
\def\emntau{e^{-\nu\tau}}
\def\tent{\f{1-\emnt}{\nu}}
\def\tentau{\f{1-\emntau}{\nu}}
\begin{document}

	\title[Eventual Collisional Relaxation After Large Echo-Driven Amplification]{Eventual Collisional Relaxation After Large Echo-Driven Amplification in a Weakly Collisional Vlasov-Poisson-Fokker-Planck System}
	
	\author{Yue Luo}
	
	\address{Tsinghua University}
	\email{luo-y21@tsinghua.org.cn}
\makeatletter
\g@addto@macro\@setauthors{%
	\vspace{-0.5em}%
	\begin{center}
		\footnotesize
		Tsinghua University\\[2pt]
		\href{mailto:luo-y21@tsinghua.org.cn}{\texttt{luo-y21@tsinghua.org.cn}}
	\end{center}%
}
\let\@setaddresses\@empty
\makeatother

	\begin{abstract}
		We study the weakly collisional Vlasov-Poisson-Fokker-Planck system near a Maxwellian equilibrium and construct global classical solutions exhibiting large transient growth before eventual collisional relaxation. The main mechanism is a nonlinear plasma-echo cascade which remains effective up to the enhanced-dissipation time scale $t\sim\nu^{-\frac13}$, producing norm inflation of size larger than $\nu^{-1}$ from arbitrarily small $\nu$-dependent perturbations. At the same time, the nonzero spatial modes decay at the enhanced-dissipation scale $e^{-\mathsf{r}\nu^{\frac13}t}$, while the zero mode relaxes at the thermal scale $e^{-\nu t}$. The estimates also reveal a sharp contrast in the smoothing effects before and after the collisional time scale $t\sim \nu^{-1}$. As a consequence, our construction gives explicit upper bounds on the Sobolev and Gevrey stability thresholds: perturbations of the corresponding $\nu$-dependent sizes can still trigger severe echo-driven amplification before collisional dissipation becomes dominant at the critical time scale $t\sim\nu^{-\frac13}$.
	\end{abstract}
	
	\maketitle

	


	\setcounter{tocdepth}{1}
	\tableofcontents

	\setcounter{tocdepth}{1}
\section{Introduction}
In collisionless plasmas, phase mixing transfers spatial oscillations to increasingly fine velocity scales, leading to the decay of the self-consistent electric field through Landau damping, even though the distribution function need not converge strongly to equilibrium. However, nonlinear resonances can interfere with this damping mechanism, potentially leading to the persistence or regeneration of electric field oscillations. The presence of weak collisions introduces dissipation at fine velocity scales, which can suppress nonlinear resonances and enhance the relaxation toward equilibrium. Understanding the competition between nonlinear resonances and collisional dissipation requires quantitative control of both the size and the regularity of the initial perturbation.

We study this competition for the single-species Vlasov-Poisson-Fokker-Planck system with a neutralizing background on $\T_x^d\times\R_v^d$, with each spatial period equal to $2\pi$ and $0<\nu\ll1$:
	\beq\label{VPF}
	\left\{\ba
	&\pa_tF+v\cdot\na_xF+E\cdot\na_vF=\nu\rr{\Delta_vF+\na_v\cdot(vF)}, \quad(x,v)\in\T^d\times\R^d,\\
	&E(t,x)=-\na_x(-\Delta_x)^{-1}\rr{\varrho(t,x)-\f{1}{|\T^d|}\int_{\T^d}\varrho dx}, \ \varrho(t,x)=\int_{\R^d} Fdv.
	\ea\right.
	\eeq

Here $F(t,x,v)$ is the particle distribution, $\varrho(t,x)$ is its spatial density, and $E(t,x)$ is the self-consistent electrostatic field. The linear Fokker-Planck operator describes velocity diffusion and friction against a thermal bath; its strength is measured by the collision frequency $\nu$. In particular, the collision operator in \eqref{VPF} fixes the temperature and mean velocity of the bath, rather than conserving the momentum and kinetic energy of the particles.

The total mass is conserved:
\[
\frac{d}{dt}\iint F(t,x,v)\,dv\,dx=0.
\]
We impose $\iint F(0,x,v)\,dv\,dx=(2\pi)^d$ and consider perturbations of the spatially homogeneous Maxwellian
\[
\sM(v)=\frac{1}{(2\pi)^{d/2}}e^{-|v|^2/2}.
\]
For sufficiently regular positive solutions, the relative entropy together with the electrostatic energy satisfies
	\beq\label{dee}\f{d}{dt}\rr{\iint F\log \f F\sM dvdx+\f12\int E^2dx}+4\nu\iint\nr{\na_v\sqrt{\f F\sM}}^2\sM dvdx=0,\eeq

The functional in \eqref{dee} is minimized by $\sM$ at the prescribed mass. This dissipation identity identifies the equilibrium and constrains the evolution, but it does not control the growth of high derivatives uniformly as $\nu\to0$. Large transient amplification in strong norms can therefore coexist with eventual relaxation.

The relevant time scales already appear in the linear transport-diffusion dynamics. For $\nu t\ll1$, a nonzero spatial Fourier mode $k$ experiences a damping factor of the form $\exp(-c\nu|k|^2t^3)$, reflecting the accumulation of velocity diffusion along mixing characteristics. Thus the lowest nonzero modes dissipate on the enhanced-dissipation scale $t\sim\nu^{-1/3}$, well before the thermal relaxation scale $t\sim\nu^{-1}$. The spatially homogeneous mode does not benefit from this mixing mechanism. The central question of this paper is whether a nonlinear plasma-echo cascade can produce large amplification up to the enhanced-dissipation scale while the same solution remains global and subsequently relaxes to $\sM$.

\underline{Collisionless damping and plasma echoes.}
Landau's linear analysis \cite{Landau} predicts damping of the electric field through phase mixing, and Penrose \cite{Penrose} formulated a spectral stability criterion for homogeneous equilibria. The nonlinear problem is more delicate because information stored in fine velocity oscillations can return to the density through resonant interactions. The plasma-wave echo was described by Gould, O'Neil, and Malmberg \cite{GOM} and observed experimentally by Malmberg, Wharton, Gould, and O'Neil \cite{MWGO}. Iterated echoes provide a mechanism by which nonlinear effects can accumulate even when individual field modes undergo linear damping.

Mouhot and Villani \cite{MV} established nonlinear Landau damping near stable homogeneous equilibria for sufficiently regular small perturbations. Bedrossian, Masmoudi, and Mouhot \cite{BMM} subsequently obtained damping for Fourier weights $e^{\lambda\langle k,\xi\rangle^s}$ with $s>1/3$, using estimates that track the echo interactions. Grenier, Nguyen, and Rodnianski \cite{GNR} gave a further proof in the analytic and Gevrey regimes. The endpoint $s=1/3$, corresponding to Gevrey-3 regularity, was established in suitable weighted spaces by Ionescu, Pausader, Wang, and Widmayer \cite{waveop}, together with the construction of nonlinear scattering operators. These results show how sufficient regularity limits the cumulative effect of resonances.

Below this regularity, echo constructions reveal obstructions to uniform damping estimates. Bedrossian \cite{echoes} constructed nonlinear echo cascades in Sobolev spaces; in the electrostatic case, the construction involves a sequence of small backgrounds and still smaller perturbations. Zillinger \cite{Zillinger} exhibited the echo-chain mechanism in equations linearized around nonhomogeneous states and identified Gevrey-3 as a critical class for that problem. Grenier, Nguyen, and Rodnianski \cite{GNRecho} constructed echo solutions near Penrose-stable data which nevertheless exhibit damping at large times. These works also emphasize that damping of the field and control of strong norms of the distribution are distinct questions. At an amplitude-dependent Gevrey-3 radius, Luo \cite{LuoEcho2026} showed that infinitely many echo cascades can drive divergence in a high Sobolev norm while the density and electric field decay. We refer to \cite{review} for a broader account of phase mixing and nonlinear echoes.

\underline{Weak collisions and stability thresholds.}
For kinetic Fokker-Planck equations, hypocoercivity describes how transport and degenerate velocity diffusion combine to produce convergence to equilibrium; see Villani \cite{Villani}. In the weakly collisional limit, however, one must additionally determine whether echoes can develop before mixing-enhanced diffusion suppresses them. Bedrossian \cite{Bj} proved nonlinear Landau damping and enhanced dissipation for a Vlasov-Poisson equation with nonlinear Fokker-Planck collisions when sufficiently regular Sobolev perturbations have size at most a small multiple of $\nu^{1/3}$. Chaturvedi, Luk, and Nguyen \cite{CLT} obtained analogous weak-collision estimates at the $\nu^{1/3}$ scale for the Vlasov-Poisson-Landau system with Coulomb interactions, where the collision operator requires a different analytical framework.

The relation between regularity and perturbation size was developed further by Luo \cite{LY} for the linear Fokker-Planck operator in \eqref{VPF}, and by Bedrossian, Zhao, and Zi \cite{BZZ} for nonlinear Fokker-Planck collisions. In particular, \cite{LY} proves stability and enhanced dissipation for sufficiently small data independently of $\nu$ when $1/3\leq s\leq1$, and for data of size $O(\nu^{(1-3s)/(3-3s)})$ when $0<s<1/3$, with the velocity weights and other hypotheses specified there. The results of \cite{BZZ} establish a related size-regularity tradeoff for nonlinear collisions and also address the collisionless limit. Bedrossian, Zhao, and Zi \cite{BZZ1} established uniform-in-collision-frequency Landau damping and enhanced dissipation for the Vlasov–Poisson–Landau system in a Gevrey regularity class, while justifying the long-time collisionless limit and controlling nonlinear plasma echoes. More recently, Zhao and Zi \cite{ZZlow} proved stability at size $\nu^{1/2}$ in a weighted space with logarithmic spatial regularity and no velocity derivatives, and constructed larger perturbations with transient growth obstructing uniform enhanced-dissipation bounds in that setting. Their lower-regularity thresholds concern different norms from the high Sobolev and Gevrey norms considered here.

\underline{The present work.}
We construct positive global classical solutions of \eqref{VPF} near the Maxwellian for which echo-driven amplification persists until a time comparable to $\nu^{-1/3}$. The initial perturbations are small in arbitrarily high Sobolev spaces, or in Gevrey classes with Fourier exponent $0<s<1/3$, yet suitable norms measured along transport characteristics become larger than $\nu^{-1}$. For these same solutions, lower-order norms of the nonzero spatial modes decay at rate $e^{-\mathsf r\nu^{1/3}t}$, while the homogeneous mode relaxes at rate $e^{-\nu t}$. Thus the construction quantifies a transient loss of perturbative control in strong norms together with global collisional relaxation.

The Sobolev construction has initial size $O_N(\nu^{1/3}(\log(1/\nu))^{7/3})$, and the Gevrey construction has initial size $O(\nu^{(1-3s)/(3-2s)})$. These give obstructions to extending uniform stability estimates to those sizes in the norms specified below. In particular, the Gevrey construction does not close the gap with the sufficient stability scale in \cite{LY}; the initial norm and the weaker norm exhibiting inflation are specified separately in Theorem~\ref{T1}. The proof must track the growth of a resonant cascade and its collisional attenuation simultaneously, and then control the exact nonlinear solution for all later times.

	\subsection{Main results}
	We set $F(t)=\sM+f(t)$ and $\varrho(t)=1+\rho$. Then $f(t)$ satisfies the following equation:
	\beq\label{1}\left\{\ba
	&\pa_tf+v\cdot\na_xf+E\cdot\na_v\sM+E\cdot\na_vf=\nu Lf,\\
	&E(t,x)=\na_x\Delta_x^{-1}\rr{\rho(t,x)},
	\ea\right.\eeq
	where $Lf=\Delta_vf+\na_v\cdot\rr{vf}$.	We recall that $\iint f(t)dvdx=0$, which implies that $\int \rho(t)dx=0$.
	\begin{thm}\label{T1}
		\begin{enumerate}
			\item For any $N\ge 10^4$, there exists $\nu(N)>0$ such that for every $0<\nu<\nu(N)$, there is a classical solution $f(t,x,v)$ to \eqref{1} on $[0,\infty)$ with $\sM(v)+f(t,x,v)>0$. The initial data satisfies
			\beq\label{initial}\nnr{\<\na\>^{N}f(0)}_{L^2_{x,v}}\ls_N\rr{\log\f{1}{\nu}}^\f73\nu^\f13;\eeq
			there exists a time $\check{T}$ with $c\nu^{-\f13}<\check{T}<\nu^{-\f13}$ such that
			\beq\label{checkT}\nnr{\<\na_x,\na_v-\check{T}\na_x\>^{304}f(\check{T})}_{L^2_{x,v}}>\f1\nu;\eeq
			and for all $t\geq0$,
			\beq\label{edted}
			\nnr{\<\na_x,\emnt\na_v-\tent\na_x\>^{296}e^{\mathsf{r}\min\{1,\nu t\}\rr{|\na_v|+|\na_x|t}} P_{\neq}f(t)}_{L^2_{x,v}}\leq e^{-\mathsf{r}\nu^\f13 t},\eeq
			\beq\label{edt0}\nnr{\<\na_v\>^{296}e^{\mathsf{r}\min\{1,\nu t\}|\na_v|} P_{0}f(t)}_{L^2_{v}}\leq e^{-\nu t},\eeq
			where $c, \mathsf{r}>0$ are universal constants.
			\item For any $0<s<\f13$, there exists $\nu(s)>0$ such that for every $0<\nu<\nu(s)$, there is a classical solution $f(t,x,v)$ to \eqref{1} on $[0,\infty)$ with $\sM(v)+f(t,x,v)>0$. The initial data satisfies
			\beq\label{initial2}\nnr{e^{\<\na\>^{s}}f(0)}_{L^2_{x,v}}\ls\nu^{\f{1-3s}{3-2s}};\eeq
			there exists a time $\check{T}$ with $c\nu^{-\f13}<\check{T}<\nu^{-\f13}$ such that
			\beq\label{checkT2}\nnr{e^{10c^{-1}\<\na_x,\na_v-\check{T}\na_x\>^{\f{2s}{3}}\log\<\na_x,\na_v-\check{T}\na_x\>}f(\check{T})}_{L^2_{x,v}}>\f1\nu;\eeq
			and for all $t\geq0$,
			\beq
			\nnr{\<\na_x,\emnt\na_v-\tent\na_x\>^{296}e^{\mathsf{r}\min\{1,\nu t\}\rr{|\na_v|+|\na_x|t}} P_{\neq}f(t)}_{L^2_{x,v}}\leq e^{-\mathsf{r}\nu^\f13 t},\eeq
			\beq\nnr{\<\na_v\>^{296}e^{\mathsf{r}\min\{1,\nu t\}|\na_v|} P_{0}f(t)}_{L^2_{v}}\leq e^{-\nu t},\eeq
			where $c, \mathsf{r}>0$ are universal constants.
		\end{enumerate}
	\end{thm}
    \begin{rmk}
    	 The vector field derivative $\na_v-t\na_x$ is adopted in \eqref{checkT} and \eqref{checkT2} instead of $\na_v$ in order to eliminate the growth caused by the transport effect. The vector field derivative  $\emnt\na_v-\tent\na_x$ appearing in \eqref{edted} originates from the coordinate shift  \eqref{coordinateshift}. We also observe that $\<\na_x,\emnt\na_v-\tent\na_x\>\sim\<\na_x,\na_v-t\na_x\>$ whenever $0\leq t\leq\nu^{-\f12}$.
    \end{rmk}

	\subsection{Notations and conventions}
\subsubsection{Fourier transform}
For $f(x,v)$ defined on $\T^d\times\R^d$ we define the Fourier transform via
\[\wh{f}(k,\xi)=\f{1}{(2\pi)^d}\int_{\T^d\times\R^d}f(x,v)e^{-ik\cdot x-i\xi\cdot v}dxdv.\]
For $E(x)$ defined on $\T^d$ we define the Fourier transform via
\[\wh{E}(k)=\f{1}{(2\pi)^d}\int_{\T^d}E(x)e^{-ik\cdot x}dx.\]
We also use $\sF$ to denote the Fourier transform. 
For $f(x,v)$ we define the following components:
\[P_{0}f:=\f{1}{(2\pi)^d}\int_{\T^d}f(x,v)dx,\quad P_kf:=\sF_xf(k,v)e^{ik\cdot x},\ k\in\Z^d, \quad P_{\neq}f:=f-P_0f.\]
Here, $P_0f$ denotes the zero spatial Fourier mode of $f$, while $P_{\neq}f$ represents the nonzero spatial Fourier mode of $f$.
\subsubsection{Fourier multipliers}
For $A(t,k,\xi)$ defined on $I\times\Z^d\times\R^d$, $f(t,x,v)$ defined on $I\times\T^d\times\R^d$, $E(t,x)$ defined on $I\times\T^d$, where $I\subset\R$, we define
\beq\label{multiplierE}\ba&A(t,\na)f(t,x,v)=\sF^{-1}\rr{A(t,k,\xi)\wh{f}(t,k,\xi)},\\ &A(t,\na)E(t,x)=\sF^{-1}\rr{A(t,k,k \tent)\wh{E}(t,k)}.\ea\eeq
\subsubsection{Miscellaneous}
\begin{itemize}
	\item We use $\ell^2_\alpha$ to denote $\ell^2(|\alpha|\leq\llcorner\f d2\lrcorner+2, \alpha\in\Z_{\geq0}^d)$.
	\item We use $\th$, $\tg$, $\tQ$ to denote the coordinate shift \eqref{coordinateshift} of $h$, $g$, $Q$, separately.
	\item $\vf(v)$ is a smooth function satisfying $\mathrm{supp} \wh{\vf}\subset[-0.9,0.9]$, $\wh{\vf}(\xi)\geq0$, $\int\wh{\vf}(\xi)d\xi=1$, and 
	\beq\label{vfdecay}\nr{\vf(v)}\leq C\exp(-|v|^{0.99}).\eeq
	\item We set $N'=6$, $N_1=300$. $D_1, D_2, D_3$ are uniform constants satisfying $1\ll D_1\ll D_2\ll D_3$ which will be determined in the sequel.
\end{itemize}
	\subsection{Main idea of the proof}
Our proof follows the strategy of Deng and Masmoudi \cite{DN}, combining linearization around a carefully chosen background solution with a detailed analysis of resonant cascades and nonlinear error estimates.
To capture the mechanism responsible for plasma echoes, we construct a solution consisting of a low-frequency background and a much smaller high-frequency perturbation. The background is chosen to be an exact time-dependent solution with vanishing electric field, while the high-frequency perturbation serves as a seed for the resonant interactions. Linearization around this time-dependent background retains the principal echo interactions. We then follow the resulting resonant cascade, estimate its nonlinear corrections by a finite expansion, and control the remainder through and beyond the thermal time scale. Proposition~\ref{P1} contains the quantitative lower bound and the global upper bounds needed to deduce Theorem~\ref{T1}.

	\subsubsection{Coordinate shift}
	In this subsection we introduce, in a similar way as \cite{Bj,BZZ}, a change of variables that reduces \eqref{1} to a system whose structure is more transparent. For $f(t,x,v)$ defined on $\R\times\T^d\times\R^d$, we set 
	\beq\label{coordinateshift}\tf(t,z,v)=f(t,z+\tent v,\emnt v)e^{-d\nu t}.\eeq
	Then if $f(t,x,v)$ solves \eqref{1} then $\tf(t,z,v)$ solves
	\beq\label{cs1}
	\left\{\ba
	&\pa_t\tf+E(t,z+\tent v)\cdot\vect\rr{\sM(\emnt v)e^{-d\nu t}+\tf}=-\nu{\vect}^2\tf,\\
	&E(t,z)=\na_z\Delta_z^{-1}\rr{\rho(t,z)},\ \rho(t,z)=\int\tf(t,z-\f{\ent-1}{\nu}v,\ent v)e^{d\nu t}dv,
	\ea\right.\eeq
	where 
	\beq\label{vect}\vect=\ent\na_v-\f{\ent-1}{\nu}\na_z.\eeq
	By Fourier transform we get
	\beq\label{csf1}
	\left\{\ba&\pa_t\wh{\tf}(t,k,\xi)+\wh{\rho}(t,k)\f{k\cdot\bar{\eta}}{|k|^2}\wh{\sM}(\bar{\eta})+\sum_\ell\wh{\rho}(t,\ell)\f{\ell\cdot\bar{\eta}}{|\ell|^2}\wh{\tf}(t,k-\ell,\xi-\ell\tent)=-\nu|\bar{\eta}|^2\wh{\tf}(t,k,\xi),\\
	&\wh{\rho}(t,k)=\wh{\tf}(t,k,k\tent),
	\ea\right.\eeq
	where
	\beq\label{bareta}\bar{\eta}(t,k,\xi)=\ent\xi-\f{\ent-1}{\nu}k.\eeq
	Following \cite{Bj,BZZ} we define
	\beq\label{deS}S(t,\tau,k,\xi):=\exp\rr{-\nu\int_\tau^t\nr{\bar{\eta}(\tau',k,\xi)}^2d\tau'}\eeq
	and
	\beq\label{deSt}S(t-\tau,k):=S(t,\tau,k,k\tent)=\exp\rr{-\nu\int_0^{t-\tau}\nr{k\f{1-e^{-\nu \tau'}}{\nu} }^2d\tau'}.\eeq

	Thanks to  Duhamel's principle, we have
	\[\ba
	\wh{\tf}(t,k,\xi)=&\wh{\tf}(0,k,\xi)S(t,0,k,\xi)-\int_0^t\wh{\rho}(\tau,k)\f{k\cdot\bar{\eta}(\tau,k,\xi)}{|k|^2}\wh{\sM}(\bar{\eta}(\tau,k,\xi))S(t,\tau,k,\xi)d\tau\\
	&-\int_0^t\sum_\ell\wh{\rho}(\tau,\ell)\f{\ell\cdot\bar{\eta}(\tau,k,\xi)}{|\ell|^2}\wh{\tf}(\tau,k-\ell,\xi-\ell\tentau)S(t,\tau,k,\xi)d\tau.
	\ea\]
	By setting $\xi:=k\tent$ we get
	\beq\\label{rho}\ba
	\wh{\rho}(t,k)=\wh{\tf}(0,k,k\tent)S(t,k)-\int_0^t\wh{\rho}(\tau,k)\f{1-e^{-\nu(t-\tau)}}{\nu}\wh{\sM}(k\f{1-e^{-\nu(t-\tau)}}{\nu})S(t-\tau,k)d\tau\notag\\
	-\sum_\ell\int_0^t\wh{\rho}(\tau,\ell)\f{\ell\cdot k}{|\ell|^2}\f{1-e^{-\nu(t-\tau)}}{\nu}\wh{\tf}(\tau,k-\ell,k\tent-\ell\tentau)S(t-\tau,k)d\tau.
	\ea\eeq
	We remark that by performing inverse Fourier transform and coordinate shift we can get the fundamental solution of the inhomogeneous Fokker-Planck equation $\pa_tf+v\cdot\na_xf=\nu\na_v\cdot\rr{\na_v f+vf}$
	\[\ba f(t)=&\sum_{k}\int\sF_xf(0,k,\ent(v-u))e^{d\nu t}\f{1}{(2\pi)^\f d2(1-e^{-2\nu t})^\f d2}\exp\rr{-\f{u^2}{2(1-e^{-2\nu t})}}\\
	&\times\exp\left\{i\rr{x-\f{\ent-1}{\nu}v+\f{\ent-1}{\nu(1+\emnt)}u }\cdot k\right\}du\ \exp\left\{\rr{\f{2}{\nu^2}\f{1-\emnt}{1+\emnt}-\f t\nu}k^2\right\}.\ea\]
	We set
	\[\cN(g,h)(t,x,v)=-E_g(t,x)\cdot\na_vh,\]
	where $E_g(t,x)=\na_x\Delta_x^{-1}\rho_g$, $\rho_g(t,x)=\int_{\R^d}gdv$, and 
	\[\cN^\dagger(\tg,\th)(t,z,v)=\cN(g,h)(t,z+\tent v,\emnt v)e^{-d\nu t}=-E_{\tg}(t,z+\tent v)\cdot\vect \th\]
	where $\tg$ ($\th$) is the coordinate shift \eqref{coordinateshift} of $g$ ($h$), $E_{\tg}(t,z)=E_{g}(t,z)$.
	
\subsubsection{Reduction to 1-$d$ case}
In this subsection we reduce the $d$-dimensional problem to the 1-dimensional problem. For $v\in\R^d$, we denote $v=(v_1,v_2,\cdots,v_d)\in\R^d$. In this subsection we set
\[\sM(v)=\f{1}{(2\pi)^\f d2}e^{-\f{|v|^2}{2}}\]
to be $d$-dimensional Maxwellian, while
\[\sM_1(v_1)=\f{1}{(2\pi)^\f 12}e^{-\f{|v_1|^2}{2}}\]
to be 1-dimensional Maxwellian. We note that $\sM(v)=\sM_1^{\otimes^{d}}=\sM_1(v_1)\cdots\sM_1(v_d)$. Moreover if $f_1(t,x_1,v_1)$ solves the 1-dimensional equation \eqref{1}, then
	\[\ba
	&f(t,x,v)=f_1(t,x_1,v_1)\otimes\sM_1^{\otimes^{d-1}}\\
	&=f_1(t,x_1,v_1)\sM_2(v_1)\cdots\sM_d(v_d)
	\ea\]
	solves the $d$-dimensional equation \eqref{1}. Furthermore if $f_1(t,x_1,v_1)$ satisfies the estimates in Theorem 1, then $f(t,x,v)=f_1(t,x_1,v_1)\otimes\sM_1^{\otimes^{d-1}}$ satisfies the estimates in the $d$-dimensional case. Thus we only need to prove Theorem 1 in the 1-$d$ case. In the following discussion we will restrict ourselves to the 1-$d$ case.

\subsubsection{Construction of the low-frequency background and the high-frequency perturbation}
We set $\sP=\f2\pi\f{1}{(1+v^2)^2}$, which satisfies $\int_\R\sP(v)dv=1$. We will construct the solution on time interval $[1,\infty)$, and finally by a time shifting argument we get the solution on $[0,\infty)$.

We first construct a low-frequency background solution $\udl{h}(t)$ on $[1,\infty)$ by setting 
\beq\label{udlh}\ba
\udl{h}(t)=&2\udl{\eps}\int\vf(\ent(v-u))e^{\nu t}\f{1}{(2\pi)^\f12(1-e^{-2\nu t})^\f12}\exp\rr{-\f{u^2}{2(1-e^{-2\nu t})}}\cos\left(x-\f{\ent-1}{\nu}v\right.\\
&\left.+\f{\ent-1}{\nu(1+\emnt)}u \right)du\ \exp\rr{\f{2}{\nu^2}\f{1-\emnt}{1+\emnt}-\f t\nu}+e^{\nu(t-1)L}\rr{\ueps\sP(v)-\ueps\sM(v)}.
\ea\eeq
We claim that $\udl{h}(t)$ is an exact solution of \eqref{1} with identically vanishing electric field. To verify this, we decompose
\[
\udl{h}=P_0\udl{h}+P_{\neq}\udl{h}.
\]
The spatially homogeneous component satisfies
\beq\label{957}P_0\uh(t)=e^{\nu(t-1)L}\rr{\ueps\sP(v)-\ueps\sM(v)},\eeq
and hence solves the homogeneous Fokker-Planck equation. The polynomial tail in this component will play a crucial role in ensuring the positivity of the full distribution function \(\sM+f\). For the nonzero spatial modes, applying the coordinate transformation \eqref{coordinateshift} and taking the Fourier transform yield
\beq\label{963}\wh{P_{\neq}\uth}(t,k,\xi)=\ueps\wh{\vf}(\xi)S(t,0,1,\xi)\vv{1}_{k=1}+\ueps\wh{\vf}(\xi)S(t,0,-1,\xi)\vv{1}_{k=-1}.\eeq
In particular, the nonzero spatial component solves the linear transport–Fokker–Planck equation. Its associated density is given by
\[\wh{\udl{\rho}}(t,k)=\ueps\wh{\vf}(\tent)S(t,1)\vv{1}_{k=1}+\ueps\wh{\vf}(-\tent)S(t,-1)\vv{1}_{k=-1}.\]
By the assumption $\mathrm{supp} \wh{\vf}\subset[-0.9,0.9]$ we have $\udl{\rho}(t)=0$ on $[1,\infty)$. Consequently, the self-consistent electric field generated by $\udl{h}$ vanishes identically. Therefore, $\udl{h}$ is an exact solution of \eqref{1} on $[1,\infty)$.

We next construct a high-frequency perturbation that will serve as the initial seed for the plasma echo cascade. To this end, we introduce the spatial frequency $K$, the velocity frequency $A$, and the time $T$ by
\beq\label{KAT}K=\llcorner\rr{\f{\ueps}{\nu^\f13}}^\f37\lrcorner,\quad A=\f{1.01K^3}{\ueps},\quad T=\f{\log\rr{1-\f{\nu A}{K+\f12}}}{-\nu}.\eeq
On the time interval $[1,T]$, we define the high-frequency perturbation by
\beq\label{h1T}\ba
h(t)=&2\eps\int\vf(\ent(v-u))e^{\nu t}\f{1}{(2\pi)^\f12(1-e^{-2\nu t})^\f12}\exp\rr{-\f{u^2}{2(1-e^{-2\nu t})}}\cos\left\{K\rr{x-\f{\ent-1}{\nu}v\right.\right.\\
	&\left.\left.+\f{\ent-1}{\nu(1+\emnt)}u}+A\ent(v-u) \right\}du\ \exp\{\rr{\f{2}{\nu^2}\f{1-\emnt}{1+\emnt}-\f t\nu}K^2\}.
\ea\eeq
We perform coordinate shift \eqref{coordinateshift} and Fourier transform and get
\beq\label{Vb}\wh{P_{\neq}\th}(t,k,\xi)=\eps\wh{\vf}(\xi-A)S(t,0,K,\xi)\vv{1}_{k=K}+\eps\wh{\vf}(\xi+A)S(t,0,-K,\xi)\vv{1}_{k=-K}.\eeq
Consequently, the associated density satisfies
\[\wh{\rho}(t,k)=\eps\wh{\vf}(K\tent-A)S(t,K)\vv{1}_{k=K}+\eps\wh{\vf}(-K\tent+A)S(t,-K)\vv{1}_{k=-K}.\]
By the definition of $T$ and the assumption $\mathrm{supp} \wh{\vf}\subset[-0.9,0.9]$ we have $\rho(t)=0$ on $[1,T]$. Hence the electric field generated by $h$ vanishes throughout this interval. Since both $\uh$ and $h$ solve the linear transport–Fokker–Planck equation, their sum $\uh+h$ is an exact solution of \eqref{1} on $[1,T]$. In particular, the nonlinear evolution on this interval is explicitly determined, allowing us to use $\uh(T)+h(T)$ as the initial data for the subsequent echo analysis.

In the remainder of the paper, we study the evolution of the mean-zero perturbation \(h(t)\) around the background solution \(\uh(t)\) for \(t\geq T\).
 We define the following linearized operator 
\[\cL h=-v\pa_xh-E_h(t,x)\pa_v\rr{\sM+\uh}+\nu Lh.\]
Then $h(t)$ solves
\beq\label{2}\left\{\ba
&\pa_th=\cL h+\cN(h,h),\\
&\wh{\th}(T,k,\xi)=\eps\wh{\vf}(\xi-A)S(T,0,K,\xi)\vv{1}_{k=K}+\eps\wh{\vf}(\xi+A)S(T,0,-K,\xi)\vv{1}_{k=-K}.
\ea\right.\eeq

The following proposition establishes the quantitative growth of the resonant cascade together with global bounds for the solution of \eqref{2}.

\begin{prop}\label{P1}
	For any $10^7\nu^\f13\rr{\log\f1\nu}^\f73\leq\ueps\ll1$, there exists $A^{-N_1}e^{-3K}<\eps<A^{-N_1}e^{-1.5K}$ ($K,A$ are defined in \eqref{KAT}) such that the mean-zero solution to \eqref{2} satisfies the following conditions:
	\begin{enumerate}
		\item There exists a time $\check{T}$ with $c\nu^{-\f13}<\check{T}<\nu^{-\f13}$ such that 
		\beq\label{Zx}\rr{\sum_k\int_{|\xi-A|\leq20D_2K}\nr{\wh{P_{\neq}\th}(\check{T},k,\xi)}^2d\xi}^\f12\geq K^{-1}\exp\rr{-1.02c^{-1}K^\f23\log K}A^{-N_1}. \eeq
		\item For any $t\geq1$, it holds that
		\beq\label{zx}\nnr{\<\na_x,\emnt\na_v-\tent\na_x\>^{N_1-4}e^{\mathsf{r}\min\{1,\nu t\}\rr{|\na_v|+|\na_x|t}} P_{\neq}h(t)}_{L^2_{x,v}}e^{\mathsf{r}\nu^\f13 t}\ls K^{-1},\eeq
		\beq\label{zX}\nnr{\<\na_v\>^{N_1-4}e^{\mathsf{r}\min\{1,\nu t\}|\na_v|} P_{0}h(t)}_{L^2_{v}}\ent\ls K^{-1}.\eeq
	\end{enumerate}
\end{prop}

The proof of Proposition~\ref{P1} constitutes the main part of our analysis. Our strategy is to construct the nonlinear solution through a finite Taylor expansion. In the following subsections, we outline the main steps of the proof. We first describe the growth of the resonant echo cascade arising from the linearized equation in Subsection 1.3.4. We then explain how to control the additional resonances generated by higher-order nonlinear interactions in Subsection 1.3.5. We also discuss the estimates needed to control the remainder of the expansion and establish global bounds for the solution. Finally, in Subsection 1.3.7, we explain how Theorem~\ref{T1} follows from Proposition~\ref{P1} through a suitable choice of parameters.


\subsubsection{Linear analysis}\label{s134}
We first study the linearized equation $\pa_t h=\cL h$ with the high-frequency initial data in \eqref{2}. The purpose is to show that the interaction between this perturbation and the low-frequency background produces a sequence of plasma echoes, with a quantitative lower bound at a time comparable to $\nu^{-1/3}$.

In the shifted coordinates, the initial perturbation consists of two parts, supported near the Fourier frequencies \((K,A)\) and \((-K,-A)\). By symmetry, it suffices to analyze the part with spatial mode \(k=K\) and velocity frequency \(\xi\) close to \(A\). Its density becomes significant when the curve $\xi=k(1-e^{-\nu t})/\nu$ passes through the region near $A$. The resulting electric field interacts with the spatial mode $-1$ of the background and generates a perturbation in mode $k-1$. This perturbation subsequently produces its own density response, allowing the same process to repeat. Thus the perturbation is transferred from spatial mode $K$ to successively lower modes, while its velocity Fourier frequency remains close to $A$.

To follow these interactions, we divide time into intervals with endpoints
\[
t_k=-\frac{1}{\nu}\log\left(1-\frac{\nu A}{k+\frac12}\right),
\]
as in \eqref{4.1}. The key estimate is a one-step transfer formula. Up to the controlled errors in \eqref{hktke}, the leading contribution has the form
\[
\wh{\th}(t_{k-1},k-1,\xi)
\simeq
-\frac{\ueps A}{k^3}
e^{-\nu A^3/(3k^3)}
\int\wh{\th}(t_k,k,\zeta)\wh{\vf}(\xi-\zeta)\,d\zeta.
\]
The factor $\ueps A/k^3$ is the echo gain, whereas $e^{-\nu A^3/(3k^3)}$ measures the collisional decay of the low-frequency background at the resonance time. The nonnegativity of $\wh{\vf}$ and control of the transfer errors allow a lower bound to propagate through the chain despite the alternating sign. The Maxwellian response is handled using the resolvent of the density Volterra equation. Frequency-localized energy estimates control the nonresonant terms and show that the contribution outside the region near $\xi=A$ remains small.

By iterating the one-step transfer estimate, we obtain a sequence of resonant interactions that transfer the initial high-frequency perturbation from mode \(K\) to successively lower spatial modes. At each step, the resonance produces amplification, while collisions weaken the transfer through an exponential damping factor. This gives the cumulative factor
\[ \prod_{j=k+1}^{K} \frac{\ueps A}{j^3}e^{-\nu A^3/(3j^3)}. \]
We choose the parameters \(K\) and \(A\) so that the resonant amplification remains effective despite collisional dissipation. In particular, the cascade reaches spatial modes of order \((\nu A^3)^{1/3}\) at times comparable to \(\nu^{-1/3}\). This shows that plasma echoes can accumulate significant growth before collisional dissipation suppresses the resonant interactions.

To organize the estimates, we introduce a nondecreasing control function $B(t)$ that measures the cumulative amplification. We choose the initial amplitude $\eps$ so that $B(\infty)=A^{-N_1}$. The same function controls upper estimates for the linear evolution and lower estimates for its contribution near $\xi=A$ at the selected resonance times. In particular, \eqref{cnu13} and \eqref{lowerb} yield a contribution of size at least
\[
K^{-1/2}A^{-N_1}
\exp\left(-\frac{1.02}{c}K^{2/3}\log K\right)
\]
in the region near $\xi=A$ at a time comparable to $\nu^{-1/3}$. Consequently, applying the transport-adapted derivatives or Gevrey weights introduces large factors depending on $A$. With a suitable choice of parameters, these factors make the norms in \eqref{checkT} and \eqref{checkT2} larger than $\nu^{-1}$, even though the unweighted contribution remains small.

\subsubsection{Nonlinear analysis}
We next explain how the linear growth persists for the nonlinear equation \eqref{2}. We construct an approximate solution by a finite expansion in the amplitude $\eps$. The first term solves the linearized equation:
\beq\label{hone}
\left\{\ba
&\pa_t h^{(1)}=\cL h^{(1)},\\
&\wh{\th^{(1)}}(T,k,\xi)
=\eps\wh{\vf}(\xi-A)S(T,0,K,\xi)\vv{1}_{k=K}
+\eps\wh{\vf}(\xi+A)S(T,0,-K,\xi)\vv{1}_{k=-K}.
\ea\right.
\eeq
For $2\leq n\leq A^{N'}$, the higher-order terms are defined recursively by
\beq\label{hn}
\pa_t h^{(n)}=\cL h^{(n)}
+\sum_{\substack{n_1+n_2=n\\n_1,n_2\geq1}}
\cN(h^{(n_1)},h^{(n_2)}),
\qquad h^{(n)}(T)=0.
\eeq
The main difficulty is that the forcing in each of these equations can generate additional echoes. Smallness of the forcing alone is therefore insufficient: we must also control its subsequent resonant amplification.

For this purpose, we distinguish contributions arising from the two initial velocity frequencies $+A$ and $-A$. We write
\[
h^{(n)}=\sum_{\substack{|q|\leq n\\2\mid n-q}}h^{(n,q)}.
\]
The two first-order components satisfy
\beq\label{fonemp}
\left\{\ba
&\pa_t h^{(1,\pm1)}=\cL h^{(1,\pm1)},\\
&\wh{\th^{(1,\pm1)}}(T,k,\xi)
=\eps\wh{\vf}(\xi\mp A)S(T,0,\pm K,\xi)\vv{1}_{k=\pm K},
\ea\right.
\eeq
and, for $n\geq2$,
\beq\label{hnq}
\pa_t h^{(n,q)}=\cL h^{(n,q)}
+\sum_{\substack{n_1+n_2=n\\q_1+q_2=q}}
\cN(h^{(n_1,q_1)},h^{(n_2,q_2)}),
\qquad h^{(n,q)}(T)=0.
\eeq
Here the sum is over positive $n_i$ and admissible $q_i$, namely $|q_i|\leq n_i$ and $2\mid n_i-q_i$ for $i=1,2$. Reality of the initial perturbation gives
\beq\label{qmq}
h^{(n,-q)}=\ov{h^{(n,q)}}.
\eeq
The index \(n\) records the order of the nonlinear interaction, while \(q\) identifies its principal velocity frequency. More precisely, weighted Fourier estimates show that \(h^{(n,q)}\) is concentrated near \(\xi=qA\) in the shifted coordinates. This localization allows us to estimate the additional echoes generated by each component according to its frequency.

We apply the inhomogeneous linear estimates to each frequency-localized component. The analysis distinguishes frequencies that have already passed resonance, frequencies that can generate further echoes, and high frequencies for which diffusion provides stronger control. The key comparison \eqref{pbqhb} bounds the amplification from these additional echoes in terms of \(B(t)^n\). This allows us to retain the small factor \(B(t)^n\) throughout the induction, with polynomial losses in \(A\), and obtain \eqref{hnqu}. Consequently, the sum of the terms with $n\geq2$ is smaller than the leading linear contribution at the observation time.

To pass from the expansion to an exact solution, we set
\[
h^{\#}=\sum_{n=1}^{A^{N'}}h^{(n)},
\qquad \wt h=h-h^{\#}.
\]
The remainder starts from zero at $T$. Its equation contains terms involving $\wt h$ and a forcing formed from interactions whose total order exceeds $A^{N'}$; see \eqref{wtheq}--\eqref{wtQ}. The high truncation order makes this forcing sufficiently small to compensate for the losses in an energy estimate on $[T,\nu^{-1}]$. For $t\geq\nu^{-1}$, we use hypocoercive estimates for the nonzero spatial modes and the homogeneous Fokker-Planck semigroup for the zero mode to close a global bootstrap. These estimates preserve the echo lower bound for the exact solution and establish its subsequent relaxation.

\subsubsection{Smoothing estimate}\label{s136}
The global bounds must also account for the regularization produced by collisions.We use the Fourier multiplier \(\fM\) to measure the regularity gained from velocity diffusion along the shifted characteristics. It is defined by
\beq\label{defM} \fM(t,k,\xi)=\exp\rr{\mathrm{r}\nu\int_0^t\rr{e^{-2\nu(t-s)}+|\bar{\eta}(s,k,\xi)|^2}^\f12ds}.\eeq
We set $\fM_\sigma(t,k,\xi)=\<k,\xi\>^\sigma \fM(t,k,\xi)$. For sufficiently small \(\rm\), the frequency-dependent contribution from differentiating the multiplier in the weighted energy estimate is absorbed by velocity diffusion, leaving a lower-order term of size \(O(\rm\nu)\) times the weighted energy. The product and commutator estimates for \(\fM\) are used separately to control the nonlinear terms. These properties, together with the corresponding semigroup estimates, are collected in Appendix B.

After returning to the original variables, \eqref{Mxj} shows that control of $\fM$ provides an exponential Fourier weight of the form
\[
\exp\left\{\mathsf r\min\{1,\nu t\}
\big(|\na_v|+t|\na_x|\big)\right\},
\]
for a sufficiently small $\mathsf r>0$. Thus the velocity regularity radius supplied by this estimate grows like $\nu t$ before the thermal time $t\sim\nu^{-1}$ and remains uniformly positive afterwards.

Combining these weighted estimates with enhanced dissipation for the nonzero spatial modes and thermal relaxation for the homogeneous mode yields the global bounds in Proposition~\ref{P1}. They control lower-order norms and therefore remain compatible with the growth in the norms appearing in \eqref{checkT} and \eqref{checkT2}.

\subsubsection{Choice of parameters and completion of the construction}\label{s137}
It remains to choose the background amplitude $\ueps$ so that Proposition~\ref{P1} gives the two conclusions of Theorem~\ref{T1}. The associated frequencies $K$ and $A$ are then determined by \eqref{KAT}. The main idea is to make the initial perturbation sufficiently small while ensuring that the resonant cascade produces large growth in stronger norms.

In the Sobolev case, we take \(K\sim N\log A\), which allows the initial perturbation to remain small in \(H^N\). The parameter relations
\[
\ueps\sim\nu^{1/3}K^{7/3},
\qquad
A\sim\nu^{-1/3}K^{2/3}
\]
then yield
\[
\ueps\sim_N\nu^{1/3}(\log(1/\nu))^{7/3}.
\]
The lower bound for the resonant cascade, together with the large velocity frequency \(A\), gives the desired growth of the higher Sobolev norm at a time comparable to \(\nu^{-1/3}\).

For the Gevrey case \(0<s<1/3\), we choose \(K\sim A^s\) so that the initial perturbation is small in the prescribed Gevrey norm. This leads to
\[
\ueps\sim_s\nu^{\frac{1-3s}{3-2s}}.
\]
The lower bound for the resonant cascade then implies the required norm inflation in the weaker Gevrey norm appearing in \eqref{checkT2}.

Finally, we add the explicitly constructed background to the perturbation and translate time by one unit to obtain a solution on $[0,\infty)$. The positive velocity tail ensures positivity of the initial distribution, which is preserved by \eqref{VPF}. Tensoring the one-dimensional solution with Maxwellians in the remaining velocity variables gives the construction in dimension $d$. The global estimates obtained above then provide the relaxation bounds in Theorem~\ref{T1} for these same solutions.

\subsection{Organization of the paper.} The rest of the paper is organized as follows. Section 2 establishes the weighted energy and bilinear estimates used throughout the proof. In particular, it provides quantitative control of Fourier localization, which allows us to estimate contributions outside the resonant frequency region. These estimates are used both to isolate the leading echo interactions in the linear analysis and to control the frequency distribution of the higher-order nonlinear terms. In Section 3, we analyze the homogeneous linearized equation, derive the one-step echo transfer estimates, and quantify the cumulative amplification in the presence of collisions. Section 4 treats the inhomogeneous linearized equation and bounds the additional echoes generated by frequency-localized forcing. In Section 5, we combine these linear estimates with the localization bounds from Section 2 to estimate the terms in the finite nonlinear expansion. We choose the initial perturbation amplitude and show that the higher-order corrections remain smaller than the leading linear contribution at the observation time. Section 6 controls the remainder of the expansion, first up to the thermal time and then globally, thereby justifying the approximation and establishing the bounds needed for eventual relaxation. In Section 7, we assemble these estimates to prove Proposition~\ref{P1} and obtain Theorem~\ref{T1} through suitable parameter choices. Finally, Appendix A contains auxiliary estimates for the resonance analysis, while Appendix B collects properties of the Fourier multipliers and the Fokker-Planck semigroup used in the smoothing and decay estimates.

\section{Energy estimate of the forced linearized equation}
This section establishes weighted estimates for the forced linearized equation
\[\pa_th=\cL h+Q.\]
on $[T,\infty)$. We assume that $h(T)$ and $Q(t)$ are mean-zero. These estimates control both the size and the Fourier localization of the solution. They will be used in Sections~3 and~4 to isolate the resonant interactions, and in Section~5 to bound the nonresonant part of each nonlinear source. Proposition~\ref{prop:bilinear} gives the bilinear estimates needed to control these sources. Proposition~\ref{3.2} combines a Volterra estimate for the density with semigroup bounds for the profile, while Proposition~\ref{3.3} applies this linear estimate inductively to the expansion coefficients. The resulting localization bounds are deliberately coarse in size; their purpose is to make contributions far from the principal frequency exponentially small.

We introduce the multiplier
\beq\label{deffA}\fA(t,k,\xi)=\exp\rr{D_1\rr{1+\ueps\tent\<\ueps|k,\xi|\>^{-\f23}}^{-0.01}\<\ueps|k,\xi|\>^\f13}\eeq
for a sufficiently large constant $D_1\geq10^3$. The argument of Lemma~2.4 in \cite{waveop} gives
\[\fA(t,k+\ell,\xi+\eta)\leq \fA(t,k,\xi)\fA(t,\ell,\eta)\ \text{and}\ \nr{D^\alpha_\xi \fA(t,k,\xi)}\ls\fA(t,k,\xi)\ \text{for}\ |\alpha|\le2.\]

We begin with the bilinear term. The estimates below separate the density of the first factor from the profile of the second factor and distinguish its zero and nonzero spatial modes. This form is suited both to the linearized background terms and to the nonlinear sources in the expansion. The proof distributes the multiplier between the two factors and uses the time weights to pass from pointwise bounds to time-integrated estimates.
\begin{prop}[Weighted bilinear estimates]\label{prop:bilinear}
	For $\sigma\geq0$, we have
	\beq\label{cNk1}\ba&\int_T^t\nnr{\fM_\sigma(\tau,\na)\rr{e^{-|\alpha|\nu\tau}v^\alpha P_{\neq}\cN^\dagger(\tg,\th)}}_{\lta L^2_{z,v}}e^{\rm\nu^\f13\tau}d\tau\ls\nnr{\fM_{\sigma+2}(\tau,\na)\rho_\tg(\tau)e^{\rm\nu^\f13\tau}}_{L^2([T,t],L^2_z)}\\
	&\times\sup_{\tau\in[T,t]}\left\{\nnr{\fM_{\sigma+1}(\tau,\na)\rr{e^{-|\alpha|\nu\tau}v^\alpha P_{\neq}\th }}_{\lta L^2_{z,v}}e^{\f32\nu\tau}+\nnr{v^\alpha\fM_{\sigma}(\tau,0,\emntau\na_v)\pa_v P_{0}h }_{\lta L^2_{v}}\right\}.
	\ea\eeq
	\beq\label{cNk2}\ba&\int_T^t\nnr{\fM_\sigma(\tau,\na)\rr{e^{-|\alpha|\nu\tau}v^\alpha P_{\neq}\cN^\dagger(\tg,\th)}}_{\lta L^2_{z,v}}e^{\rm\nu^\f13\tau+\f12\nu\tau}d\tau\ls\nnr{\fM_{\sigma+2}(\tau,\na)\rho_\tg(\tau)e^{\rm\nu^\f13\tau}}_{L^2([T,t],L^2_z)}\\
	&\times\sup_{\tau\in[T,t]}\left\{\nnr{\fM_{\sigma+1}(\tau,\na)\rr{e^{-|\alpha|\nu\tau}v^\alpha P_{\neq}\th }}_{\lta L^2_{z,v}}e^{2\nu\tau}+\nnr{v^\alpha\fM_{\sigma}(\tau,0,\emntau\na_v)\pa_v P_{0}h }_{\lta L^2_{v}}e^{\f12\nu\tau}\right\}.
	\ea\eeq
	\beq\label{cN0}\ba&\int_T^t\nnr{v^\alpha\<\na_v\>^\sigma\fM(\tau,0,\emntau\na_v)P_{0}\cN(g,h)}_{\lta L^2_{v}}e^{\nu\tau}d\tau\ls\nnr{\fM_{\sigma+2}(\tau,\na)\rho_g(\tau)e^{(\sigma+3)\nu\tau}}_{L^2([T,t],L^2_z)}\\
	&\times\sup_{\tau\in[T,t]}\nnr{\fM_{\sigma+1}(\tau,\na)\rr{e^{-|\alpha|\nu\tau}v^\alpha P_{\neq}\th }}_{\lta L^2_{z,v}}.\ea\eeq	
\end{prop}
\begin{proof}
We estimate the nonzero and zero spatial outputs separately. For the nonzero output, we split the second factor into its nonzero and zero modes, distribute the multiplier using \eqref{M1}, and apply convolution estimates. Integrating in time then gives the two time-weighted versions. The zero output is treated in the original velocity variable.

    We first prove \eqref{cNk1}. We split $P_{\neq}\cN^\dagger(\tg,\th)$ into two parts.
    \[P_{\neq}\cN^\dagger(\tg,\th)=P_{\neq}\cN^\dagger(\tg,P_{\neq}\th)+P_{\neq}\cN^\dagger(\tg,P_0\th).\] 
    For the first term, fix $|\alpha|\leq2$. Then
    \[\ba
    &\nnr{\fM_\sigma(\tau,\na)\rr{e^{-|\alpha|\nu\tau}v^\alpha P_{\neq}\cN^\dagger(\tg,P_{\neq}\th)}}_{L^2_{z,v}}\ls\sum_{\alpha_1+\alpha_2=\alpha}\nnr{\sum_\ell\nr{\sF\rr{\fM_\sigma(\tau,\na)\rho_\tg}(\tau,\ell)}e^{-|\alpha_1|\nu\tau}\right.\\
	&\quad \left.\f{|\pa_\xi^{\alpha_1}\bar{\eta}(\tau,k,\xi)|}{|\ell|}\nr{\sF\rr{\fM_{\sigma}(\tau,\na)\rr{e^{-|\alpha_2|\nu\tau}v^{\alpha_2} P_{\neq}\th }}(\tau,k-\ell,\xi-\ell\tentau)} }_{L^2_{k,\xi}}\\
    \ls&\sum_{\alpha_1+\alpha_2=\alpha}\nnr{\sum_\ell\nr{\sF\rr{\fM_{\sigma+1}(\tau,\na)\rho_\tg}(\tau,\ell)}\f{\entau}{|\ell|}\right.\\
	&\quad \left.\nr{\sF\rr{\fM_{\sigma+1}(\tau,\na)\rr{e^{-|\alpha_2|\nu\tau}v^{\alpha_2} P_{\neq}\th }}(\tau,k-\ell,\xi-\ell\tentau)}}_{L^2_{k,\xi}}\\
    \ls&\nnr{\fM_{\sigma+1}(\tau,\na)\rho_\tg(\tau)}_{L^2_z}\nnr{\fM_{\sigma+1}(\tau,\na)\rr{e^{-|\alpha|\nu\tau}v^\alpha P_{\neq}\th }}_{\lta L^2_{z,v}}e^{\nu\tau}.
    \ea\]
    The second inequality uses $e^{-|\alpha_1|\nu\tau}|\pa_\xi^{\alpha_1}\bar{\eta}(\tau,k,\xi)|\leq\entau\<k-\ell,\xi-\ell\tentau\>\<\tentau\>$. The third inequality uses Young's inequality. For the second term, fix $|\alpha|\leq2$. Then
    \[\ba
    &\nnr{\fM_\sigma(\tau,\na)\rr{e^{-|\alpha|\nu\tau}v^\alpha P_{\neq}\cN^\dagger(\tg,P_{0}\th)}}_{L^2_{z,v}}\ls\nnr{\nr{\sF\rr{\fM_\sigma(\tau,\na)\rho_\tg}(\tau,k)}\f{1}{|k|}\left|\sF\rr{\fM_{\sigma}(\tau,0,e^{-\nu\tau}\na_v)v^{\alpha} \pa_vP_{0}h}\right.\right.\\
    &\left.\left.(\tau,0,\bar{\eta}(\tau,k,\xi))\right|}_{L^2_{k,\xi}}
    \ls\nnr{\fM_{\sigma}(\tau,\na)\rho_\tg(\tau)}_{L^2_z}\nnr{v^{\alpha}\fM_{\sigma}(\tau,0,\emntau\na_v) \pa_vP_{0}h }_{\lta L^2_{v}}e^{-\f12\nu\tau}.
    \ea\]
    Thus
    \[\ba
    &\int_T^t\nnr{\fM_\sigma(\tau,\na)\rr{e^{-|\alpha|\nu\tau}v^\alpha P_{\neq}\cN^\dagger(\tg,\th)}}_{\lta L^2_{z,v}}e^{\rm\nu^\f13\tau}d\tau\ls\int^t_T\nnr{\fM_{\sigma+1}(\tau,\na)\rho_\tg(\tau)}_{L^2_z}e^{\rm\nu^\f13\tau-\f12\nu\tau}d\tau\\ &\times\sup_{\tau\in[T,t]}\rr{\nnr{\fM_{\sigma+1}(\tau,\na)\rr{e^{-|\alpha|\nu\tau}v^\alpha P_{\neq}\th }}_{\lta L^2_{z,v}}e^{\f32\nu\tau}+\nnr{v^\alpha\fM_{\sigma}(\tau,0,\emntau\na_v)\pa_v P_{0}h }_{\lta L^2_{v}}}\\
    \ls&\nnr{\fM_{\sigma+2}(\tau,\na)\rho_\tg(\tau)e^{\rm\nu^\f13\tau}}_{L^2([T,t],L^2_z)}
    \sup_{\tau\in[T,t]}\left(\nnr{\fM_{\sigma+1}(\tau,\na)\rr{e^{-|\alpha|\nu\tau}v^\alpha P_{\neq}\th }}_{\lta L^2_{z,v}}e^{\f32\nu\tau}\right.\\
    &\left.+\nnr{v^\alpha\fM_{\sigma}(\tau,0,\emntau\na_v)\pa_v P_{0}h }_{\lta L^2_{v}}\right).
    \ea\]
    The second inequality uses
    \[\ba
    &\int^t_T\nnr{\fM_{\sigma+1}(\tau,\na)\rho_\tg(\tau)}_{L^2_z}e^{\rm\nu^\f13\tau-\f12\nu\tau}d\tau\ls\nnr{\fM_{\sigma+2}(\tau,\na)\rho_\tg(\tau)e^{\rm\nu^\f13\tau}}_{L^2([T,t],L^2_z)}\\
    &\times\rr{\int_T^t\<\tentau\>^{-2}\emntau d\tau}^\f12
    \ls\nnr{\fM_{\sigma+2}(\tau,\na)\rho_\tg(\tau)e^{\rm\nu^\f13\tau}}_{L^2([T,t],L^2_z)}.
    \ea\]
    The same estimates with the additional time weight give \eqref{cNk2}. We next prove \eqref{cN0}. 
    \[\ba
    &\nnr{v^\alpha\<\na_v\>^\sigma\fM(\tau,0,\emntau\na_v)P_{0}\cN(g,h)}_{\lta L^2_{v}}\ls\nnr{\sum_\ell\nr{\sF\rr{\fM_{\sigma+1}(\tau,\na)\rho_g}(\tau,\ell)}\right.\\
	&\quad \left.\left|\sF\rr{\fM_{\sigma+1}(\tau,\na)\rr{e^{-|\alpha|\nu\tau}v^{\alpha} P_{\neq}\th }}(\tau,-\ell,\emntau\xi-\ell\tentau)\right|e^{(\sigma+1)\nu\tau} }_{\lta L^2_{\xi}}\\
	&\ls\nnr{\fM_{\sigma+1}(\tau,\na)\rho_\tg(\tau)}_{L^2_z}\nnr{\fM_{\sigma+1}(\tau,\na)\rr{e^{-|\alpha|\nu\tau}v^\alpha P_{\neq}\th }}_{\lta L^2_{z,v}}e^{(\sigma+\f32)\nu\tau}.
    \ea\]
    The first inequality uses $\<\xi\>\leq\entau\<\emntau\xi\>\ls\entau\<-\ell,\emntau\xi-\ell\tentau\>\<\tentau\>$. Thus
    \[\ba&
    \int_T^t\nnr{v^\alpha\fM_\sigma(\tau,0,\emntau\na_v)P_{0}\cN(g,h)}_{\lta L^2_{v}}e^{\nu\tau}d\tau\ls\int_T^t\nnr{\fM_{\sigma+1}(\tau,\na)\rho_g(\tau)e^{(\sigma+3)\nu\tau}}_{L^2_z}e^{-\f12\nu\tau}d\tau\\
	&\quad \sup_{\tau\in[T,t]}\left\|\fM_{\sigma+1}(\tau,\na)\rr{e^{-|\alpha|\nu\tau}v^\alpha P_{\neq}\th }\right\|_{\lta L^2_{z,v}}\\
	&\ls\nnr{\fM_{\sigma+2}(\tau,\na)\rho_g(\tau)e^{(\sigma+3)\nu\tau}}_{L^2([T,t],L^2_z)}
    \sup_{\tau\in[T,t]}\nnr{\fM_{\sigma+1}(\tau,\na)\rr{e^{-|\alpha|\nu\tau}v^\alpha P_{\neq}\th }}_{\lta L^2_{z,v}}.
    \ea\]
\end{proof}

The bilinear bounds provide the input for the weighted linear estimate. We retain an exponential weight centered at an arbitrary frequency $(K',A')$, so that the same estimate can track packets at different orders of the expansion. The density is controlled first through its Volterra resolvent; Duhamel's formula then gives the nonzero profile, and the homogeneous Fokker--Planck semigroup controls the zero mode.

\begin{prop}[Weighted linear estimates]\label{3.2}
	Let $Z(t)$ be a nonnegative nondecreasing function on $[T,\infty)$.	There exists $0<r_1<0.1$ such that for any $0\leq r_2\leq r_1$,  $(K',A')\in\Z\times\R$, $\mathrm{K}>0$, the following holds:
	\begin{itemize}
		\item[(1)] Fix $\sigma\geq2$. If 
		\[\ba
		&\int_T^t\nnr{e^{r_2\<\na_z-K',\na_v-A'\>}\fA\fM_\sigma(\tau,\na)\rr{e^{-|\alpha|\nu\tau}v^\alpha P_{\neq}\tQ(\tau)}}_{\lta L^2_{z,v}}e^{\rm\nu^\f13\tau+\f12\nu\tau}d\tau\leq Z(t),\quad P_0Q(t)=0,\\
		&\nnr{e^{r_2\<\na_z-K',\na_v-A'\>}\fA\fM_\sigma(T,\na)\rr{e^{-|\alpha|\nu T}v^\alpha P_{\neq}\th(T)}}_{\lta L^2_{z,v}}e^{\rm\nu^\f13 T+\f12\nu T}\leq \mathrm{K},\quad \quad P_0h(T)=0,\\
		\ea\]
		then
		\beq\label{nyb}\ba
		&\nnr{e^{r_2\<\na_z-K',\na_v-A'\>}\fA\fM_\sigma(\tau,\na)\rho(\tau)e^{\rm\nu^\f13\tau}}_{L^2([T,t], L^2_{z})}\\
	&+\left\|e^{r_2\<\na_z-K',\na_v-A'\>}\fA\fM_{\sigma-2}(t,\na)\rr{e^{-|\alpha|\nu t}v^\alpha P_{\neq}\th(t)} \right\|_{\lta L^2_{z,v}}\\
		&\times e^{\rm\nu^\f13t}
		+\nnr{v^\alpha\<\na_v\>^{\sigma-2} e^{r_2\<K',\emnt\na_v-A'\>}\fA\fM(t,0,\emnt\na_v)P_{0}h(t)}_{\lta L^2_{v}}e^{\nu t}\ls Z(t)+\mathrm{K}.
		\ea\eeq
		\item[(2)] Fix $\sigma\geq0$. If 
		\[\ba
		&\int^t_T\nnr{v^\alpha\<\na_v\>^\sigma e^{r_2\<K',\emntau\na_v-A'\>}\fA\fM(\tau,0,\emntau\na_v)P_{0}Q(\tau)}_{\lta L^2_{v}}e^{\nu \tau}d\tau\leq Z(t),\quad P_{\neq}\tQ(t)=0,\\
	    &\nnr{v^\alpha\<\na_v\>^\sigma e^{r_2\<K',e^{-\nu T}\na_v-A'\>}\fA\fM(T,0,e^{-\nu T}\na_v)P_{0}h(T)}_{\lta L^2_{v}}e^{\nu T}\leq \mathrm{K},\quad P_{\neq}\th(T)=0,
		\ea\]
		Then 
		\beq\label{myb}\ba
		\nnr{v^\alpha\<\na_v\>^\sigma e^{r_2\<K',\emnt\na_v-A'\>}\fA\fM(t,0,\emnt\na_v)P_{0}h(t)}_{\lta L^2_{v}}e^{\nu t}\ls Z(t)+\mathrm{K},\quad P_{\neq}\th(t)=0.
		\ea\eeq
	\end{itemize}
\end{prop}
\begin{proof}
The main point is to close the density estimate before estimating the profile. The resolvent controls the Maxwellian response, the small background zero mode is absorbed perturbatively, and Lemma~\ref{resonance} controls the echo kernel through the time dependence of $\fA$. Once the density is bounded, the profile estimates follow from Duhamel's formula and the bilinear estimates.

	We define the following multiplier
	\[\fB_\sigma(t,k,\xi)=e^{r_2\<k-K',\xi-A'\>}\fA(t,k,\xi)\fM_\sigma(t,k,\xi)\]
	for $\sigma\geq0$.
	
	\underline{Proof of (1).} \uwave{Step 1: Estimate of $\rho$.} Using the change of variables \eqref{coordinateshift}, taking the Fourier transform, and applying Duhamel's formula, we obtain
	\beq\label{duhah}\ba
	\wh{\th}(t,k,\xi)&=\wh{\th}(T,k,\xi)S(t,T,k,\xi)+\int_T^t\wh{\tQ}(\tau,k,\xi)S(t,\tau,k,\xi)d\tau-\int_T^t\wh{\rho}(\tau,k)\f{\bar{\eta}(\tau,k,\xi)}{k}\\
	&\times\wh{\sM}(\bar{\eta}(\tau,k,\xi))S(t,\tau,k,\xi)d\tau
	-\sum_\ell\int_T^t\wh{\rho}(\tau,\ell)\f{\bar{\eta}(\tau,k,\xi)}{\ell}\wh{\udl{\th}}(\tau,k-\ell,\xi-\ell\tentau)S(t,\tau,k,\xi)d\tau
	\ea\eeq
	For $k>0$ we have 
	\beq\label{rba}\ba
	\wh{\rho}(t,k)&=\wh{\th}(T,k,k\tent)S(t-T,k)+\int_T^t\wh{\tQ}(\tau,k,k\tent)S(t-\tau,k)d\tau-\int_T^t\wh{\rho}(\tau,k)\f{1-e^{-\nu(t-\tau)}}{\nu}\\
	&\times\wh{\sM}(k\f{1-e^{-\nu(t-\tau)}}{\nu})S(t-\tau,k)d\tau
	-\int_T^t\wh{\rho}(\tau,k)\f{1-e^{-\nu(t-\tau)}}{\nu}\wh{\uh}(\tau,0,k\f{1-e^{-\nu(t-\tau)}}{\nu})S(t-\tau,k)d\tau\\
	&-\int_T^t\wh{\rho}(\tau,k+1)\f{k}{k+1}\f{1-e^{-\nu(t-\tau)}}{\nu}\wh{\uth}(\tau,-1,k\tent-(k+1)\tentau)S(t-\tau,k)d\tau.
	\ea\eeq
	For $k<0$ we have
	\[\ba
	\wh{\rho}(t,k)&=\wh{\th}(T,k,k\tent)S(t-T,k)+\int_T^t\wh{\tQ}(\tau,k,k\tent)S(t-\tau,k)d\tau-\int_T^t\wh{\rho}(\tau,k)\f{1-e^{-\nu(t-\tau)}}{\nu}\\
	&\times\wh{\sM}(k\f{1-e^{-\nu(t-\tau)}}{\nu})S(t-\tau,k)d\tau
	-\int_T^t\wh{\rho}(\tau,k)\f{1-e^{-\nu(t-\tau)}}{\nu}\wh{\uh}(\tau,0,k\f{1-e^{-\nu(t-\tau)}}{\nu})S(t-\tau,k)d\tau\\
	&-\int_T^t\wh{\rho}(\tau,k-1)\f{k}{k-1}\f{1-e^{-\nu(t-\tau)}}{\nu}\wh{\uth}(\tau,1,k\tent-(k-1)\tentau)S(t-\tau,k)d\tau,
	\ea\]
	where we use \eqref{963} and $\mathrm{supp} \wh{\vf}\subset[-0.9,0.9]$. By symmetry we only need to consider the $k>0$ case.	For $k\neq0$ we introduce the following resolvent
	\beq\label{definereso}\wh{R}(t,k)=\tent\wh{\sM}(k\tent)S(t,k)-\int_0^t\wh{R}(\tau,k)\f{1-e^{-\nu(t-\tau)}}{\nu}\wh{\sM}(k\f{1-e^{-\nu(t-\tau)}}{\nu})S(t-\tau,k)d\tau.\eeq
	By Proposition~2.1 of \cite{LY}, there exists $r_1>0$ such that
	\beq\label{resolvent}\nr{\wh{R}(t,k)}\ls\f{1}{|k|}e^{-10r_1|k|t}.\eeq
	Then we can rewrite the density equation in the following way
	\[\ba
	&\wh{\rho}(t,k)=\wh{H}(t,k)-\int_T^t\wh{H}(\tau,k)\wh{R}(t-\tau,k)d\tau,\quad 
	\wh{H}(t,k)=\wh{\th}(T,k,k\tent)S(t-T,k)\\
	+&\int_T^t\wh{\tQ}(\tau,k,k\tent)S(t-\tau,k)d\tau-\int_T^t\wh{\rho}(\tau,k)\f{1-e^{-\nu(t-\tau)}}{\nu}\wh{\uh}(\tau,0,k\f{1-e^{-\nu(t-\tau)}}{\nu})S(t-\tau,k)d\tau\\
	-&\int_T^t\wh{\rho}(\tau,k+1)\f{k}{k+1}\f{1-e^{-\nu(t-\tau)}}{\nu}\wh{\uth}(\tau,-1,k\tent-(k+1)\tentau)S(t-\tau,k)d\tau.
	\ea\]
	By Young's inequality and \eqref{resolvent} we have
	\beq\label{egz}\ba
	&\nnr{\fB_\sigma(\tau,\na)\rho(\tau)e^{\rm\nu^\f13\tau}}_{L^2([T,t], L^2_{z})}\leq\rr{1+\sup_k\int_0^\infty e^{r_2|kt|+|kt|^\f13+3\rm|kt|+10}|\wh{R}(t,k)|dt}\\
	&\times\nnr{\fB_\sigma(\tau,\na)H(\tau)e^{\rm\nu^\f13\tau}}_{L^2([T,t], L^2_{z})}
	\ls\nnr{\fB_\sigma(\tau,\na)H(\tau)e^{\rm\nu^\f13\tau}}_{L^2([T,t], L^2_{z})},
	\ea\eeq
	where we use 
	\[\ba
	&\fM(t,k,k\tent)\leq\fM(\tau,k,k\tent)e^{\rm|k|(t-\tau)+1}
	\leq\fM(\tau,k,k\tentau)\fM(\tau,0,k\f{\emntau-\emnt}{\nu})e^{\rm|k|(t-\tau)+1}\\
	\leq& \fM(\tau,k,k\tentau)e^{2\rm|k|(t-\tau)+1},\\
	&\quad\<k\tent\>\leq\<k\tentau\>+\<k\f{\emntau-\emnt}{\nu}\>\leq\<k\tentau\>+\<k(t-\tau)\>.
	\ea\]
	The first inequality uses \eqref{MSkt}. The second inequality uses \eqref{M1}. The third inequality uses \eqref{M0}. Thus we only need to estimate $\nnr{\fB_\sigma(\tau,\na)H(\tau)e^{\rm\nu^\f13\tau}}_{L^2([T,t], L^2_{z})}$. Denote the four terms in $\wh{H}(t,k)$ by $\wh{\aleph_1}, \wh{\aleph_2}, \wh{\aleph_3}, \wh{\aleph_4}$.
	\beq\label{egx}\ba
	&\nnr{\fB_\sigma(t,\na)\aleph_1(t)e^{\rm\nu^\f13t}}_{L^2([T,\tilde{T}], L^2_{z})}^2\ls
	\sum_{k\neq0}\int_T^{\tilde{T}}\nr{\sF\rr{\fB_\sigma(t,\na)\th}(T,k,k\tent)S(t-T,k)e^{\rm\nu^\f13t}}^2dt\\
	\ls&\sum_{k\neq0}\int_T^{\tilde{T}}\nr{\sF\rr{\fB_\sigma(T,\na)\th}(T,k,k\tent)}^2e^{-\nu t}dt\ e^{2\rm\nu^\f13T+\nu T}\leq \nnr{\fB_\sigma(T,\na)P_{\neq}\th(T)}_{L^2_{z,v}}^2e^{2\rm\nu^\f13T+\nu T}.
	\ea\eeq
	The second inequality uses the following inequality
	\[\ba
	&\fM(t,k,k\tent)S^{0.1}(t-T,k)e^{\rm\nu^\f13t}\ls\fM(T,k,k\tent)e^{\rm\nu^\f13T}e^{\f12\nu(T-t)}e^{\rm\min\{\nu T,1\}+\rm\nu},
	\ea\]
	which can be deduced from \eqref{MSk} and \eqref{Snu13}. For the third step we make the change of variables $t'=\f{1-e^{-\nu t}}{\nu}$. Similarly for $\aleph_2$ we have
	\beq\label{egc}\ba
	&\nnr{\fB_\sigma(t,\na)\aleph_2(t)e^{\rm\nu^\f13t}}_{L^2([T,\tilde{T}], L^2_{z})}\ls\nnr{\int_T^t\sF\rr{\fB_\sigma(t,\na)P_{\neq}\tQ}(\tau,k,k\tent)S(t-\tau,k)e^{\rm\nu^\f13t}d\tau}_{L^2([T,\tilde{T}], L^2_{k})}\\
	\ls&\int_T^{\tilde{T}}\nnr{\sF\rr{\fB_\sigma(\tau,\na)P_{\neq}\tQ}(\tau,k,k\tent)e^{-\f12\nu t}}_{L^2([T,\tilde{T}], L^2_{k})}e^{\rm\nu^\f13\tau+\f12\nu\tau}d\tau\\
	&\ls\int_T^{\tilde{T}}\nnr{\fB_\sigma(\tau,\na)P_{\neq}\tQ(\tau)}_{L^2_{z,v}}e^{\rm\nu^\f13\tau+\f12\nu\tau}d\tau.
	\ea\eeq
	The second inequality uses Minkowski's inequality. For $\aleph_3$,
	\[\ba
	&\nr{\sF\rr{\fB_\sigma(t,\na)\aleph_3}(t,k)}e^{\rm\nu^\f13t}\ls\int_T^t\nr{\sF\rr{\fB_\sigma(\tau,\na)\rho}(\tau,k)}\nr{\sF\rr{e^{r_2\<\emntau\na_v\>}\fA\fM_{\sigma}(\tau,0,\emntau\na_v)\pa_vP_0\uh}\right.\\
		&\left.\rr{\tau,0,k\f{1-e^{-\nu(t-\tau)}}{\nu}}}
	e^{\rm\nu^\f13\tau}S^{0.9}(t-\tau,k)d\tau
	\ls\int_T^t\nr{\sF\rr{\fB_\sigma(\tau,\na)\rho}(\tau,k)}e^{\rm\nu^\f13\tau}\<\f{1-e^{-\nu(t-\tau)}}{\nu}\>^{-2}e^{-\nu(t-\tau)}d\tau\\
	&\times\sup_{\tau\in[T,t]} \nnr{v^\alpha\<\na_v\>^2 e^{r_2\<\emntau\na_v\>}\fA\fM_{\sigma}(\tau,0,\emntau\na_v)\pa_vP_0\uh(\tau)}_{\lta L^2_v}\ls\ueps\int_T^t\nr{\sF\rr{\fB_\sigma(\tau,\na)\rho}(\tau,k)}e^{\rm\nu^\f13\tau}\\
	&\times\<\f{1-e^{-\nu(t-\tau)}}{\nu}\>^{-2}e^{-\nu(t-\tau)}d\tau.
	\ea\]
	The last inequality uses
	\beq\label{beo}\ba
	&\nnr{v^\alpha\<\na_v\>^3 e^{r_2\<\emntau\na_v\>}\fA\fM_{\sigma}(\tau,0,\emntau\na_v)e^{\nu(\tau-1)L}P_0\uh(1)}_{\lta L^2_v}=\nnr{v^\alpha\<\na_v\>^3 \fM_{\sigma}(\tau,0,\emntau\na_v)\right.\\
	&\left.\times e^{\nu(\tau-1)L}e^{r_2\<e^{-\nu}\na_v\>}\fA(\tau,0,e^{-\nu}\na_v)P_0\uh(1)}_{\lta L^2_v}
	\ls\nnr{v^\alpha\<\na_v\>^3 \fA\fM_{\sigma}(1,0,e^{-\nu}\na_v)e^{r_2\<e^{-\nu}\na_v\>}P_0\uh(1)}_{\lta L^2_v}\ls\ueps,
	\ea\eeq
	which follows from \eqref{957} and Proposition~\ref{etl}. Thus by Schur's test 
	\beq\label{egv}
	\nnr{\fB_\sigma(t,\na)\aleph_3(t)e^{\rm\nu^\f13t}}_{L^2([T,\tilde{T}], L^2_{z})}\ls\ueps\nnr{\fB_\sigma(t,\na)\rho(t)e^{\rm\nu^\f13t}}_{L^2([T,\tilde{T}],L^2_z)}.
	\eeq
	For $\aleph_4$,
	\[\ba
	&\nr{\sF\rr{\fB_\sigma(t,\na)\aleph_4}(t,k)}e^{\rm\nu^\f13t}\ls\int_T^t\nr{\sF\rr{\fB_\sigma(\tau,\na)\rho}(\tau,k+1)}\f{k(1-e^{-\nu(t-\tau)})}{(k+1)\nu}\nr{\sF\rr{e^{r_2\<\na\>}\fA\fM_{\sigma}(\tau,\na)\uth}\right.\\&\left.\rr{\tau,-1,k\tent-(k+1)\tentau}}
	\f{\fA(t,k,k\tent)}{\fA(\tau,k,k\tent)}e^{\rm\nu^\f13\tau}S^{0.9}(t-\tau,k)d\tau\\
	\ls&\ueps\int_T^t\nr{\sF\rr{\fB_\sigma(\tau,\na)\rho}(\tau,k+1)}e^{\rm\nu^\f13\tau}\f{1-\emnt}{\nu k}\wh{\vf}(k\tent-(k+1)\tentau)\f{\fA(t,k,k\tent)}{\fA(\tau,k,k\tent)}e^{-\nu t}d\tau\\
	\ls&D_1^{-\f14}\int_T^t\nr{\sF\rr{\fB_\sigma(\tau,\na)\rho}(\tau,k+1)}e^{\rm\nu^\f13\tau}\wh{\vf}(k\tent-(k+1)\tentau)|k|e^{-\nu t}d\tau
	\ea\]
	The second inequality uses 
	\beq\label{hudlneq}\ba
	&\nr{\sF\rr{e^{r_2\<\na\>}\fA\fM_{\sigma}(\tau,\na)\uth}\rr{\tau,-1,\xi}}\ls\ueps\wh{\vf}(\xi)S(\tau,0,-1,\xi)\fM(\tau,-1,\xi)\\
	&\ls\ueps\wh{\vf}(\xi)S^{0.9}(\tau,0,-1,\xi)\ls\ueps\wh{\vf}(\xi) e^{-0.009\nu^\f13\tau},
	\ea\eeq
	which follows from \eqref{963}, $\mathrm{supp} \wh{\vf}(\xi)\subset[-0.9,0.9]$, \eqref{MSk} and \eqref{Snu13}. Moreover we have $\f{k(1-e^{-\nu(t-\tau)})}{\nu}\leq\entau\<-1,k\tent-(k+1)\tentau\>\<\tentau\>$. The third inequality uses Lemma~\ref{resonance}. By Schur's test we have
	\beq\label{egn}\nnr{\fB_\sigma(t,\na)\aleph_4(t)e^{\rm\nu^\f13t}}_{L^2([T,\tilde{T}], L^2_{z})}\ls D_1^{-\f14}\nnr{\fB_\sigma(t,\na)\rho(t)e^{\rm\nu^\f13t}}_{L^2([T,\tilde{T}],L^2_z)}.\eeq
	Combining \eqref{egz}, \eqref{egx}, \eqref{egc}, \eqref{egv}, and \eqref{egn} and choosing $D_1$ sufficiently large we obtain 
	\beq\label{rhooo}\ba&\nnr{\fB_\sigma(t,\na)\rho(t)e^{\rm\nu^\f13t}}_{L^2([T,\tilde{T}],L^2_z)}\ls\nnr{\fB_\sigma(T,\na)P_{\neq}\th(T)}_{L^2_{z,v}}e^{\rm\nu^\f13T+\f12\nu T}\\
	&+\int_T^{\tilde{T}}\nnr{\fB_\sigma(\tau,\na)P_{\neq}\tQ(\tau)}_{L^2_{z,v}}e^{\rm\nu^\f13\tau+\f12\nu\tau}d\tau\ls Z(\tilde{T})+\mathrm{K}. \ea\eeq
	
	\uwave{Step 2: Estimate of $P_{\neq}\th(t)$.} By \eqref{duhah} we have 
	\beq\label{369}\ba&\wh{P_{\neq}\th}(t,k,\xi)=\wh{P_{\neq}\th}(T,k,\xi)S(t,T,k,\xi)+\int_T^t\wh{P_{\neq}\tQ}(\tau,k,\xi)S(t,\tau,k,\xi)d\tau\\
	&+\int_T^t\sF\rr{P_{\neq}\cN^\dagger\rr{\th,\sM(\emntau v)\emntau}}(\tau,k,\xi)S(t,\tau,k,\xi)d\tau+\int_T^t\sF\rr{P_{\neq}\cN^\dagger\rr{\th,\udl{\th}}}(\tau,k,\xi)S(t,\tau,k,\xi)d\tau.\ea\eeq
	Denote the four terms on the right-hand side by $\wh{\beth_1}$, $\wh{\beth_2}$, $\wh{\beth_3}$, and $\wh{\beth_4}$. We first estimate $\beth_1$.
	\beq\label{beth11}\ba&\nnr{\fB_{\sigma-2}(t,\na)\rr{e^{-\nu|\alpha|t}v^\alpha\beth_1}}_{\lta L^2_{z,v}}e^{\rm\nu^\f13 t}\ls\sum_{|\alpha_1|+|\alpha_2|\leq2}\nnr{\fB_{\sigma-2}(t,k,\xi)D_\xi^{\alpha_1}\wh{P_{\neq}\th}(T,k,\xi)D_\xi^{\alpha_2}S(t,T,k,\xi)}_{L^2_{k,\xi}}\\
	&\times e^{\rm\nu^\f13 t}e^{-\nu\rr{|\alpha_1|+|\alpha_2|}t}
	\ls\nnr{\fB_{\sigma-2}(t,k,\xi)D_\xi^{\alpha}\wh{P_{\neq}\th}(T,k,\xi)S^{0.9}(t,T,k,\xi)e^{-\nu|\alpha|t}}_{\lta L^2_{k,\xi}}e^{\rm\nu^\f13 t}\\
	\ls&\nnr{\fB_{\sigma-2}(T,k,\xi)D_\xi^{\alpha}\wh{P_{\neq}\th}(T,k,\xi)e^{-\nu|\alpha|T}}_{\lta L^2_{k,\xi}}e^{\rm\nu^\f13 T}\ls\nnr{\fB_{\sigma-2}(T,\na)\rr{e^{-\nu|\alpha|T}v^\alpha P_{\neq}\th(T)}}_{\lta L^2_{z,v}}e^{\rm\nu^\f13 T}.
	\ea\eeq
	The third inequality uses \eqref{Sds}. The fourth inequality uses $\fM(t,k,\xi)S^{0.9}(t,T,k,\xi)e^{\rm\nu^\f13 t}\ls\fM(T,k,\xi)e^{\rm\nu^\f13 T}$, which follows from \eqref{MSk} and \eqref{Snu13}. Similarly for $\beth_2$ we have
	\beq\label{beth21}\nnr{\fB_{\sigma-2}(t,\na)\rr{e^{-\nu|\alpha|t}v^\alpha\beth_2}}_{\lta L^2_{z,v}}e^{\rm\nu^\f13 t}\ls\int_T^t\nnr{\fB_{\sigma-2}(\tau,\na)\rr{e^{-\nu|\alpha|\tau}v^\alpha P_{\neq}\tQ(\tau)}}_{\lta L^2_{z,v}}e^{\rm\nu^\f13 \tau}d\tau.
	\eeq
	We next estimate $\beth_4$.
	\beq\label{beth41}\ba
	&\nnr{\fB_{\sigma-2}(t,\na)\rr{e^{-\nu|\alpha|t}v^\alpha\beth_4}}_{\lta L^2_{z,v}}e^{\rm\nu^\f13 t}\ls\int_T^t\nnr{\fB_{\sigma-2}(\tau,\na)\rr{e^{-\nu|\alpha|\tau}v^\alpha P_{\neq}\cN^\dagger\rr{\th,\udl{\th}}(\tau) }}_{\lta L^2_{z,v}}e^{\rm\nu^\f13 \tau}d\tau\\
	&\ls\nnr{\fB_\sigma(\tau,\na)\rho(\tau)e^{\rm\nu^\f13\tau}}_{L^2([T,t],L^2_z)}\sup_{\tau\in[T,t]}\left\{\nnr{e^{r_2\<\na\>}\fA\fM_{\sigma-1}(\tau,\na)\rr{e^{-\nu|\alpha|\tau}v^\alpha P_{\neq}\udl{\th}(\tau)}}_{\lta L^2_{z,v}}e^{\f32\nu\tau}\right.\\
	&\left.+\nnr{v^\alpha e^{r_2\<\emntau\na_v\>}\fA\fM_{\sigma-2}(\tau,0,\emntau\na_v)\pa_v P_{0}\udl{h}(\tau)}_{\lta L^2_{v}}\right\}
	\ls\ueps\nnr{\fB_\sigma(\tau,\na)\rho(\tau)e^{\rm\nu^\f13\tau}}_{L^2([T,t],L^2_z)}.
	\ea\eeq
	The third inequality follows from \eqref{cNk1} after checking that the multipliers $e^{r_2\<\na_z-K',\na_v-A'\>}$ and $\fA(t,\na)$ can be distributed between $\rho$ and $\udl{h}$. The last inequality uses \eqref{beo} and \eqref{hudlneq}. Similarly we can estimate $\beth_3$.
	\beq\label{beth31}\ba
	&\nnr{\fB_{\sigma-2}(t,\na)\rr{e^{-\nu|\alpha|t}v^\alpha\beth_3}}_{\lta L^2_{z,v}}e^{\rm\nu^\f13 t}\ls\nnr{\fB_\sigma(\tau,\na)\rho(\tau)e^{\rm\nu^\f13\tau}}_{L^2([T,t],L^2_z)}\\
	&\times\sup_{\tau\in[T,t]}\nnr{v^\alpha e^{r_2\<\emntau\na_v\>}\fA\fM_{\sigma-2}(\tau,0,\emntau\na_v)\pa_v \sM(v)}_{\lta L^2_{v}}\ls\nnr{\fB_\sigma(\tau,\na)\rho(\tau)e^{\rm\nu^\f13\tau}}_{L^2([T,t],L^2_z)}.
	\ea\eeq
	Combining \eqref{beth11}, \eqref{beth21}, \eqref{beth31}, and \eqref{beth41} and using \eqref{rhooo} we obtain
	\beq\label{neqneqh}\nnr{\fB_{\sigma-2}(t,\na)\rr{e^{-\nu|\alpha|t}v^\alpha P_{\neq}\th(t)}}_{\lta L^2_{z,v}}e^{\rm\nu^\f13 t}\ls Z(t)+\mathrm{K}.\eeq
	
	\uwave{Step 3: Estimate of $P_0h(t)$.}  By Duhamel's formula we have
	\[P_0h(t)=e^{\nu(t-T)L}P_0h(T)+\int_T^te^{\nu(t-\tau)L}P_0Q(\tau)d\tau+\int_T^te^{\nu(t-\tau)L}P_0\cN(h,\uh)d\tau=\int_T^te^{\nu(t-\tau)L}P_0\cN(h,\uh)d\tau.\]
	By Proposition~\ref{etl} we have
	\beq\label{p0poh}\ba
	&\nnr{v^\alpha \<\na_v\>^{\sigma-2}\fB_0(t,0,\emnt\na_v)P_0h(t)}_{\lta L^2_{v}}\ent\leq\int_T^t\nnr{v^\alpha \<\na_v\>^{\sigma-2}\fM(t,0,\emnt\na_v)e^{\nu(t-\tau)L}\fA(t,0,\emntau\na_v)\right.\\
	&\left.e^{r_2\<K',\emntau\na_v-A'\>}P_0\cN(h,\uh)(\tau)}_{\lta L^2_{v}}\ent d\tau
	\ls\int_T^t\nnr{v^\alpha \<\na_v\>^{\sigma-2}\fB_0(\tau,0,\emntau\na_v)P_0\cN(h,\uh)(\tau)}_{\lta L^2_{v}}\entau d\tau\\
	\ls&\nnr{\fB_\sigma(\tau,\na)\rho(\tau)e^{\rm\nu^\f13\tau}}_{L^2([T,t],L^2_z)}\sup_{\tau\in[T,t]}\nnr{e^{r_2\<\na\>}\fA\fM_{\sigma-1}(\tau,\na)\rr{e^{-\nu|\alpha|\tau}v^\alpha P_{\neq}\udl{\th}(\tau)}}_{\lta L^2_{z,v}}\\
	\ls&\ueps\nnr{\fB_\sigma(\tau,\na)\rho(\tau)e^{\rm\nu^\f13\tau}}_{L^2([T,t],L^2_z)}.
	\ea\eeq
	The third inequality follows from the proof of \eqref{cN0}, with the additional weights distributed in the same way. Combining \eqref{rhooo}, \eqref{neqneqh}, and \eqref{p0poh} we obtain \eqref{nyb}.
	
	\underline{Proof of (2).} By \eqref{rba} and \eqref{duhah} we obtain $\rho(t)=0$, $P_{\neq}\th(t)=0$. By Duhamel's formula we obtain
	\[P_0h(t)=e^{\nu(t-T)L}P_0h(T)+\int_T^te^{\nu(t-\tau)L}P_0Q(\tau)d\tau.\]
	By Proposition~\ref{etl} we have  
	\[\ba
	&\nnr{v^\alpha \<\na_v\>^{\sigma}\fB_0(t,0,\emnt\na_v)e^{\nu(t-T)L}P_0h(T)}_{\lta L^2_{v}}\ent=\nnr{v^\alpha \<\na_v\>^{\sigma}\fM(t,0,\emnt\na_v)e^{\nu(t-T)L}\fA(t,0,e^{-\nu T}\na_v)\right.\\
		&\left.e^{r_2\<K',e^{-\nu T}\na_v-A'\>}P_0h(T)}_{\lta L^2_{v}}\ent\ls\nnr{v^\alpha \<\na_v\>^{\sigma}\fB_0(T,0,e^{-\nu T}\na_v)P_0h(T)}_{\lta L^2_{v}}e^{\nu T}.
	\ea\]
	Similarly we have
	\[\ba
	&\nnr{v^\alpha \<\na_v\>^{\sigma}\fB_0(t,0,\emnt\na_v)\int_T^te^{\nu(t-\tau)L}P_0Q(\tau)d\tau}_{\lta L^2_{v}}\ent\\
	&\ls\int_T^t\nnr{v^\alpha \<\na_v\>^{\sigma}\fB_0(\tau,0,\emntau\na_v)P_0Q(\tau)}_{\lta L^2_{v}}\entau d\tau.
	\ea\]
	This proves \eqref{myb}.
\end{proof}

We now apply Proposition~\ref{3.2} to the coefficients $h^{(n,q)}$ defined in \eqref{hnq}. Their principal frequencies add under bilinear interactions. A small decrease in the localization radius at each order absorbs the derivative losses, allowing an induction over $n$. The next proposition supplies the Fourier tail estimates used throughout the resonant analysis; the sharper size bound in terms of the cascade envelope will be proved in Section~5.

\begin{prop}[Fourier localization of the expansion coefficients]\label{3.3}
	For $1\leq n\leq A^{N'}$, we set $\tilde{r}_n=r_1(1-\f14nA^{-N'})$. The following holds for $1\leq n\leq A^{N'}$, $|q|\leq n$, $2\mid n-q$, $t\geq T$
	\beq\label{hnqenergy}\ba
	&\nnr{e^{\tilde{r}_n\<\na_z-qK,\na_v-qA\>}\fA\fM_4(\tau,\na)\rho^{(n,q)}(\tau)e^{\rm\nu^\f13\tau}}_{L^2([T,t], L^2_{z})}+\left\|e^{\tilde{r}_n\<\na_z-qK,\na_v-qA\>}\fA\fM_{2}(t,\na)\rr{e^{-|\alpha|\nu t}v^\alpha\right.\right.\\
	&\left.\left. P_{\neq}\th^{(n,q)}(t)} \right\|_{\lta L^2_{z,v}}e^{\rm\nu^\f13t}
	+\nnr{v^\alpha\<\na_v\> e^{\tilde{r}_n\<qK,\emnt\na_v-qA\>}\fA\fM(t,0,\emnt\na_v)P_{0}h^{(n,q)}(t)}_{\lta L^2_{v}}e^{\nu t}\leq e^{10(2n-1)D_1K}.
	\ea\eeq
\end{prop}
\begin{proof}
We argue by induction on the interaction order. The initial packet gives the base case. At higher orders, the gap between successive localization radii absorbs the extra polynomial frequency weights in the bilinear sources, after which Proposition~\ref{3.2} propagates the estimate.

	We set $\fB_\sigma^{(n,q)}(t,k,\xi)=e^{\tilde{r}_n\<k-qK,\xi-qA\>}\fA\fM_\sigma(t,k,\xi)$. For the base case $n=1$, recall \eqref{fonemp}. By Proposition~\ref{3.2} we have
	\[\ba
	&\nnr{\fB_4^{(1,1)}(\tau,\na)\rho^{(1,1)}(\tau)e^{\rm\nu^\f13\tau}}_{L^2([T,t], L^2_{z})}+\left\|\fB_2^{(1,1)}(t,\na)\rr{e^{-|\alpha|\nu t}v^\alpha P_{\neq}\th^{(1,1)}(t)} \right\|_{\lta L^2_{z,v}}e^{\rm\nu^\f13t}\\
	&+\nnr{v^\alpha\<\na_v\>^2\fB_0^{(1,1)}(t,0,\emnt\na_v)P_{0}h^{(1,1)}(t)}_{\lta L^2_{v}}e^{\nu t}\\
	&\ls\left\|\fB_4^{(1,1)}(T,\na)\rr{e^{-|\alpha|\nu T}v^\alpha P_{\neq}\th^{(1,1)}(T)} \right\|_{\lta L^2_{z,v}}e^{\rm\nu^\f13T+\f12\nu T}\\
	&+\nnr{v^\alpha\<\na_v\>^2\fB_0^{(1,1)}(T,0,e^{-\nu T}\na_v)P_{0}h^{(1,1)}(T)}_{\lta L^2_{v}}e^{\nu T}
	\ls\eps\nnr{\<K,\xi\>^4\fA\fM(T,K,\xi)D^\alpha_\xi\wh{\vf}(\xi-A)S^{0.9}(T,0,K,\xi)}_{\lta L^2_\xi}\\
	&\ls\eps A^4e^{1.1D_1K}\nnr{D^\alpha_\xi\wh{\vf}(\xi)}_{\lta L^2_\xi}\ls\eps A^4e^{1.1D_1K},
	\ea\]
	where we use $\log \fA(T,K,\xi)\leq D_1\<\ueps|K,\xi|\>^\f13\leq1.1D_1K$ on the support of $\wh{\th^{(1,1)}}(T)$.
	
	Suppose \eqref{hnqenergy} holds through order $n-1$, where $n\geq2$. We prove it at order $n$. Recall \eqref{hnq}. Fix $(n,q)=(n_1,q_1)+(n_2,q_2)$. We next estimate the bilinear terms.
	\beq\label{ucm}\ba
	&\int_T^t\left\|\fB_4^{(n,q)}(\tau,\na)\rr{e^{-|\alpha|\nu \tau}v^\alpha P_{\neq}\cN^\dagger(\th^{(n_1,q_1)},\th^{(n_2,q_2)})(\tau)} \right\|_{\lta L^2_{z,v}}e^{\rm\nu^\f13\tau+\f12\nu\tau}d\tau\\
	&\ls\int_T^tA^{4(N'+1)}\left\|\fB_0^{(n-1,q)}(\tau,\na)\right.\\
	&\left.\rr{e^{-|\alpha|\nu \tau}v^\alpha P_{\neq}\cN^\dagger(\th^{(n_1,q_1)},\th^{(n_2,q_2)})(\tau)} \right\|_{\lta L^2_{z,v}}e^{\rm\nu^\f13\tau+\f12\nu\tau}d\tau\\
	&\ls A^{4(N'+1)}\nnr{\fB_2^{(n_1,q_1)}(\tau,\na)\rho^{(n_1,q_1)}(\tau)e^{\rm\nu^\f13\tau}}_{L^2([T,t], L^2_{z})}\\
    &\times\sup_{\tau\in[T,t]}\rr{\left\|\fB_1^{(n_2,q_2)}(\tau,\na)\rr{e^{-|\alpha|\nu \tau}v^\alpha P_{\neq}\th^{(n_2,q_2)}(\tau)} \right\|_{\lta L^2_{z,v}}e^{2\nu \tau}\right.\\
    &\left.+\nnr{v^\alpha\<\na_v\>\fB_0^{(n_2,q_2)}(\tau,0,\emntau\na_v)P_{0}h^{(n_2,q_2)}(\tau)}_{\lta L^2_{v}}e^{\f12\nu \tau}}\ls A^{4(N'+1)}e^{10(2n-2)D_1K}.
	\ea\eeq
	The first inequality uses 
	\[\<k,\xi\>^4e^{-\f{r_1}{4}A^{-N'}\<k-qK,\xi-qA\>}\ls\rr{\<k-qK,\xi-qA\>^4+\<qK,qA\>^4}e^{-\f{r_1}{4}A^{-N'}\<k-qK,\xi-qA\>}\ls A^{4N'}+A^{4(N'+1)}.\]
	The second inequality follows from the weighted version of the bilinear estimate. By \eqref{cN0}
	\beq\label{ucl}\ba
	&\int_T^t\nnr{v^\alpha\<\na_v\>\fB_0^{(n,q)}(\tau,0,\emntau\na_v)P_{0}\cN(h^{(n_1,q_1)},h^{(n_2,q_2)})(\tau)}_{\lta L^2_{v}}\entau d\tau\\
	&\ls\nnr{\fB_3^{(n_1,q_1)}(\tau,\na)\rho^{(n_1,q_1)}(\tau)e^{\rm\nu^\f13\tau}}_{L^2([T,t], L^2_{z})}\\
	&\times\sup_{\tau\in[T,t]}\left\|\fB_2^{(n_2,q_2)}(\tau,\na)\rr{e^{-|\alpha|\nu \tau}v^\alpha P_{\neq}\th^{(n_2,q_2)}(\tau)} \right\|_{\lta L^2_{z,v}}\ls e^{10(2n-2)D_1K}.
	\ea\eeq
	By \eqref{hnq}, $h^{(n,q)}=h^{(n,q)}_1+h^{(n,q)}_2$, where $h^{(n,q)}_1$ solves \[\pa_t h_1^{(n,q)}=\cL h_1^{(n,q)}+\sum_{\substack{n_1+n_2=n\\q_1+q_2=q}}
	P_{\neq}\cN(h^{(n_1,q_1)},h^{(n_2,q_2)}),\quad h_1^{(n,q)}(T)=0\] and $h_2^{(n,q)}$ solves \[\pa_t h_2^{(n,q)}=\cL h_2^{(n,q)}+\sum_{\substack{n_1+n_2=n\\q_1+q_2=q}}
	P_{0}\cN(h^{(n_1,q_1)},h^{(n_2,q_2)}),\quad h_2^{(n,q)}(T)=0.\]
	Combining \eqref{ucl}, \eqref{ucm} and Proposition~\ref{3.2} we obtain
	\[\ba
	&\nnr{\fB_4^{(n,q)}(\tau,\na)\rho^{(n,q)}(\tau)e^{\rm\nu^\f13\tau}}_{L^2([T,t], L^2_{z})}+\left\|\fB_2^{(n,q)}(t,\na)\rr{e^{-|\alpha|\nu t}v^\alpha P_{\neq}\th^{(n,q)}(t)} \right\|_{\lta L^2_{z,v}}e^{\rm\nu^\f13t}\\
	&+\nnr{v^\alpha\<\na_v\>\fB_0^{(n,q)}(t,0,\emnt\na_v)P_{0}h^{(n,q)}(t)}_{\lta L^2_{v}}e^{\nu t}\ls n^2A^{4(N'+1)}e^{10(2n-2)D_1K}\ls A^{6N'+4}e^{10(2n-2)D_1K},
	\ea\]
    where we use $\#\{(n_1,q_1,n_2,q_2):(n,q)=(n_1,q_1)+(n_2,q_2), 2\mid n_1-q_1, 2\mid n_2-q_2\}\ls n^2$. Note that $A^{6N'+5}\leq e^{10D_1K}$. This proves \eqref{hnqenergy} at order $n$ and closes the induction.
\end{proof}

\section{Linear analysis: the homogeneous case}
This section studies the homogeneous linearized equation
\[\left\{\ba&\partial_t h=\mathcal{L} h, \\ &\wh{\th}(T,k,\xi)=\eps\wh{\vf}(\xi- A)S(T,0, K,\xi)\vv{1}_{k= K}\ea\right.\]
on $[T,\infty)$. This is the equation for the leading coefficient $h^{(1,1)}$. Our goal is to quantify its amplification along the echo chain while retaining the attenuation caused by collisions. We first use Section~2 to isolate one resonant interaction on each critical interval. We then estimate the density response, derive the one-step transfer in Proposition~\ref{4.3}, and iterate it in Proposition~\ref{4.4}. Positivity of the leading convolution gives a lower bound, while Proposition~\ref{prop:homogeneous-upper} converts the discrete upper bounds into norms on the full time interval. These results will determine the cascade envelope and preserve the leading contribution in the nonlinear construction.

For $0\leq k\leq K$, set
\beq\label{4.1}t_k=-\f1\nu\log\rr{1-\f{\nu A}{k+\f12}}.\eeq
Set $D_2=\f{D_1}{r_1}$, where $r_1$ is defined in Proposition~\ref{3.2}. The critical times order the successive resonances from mode $K$ down to mode $1$. The next proposition reduces the density equation and the transfer from $k$ to $k-1$ to their resonant contributions. Its exponentially small remainders follow from the Fourier localization bound \eqref{hnqenergy}.

\begin{prop}[Reduction to a single resonant interaction]\label{prop:resonant-reduction}
	For $1\leq k\leq K$ and $t\in[t_k,t_{k-1}]$, we have
	\beq\label{rhoktk}\ba
	\wh{\rho}(t,k)&=\wh{\th}(t_k,k,k\tent)S(t-t_k,k)-\int_{t_k}^t\wh{\rho}(\tau,k)\f{1-e^{-\nu(t-\tau)}}{\nu}\wh{\sM}(k\f{1-e^{-\nu(t-\tau)}}{\nu})S(t-\tau,k)d\tau\\
	&-\int_{t_k}^t\wh{\rho}(\tau,k)\f{1-e^{-\nu(t-\tau)}}{\nu}\wh{\uh}(\tau,0,k\f{1-e^{-\nu(t-\tau)}}{\nu})S(t-\tau,k)d\tau+\wh{\sR_1}(t,k).	
	\ea\eeq
	For $2\leq k\leq K$, we have
	\beq\label{hktk}\ba\wh{\th}(t_{k-1},k-1,\xi)=&-\int_T^{t_{k-1}}\vv{1}_{|k\tentau-A|\leq40D_2K}\wh{\rho}(\tau,k)\f{\bar{\eta}(\tau,k-1,\xi)}{k}\wh{\udl{\th}}(\tau,-1,\xi-k\tentau)\\
	&\times S(t_{k-1},\tau,k-1,\xi)d\tau+\wh{\sR_2}(k-1,\xi),\ea\eeq
    where $\sR_1$ and $\sR_2$ satisfy
    \beq\label{341}\ba\nnr{\sF\rr{\fM_4(\tau,\na)\sR_1}(\tau,k)e^{\rm\nu^\f13\tau}}_{L^2[t_k,t_{k-1}]}+ \nnr{\sF\rr{\fM_2(t_{k-1},\na)\sR_2}(k-1,\xi)}_{L^2_\xi}e^{\rm\nu^\f13t_{k-1}}\leq e^{-9D_1K}.
    \ea\eeq
\end{prop}
\begin{proof}
The localization estimate makes every density contribution outside the principal frequency region exponentially small. On a critical interval, this leaves the active density mode and its transfer to the next lower mode. We estimate the discarded terms using the bilinear bounds from Section~2.

	The localization estimate \eqref{hnqenergy} gives
	\beq\label{rhoexsm}\nnr{\vv{1}_{\<\na_z-K,\na_v-A\>\geq40D_2K}\fM_4(\tau,\na)\rho(\tau)e^{\rm\nu^\f13\tau}}_{L^2([T,t],L^2_z)}\leq e^{-10D_1K}.\eeq
	Note that for $t\geq T$,
	\beq\label{rmkk}|k\tent-A|\leq40D_2K\quad\text{implies}\quad1\leq k\leq K,\ t\in[t_k,t_{k-1}].\eeq
    By Duhamel's formula (see \eqref{rba}), \eqref{rhoktk} holds with 
    \[\wh{\sR_1}(t,k)=\int_{t_k}^t\sF\rr{\cN^\dagger(P_{k+1}\th,P_{-1}\underline{\th})}(\tau,k,k\tent)S(t-\tau,k)d\tau.\]
    Consequently,
    \beq\label{RRR111}\ba
    &\nnr{\sF\rr{\fM_4(\tau,\na)\sR_1}(\tau,k)e^{\rm\nu^\f13\tau}}_{L^2[t_k,t_{k-1}]}\ls\int_{t_k}^{t_{k-1}}\nnr{\fM_4(\tau,\na)\cN^\dagger(P_{k+1}\th,P_{-1}\udl{\th})(\tau)}_{L^2_{z,v}}e^{\rm\nu^\f13\tau+\f12\nu\tau}d\tau\\
    \ls&\nnr{\sF\rr{\fM_6(\tau,\na)\rho}(\tau,k+1)e^{\rm\nu^\f13\tau}}_{L^2[t_k,t_{k-1}]}\sup_{\tau\in[t_k,t_{k-1}]}\nnr{\fM_5(\tau,\na)P_{\neq}\udl{\th}(\tau)}_{L^2_{z,v}}e^{2\nu\tau}\\
    \ls&\ueps A^2\nnr{\sF\rr{\vv{1}_{\<\na_z-K,\na_v-A\>\geq40D_2K}e^{\f{\tilde{r}_1}4\<\na_z-K,\na_v-A\>}\fM_4(\tau,\na)\rho}(\tau,k+1)e^{\rm\nu^\f13\tau}}_{L^2[t_k,t_{k-1}]}\ls e^{-9.5D_1K}.
    \ea\eeq
    The first inequality follows by the argument of \eqref{egc}. The second inequality uses \eqref{cNk2}. The third inequality uses $\<k,\xi\>^2e^{-\f{\tilde{r}_1}4\<k-K,\xi-A\>}\ls A^2$. The fourth inequality uses \eqref{hnqenergy} and \eqref{rmkk}. By Duhamel's formula (see \eqref{369}), \eqref{hktk} holds with
    \[\ba
    &\wh{\sR_2}(k-1,\xi)=-\int_T^{t_{k-1}}\vv{1}_{|k\tentau-A|\geq40D_2K}\wh{\rho}(\tau,k)\f{\bar{\eta}(\tau,k-1,\xi)}{k}\wh{\udl{\th}}(\tau,-1,\xi-k\tentau)S(t_{k-1},\tau,k-1,\xi)d\tau\\
    &+\int_T^{t_{k-1}}\sF\rr{\cN^\dagger(P_{k-1}\th,\sM(\emntau v)\emntau+P_0\udl{\th})}(\tau,k-1,\xi)S(t_{k-1},\tau,k-1,\xi)d\tau\\
    &+\int_T^{t_{k-1}}\sF\rr{\cN^\dagger(P_{k-2}\th,P_{1}\udl{\th})}(\tau,k-1,\xi)S(t_{k-1},\tau,k-1,\xi)d\tau.
    \ea\]
    Denote the three terms on the right-hand side by $\wh{\beth'_1}$, $\wh{\beth'_2}$, $\wh{\beth'_3}$. We first estimate $\beth'_3$. As in \eqref{beth41} we have
    \beq\label{beth233}\ba
    \nnr{\sF\rr{\fM_2(t_{k-1},\na)\beth'_3}(k-1,\xi)}_{L^2_\xi}e^{\rm\nu^\f13t_{k-1}}\ls\int_T^{t_{k-1}}\nnr{\fM_2(\tau,\na)\cN^\dagger(P_{k-2}\th,P_{1}\udl{\th})}_{L^2_{z,v}}e^{\rm\nu^\f13\tau}d\tau\\
    \ls\nnr{\sF\rr{\fM_4(\tau,\na)\rho}(\tau,k-2)e^{\rm\nu^\f13\tau}}_{L^2[T,t_{k-1}]}\sup_{\tau\in[T,t_{k-1}]}\nnr{\fM_3(\tau,\na)P_{\neq}\udl{\th}(\tau)}_{L^2_{z,v}}e^{\f32\nu\tau}\ls\ueps e^{-10D_1K}.
    \ea\eeq
    The last inequality uses \eqref{rhoexsm} and \eqref{rmkk}. Similarly we have
    \beq\label{beth222}\ba
    &\nnr{\sF\rr{\fM_2(t_{k-1},\na)\beth'_2}(k-1,\xi)}_{L^2_\xi}e^{\rm\nu^\f13t_{k-1}}\ls\nnr{\sF\rr{\fM_4(\tau,\na)\rho}(\tau,k-1)e^{\rm\nu^\f13\tau}}_{L^2[T,t_{k-1}]}\\
    &\times\sup_{\tau\in[T,t_{k-1}]}\nnr{\fM_2(\tau,0,\emntau\na_v)\rr{\pa_v\sM(v)+\pa_vP_0\uh(\tau)}}_{L^2_v}\ls e^{-10D_1K}.
    \ea\eeq
    \beq\label{beth211}\ba
    &\nnr{\sF\rr{\fM_2(t_{k-1},\na)\beth'_1}(k-1,\xi)}_{L^2_\xi}e^{\rm\nu^\f13t_{k-1}}\ls\nnr{\sF\rr{\vv{1}_{|\na_v-A|\geq40D_2K}\fM_4(\tau,\na)\rho}(\tau,k)e^{\rm\nu^\f13\tau}}_{L^2[T,t_{k-1}]}\\
    &\times\sup_{\tau\in[T,t_{k-1}]}\nnr{\fM_3(\tau,\na)P_{\neq}\udl{\th}(\tau)}_{L^2_{z,v}}e^{\f32\nu\tau}\ls\ueps e^{-10D_1K}.
    \ea\eeq
    Combining \eqref{RRR111}, \eqref{beth211}, \eqref{beth222}, and \eqref{beth233} we obtain \eqref{341}.
\end{proof}

The reduction above leaves a scalar Volterra equation for the density on each critical interval. The next proposition compares this density with the free trace corrected by the Maxwellian resolvent, and gives a weighted $L^2$ bound with the explicit factor $k^{-1/2}$. The comparison controls the error in the leading transfer, while the upper bound is used to propagate weighted norms along the chain.

\begin{prop}[Density estimates on a critical interval]\label{prop:critical-density}
	For $1\leq k\leq K$, we have
	\beq\label{rhoklower}\ba&\nnr{\wh{\rho}(t,k)\emnt-\int_{t_k}^t\wh{\th}(t_k,k,k\tentau)S(\tau-t_k,k)\rr{\delta_0(t-\tau)-\wh{R}(t-\tau,k)}\emnt d\tau}_{L^1[t_k,t_{k-1}]}\\
	&\ls\f{\ueps}{k^3}\nnr{\wh{\th}(t_k,k,\xi)}_{L^1_\xi}+e^{-9D_1K},\ea\eeq
	\beq\label{rhokupper}\ba
	&\nnr{\sF\rr{\fM(\tau,\na)\rho}(\tau,k)e^{\rm\nu^\f13\tau}}_{L^2[t_k,t_{k-1}]}\\
	&\leq k^{-\f12}\nnr{\sF\rr{\fM(t_k,\na)\th}(t_k,k,\xi)}_{L^2_\xi}e^{\nu t_k+\rm\nu^\f13t_{k-1}}\rr{1+\f{C_0}{k^2}}+e^{-8.5D_1K},
	\ea\eeq
	where $\wh{R}(t,k)$ is defined in \eqref{definereso}, and $C_0$ is a universal constant.
\end{prop}
\begin{proof}
We invert the Maxwellian Volterra term using the resolvent. The background zero-mode term is a small perturbation. An $L^1$ estimate gives the comparison formula, while a weighted $L^2$ estimate retains the trace factor $k^{-1/2}$ needed for the one-step upper bound.

	By \eqref{rhoktk} and \eqref{definereso} we have
	\[\ba
	&\wh{\rho}(t,k)=\wh{H}(t,k)-\int_{t_k}^t\wh{H}(\tau,k)\wh{R}(t-\tau,k)d\tau,\quad \wh{H}(t,k)=\wh{\th}(t_k,k,k\tent)S(t-t_k,k)\\
	&-\int_{t_k}^t\wh{\rho}(\tau,k)\f{1-e^{-\nu(t-\tau)}}{\nu}\wh{\uh}(\tau,0,k\f{1-e^{-\nu(t-\tau)}}{\nu})S(t-\tau,k)d\tau+\wh{\sR_1}(t,k).
	\ea\]
	Denote the first two terms in $\wh{H}(t,k)$ by $\wh{\aleph'_1}$, $\wh{\aleph'_2}$. We have
	\[\ba&\nnr{\wh{\aleph'_1}(t,k)\emnt}_{L^1[t_k,t_{k-1}]}\leq\int_{t_k}^{t_{k-1}}\nr{\wh{\th}(t_k,k,k\tent)}\emnt dt\leq k^{-1}\nnr{\wh{\th}(t_k,k,\xi)}_{L^1_\xi},\\
	&\nnr{\wh{\sR_1}(t,k)\emnt}_{L^1[t_k,t_{k-1}]}\ls\nnr{\wh{\<\na\>\sR_1}(t,k)}_{L^2[t_k,t_{k-1}]}\rr{\int_{t_k}^{t_{k-1}}\<\tent\>^{-2}\emnt dt}^\f12\ls e^{-10D_1K}.\\
	&\nnr{\wh{\aleph'_2}(t,k)\emnt}_{L^1[t_k,t_{k-1}]}\ls\int_{[t_k,t_{k-1}]^2}\nr{\wh{\rho}(\tau,k)}\emntau\f{1-e^{-\nu(t-\tau)}}{\nu}\nr{\wh{\uh}(\tau,0,k\f{1-e^{-\nu(t-\tau)}}{\nu})}\vv{1}_{\tau\le t}e^{-\nu(t-\tau)}d\tau dt\\
	&\leq k^{-2}\nnr{\wh{\rho}(t,k)\emnt}_{L^1[t_k,t_{k-1}]}\sup_{\tau\in[t_k,t_{k-1}]}\nnr{\xi\wh{\udl{h}}(\tau,0,\xi)}_{L^1_\xi}\ls\f{\ueps}{k^2}\nnr{\wh{\rho}(t,k)\emnt}_{L^1[t_k,t_{k-1}]}.
	\ea\]
	The third inequality uses \eqref{Snu13}, Fubini's theorem and the coordinate shift $\xi=k\f{1-e^{-\nu(t-\tau)}}{\nu}$. By \eqref{resolvent} we have 
	\[\ba
	&\nnr{\int_{t_k}^{t}\wh{\aleph'_1}(\tau,k)\wh{R}(t-\tau,k)e^{-\nu t}d\tau}_{L^1[t_k,t_{k-1}]}\leq\nnr{\wh{\aleph'_1}(\tau,k)\emntau}_{L^1[t_k,t_{k-1}]}\int_0^\infty\nr{\wh{R}(t,k)}dt\\
	&\ls k^{-2}\nnr{\wh{\aleph'_1}(\tau,k)\emntau}_{L^1[t_k,t_{k-1}]}.
	\ea\]
	Consequently,
	\[\ba
	&\nnr{\wh{\rho}(t,k)\emnt-\int_{t_k}^t\wh{\aleph'_1}(\tau,k)\rr{\delta_0(t-\tau)-\wh{R}(t-\tau,k)}\emnt d\tau}_{L^1[t_k,t_{k-1}]}\\
	&=\nnr{\int_{t_k}^t\left\{\wh{\aleph'_2}(\tau,k)+\wh{\sR_1}(\tau,k)\right\}\rr{\delta_0(t-\tau)-\wh{R}(t-\tau,k)}\emnt d\tau}_{L^1[t_k,t_{k-1}]}\\
	&\ls\f{\ueps}{k^2}\nnr{\wh{\rho}(t,k)\emnt}_{L^1[t_k,t_{k-1}]}+e^{-10D_1K}\leq\f{\ueps}{k^2}\left\{\nnr{\wh{\rho}(t,k)\emnt-\int_{t_k}^t\wh{\aleph'_1}(\tau,k)\right.\right.\\
&\left.\left.\times\rr{\delta_0(t-\tau)-\wh{R}(t-\tau,k)}\emnt d\tau}_{L^1[t_k,t_{k-1}]}+k^{-1}\nnr{\wh{\th}(t_k,k,\xi)}_{L^1_\xi}\rr{1+\nnr{\wh{R}(t,k)}_{L^1_t}}\right\}+e^{-10D_1K}.
	\ea\]
	We obtain \eqref{rhoklower}. We next prove \eqref{rhokupper}. Let $C_0=2\sup_kk^2\nnr{e^{3\rm|kt|+10}\wh{R}(t,k)}_{L^1_t}$. We first estimate $\aleph_1'$. As in \eqref{egx} we have
	\[\ba
	&\nnr{\sF\rr{\fM(\tau,\na)\aleph_1'}(\tau,k)e^{\rm\nu^\f13\tau}}_{L^2[t_k,t_{k-1}]}^2\leq\int_{t_k}^{t_{k-1}}\nr{\sF\rr{\fM(t_k,\na)\th}(t_k,k,k\tent)}^2\emnt dt e^{\nu t_k+2\rm\nu t_k+2\rm\nu^\f13t_{k-1}}\\
	&\leq\nnr{\sF\rr{\fM(t_k,\na)\th}(t_k,k,\xi)}_{L^2_\xi}^2e^{2\nu t_k+2\rm\nu^\f13t_{k-1}}k^{-1}.	
	\ea\]
	As in \eqref{egv} we have
	\[\ba
	&\nr{\sF\rr{\fM(t,\na)\aleph_2'}(t,k)e^{\rm\nu^\f13t}}\ls\int_{t_k}^t\nr{\sF\rr{\fM(\tau,\na)\rho}(\tau,k)}e^{\rm\nu^\f13\tau+\rm\nu t_k}\f{1}k\nr{\sF\rr{\fM(\tau,0,\emntau\na_v)\pa_vP_0\uh}\right.\\
	&\left.\quad(\tau,0,k\f{1-e^{-\nu(t-\tau)} }{\nu}) }e^{-\nu(t-\tau)}d\tau\\
	&\ls\sup_{\tau\in[t_k,t]}\nnr{v^\alpha\<\na_v\>^3\fM(\tau,0,\emntau\na_v)P_0\uh(\tau)}_{\lta L^2_v}\f1k\int_{t_k}^t\nr{\sF\rr{\fM(\tau,\na)\rho}(\tau,k)}e^{\rm\nu^\f13\tau}\\
    &\times \<k\f{1-e^{-\nu(t-\tau)} }{\nu}\>^{-2}e^{-\nu(t-\tau)}d\tau
	\ls\f\ueps {k}\int_{t_k}^t\nr{\sF\rr{\fM(\tau,\na)\rho}(\tau,k)}e^{\rm\nu^\f13\tau}\<k\f{1-e^{-\nu(t-\tau)} }{\nu}\>^{-2}e^{-\nu(t-\tau)}d\tau.
	\ea\]
	By Schur's test we have
	\[\nnr{\sF\rr{\fM(\tau,\na)\aleph_2'}(\tau,k)e^{\rm\nu^\f13\tau}}_{L^2[t_k,t_{k-1}]}\ls\f{\ueps}{k^2}\nnr{\sF\rr{\fM(\tau,\na)\rho}(\tau,k)e^{\rm\nu^\f13\tau}}_{L^2[t_k,t_{k-1}]}.\]
	Consequently,
	\[\ba&
	\nnr{\sF\rr{\fM(\tau,\na)\rho}(\tau,k)e^{\rm\nu^\f13\tau}}_{L^2[t_k,t_{k-1}]}\\
	&\leq\rr{1+\int_0^\infty e^{3\rm|kt|+10}\nr{\wh{R}(t,k)}dt}\nnr{\sF\rr{\fM(\tau,\na)H}(\tau,k)e^{\rm\nu^\f13\tau}}_{L^2[t_k,t_{k-1}]}\\&
	\leq(1+\f{C_0}{2k^2})k^{-\f12}\nnr{\sF\rr{\fM(t_k,\na)\th}(t_k,k,\xi)}_{L^2_\xi}e^{\rm\nu^\f13t_{k-1}+\nu t_k}\\
	&+\f{\ueps^{0.9}}{k^2}\nnr{\sF\rr{\fM(\tau,\na)\rho}(\tau,k)e^{\rm\nu^\f13\tau}}_{L^2[t_k,t_{k-1}]}+e^{-8.9D_1K}.
	\ea\]
	Thus we have \eqref{rhokupper}.
\end{proof}

We now substitute the density estimates into the resonant transfer formula \eqref{hktk}. Near resonance, freezing the interaction coefficient produces a convolution with $\wh{\vf}$ and the factor $\ueps A/k^3$, together with the collisional attenuation $e^{-\nu A^3/(3k^3)}$. Proposition~\ref{4.3} quantifies this approximation in $L^1$ and gives the corresponding weighted $L^2$ upper bound. The relative errors are kept explicit because they must remain controlled after iteration.

\begin{prop}[One-step echo transfer]\label{4.3}
	There exists $D_3\geq10^3$ such that for $\llcorner\nu^\f13A\lrcorner+1\leq k\leq K$, we have
	\beq\label{hktke}\ba
	&\nnr{\wh{\th}(t_{k-1},k-1,\xi)+\f{\ueps A}{k^3}e^{-\f{\nu A^3}{3k^3}}\int\wh{\th}(t_k,k,\zeta)\wh{\vf}(\xi-\zeta)d\zeta}_{L^1_\xi}\\
	&\leq\f{\ueps A}{k^3}e^{-\f{\nu A^3}{3k^3}}D_3\rr{\f{1}{k^2}+\f{K}{A}+\f{\nu A^3}{k^4}}\nnr{\wh{\th}(t_k,k,\zeta)}_{L^1_\zeta}+e^{-8D_1K}.
	\ea\eeq
	For $2\leq k\leq K$, we have
	\beq\label{bktku}\ba&\nnr{\sF\rr{\fM(t_{k-1},\na)\th}(t_{k-1},k-1,\xi)}_{L^2_\xi}e^{\rm\nu^\f13t_{k-1}}\leq\f{\ueps A}{k^3}e^{-\f{\nu A^3}{3k^3}+2\rm(t_{k-1}-t_{k})}\\&\times\left\{1+D_3\rr{\f{1}{k^2}+\f{K}{A}}\right\}\nnr{\sF\rr{\fM(t_{k},\na)\th}(t_{k},k,\xi)}_{L^2_\xi}e^{\rm\nu^\f13t_{k}}+e^{-8D_1K}.\ea
	\eeq
\end{prop}
\begin{proof}
We first estimate the semigroup on the resonant support. We then compare the resonant integral with the frozen convolution, separating the error in its coefficient from the error in the density approximation. Finally, we use the weighted density bound to obtain the $L^2$ transfer estimate.

	\underline{Step 1: Some auxiliary inequalities.} We first state and prove some auxiliary inequalities. For $|\xi-A|\leq40D_2K+1$, $2\leq k\leq K$, we have
	\beq\label{auxi1}-\log S(t_{k-1},t_k,k-1,\xi)\leq\f{5\nu A^3}{k^4},\quad-\log S(t_{k-1}-t_k,k)\leq\f{2\nu A^3}{k^4}.\eeq
	For $|\xi|\leq1$, $T\leq t\leq \nu^{-\f12}$, we have
	\beq\label{auxi2}\nr{\log S(t,0,-1,\xi)+\f13\nu t^3}\leq2\nu t^2.\eeq
	We first prove \eqref{auxi1}.
	\[\ba
	&-\log S(t_{k-1},t_k,k-1,\xi)=\nu\int_{t_k}^{t_{k-1}}e^{3\nu s}\nr{\xi-(k-1)\f{1-e^{-\nu s}}{\nu} }^2e^{-\nu s}ds\\
	\leq&1.1\nu\int_{\f{A}{k+\f12}}^{\f{A}{k-\f12}}\nr{\xi-(k-1)s}^2ds	\leq1.1\nu\int_{\f{A}{k+\f12}}^{\f{A}{k-\f12}}\nr{\f{\f32A}{k+\f12}+40D_2K+1}^2ds\leq\f{5\nu A^3}{k^4}.
	\ea\]
	The second inequality uses the change of variables $s'=\f{1-e^{-\nu s}}{\nu}$. By \eqref{St} we have
	\[-\log S(t_{k-1}-t_k,k)\leq\nu k^2(t_{k-1}-t_k)^3\leq\f{2\nu A^3}{k^4}.\]
	From \eqref{S1} we obtain \eqref{auxi2}, where we use the following Taylor expansion on $\nu t\leq1$
	\[\f{2}{\nu^2}\f{\ent-1}{\ent+1}-\f{t}{\nu}=-\f{1}{12}\nu t^3+\mathcal{O}(\nu^2t^4),\quad \f{\ent-1}{\nu(\ent+1)}=\f12 t-\f{1}{24}\nu^2t^3+\mathcal{O}(\nu^3t^4). \]
	
	\underline{Step 2: Proof of \eqref{hktke}.} Equation~\eqref{hktk} gives
	\[\wh{\th}(t_{k-1},k-1,\xi)=\wh{\gimel_1}+\wh{\gimel_2}-\f{\ueps A}{k^3}e^{-\f{\nu A^3}{3k^3}}\int\wh{\th}(t_k,k,\zeta)\wh{\vf}(\xi-\zeta)d\zeta+\wh{\sR_2},\]
	where
	\[\ba
	\wh{\gimel_1}=&-\int_T^{t_{k-1}}\vv{1}_{|k\tentau-A|\leq40D_2K}\wh{\th}(t_k,k,k\tentau)S(\tau-t_k,k)\f{\bar{\eta}(\tau,k-1,\xi)}{k}\wh{\udl{\th}}(\tau,-1,\xi-k\tentau)\\
	&\times S(t_{k-1},\tau,k-1,\xi)d\tau+\f{\ueps A}{k^3}e^{-\f{\nu A^3}{3k^3}}\int\wh{\th}(t_k,k,\zeta)\wh{\vf}(\xi-\zeta)d\zeta,\\
	\wh{\gimel_2}=&-\int_T^{t_{k-1}}\vv{1}_{|k\tentau-A|\leq40D_2K}\left\{\wh{\rho}(\tau,k)-\wh{\th}(t_k,k,k\tentau)S(\tau-t_k,k)\right\}\f{\bar{\eta}(\tau,k-1,\xi)}{k}\\
	&\times\wh{\udl{\th}}(\tau,-1,\xi-k\tentau)S(t_{k-1},\tau,k-1,\xi)d\tau.
	\ea\]
	
	We first estimate $\wh{\gimel_1}$. Note that $\wh{\gimel_1}=\wh{\gimel_{11}}+\wh{\gimel_{12}}$, where
	\[\ba\wh{\gimel_{11}}=&-\ueps\int_T^{t_{k-1}}\vv{1}_{|k\tentau-A|\leq40D_2K}\wh{\th}(t_k,k,k\tentau)\wh{\vf}(\xi-k\tentau)\f{\emntau}{k}\left\{S(\tau-t_k,k)\right.\\
	&\left.\bar{\eta}(\tau,k-1,\xi)S(\tau,0,-1,\xi-k\tentau)S(t_{k-1},\tau,k-1,\xi)\entau-\f Ake^{-\f{\nu A^3}{3k^3}}\right\}d\tau,\\
	\wh{\gimel_{11}}=&\f{\ueps A}{k^3}e^{-\f{\nu A^3}{3k^3}}\int\wh{\th}(t_k,k,\zeta)\vv{1}_{|\zeta-A|\geq40D_2K}\wh{\vf}(\xi-\zeta)d\zeta,
	\ea\]
	where we use the change of variables $\zeta=k\tentau$ and \eqref{963}.	We first estimate $\wh{\gimel_{11}}$. By \eqref{auxi1} and \eqref{auxi2} we obtain the following bounds when $|k\tentau-A|\leq40D_2K$, $|\xi-k\tentau|\leq1$, $\llcorner\nu^\f13A\lrcorner+1\leq k\leq K$, we have
	\[\ba
	&\nr{\emntau\bar{\eta}(\tau,k-1,\xi)-\f Ak}\leq\nr{\xi-k\tentau}+\f1k\nr{k\tentau-A}\leq1+\f{40D_2K}{k}\leq\f Ak\f{41D_2K}{A}.\\
	&\nr{e^{\f{\nu A^3}{3k^3}}S(\tau,0,-1,\xi-k\tentau)e^{2\nu\tau}S(\tau-t_k,k)S(t_{k-1},\tau,k-1,\xi)-1}\leq\f{50\nu A^3}{k^4}.
	\ea\]
	Note that $\llcorner\nu^\f13A\lrcorner+1\leq k\leq K$ implies $\f{\nu A^3}{k^4}\ll1$. To prove the second inequality we use  
	\beq\label{key1}\nr{\log S(\tau,0,-1,\xi-k\tentau)+\f{\nu A^3}{3k^3}}\leq\f{\nu}{3}\nr{\tau^3-\f{A^3}{k^3}}+2\nu\tau^2\leq10\nu\f{A^2}{k^2}\f{40D_2K}{k},\eeq
	which can be deduced from \eqref{auxi2}. Thus
	\beq\label{key}\ba
	&\nr{S(\tau-t_k,k)\bar{\eta}(\tau,k-1,\xi)S(\tau,0,-1,\xi-k\tentau)S(t_{k-1},\tau,k-1,\xi)\entau-\f Ake^{-\f{\nu A^3}{3k^3}}}\\
	&\leq\f Ake^{-\f{\nu A^3}{3k^3}}50D_2\rr{\f{K}{A}+\f{\nu A^3}{k^4}}.
	\ea\eeq
	Thus by the coordinate shift $\zeta=k\tentau$ we obtain
	\beq\label{gimel11}\nnr{\wh{\gimel_{11}}}_{L^1_\xi}\leq\f{\ueps A}{k^3}e^{-\f{\nu A^3}{3k^3}}50D_2\rr{\f{K}{A}+\f{\nu A^3}{k^4}}\nnr{\wh{\th}(t_k,k,\zeta)}_{L^1_\zeta}\nnr{\wh{\vf}}_{L^1_\xi}.\eeq
	We recall that $\wh{\vf}(\xi)\geq0$, $\nnr{\wh{\vf}}_{L^1_\xi}=1$. We next estimate $\wh{\gimel_{12}}$.
	\beq\label{gimel12}\ba\nnr{\wh{\gimel_{12}}}_{L^1_\xi}\leq&\f{\ueps A}{k^3}e^{-\f{\nu A^3}{3k^3}}\nnr{\wh{\th}(t_k,k,\zeta)\vv{1}_{|\zeta-A|\geq40D_2K}}_{L^1_\zeta}\nnr{\wh{\vf}}_{L^1_\xi}\ls\f{\ueps A}{k^3}e^{-\f{\nu A^3}{3k^3}}\\
	&\times\nnr{\sF\rr{e^{\tilde{r}_1\<\na_z-K,\na_v-A\>}\<\na\>\th}(t_k,k,\xi)}_{L^2_\xi}e^{-40D_1K}\ls e^{-20D_1K}.\ea\eeq
	The last inequality uses \eqref{hnqenergy}. 
	
	We next estimate $\wh{\gimel_2}$. By \eqref{963} we have
	\[\ba\wh{\gimel_2}=&-\ueps\int_T^{t_{k-1}}\vv{1}_{|k\tentau-A|\leq40D_2K}\left\{\wh{\rho}(\tau,k)-\wh{\th}(t_k,k,k\tentau)S(\tau-t_k,k)\right\}\wh{\vf}(\xi-k\tentau)\f{\emntau}{k}\\&\times \bar{\eta}(\tau,k-1,\xi)S(\tau,0,-1,\xi-k\tentau)S(t_{k-1},\tau,k-1,\xi)\entau d\tau.\ea\]
	As in \eqref{key} we have
	\[\nr{\bar{\eta}(\tau,k-1,\xi)S(\tau,0,-1,\xi-k\tentau)S(t_{k-1},\tau,k-1,\xi)\entau}\ls\f{ A}{k}e^{-\f{\nu A^3}{3k^3}}.\]
	Consequently,
	\beq\label{gimel2}\ba
	\nnr{\wh{\gimel_{2}}}_{L^1_\xi}\ls&\f{\ueps A}{k^2}e^{-\f{\nu A^3}{3k^3}}\nnr{\wh{\rho}(\tau,k)\emntau-\wh{\th}(t_k,k,k\tentau)S(\tau-t_k,k)\emntau}_{L^1[t_k,t_{k-1}]}\\
	\ls&\f1{k^2}\f{\ueps A}{k^3}e^{-\f{\nu A^3}{3k^3}}\nnr{\wh{\th}(t_k,k,\xi)}_{L^1_\xi}+e^{-8.9D_1K},
	\ea\eeq
	where we use 
	\[\nnr{\wh{\rho}(\tau,k)\emntau-\wh{\th}(t_k,k,k\tentau)S(\tau-t_k,k)\emntau}_{L^1[t_k,t_{k-1}]}\ls\f{1}{k^3}\nnr{\wh{\th}(t_k,k,\xi)}_{L^1_\xi}+e^{-9D_1K},\]
	which can be deduced from \eqref{rhoklower} and \eqref{resolvent}. Combining \eqref{gimel11}, \eqref{gimel12}, \eqref{gimel2}, and \eqref{341} we obtain \eqref{hktke}.
	
	\underline{Step 3: Proof of \eqref{bktku}.} By \eqref{hktk} and \eqref{963} we have 
	\[\ba\wh{\th}(t_{k-1},k-1,\xi)=&-\ueps\int_T^{t_{k-1}}\vv{1}_{|k\tentau-A|\leq40D_2K}\wh{\rho}(\tau,k)\wh{\vf}(\xi-k\tentau)\f{e^{-2\nu\tau}}{k}\bar{\eta}(\tau,k-1,\xi)\\
	&\times S(\tau,0,-1,\xi-k\tentau)S(t_{k-1},\tau,k-1,\xi)e^{2\nu\tau} d\tau+\wh{\sR_2}.\ea\]
	For $|k\tentau-A|\leq40D_2K$, $|\xi-k\tentau|\leq1$, $2\leq k\leq K$, we have
	\beq\label{mtp}\ba
	&e^{3\nu\tau}\fM(t_{k-1},k-1,\xi)S(t_{k-1},\tau,k-1,\xi)\leq e^{3\nu t_{k-1}+2\rm\nu t_{k-1}}\fM(\tau,k-1,\xi)\\
	\leq& e^{4\nu t_{k-1}}\fM(\tau,k,k\tentau)\fM(\tau,-1,\xi-k\tentau)\leq e^{4\nu t_{k-1}+\nu t_{k-1}^2}\fM(\tau,k,k\tentau).
	\ea\eeq
	The first inequality uses \eqref{MSk} and \eqref{Snu13}. The second inequality uses \eqref{M1}. The third inequality uses \eqref{Mxj}. (Note that $|\xi-k\tentau|\leq1$.) By \eqref{key1} we have
	\[S(\tau,0,-1,\xi-k\tentau)e^{\f{\nu A^3}{3k^3}}\leq e^{400D_2K\f{\nu A^2}{k^3}}\leq1+500D_2\f{\nu A^2K}{k^3}.\]
	Note that by \eqref{KAT} we have $\nu A^2K\sim\ueps^\f37\ll1$. These together with $\nr{\emntau\bar{\eta}(\tau,k-1,\xi)}\leq\f Ak\rr{1+\f{41D_2K}{A}}$ yield 
	\[\ba
	&\nnr{\sF\rr{\fM(t_{k-1},\na)\th}(t_{k-1},k-1,\xi)}_{L^2_\xi}e^{\rm\nu^\f13t_{k-1}}\leq\f{\ueps A}{k^2}e^{-\f{\nu A^3}{3k^3}+\rm\nu^\f13(t_{k-1}-t_k)}\left\{1+10^3D_2\rr{\f KA+\f{\nu A^2K}{k^3}}\right\}e^{-\f32\nu t_{k}}\\
	&\times \nnr{\int_{t_k}^{t_{k-1}}\nr{\sF\rr{\fM(\tau,\na)\rho}(\tau,k)}e^{\rm\nu^\f13\tau}\wh{\vf}(\xi-k\tentau)e^{-\f{\nu\tau}{2}}d\tau }_{L^2_\xi}+e^{-8.5D_1K}
	\leq\f{\ueps A}{k^{2.5}}e^{-\f{\nu A^3}{3k^3}}\\
	&\times e^{\rm\nu^\f13(t_{k-1}-t_k)}\left\{1+10^3D_2\rr{\f KA+\f{\nu A^2K}{k^3}}\right\}e^{-\f32\nu t_{k}}\nnr{\sF\rr{\fM(\tau,\na)\rho}(\tau,k)e^{\rm\nu^\f13\tau}}_{L^2[t_k,t_{k-1}]}\nnr{\wh{\vf}}_{L^1_\xi}\\
	&+e^{-8.5D_1K}
	\leq\f{\ueps A}{k^3}e^{-\f{\nu A^3}{3k^3}+2\rm\nu^\f13(t_{k-1}-t_k)}\left\{1+D_3\rr{\f KA+\f1{k^2}}\right\}\nnr{\sF\rr{\fM(t_{k},\na)\th}(t_{k},k,\xi)}_{L^2_\xi}e^{\rm\nu^\f13t_k}+e^{-8.4D_1K}.
	\ea\]
	The second inequality uses Schur's test. The third inequality uses \eqref{rhokupper}.
\end{proof}

The one-step estimate can now be iterated from the prescribed packet at mode $K$. Proposition~\ref{4.4} records both the leading convolution profile and the accumulated error, as well as the product bound for the weighted norm. The lower-bound range is chosen so that the sum of the relative errors remains small.

\begin{prop}[Iteration along the echo chain]\label{4.4}
	For $\llcorner\rr{\nu A^3D_3}^\f13\lrcorner+1\leq k\leq K$, we have
	\beq\label{hktkl}\ba
	&\nnr{\wh{\th}(t_k,k,\xi)-(-1)^{K-k}\prod_{i=k+1}^{K}\f{\ueps A}{i^3}e^{-\f{\nu A^3}{3i^3}}\wh{\th}(T,K,\cdot)*\rr{\wh{\vf}^{*(K-k)}}}_{L^1_\xi}\leq2D_3\left\{\sum_{i=k+1}^{K}\rr{\f1{i^2}+\f KA+\f{\nu A^3}{i^4}}\right\}\\
	&\times\left\{\prod_{i=k+1}^{K}\f{\ueps A}{i^3}e^{-\f{\nu A^3}{3i^3}}\right\}\nnr{\wh{\th}(T,K,\xi)}_{L^1_\xi}+\left\{\prod_{i=k+1}^{K}\f{3\ueps A}{i^3}\right\}e^{-8D_1K},
	\ea\eeq
	where $\vf^{*j}:=\underbrace{\vf*\vf*\cdots*\vf}_{j}$, $\vf^{*0}=\delta$.
	For $1\leq k\leq K$, we have
	\beq\label{hktku}\ba
	&\nnr{\sF\rr{\fM(t_{k},\na)\th}(t_{k},k,\xi)}_{L^2_\xi}e^{\rm\nu^\f13t_k}\leq\left\{\prod_{i=k+1}^{K}\f{\ueps A}{i^3}e^{-\f{\nu A^3}{3i^3}+2\rm\nu^\f13(t_{i-1}-t_{i})}\rr{1+D_3\rr{\f1{i^2}+\f KA}}\right\}\\
	&\times\nnr{\sF\rr{\fM(T,\na)\th}(T,K,\xi)}_{L^2_\xi}e^{\rm\nu^\f13T}+\left\{\prod_{i=k+1}^{K}\f{3\ueps A}{i^3}\rr{1+\f{D_3}{i^2}}\right\}e^{-8D_1K}.
	\ea\eeq
\end{prop}
\begin{proof}
Both estimates follow by descending induction on the spatial mode. For the approximation in $L^1$, convolution with $\wh{\vf}$ does not increase the error, and the one-step relative errors accumulate additively. For the weighted upper bound, we iterate \eqref{bktku} and retain the corresponding multiplicative factors.

	We first prove \eqref{hktkl} by descending induction. The case $k'=K$ is immediate. Suppose the estimate holds for $k\leq k'\leq K$; we prove it for $k'=k-1$. By the inductive hypothesis for $k'=k$ we have
	\beq\label{key3}\ba
	&\nnr{\wh{\th}(t_k,k,\xi)}_{L^1_\xi}\leq\left\{1+2D_3\sum_{i=k+1}^{K}\rr{\f1{i^2}+\f KA+\f{\nu A^3}{i^4}}\right\}\\
	&\quad \left\{\prod_{i=k+1}^{K}\f{\ueps A}{i^3}e^{-\f{\nu A^3}{3i^3}}\right\}\nnr{\wh{\th}(T,K,\xi)}_{L^1_\xi}+\left\{\prod_{i=k+1}^{K}\f{3\ueps A}{i^3}\right\}e^{-8D_1K}.
	\ea\eeq
	By \eqref{hktke} we have
	\[\ba
	&\nnr{\wh{\th}(t_{k-1},k-1,\xi)-(-1)^{K-k+1}\prod_{i=k}^{K}\f{\ueps A}{i^3}e^{-\f{\nu A^3}{3i^3}}\wh{\th}(T,K,\cdot)*\rr{\wh{\vf}^{*(K-k+1)}}}_{L^1_\xi}\leq\f{\ueps A}{k^3}e^{-\f{\nu A^3}{3k^3}}\nnr{\wh{\th}(t_k,k,\xi)-(-1)^{K-k}\right.\\
		&\left.\times\prod_{i=k+1}^{K}\f{\ueps A}{i^3}e^{-\f{\nu A^3}{3i^3}}\wh{\th}(T,K,\cdot)*\rr{\wh{\vf}^{*(K-k)}}}_{L^1_\xi}+\f{\ueps A}{k^3}e^{-\f{\nu A^3}{3k^3}}D_3\rr{\f1{k^2}+\f KA+\f{\nu A^3}{k^4}}\nnr{\wh{\th}(t_k,k,\xi)}_{L^1_\xi}+e^{-8D_1K}\\
	\leq&\prod_{i=k}^{K}\f{\ueps A}{i^3}e^{-\f{\nu A^3}{3i^3}}\nnr{\wh{\th}(T,K,\xi)}_{L^1_\xi}\left\{2D_3\sum_{i=k+1}^{K}\rr{\f1{i^2}+\f KA+\f{\nu A^3}{i^4}}+D_3\rr{\f1{k^2}+\f KA+\f{\nu A^3}{k^4}}\left[1+2D_3\right.\right.\\
	&\left.\left.\times\sum_{i=k+1}^{K}\rr{\f1{i^2}+\f KA+\f{\nu A^3}{i^4}}\right]\right\}\\
	&+\f{\ueps A}{k^3}e^{-\f{\nu A^3}{3k^3}}\left\{1+D_3\rr{\f1{k^2}+\f KA+\f{\nu A^3}{k^4}}\right\}\left\{\prod_{i=k+1}^{K}\f{3\ueps A}{i^3}\right\}e^{-8D_1K}+e^{-8D_1K}\\
	\leq&\prod_{i=k}^{K}\f{\ueps A}{i^3}e^{-\f{\nu A^3}{3i^3}}\nnr{\wh{\th}(T,K,\xi)}_{L^1_\xi}2D_3\sum_{i=k}^{K}\rr{\f1{i^2}+\f KA+\f{\nu A^3}{i^4}}+\left\{\prod_{i=k}^{K}\f{3\ueps A}{i^3}\right\}e^{-8D_1K}.
	\ea\]
	The first inequality uses $\nnr{\wh{\vf}}_{L^1_\xi}=1$, which implies that $\nnr{\wh{g}*\wh{\vf}}_{L^1_\xi}\leq\nnr{\wh{g}}_{L^1_\xi}$. The second inequality uses the inductive hypothesis for $k'=k$ and \eqref{key3}. The third inequality uses the fact that
	\[2D_3\sum_{i=k+1}^{K}\rr{\f1{i^2}+\f KA+\f{\nu A^3}{i^4}}\leq0.9\]
	for $\llcorner\rr{\nu A^3D_3}^\f13\lrcorner+1\leq k\leq K$. This proves that \eqref{hktkl} holds for $k'=k-1$.
	
	We prove \eqref{hktku} by the same descending induction. \eqref{hktku} holds immediately for $k'=K$. We assume \eqref{hktku} holds for $k\leq k'\leq K$. We next show that \eqref{hktku} holds for $k'=k-1$. By \eqref{bktku} we have
	\[\ba
	&\nnr{\sF\rr{\fM(t_{k-1},\na)\th}(t_{k-1},k-1,\xi)}_{L^2_\xi}e^{\rm\nu^\f13t_{k-1}}\leq  \left\{\prod_{i=k}^{K}\f{\ueps A}{i^3}e^{-\f{\nu A^3}{3i^3}+2\rm\nu^\f13(t_{i-1}-t_{i})}\rr{1+D_3\rr{\f1{i^2}+\f KA}}\right\}e^{\rm\nu^\f13T}\\
	&\times\nnr{\sF\rr{\fM(T,\na)\th}(T,K,\xi)}_{L^2_\xi}
	+\f{\ueps A}{k^3}e^{-\f{\nu A^3}{3k^3}+2\rm\nu^\f13(t_{k-1}-t_{k})}\rr{1+D_3\rr{\f1{k^2}+\f KA}}\left\{\prod_{i=k+1}^{K}\f{3\ueps A}{i^3}\rr{1+\f{D_3}{i^2}}\right\}\\
	&\times e^{-8D_1K}+e^{-8D_1K}\leq\left\{\prod_{i=k}^{K}\f{\ueps A}{i^3}e^{-\f{\nu A^3}{3i^3}+2\rm\nu^\f13(t_{i-1}-t_{i})}\rr{1+D_3\rr{\f1{i^2}+\f KA}}\right\}\nnr{\sF\rr{\fM(T,\na)\th}(T,K,\xi)}_{L^2_\xi}\\
	&\times e^{\rm\nu^\f13T}+\left\{\prod_{i=k}^{K}\f{3\ueps A}{i^3}\rr{1+\f{D_3}{i^2}}\right\}e^{-8D_1K}.
	\ea\]
	Thus \eqref{hktku} holds for $k'=k-1$.
\end{proof}

The next corollary extracts the amplification bounds from Proposition~\ref{4.4}. Nonnegativity of the initial Fourier packet and of $\wh{\vf}$ prevents cancellation in the leading convolution, giving the lower bound. Summability of the weighted relative errors gives the upper bound. These are the discrete estimates used to construct $B(t)$ in Section~5.

\begin{cor}
	For $\llcorner\rr{\nu A^3D_3}^\f13\lrcorner+1\leq k\leq K$, we have
	\beq\label{hklow}\nnr{\wh{\th}(t_k,k,\xi)}_{L^1_\xi}\geq0.1\left\{\prod_{i=k+1}^{K}\f{\ueps A}{i^3}e^{-\f{\nu A^3}{3i^3}}\right\}\nnr{\wh{\th}(T,K,\xi)}_{L^1_\xi}-e^{-7D_1K}.\eeq
	For $1\leq k\leq K$, we have
	\beq\label{hkupp}\ba
	&\nnr{\sF\rr{\fM(t_{k},\na)\th}(t_{k},k,\xi)}_{L^2_\xi}e^{\rm\nu^\f13t_{k}}\\
	&\ls\left\{\prod_{i=k+1}^{K}\f{\ueps A}{i^3}e^{-\f{\nu A^3}{3i^3}+2\rm\nu^\f13(t_{i-1}-t_{i})}\right\}\nnr{\sF\rr{\fM(T,\na)\th}(T,K,\xi)}_{L^2_\xi}e^{\rm\nu^\f13T}+e^{-7D_1K}.
	\ea\eeq
\end{cor}
\begin{proof}
	Equation~\eqref{hktkl} gives
	\[\ba
	&\nnr{\wh{\th}(t_k,k,\xi)}_{L^1_\xi}\geq\nnr{\prod_{i=k+1}^{K}\f{\ueps A}{i^3}e^{-\f{\nu A^3}{3i^3}}\wh{\th}(T,K,\cdot)*\rr{\wh{\vf}^{*(K-k)}}}_{L^1_\xi}-\nnr{\wh{\th}(t_k,k,\xi)-(-1)^{K-k}\prod_{i=k+1}^{K}\f{\ueps A}{i^3}e^{-\f{\nu A^3}{3i^3}}\right.\\
	&\left.\wh{\th}(T,K,\cdot)*\rr{\wh{\vf}^{*(K-k)}}}_{L^1_\xi}
	\geq\left\{\prod_{i=k+1}^{K}\f{\ueps A}{i^3}e^{-\f{\nu A^3}{3i^3}}\right\}\nnr{\wh{\th}(T,K,\xi)}_{L^1_\xi}\left\{1-2D_3\sum_{i=k+1}^{K}\rr{\f1{i^2}+\f KA+\f{\nu A^3}{i^4}}\right\}\\
	&-e^{-7D_1K}\geq0.1\left\{\prod_{i=k+1}^{K}\f{\ueps A}{i^3}e^{-\f{\nu A^3}{3i^3}}\right\}\nnr{\wh{\th}(T,K,\xi)}_{L^1_\xi}-e^{-7D_1K}.
	\ea\]
	The second inequality uses $\wh{\th}(T,K,\xi)\geq0$, $\wh{\vf}(\xi)\geq0$, which implies that $\nnr{\wh{\th}(T,K,\cdot)*\rr{\wh{\vf}^{*(K-k)}}}_{L^1_\xi}=\nnr{\wh{\th}(T,K,\xi)}_{L^1_\xi}\nnr{\wh{\vf}(\xi)}_{L^1_\xi}^{(K-k)}=\nnr{\wh{\th}(T,K,\xi)}_{L^1_\xi}$. The third inequality uses the fact that
	\[2D_3\sum_{i=k+1}^{K}\rr{\f1{i^2}+\f KA+\f{\nu A^3}{i^4}}\leq0.9\]
	for $\llcorner\rr{\nu A^3D_3}^\f13\lrcorner+1\leq k\leq K$. Note that by Stirling's formula we have $\prod_{i=k+1}^{K}\f{3\ueps A}{i^3}\leq e^{10K}$. By \eqref{hktku}, to show \eqref{hkupp}, we only need to prove
	\[\prod_{i=k+1}^{K}\rr{1+D_3\rr{\f1{i^2}+\f KA}}\leq\left\{\prod_{i=k+1}^{K}\rr{1+\f{D_3}{i^2}}\right\}\rr{1+\f{D_3K}{A}}^K\ls\rr{1+\f{D_3K}{A}}^K\ls1.\]
	The last step follows from the fact that $\f{K^2}{A}\ll1$.
\end{proof}

To use the discrete estimates in the nonlinear argument, we need bounds for the density and profile at arbitrary times. The next proposition sums the density estimates over the critical intervals and then applies Duhamel's formula to the nonzero modes and the Fokker--Planck semigroup to the zero mode. Thus the entire homogeneous response is controlled by the largest critical-time amplitude reached so far.

\begin{prop}[Global upper bounds for the homogeneous response]\label{prop:homogeneous-upper}
	For $t\geq T$, we have
	\beq\label{h11upper}\ba
	&\nnr{\fM_4(\tau,\na)\rho(\tau)e^{\rm\nu^\f13\tau}}_{L^2([T,t], L^2_{z})}+\left\|\fM_{2}(t,\na)\rr{e^{-|\alpha|\nu t}v^\alpha P_{\neq}\th(t)} \right\|_{\lta L^2_{z,v}}e^{\rm\nu^\f13t}\\
	&+\nnr{v^\alpha\<\na_v\>\fM(t,0,\emnt\na_v)P_{0}h(t)}_{\lta L^2_{v}}\\
	&\times e^{\nu t}
	\ls A^4K\sup_{k:T\leq t_k\leq t}\nnr{\sF\rr{\fM(t_{k},\na)\th}(t_{k},k,\xi)}_{L^2_\xi}e^{\rm\nu^\f13t_{k}}\\
	&+\left\|\fM_{2}(T,\na)\rr{e^{-|\alpha|\nu T}v^\alpha P_{\neq}\th(T)} \right\|_{\lta L^2_{z,v}}e^{\rm\nu^\f13T}+e^{-8D_1K}.
	\ea\eeq
\end{prop}
\begin{proof}
We split the density into the localized frequency region and its complement. The complement is exponentially small, while the localized part is bounded by summing \eqref{rhokupper} over critical intervals. The profile estimates then follow from the density bound by the same Duhamel and zero-mode semigroup arguments as in Proposition~\ref{3.2}.

	We first estimate $\rho(t)$.
	\[\ba
	&\nnr{\fM_4(\tau,\na)\rho(\tau)e^{\rm\nu^\f13\tau}}_{L^2([T,t], L^2_{z})}\leq\nnr{\vv{1}_{\<\na_z-K,\na_v-A\>\leq40D_2K}\fM_4(\tau,\na)\rho(\tau)e^{\rm\nu^\f13\tau}}_{L^2([T,t], L^2_{z})}\\
	&+\nnr{\vv{1}_{\<\na_z-K,\na_v-A\>\geq40D_2K}\right.\\
	&\left.\times\fM_4(\tau,\na)\rho(\tau)e^{\rm\nu^\f13\tau}}_{L^2([T,t], L^2_{z})}
	\ls A^4\sum_{k:T\leq t_k\leq t}\nnr{\sF\rr{\fM(\tau,\na)\rho}(\tau,k)e^{\rm\nu^\f13\tau}}_{L^2[t_k,t_{k-1}]}+e^{-10D_1K}\\
	&\ls A^4K\sup_{k:T\leq t_k\leq t}\nnr{\sF\rr{\fM(t_{k},\na)\th}(t_{k},k,\xi)}_{L^2_\xi}e^{\rm\nu^\f13t_{k}}+e^{-8D_1K}.
	\ea\]
	The second inequality uses \eqref{rmkk}, \eqref{rhoexsm} and $\<k,\xi\>^4\vv{1}_{\<k-K,\xi-A\>\leq40D_2K}\ls A^4$. The third inequality uses \eqref{rhokupper}. We next estimate $h(t)$. Arguing as in the proof of \eqref{neqneqh} and \eqref{p0poh} we have
	\[\ba
	&\left\|\fM_{2}(t,\na)\rr{e^{-|\alpha|\nu t}v^\alpha P_{\neq}\th(t)} \right\|_{\lta L^2_{z,v}}e^{\rm\nu^\f13t}
	+\nnr{v^\alpha\<\na_v\>\fM(t,0,\emnt\na_v)P_{0}h(t)}_{\lta L^2_{v}}\ent\\
	\ls&\left\|\fM_{2}(T,\na)\rr{e^{-|\alpha|\nu T}v^\alpha P_{\neq}\th(T)} \right\|_{\lta L^2_{z,v}}e^{\rm\nu^\f13T}+\nnr{\fM_4(\tau,\na)\rho(\tau)e^{\rm\nu^\f13\tau}}_{L^2([T,t], L^2_{z})}.
	\ea\]
	This proves \eqref{h11upper}.
\end{proof}

\section{Linear analysis: the inhomogeneous case}
This section extends the preceding analysis to a source localized near a single frequency. We consider
\[\left\{\ba&\partial_t h=\mathcal{L} h+Q, \\ &\th(T)=0\ea\right.\]
on $[T,\infty)$, with $Q(\tau)$ satisfying $\mathrm{supp} \wh{\tQ}(t)\subset\{K'\}\times[A'-1,A'+1]$ for $A'\geq0$, $K'\neq0$, $\wh{Q}(t,0,0)=0$ and
\[\int^t_T\nnr{\fM(\tau,\na)\rr{e^{-\nu|\alpha|\tau}v^\alpha \tQ(\tau)}}_{\lta L^2_{z,v}}e^{\rm\nu^\f13\tau+\f12\nu\tau} d\tau\leq Z(t).\]
Here $Z(t)$ is nondecreasing on $[T,\infty)$ and bounds the accumulated source. Our aim is to express the response in terms of $Z$, polynomial frequency factors, and any additional echo amplification. In Section~5, these bounds will be applied to a frequency partition of each nonlinear source.

We distinguish three regimes according to the resonance geometry and the strength of diffusion. In Case~1, the source has no future density resonance, so only its direct semigroup response remains. In Case~2, a future echo chain is possible and $A'\leq\nu^{-1/2}/3$; Propositions~\ref{6.1} and~\ref{6.2} provide the recurrence and its starting bound, Proposition~\ref{prop:forced-iteration} iterates the recurrence, and Proposition~\ref{prop:forced-global} recovers full-time norms. In Case~3, $A'\geq\nu^{-1/2}/3$, and collisional decay makes the echo interaction perturbative, giving Proposition~\ref{prop:forced-diffusive}.

\subsection{Case 1: $K'<0$ or $K'\f{1-e^{-\nu T}}{\nu}\geq A'+1$.}
The source support does not meet a future density trace in this regime. Equation~\eqref{rba} therefore gives $\rho(t)=0$ on $[T,\infty)$. By \eqref{duhah} we have
	\[\wh{\th}(t,k,\xi)=\int_T^t\wh{\tQ}(\tau,k,\xi)S(t,\tau,k,\xi)d\tau.\]
	Consequently,
\beq\label{Case1}\ba&\nnr{\fM_2(t,\na)\rr{e^{-\nu|\alpha|t}v^\alpha P_{\neq}\th(t)}}_{\lta L^2_{z,v}}e^{\rm\nu^\f13t}\ls\<K',A'\>^2\int^t_T\nnr{\fM(\tau,\na)\rr{e^{-\nu|\alpha|\tau}v^\alpha \tQ(\tau)}}_{\lta L^2_{z,v}}e^{\rm\nu^\f13\tau} d\tau\\
&\leq\<K',A'\>^2Z(t),\quad \rho(t)=0,\quad P_0h(t)=0.\ea\eeq

\subsection{Case 2: $K'>0$ and $K'\f{1-e^{-\nu T}}{\nu}\leq A'+1$ and $A'\leq\f{\nu^{-\f12}}{3}$.} For $0\leq k\leq K'$ we set
\beq\label{401}t'_k=\f{-1}{\nu}\log\rr{1-\f{\nu A'}{k+\f12}}.\eeq
These critical times organize the possible cascade generated by the source. We first derive the density and transfer estimates for modes below $K'$ in Proposition~\ref{6.1}. Fourier localization makes the forcing and off-resonant interactions exponentially small on these later intervals. Proposition~\ref{6.2} then estimates the direct response at mode $K'$ and supplies the initial amplitude for the recurrence.

\begin{prop}[Forced one-step transfer estimates]\label{6.1}
	Recall $C_0=2\sup_kk^2\nnr{e^{3\rm|kt|+10}\wh{R}(t,k)}_{L^1_t}$. For $1\leq k\leq K'-1$, $t\geq t'_k$ we have
	\beq\label{rho31}\ba&\nnr{\sF\rr{\fM(\tau,\na)\rho}(\tau,k)e^{\rm\nu^\f13\tau}}_{L^2[t_k',t]}\leq k^{-\f12}\nnr{\sF\rr{\fM(t_k',\na)\th}(t'_k,k,\xi)}_{L^2_\xi}e^{\rm\nu^\f13t+\nu t_k'}(1+\f{C_0}{k^2})\\
	&+\ueps^{0.8}{A'}^2e^{-5D_1(\ueps A')^{\f13}}Z(t).\ea\eeq
	For $2\leq k\leq K'-1$ we have
	\beq\label{h31}\ba
	&\nnr{\sF\rr{\fM(t_{k-1}',\na)\th}(t_{k-1}',k-1,\xi)}_{L^2_\xi}e^{\rm\nu^\f13t_{k-1}'}\leq\f{\ueps A'}{k^3}e^{-\f{\nu A'^3}{3k^3}+\f{10D_2\nu{A'}^2\rr{\ueps A'}^\f13}{k^3}}\left\{1+D_3\rr{\f{\rr{\ueps A'}^\f13+k}{A'}+\f1{k^2} }\right\}\\
	&\times e^{2\rm\nu^\f13(t'_{k-1}-t'_k)}\nnr{\sF\rr{\fM(t_k',\na)\th}(t'_k,k,\xi)}_{L^2_\xi}e^{\rm\nu^\f13t_k'}+\ueps^{0.7}{A'}^2e^{-4D_1(\ueps A')^{\f13}}Z(t'_{k-1}).
	\ea\eeq
\end{prop}
\begin{proof}
The argument follows the homogeneous reduction, but the exponentially small remainders now carry the source envelope $Z$. We first localize the response around $A'$, then control the off-resonant terms, estimate the density through its resolvent, and substitute this estimate into the transfer integral.

	We split the proof into three steps.
	
	\underline{Step 1: Localization and resonant reduction.} We prove 
	\beq\label{rho33}\ba
	\wh{\rho}(t,k)&=\wh{\th}(t'_k,k,k\tent)S(t-t'_k,k)-\int_{t'_k}^t\wh{\rho}(\tau,k)\f{1-e^{-\nu(t-\tau)}}{\nu}\wh{\sM}(k\f{1-e^{-\nu(t-\tau)}}{\nu})S(t-\tau,k)d\tau\\
	&-\int_{t'_k}^t\wh{\rho}(\tau,k)\f{1-e^{-\nu(t-\tau)}}{\nu}\wh{\uh}(\tau,0,k\f{1-e^{-\nu(t-\tau)}}{\nu})S(t-\tau,k)d\tau+\wh{\sR'_1}(t,k)
	\ea\eeq
	for $1\leq k\leq K'-1$. 
	\beq\label{h33}\ba\wh{\th}(t'_{k-1},k-1,\xi)=&-\int_T^{t'_{k-1}}\vv{1}_{|k\tentau-A'|\leq10D_2\rr{\ueps A'}^\f13}\wh{\rho}(\tau,k)\f{\bar{\eta}(\tau,k-1,\xi)}{k}\wh{\udl{\th}}(\tau,-1,\xi-k\tentau)\\
	&\times S(t'_{k-1},\tau,k-1,\xi)d\tau+\wh{\sR'_2}(k-1,\xi)\ea\eeq
	for $2\leq k\leq K'-1$. The remainders $\sR'_1$ and $\sR'_2$ satisfy
	\beq\label{R31}\ba&\nnr{\sF\rr{\fM(\tau,\na)\sR'_1}(\tau,k)e^{\rm\nu^\f13\tau}}_{L^2[t'_k,t]}\leq\ueps^{0.9}{A'}^2e^{-5D_1(\ueps A')^{\f13}}Z(t),\\
	&\nnr{\sF\rr{\fM(t'_{k-1},\na)\sR'_2}(k-1,\xi)}_{L^2_\xi}e^{\rm\nu^\f13t'_{k-1}}\leq \ueps^{0.9}{A'}^2e^{-5D_1(\ueps A')^{\f13}}Z(t'_{k-1}).
	\ea\eeq
	By \eqref{nyb}, under the condition $K'>0$ and $K'\f{1-e^{-\nu T}}{\nu}\leq A'+1$, for $\sigma\geq2$
	\beq\label{energy3}\ba
	\nnr{e^{r_1\<\na_z-K',\na_v-A'\>}\fM_{\sigma}(\tau,\na)\rho(\tau)e^{\rm\nu^\f13\tau} }_{L^2([T,t],L^2_z)}\ls e^{D_1\<\ueps\nr{K',A'}\>^\f13}\<K',A'\>^\sigma\\
	\times\int_T^t\nnr{\fM(\tau,\na)\rr{e^{-\nu|\alpha|\tau}v^\alpha \tQ(\tau)}}_{\lta L^2_{z,v}}e^{\rm\nu^\f13\tau+\f12\nu\tau} d\tau
	\ls e^{2D_1(\ueps A')^\f13}{A'}^\sigma Z(t).
	\ea\eeq
	Note that the condition $K'>0$, $K'\f{1-e^{-\nu T}}{\nu}\leq A'+1$ implies $\ueps A'\gs\ueps T\gs K^2\gg1$, $A'\gs K'T$. Consequently,
	\beq\label{remaining3}
	\nnr{\vv{1}_{\<\na_z-K',\na_v-A'\>\geq10D_2(\ueps A')^\f13}\fM_{\sigma}(\tau,\na)\rho(\tau)e^{\rm\nu^\f13\tau} }_{L^2([T,t],L^2_z)}\leq {A'}^\sigma e^{-5D_1(\ueps A')^\f13}Z(t).
	\eeq
	Recall $D_2=\f{D_1}{r_1}$, where $r_1$ is defined in Proposition~\ref{3.2}. Note that $A'\leq\f{\nu^{-\f13}}{3}$ implies that for $1\leq k\leq K'$,
	\beq\label{remark3}\nr{k\tentau-A'}\leq10D_2\rr{\ueps A'}^\f13\quad \text{implies that}\quad t'_k\leq t\leq t'_{k-1},\eeq
	which follows from $\nr{kt_k'-A'}\sim\nr{kt_{k-1}'-A'}\sim\f{A'}{k}\gs\f{A'}{K'}\gs T\gs\f{1}{\ueps^\f17\nu^\f27}\gs\rr{\ueps A'}^\f13$. The rest of the proof of \eqref{rho33} and \eqref{h33} is essentially the same as the proof of \eqref{rhoktk} and \eqref{hktk}. We prove \eqref{rho33} as an example. By Duhamel's formula (see \eqref{rba}), \eqref{rho33} holds with 
	\[\wh{\sR'_1}(t,k)=\int_{t'_k}^t\sF\rr{\cN^\dagger(P_{k+1}\th,P_{-1}\udl\th)}(\tau,k,k\tent)S(t-\tau,k)d\tau.\]
	Consequently,
	\[\ba
	&\nnr{\sF\rr{\fM(\tau,\na)\sR'_1}(\tau,k)e^{\rm\nu^\f13\tau}}_{L^2[t'_k,t]}\ls\int_{t'_k}^{t}\nnr{\fM(\tau,\na)\cN^\dagger(P_{k+1}\th,P_{-1}\udl{\th})(\tau)}_{L^2_{z,v}}e^{\rm\nu^\f13\tau+\f12\nu\tau}d\tau\\
	\ls&\nnr{\sF\rr{\fM_2(\tau,\na)\rho}(\tau,k+1)e^{\rm\nu^\f13\tau}}_{L^2[t'_k,t]}\sup_{\tau\in[t'_k,t]}\nnr{\fM_1(\tau,\na)P_{\neq}\udl{\th}(\tau)}_{L^2_{z,v}}e^{2\nu\tau}\\
	\ls&\ueps\nnr{\sF\rr{\vv{1}_{\nr{\na_v-A'}\geq10D_2\rr{\ueps A'}^\f13}\fM_2(\tau,\na)\rho}(\tau,k+1)e^{\rm\nu^\f13\tau}}_{L^2[t'_k,t]}\ls \ueps{A'}^2e^{-5D_1(\ueps A')^{\f13}}Z(t).
	\ea\]
	\eqref{h33} follows by the same argument.
	
	\underline{Step 2: Density estimate.} Proof of \eqref{rho31}. By \eqref{rho33} and \eqref{definereso} we have
	\[\ba
	&\wh{\rho}(t,k)=\wh{H}(t,k)-\int_{t'_k}^t\wh{H}(\tau,k)\wh{R}(t-\tau,k)d\tau,\quad \wh{H}(t,k)=\wh{\th}(t'_k,k,k\tent)S(t-t'_k,k)\\
	&-\int_{t'_k}^t\wh{\rho}(\tau,k)\f{1-e^{-\nu(t-\tau)}}{\nu}\wh{\uh}(\tau,0,k\f{1-e^{-\nu(t-\tau)}}{\nu})S(t-\tau,k)d\tau+\wh{\sR_1'}(t,k).
	\ea\]
	By following the proof of \eqref{rhokupper} we obtain \eqref{rho31}.
	
	\underline{Step 3: One-step transfer.} Proof of \eqref{h31}. On the support $\nr{\xi-k\tentau}\leq1$,\linebreak $\nr{A'-k\tentau}\leq10D_2\rr{\ueps A'}^\f13$
	\[\nr{\bar{\eta}(\tau,k-1,\xi)}\leq\entau\rr{\nr{\xi-k\tentau}+\f1k\nr{k\tentau}}\leq\entau\f{A'}{k}\rr{1+\f{k+10D_2\rr{\ueps A'}^\f13}{A'}},\]
	\[\nr{\log S(\tau,0,-1,\xi-k\tentau)+\f{\nu {A'}^3}{3k^3}}\leq\f{\nu}{3}\nr{\tau^3-\f{{A'}^3}{k^3}}+2\nu\tau^2\leq\f{\nu {A'}^2}{k^2}\rr{4+\f{10D_2\rr{\ueps A'}^\f13}{k}}.\]
	Equations~\eqref{h33} and \eqref{963} give
	\[\ba&\wh{\th}(t'_{k-1},k-1,\xi)=-\ueps\int_T^{t'_{k-1}}\vv{1}_{|k\tentau-A'|\leq10D_2\rr{\ueps A'}^\f13}\wh{\rho}(\tau,k)\f{\bar{\eta}(\tau,k-1,\xi)}{k}\wh{\vf}(\xi-k\tentau)\\
	&\times S(\tau,0,-1,\xi-k\tentau)S(t'_{k-1},\tau,k-1,\xi)d\tau+\wh{\sR'_2}(k-1,\xi):=\wh{\gimel}(k-1,\xi)+\wh{\sR'_2}(k-1,\xi).\ea\]
	For $\wh{\gimel}$ we have
	\[\ba
	&\nnr{\sF\rr{\fM(t'_{k-1},\na)\gimel}(k-1,\xi)}_{L^2_\xi}e^{\rm\nu^\f13t'_{k-1}}\leq\f{\ueps A'}{k^2}e^{-\f{\nu {A'}^3}{3k^3}+\f{10D_2\nu{A'}^2\rr{\ueps A'}^\f13}{k^3}}e^{\rm\nu^\f13(t'_{k-1}-t'_k)}\\
	&\quad \left\{1+\f{100D_2\rr{\rr{\ueps A'}^\f13+k}}{A'}+\f{20\nu{A'}^2}{k^2}\right\} e^{-\f32\nu t_k'}\nnr{\int_{t'_k}^{t'_{k-1}}\nr{\sF\rr{\fM(\tau,\na)\rho}(\tau,k)}e^{\rm\nu^\f13\tau}\right.\\
	&\left.\wh{\vf}(\xi-k\tentau)e^{-\f{\nu\tau}{2}}d\tau}_{L^2_\xi}
	\leq\f{\ueps A'}{k^{2.5}}e^{-\f{\nu {A'}^3}{3k^3}+\f{10D_2\nu{A'}^2\rr{\ueps A'}^\f13}{k^3}}e^{\rm\nu^\f13(t'_{k-1}-t'_k)}\\
	&\times\left\{1+\f{100D_2\rr{\rr{\ueps A'}^\f13+k}}{A'}+\f{20\nu{A'}^2}{k^2}\right\}e^{-\f32\nu t_k'}\nnr{\sF\rr{\fM(\tau,\na)\rho}(\tau,k)e^{\rm\nu^\f13\tau}}_{L^2[t_k',t_{k-1}']}\nnr{\wh{\vf}(\xi)}_{L^1_\xi}\\
	&\leq\f{\ueps A'}{k^{3}}e^{-\f{\nu {A'}^3}{3k^3}+\f{10D_2\nu{A'}^2\rr{\ueps A'}^\f13}{k^3}}\left\{1+D_3\rr{\f{\rr{\rr{\ueps A'}^\f13+k}}{A'}+\f{1}{k^2}}\right\}\nnr{\sF\rr{\fM(t_k',\na)\th}(t'_k,k,\xi)}_{L^2_\xi}e^{\rm\nu^\f13t_k'}\\
	&\times e^{2\rm\nu^\f13(t'_{k-1}-t'_k)}+\ueps^{0.8}{A'}^2e^{-4D_1(\ueps A')^{\f13}}Z(t'_{k-1}).
	\ea\]
	We use a variant of \eqref{mtp} to deal with the multipliers. Note that $\nu{t_{k-1}'}^2\ls\f{\nu {A'}^2}{k^2}\leq1$. The last inequality uses \eqref{rho31} and the fact 
	\[\ba\f{\ueps A'}{k^{2.5}}e^{-\f{\nu {A'}^3}{3k^3}+\f{10D_2\nu{A'}^2\rr{\ueps A'}^\f13}{k^3}}\left\{1+\f{100D_2\rr{\rr{\ueps A'}^\f13+k}}{A'}+\f{20\nu{A'}^2}{k^2}\right\}\\
	\times e^{\rm\nu^\f13(t'_{k-1}-t'_k)}\ueps^{0.8}{A'}^2e^{-5D_1\rr{\ueps A'}^\f13}Z(t'_{k-1})\leq\ueps^{0.8}{A'}^2e^{-4D_1(\ueps A')^{\f13}}Z(t'_{k-1}).\ea\]
	This, together with \eqref{R31}, proves \eqref{h31}.
\end{proof}

The recurrence below the source mode still requires a starting estimate. At mode $K'$, no higher mode contributes to the density equation, so the accumulated forcing can be estimated directly through the Volterra resolvent. The next proposition bounds this density and the first packet transferred to mode $K'-1$.

\begin{prop}[Response at the source mode]\label{6.2}
	Suppose $K'\ge2$.
	For any $t\in [T,t_{K'-1}']$ we have
	\beq\label{rho32}\nnr{\sF\rr{\fM(\tau,\na)\rho}(\tau,K')e^{\rm\nu^\f13\tau}}_{L^2[T,t]}\ls Z(t).\eeq
	Moreover,
	\beq\label{h32}\nnr{\sF\rr{\fM(t_{K'-1}',\na)\th}(t_{K'-1}',K'-1,\xi)}_{L^2_\xi}e^{\rm\nu^\f13t_{K'-1}'}\leq {A'}^2Z(t_{K'-1}').\eeq
\end{prop}
\begin{proof}
We start with the source mode, where the triangular density equation has no contribution from higher modes. The source trace is bounded by its accumulated norm $Z$. Inserting this density bound into the resonant transfer integral gives the first lower-mode packet.

	By Duhamel's formula (see \eqref{rba}), we know that $\wh{\rho}(t,k)=0$ for $k>K'$, $t\in[T,\infty)$. On $[T,t_{K'-1}']$ we have
	\[\ba
	&\wh{\rho}(t,K')=\wh{H}(t,K')-\int_{t'_k}^t\wh{H}(\tau,K')\wh{R}(t-\tau,K')d\tau,\quad \wh{H}(t,K')=\int_T^t\wh{\tQ}(\tau,K',K'\tent)S(t-\tau,K')d\tau\\
	&-\int_{T}^t\wh{\rho}(\tau,K')\f{1-e^{-\nu(t-\tau)}}{\nu}\wh{\uh}(\tau,0,K'\f{1-e^{-\nu(t-\tau)}}{\nu})S(t-\tau,K')d\tau.
	\ea\]
	Following the proof of \eqref{rhokupper} we obtain \eqref{rho32}. Note that
	\[\ba&\nnr{\fM(t,K',K'\tent)\int_T^{t}\wh{\tQ}(\tau,K',K'\tent)S(t-\tau,K')d\tau e^{\rm\nu^\f13t}}_{L^2[T,\tilde{T}]}\\
	\ls&\int^{\tilde{T}}_T\nnr{\sF\rr{\fM(\tau,\na)\tQ}(\tau,K',\xi)}_{L^2_\xi}e^{\rm\nu^\f13\tau+\f12\nu\tau}d\tau\ls Z(\tilde{T})\ea\]
	for $\tilde{T}\in[T,t_{K'-1}']$. We next prove \eqref{h32}. As in \eqref{h33} we obtain
	\[\ba\wh{\th}(t'_{K'-1},K'-1,\xi)=&-\int_T^{t'_{K'-1}}\vv{1}_{|K'\tentau-A'|\leq10D_2\rr{\ueps A'}^\f13}\wh{\rho}(\tau,K')\f{\bar{\eta}(\tau,K'-1,\xi)}{K'}\wh{\udl{\th}}(\tau,-1,\xi-K'\tentau)\\
	&\times S(t'_{K'-1},\tau,K'-1,\xi)d\tau+\wh{\sR'_2}(K'-1,\xi).\ea\]
	Following the proof of \eqref{h31} we obtain
	\[\ba&\nnr{\sF\rr{\fM(t_{K'-1}',\na)\th}(t_{K'-1}',K'-1,\xi)}_{L^2_\xi}e^{\rm\nu^\f13t_{K'-1}'}\ls\ueps A'\nnr{\sF\rr{\fM(\tau,\na)\rho}(\tau,K')e^{\rm\nu^\f13\tau}}_{L^2[T,t'_{K'-1}]}\\
	&+\ueps^{0.9}{A'}^2e^{-4D_1\rr{\ueps A'}^\f13}Z(t'_{K'-1}).\ea\]
	This proves \eqref{h32}.
\end{proof}

We now combine the starting bound with the recurrence \eqref{h31}. Only an intermediate range of modes can contribute appreciable amplification: the transfer is contractive above the upper cutoff, and collisional attenuation makes it contractive below the lower cutoff. The next proposition isolates the product over this range and controls the accumulated errors. Set
\beq\label{kAot}\kAo=\llcorner\rr{2\ueps A'}^\f13\lrcorner,\quad \kAt=\llcorner\nu^\f12A'\rr{\ueps A'}^\f16\lrcorner+1.\eeq

\begin{prop}[Iteration of the forced recurrence]\label{prop:forced-iteration}
	For $1\leq k\leq K'-1$ we have
	\beq\label{h39}\ba
	&\nnr{\sF\rr{\fM(t_k',\na)\th}(t'_k,k,\xi)}_{L^2_\xi}e^{\rm\nu^\f13t_k'}\\
	&\ls\left\{\prod_{i=\max\{k+1,\kAt\}}^{\min\{K'-1,\kAo\}}\f{\ueps A'}{i^3}e^{-\f{\nu {A'}^3}{3i^3}+2\rm\nu^\f13(t'_{i-1}-t'_i)}\right\}{A'}^2Z(t'_{K'-1})+\ueps^{0.7}{A'}^2e^{-3D_1\rr{\ueps A'}^\f13}Z(t'_{k})
	\ea\eeq
\end{prop}
\begin{proof}
We first write the iterated recurrence with its accumulated inhomogeneous errors. We then show that the transfer factors are contractive outside the interval between the two cutoffs and bound the remaining product. Finally, the source-mode estimate \eqref{h32} supplies the starting value.

	To simplify the notation we set
	\[\mathbf{B}(k)=\f{\ueps A'}{k^3}e^{-\f{\nu A'^3}{3k^3}+\f{10D_2\nu{A'}^2\rr{\ueps A'}^\f13}{k^3}+2\rm\nu^\f13(t'_{k-1}-t'_k)}\left\{1+D_3\rr{\f{\rr{\ueps A'}^\f13+k}{A'}+\f1{k^2} }\right\}\]
	for $1\leq k\leq\f{A'+1}{\f{1-e^{-\nu T}}{\nu} }$. We split the proof into three steps.
	
	\underline{Step 1: Iterating the recurrence.} We prove 
	\beq\label{h98}\ba
	&\nnr{\sF\rr{\fM(t_k',\na)\th}(t'_k,k,\xi)}_{L^2_\xi}e^{\rm\nu^\f13t_k'}\leq\prod_{i=k+1}^{K'-1}\mathbf{B}(i)\nnr{\sF\rr{\fM(t_{K'-1}',\na)\th}(t_{K'-1}',K'-1,\xi)}_{L^2_\xi}e^{\rm\nu^\f13t_{K'-1}'}\\
	&+\ueps^{0.7}{A'}^2e^{-4D_1\rr{\ueps A'}^\f13}Z(t'_{k})\sum_{i=k}^{K'-2}\prod_{j=k+1}^{i}\mathbf{B}(j).\ea\eeq
	We use descending induction. The case $\tilde{k}=K'-1$ is immediate. Suppose \eqref{h98} holds for $k\leq\tilde{k}\leq K'-1$; we prove it for $\tilde{k}=k-1$. By \eqref{h31} we have
	\[\ba
	&\nnr{\sF\rr{\fM(t_{k-1}',\na)\th}(t'_{k-1},k-1,\xi)}_{L^2_\xi}e^{\rm\nu^\f13t_{k-1}'}\leq \mathbf{B}(k)\left\{\prod_{i=k+1}^{K'-1}\mathbf{B}(i)\nnr{\sF\rr{\fM(t_{K'-1}',\na)\th}(t_{K'-1}',K'-1,\xi)}_{L^2_\xi}\right.\\
	&\left.\times e^{\rm\nu^\f13t_{K'-1}'}+\ueps^{0.7}{A'}^2e^{-4D_1\rr{\ueps A'}^\f13}Z(t'_{k})\sum_{i=k}^{K'-2}\prod_{j=k+1}^{i}\mathbf{B}(j)\right\}+\ueps^{0.7}{A'}^2e^{-4D_1\rr{\ueps A'}^\f13}Z(t'_{k-1})
	\leq\prod_{i=k}^{K'-1}\mathbf{B}(i)\\
	&\times\nnr{\sF\rr{\fM(t_{K'-1}',\na)\th}(t_{K'-1}',K'-1,\xi)}_{L^2_\xi}e^{\rm\nu^\f13t_{K'-1}'}+\ueps^{0.7}{A'}^2e^{-4D_1\rr{\ueps A'}^\f13}Z(t'_{k-1})\rr{\sum_{i=k}^{K'-2}\prod_{j=k}^{i}\mathbf{B}(j)+1}\\
	&=\prod_{i=k}^{K'-1}\mathbf{B}(i)\nnr{\sF\rr{\fM(t_{K'-1}',\na)\th}(t_{K'-1}',K'-1,\xi)}_{L^2_\xi}e^{\rm\nu^\f13t_{K'-1}'}+\ueps^{0.7}{A'}^2e^{-4D_1\rr{\ueps A'}^\f13}Z(t'_{k-1})\sum_{i=k-1}^{K'-2}\prod_{j=k}^{i}\mathbf{B}(j).
	\ea\]
	
	\underline{Step 2: Restricting the amplification product.} We prove
	\beq\label{Bi}\prod_{i=k+1}^{K'-1}\mathbf{B}(i)\ls\prod_{i=\max\{k+1,\kAt\}}^{\min\{K'-1,\kAo\}}\f{\ueps A'}{i^3}e^{-\f{\nu {A'}^3}{3i^3}+2\rm\nu^\f13(t'_{i-1}-t'_i)}\ls e^{10\rr{\ueps A'}^\f13}\eeq
	for $1\leq k\leq K'-1$.
	We first prove that if $\kAo<k\leq\f{A'+1}{\f{1-e^{-\nu T}}{\nu} }$ or $1\leq k<\kAt$, then $\mathbf{B}(k)\leq\f23$. For $\kAo<k\leq\f{A'+1}{\f{1-e^{-\nu T}}{\nu} }$,
	\[\mathbf{B}(k)\leq\f{\ueps A'}{2\ueps A'}\left\{1+D_3\rr{\f{\rr{\ueps A'}^\f13+\f{A'+1}{\f{1-e^{-\nu T}}{\nu} }}{A'}+\f1{\rr{\ueps A'}^\f23}}\right\}\leq\f23.\]
	For $1\leq k<\kAt$,
	\[\mathbf{B}(k)\leq\ueps A'e^{-\f{1}{6\nu^\f12\rr{\ueps A'}^\f12}}\rr{1+2D_3}\leq\f{\ueps\nu^{-\f12}}{3}e^{-\f{1}{2\sqrt{3}}\ueps^{-\f12}\nu^{-\f16}}\rr{1+2D_3}\ll1.\]
	We next prove \eqref{Bi}.
	\[\ba
	&\prod_{i=k+1}^{K'-1}\mathbf{B}(i)\leq\prod_{i=\max\{k+1,\kAt\}}^{\min\{K'-1,\kAo\}}\mathbf{B}(i)\leq\prod_{i=\max\{k+1,\kAt\}}^{\min\{K'-1,\kAo\}}\f{\ueps A'}{i^3}e^{-\f{\nu {A'}^3}{3i^3}}\prod_{i=\kAt}^{\kAo}e^{\f{10D_2\nu{A'}^2\rr{\ueps A'}^\f13}{i^3}}\left\{1+\right.\\
	&\left.D_3\rr{\f{\rr{\ueps A'}^\f13+i}{A'}+\f1{i^2} }\right\}
	\leq\prod_{i=\max\{k+1,\kAt\}}^{\min\{K'-1,\kAo\}}\f{\ueps A'}{i^3}e^{-\f{\nu {A'}^3}{3i^3}}\exp\left\{\sum_{i=\kAt}^{\kAo}\f{10D_2\nu{A'}^2\rr{\ueps A'}^\f13}{i^3}\right.\\
	&\left.+\sum_{i=\kAt}^{\kAo}D_3\rr{\f{\rr{\ueps A'}^\f13+i}{A'}+\f1{i^2}} \right\}
	\ls\prod_{i=\max\{k+1,\kAt\}}^{\min\{K'-1,\kAo\}}\f{\ueps A'}{i^3}e^{-\f{\nu {A'}^3}{3i^3}}.
	\ea\]
	By Stirling's formula we have
	\[\prod_{i=\max\{k+1,\kAt\}}^{\min\{K'-1,\kAo\}}\f{\ueps A'}{i^3}e^{-\f{\nu {A'}^3}{3i^3}+2\rm\nu^\f13(t'_{i-1}-t'_i)}\leq\prod_{i=\max\{k+1,\kAt\}}^{\min\{K'-1,\kAo\}}\f{2\ueps A'}{i^3}\leq\prod_{i=1}^{\llcorner\rr{2\ueps A'}^\f13\lrcorner}\f{2\ueps A'}{i^3}\ls e^{10\rr{\ueps A'}^\f13}.\]
	
	\underline{Step 3: Inserting the source-mode bound.} We prove \eqref{h39}. By \eqref{h98} and \eqref{Bi} we have
	\[\ba
	&\nnr{\sF\rr{\fM(t_k',\na)\th}(t'_k,k,\xi)}_{L^2_\xi}e^{\rm\nu^\f13t_k'}\ls\prod_{i=\max\{k+1,\kAt\}}^{\min\{K'-1,\kAo\}}\f{\ueps A'}{i^3}e^{-\f{\nu {A'}^3}{3i^3}+2\rm\nu^\f13(t'_{i-1}-t'_i)}\\
	&\times\nnr{\sF\rr{\fM(t_{K'-1}',\na)\th}(t_{K'-1}',K'-1,\xi)}_{L^2_\xi}e^{\rm\nu^\f13t_{K'-1}'}+\ueps^{0.7}{A'}^2K'e^{3\rr{\ueps A'}^\f13}e^{-4D_1\rr{\ueps A'}^\f13}Z(t'_{k})\\
	&\ls\prod_{i=\max\{k+1,\kAt\}}^{\min\{K'-1,\kAo\}}\f{\ueps A'}{i^3}e^{-\f{\nu {A'}^3}{3i^3}+2\rm\nu^\f13(t'_{i-1}-t'_i)}{A'}^2Z(t'_{K'-1})+\ueps^{0.7}{A'}^2e^{-3D_1\rr{\ueps A'}^\f13}Z(t'_{k}).
	\ea\]
	The last inequality uses \eqref{h32}. Note that $K'\ls\f{A'}{T}\ll\ueps A'$.
\end{proof}

The preceding proposition controls the packets at critical times. We next recover density and profile estimates on $[T,t]$ by summing over the resonant modes and treating the direct source contribution separately. Combined with \eqref{h39}, this gives the response estimate used in the nonlinear induction.

\begin{prop}[Global bounds in the resonant regime]\label{prop:forced-global}
	For $t\in[T,\infty)$, we have
	\beq\label{hrhoh}\ba
	&\nnr{\fM_4(\tau,\na)\rho(\tau)e^{\rm\nu^\f13\tau}}_{L^2([T,t], L^2_{z})}+\left\|\fM_{2}(t,\na)\rr{e^{-|\alpha|\nu t}v^\alpha P_{\neq}\th(t)} \right\|_{\lta L^2_{z,v}}e^{\rm\nu^\f13t}\\
	&+\nnr{v^\alpha\<\na_v\>\fM(t,0,\emnt\na_v)P_{0}h(t)}_{\lta L^2_{v}}\\
	&\times e^{\nu t}
	\ls {A'}^6K'Z(t)+{A'}^4K'\sup_{\substack{k:k\leq K'-1\\T\leq t'_k\leq t}}\nnr{\sF\rr{\fM(t'_{k},\na)\th}(t'_{k},k,\xi)}_{L^2_\xi}e^{\rm\nu^\f13t'_{k}}.
	\ea\eeq
\end{prop}
\begin{proof}
We estimate the density by summing the localized critical-mode bounds and adding the direct response at $K'$. The Fourier tail is controlled by \eqref{remaining3}. With this density estimate available, the weighted profile bounds follow from Duhamel's formula and the zero-mode semigroup.

	We first estimate $\rho(t)$. By \eqref{remaining3} and \eqref{remark3} we have
	\[\ba
	&\nnr{\fM_4(\tau,\na)\rho(\tau)e^{\rm\nu^\f13\tau}}_{L^2([T,t], L^2_{z})}\leq\nnr{\vv{1}_{\nr{\na_v-A'}\leq10D_2\rr{\ueps A'}^\f13}\fM_4(\tau,\na)\rho(\tau)e^{\rm\nu^\f13\tau}}_{L^2([T,t], L^2_{z})}\\
	&+\nnr{\vv{1}_{\nr{\na_v-A'}\geq10D_2\rr{\ueps A'}^\f13}\right.\\
	&\left.\times\fM_4(\tau,\na)\rho(\tau)e^{\rm\nu^\f13\tau}}_{L^2([T,t], L^2_{z})}
	\ls {A'}^4\left\{\sum_{\substack{k:k\leq K'-1\\T\leq t'_k\leq t}}\nnr{\sF\rr{\fM(\tau,\na)\rho}(\tau,k)e^{\rm\nu^\f13\tau}}_{L^2[t'_k,t]}\right.\\
	&\left.+\nnr{\sF\rr{\fM(\tau,\na)\rho}(\tau,K')e^{\rm\nu^\f13\tau}}_{L^2[T,\min\{t'_{K'-1},t\}]}\right\}+{A'}^4e^{-5D_1\rr{\ueps A'}^\f13}Z(t)
	\ls {A'}^4e^{-5D_1\rr{\ueps A'}^\f13}Z(t)\\
	&+{A'}^4\left\{K'\sup_{\substack{k:k\leq K'-1\\T\leq t'_k\leq t}}\left\{\nnr{\sF\rr{\fM(t'_{k},\na)\th}(t'_{k},k,\xi)}_{L^2_\xi}e^{\rm\nu^\f13t'_{k}}+\ueps^{0.8}{A'}^2e^{-5D_1\rr{\ueps A'}^\f13}Z(t)\right\}+Z(t)\right\}\\
	&\ls{A'}^6K'Z(t)+{A'}^4K'\sup_{\substack{k:k\leq K'-1\\T\leq t'_k\leq t}}\nnr{\sF\rr{\fM(t'_{k},\na)\th}(t'_{k},k,\xi)}_{L^2_\xi}e^{\rm\nu^\f13t'_{k}}.
	\ea\]
	Note that $\wh{\rho}(t,k)=0$ for $k>K'$. The third inequality uses \eqref{rho31} and \eqref{rho32}. We next estimate $h(t)$. Arguing as in the proof of \eqref{neqneqh} and \eqref{p0poh} we have
	\[\ba
	&\left\|\fM_{2}(t,\na)\rr{e^{-|\alpha|\nu t}v^\alpha P_{\neq}\th(t)} \right\|_{\lta L^2_{z,v}}e^{\rm\nu^\f13t}
	+\nnr{v^\alpha\<\na_v\>\fM(t,0,\emnt\na_v)P_{0}h(t)}_{\lta L^2_{v}}\ent\\
	\ls&\int^t_T\nnr{\fM_2(\tau,\na)\rr{e^{-\nu|\alpha|\tau}v^\alpha \tQ(\tau)}}_{\lta L^2_{z,v}}e^{\rm\nu^\f13\tau} d\tau+\nnr{\fM_4(\tau,\na)\rho(\tau)e^{\rm\nu^\f13\tau}}_{L^2([T,t], L^2_{z})}\\
	\ls &{A'}^2Z(t)+\nnr{\fM_4(\tau,\na)\rho(\tau)e^{\rm\nu^\f13\tau}}_{L^2([T,t], L^2_{z})}.
	\ea\]
\end{proof}

\subsection{Case 3: $K'>0$ and $K'\f{1-e^{-\nu T}}{\nu}\leq A'+1$ and $A'\geq\f{\nu^{-\f12}}{3}$.}
At these higher velocity frequencies, the background packet has already undergone substantial collisional decay by the time it can interact resonantly with the source. The resulting echo kernel has small operator norm, so the density estimate closes by absorption rather than by a discrete iteration. This yields a bound with only a polynomial loss in $A'$, which is essential for high-frequency pieces of the nonlinear source.

\begin{prop}[High-frequency forced response]\label{prop:forced-diffusive}
	For $t\in[T,\infty)$, we have
	\beq\label{hrho3}\ba
	&\nnr{\fM_4(\tau,\na)\rho(\tau)e^{\rm\nu^\f13\tau}}_{L^2([T,t], L^2_{z})}+\left\|\fM_{2}(t,\na)\rr{e^{-|\alpha|\nu t}v^\alpha P_{\neq}\th(t)} \right\|_{\lta L^2_{z,v}}e^{\rm\nu^\f13t}\\
	&+\nnr{v^\alpha\<\na_v\>\fM(t,0,\emnt\na_v)P_{0}h(t)}_{\lta L^2_{v}}e^{\nu t}
	\ls {A'}^4Z(t).\ea\eeq
\end{prop}
\begin{proof}
We isolate the resonant density kernel and show that its operator norm is small because of collisional decay of the background. The Fourier tail is exponentially small. After absorbing the background and echo terms in the density inequality, we recover the profile bounds from Duhamel's formula.

	We split the proof into three steps.
	
	\underline{Step 1: Separating the Fourier tail.} We prove for $k\ge1$,
	\beq\label{rho96}\ba
	&\wh{\rho}(t,k)=\int_T^t\wh{\tQ}(\tau,k,k\tent)S(t-\tau,k)d\tau-\int_T^t\wh{\rho}(\tau,k)\f{1-e^{-\nu(t-\tau)}}{\nu}\wh{\sM}(k\f{1-e^{-\nu(t-\tau)}}{\nu})S(t-\tau,k)d\tau\\
	&-\int_T^t\wh{\rho}(\tau,k)\f{1-e^{-\nu(t-\tau)}}{\nu}\wh{\uh}(\tau,0,k\f{1-e^{-\nu(t-\tau)}}{\nu})S(t-\tau,k)d\tau-\int_T^t\wh{\rho}(\tau,k+1)\vv{1}_{\nr{A'-(k+1)\tentau}\leq10D_2\rr{\ueps A'}^\f13}\\
	&\times\f{k}{k+1}\f{1-e^{-\nu(t-\tau)}}{\nu}\wh{\uth}(\tau,-1,k\tent-(k+1)\tentau)S(t-\tau,k)d\tau+\wh{\sR}(t,k),
	\ea\eeq
	with $\wh{\sR}(t,k)$ satisfying
	\[\nnr{\fM_4(\tau,\na)\sR(\tau)e^{\rm\nu^\f13\tau}}_{L^2([T,t],L^2_z)}\ls\ueps {A'^6}e^{-5D_1\rr{\ueps A'}^\f13}Z(t).\]
	By \eqref{rba} we know that \eqref{rho96} holds with 
	\[\wh{\sR}(t,k)=\int_T^t\sF\rr{\cN^\dagger\rr{\vv{1}_{\nr{A'-\na_v}\geq10D_2\rr{\ueps A'}^\f13}\th,P_{-1}\uth}}(\tau,k,k\tent)S(t-\tau,k)d\tau.\]
	Arguing as in the proof of \eqref{rho33} we have 
	\[\ba
	&\nnr{\fM_4(\tau,\na)\sR(\tau)e^{\rm\nu^\f13\tau}}_{L^2([T,t],L^2_z)}\ls\int_T^t\nnr{\fM_4(\tau,\na)\cN^\dagger\rr{\vv{1}_{\nr{A'-\na_v}\geq10D_2\rr{\ueps A'}^\f13}\th,P_{-1}\uth}}_{L^2_{z,v}}e^{\rm\nu^\f13\tau+\f12\nu\tau}d\tau\\
	\ls&\nnr{\vv{1}_{\nr{A'-\na_v}\geq10D_2\rr{\ueps A'}^\f13}\fM_6(\tau,\na)\rho(\tau)e^{\rm\nu^\f13\tau}}_{L^2([T,t],L^2_z)}\\
	&\quad \sup_{\tau\in[T,t]}\nnr{\fM_5(\tau,\na)P_{-1}\uth(\tau)}_{L^2_{z,v}}e^{2\nu\tau}
	\ls\ueps{A'}^6e^{-5D_1\rr{\ueps A'}^\f13}Z(t).
	\ea\]
	The last inequality uses \eqref{remaining3}.
	
	\underline{Step 2: Absorbing the echo kernel.} We prove
	\[\nnr{\fM_4(\tau,\na)\rho(\tau)e^{\rm\nu^\f13\tau}}_{L^2([T,t], L^2_{z})}\ls {A'}^4Z(t).\]
	By \eqref{rho96} and \eqref{definereso} we obtain
	\[\ba
	&\wh{\rho}(t,k)=\wh{H}(t,k)-\int_T^t\wh{H}(\tau,k)\wh{R}(t-\tau,k)d\tau,\quad 
	\wh{H}(t,k)=\int_T^t\wh{\tQ}(\tau,k,k\tent)S(t-\tau,k)d\tau\\
	&-\int_T^t\wh{\rho}(\tau,k)\f{1-e^{-\nu(t-\tau)}}{\nu}\wh{\uh}(\tau,0,k\f{1-e^{-\nu(t-\tau)}}{\nu})S(t-\tau,k)d\tau
	-\int_T^t\wh{\rho}(\tau,k+1)\vv{1}_{\nr{A'-(k+1)\tentau}\leq10D_2\rr{\ueps A'}^\f13}\\
	&\times\f{k}{k+1}\f{1-e^{-\nu(t-\tau)}}{\nu}\wh{\uth}(\tau,-1,k\tent-(k+1)\tentau)S(t-\tau,k)d\tau+\wh{\sR}(t,k).
	\ea\]
	Denote the first three terms in $\wh{H}(t,k)$ by $\wh{\daleth_1}$, $\wh{\daleth_2}$, $\wh{\daleth_3}$. Following the proof of \eqref{rhokupper} we obtain
	\[\nnr{\fM_4(\tau,\na)\rr{\daleth_1+\daleth_2}(\tau)e^{\rm\nu^\f13\tau}}_{L^2([T,t], L^2_{z})}\ls {A'}^4Z(t)+\ueps\nnr{\fM_4(\tau,\na)\rho(\tau)e^{\rm\nu^\f13\tau}}_{L^2([T,t], L^2_{z})}.\]
	Now we estimate $\wh{\daleth_3}$. The following holds on $\nr{k\tent-(k+1)\tentau}\leq1$, $\nr{A'-(k+1)\tentau}\leq10D_2\rr{\ueps A'}^\f13$
	\[\ba
	&\nr{k\f{1-e^{-\nu(t-\tau)}}{\nu}}\leq\entau\left\{\nr{k\tent-(k+1)\tentau}+\f{1}{k+1}\nr{(k+1)\tentau}\right\}\\
	&\leq\entau\rr{1+\f{A'+10D_2\rr{\ueps A'}^\f13}{k+1}}.
	\ea\]
	Thus we have
	\[\ba
	&\nr{\sF\rr{\fM_4(t,\na)\daleth_3}(t,k)}e^{\rm\nu^\f13t}\ls\f{A'}{(k+1)^2}\int_T^t\nr{\sF\rr{\fM_4(\tau,\na)\rho}(\tau,k+1)}e^{\rm\nu^\f13\tau}\\
	&\quad \vv{1}_{\nr{A'-(k+1)\tentau}\leq10D_2\rr{\ueps A'}^\f13}e^{2\nu\tau}\\
	&\times \nr{\sF\rr{\fM_4(\tau,\na)\uth}(\tau,-1,k\tent-(k+1)\tentau)}\emnt d\tau\\
	&\ls\f{\ueps A'}{(k+1)^2}\int_T^t\nr{\sF\rr{\fM_4(\tau,\na)\rho}(\tau,k+1)}e^{\rm\nu^\f13\tau}\\
	&\times\vv{1}_{\nr{A'-(k+1)\tentau}\leq10D_2\rr{\ueps A'}^\f13}\nr{\wh{\vf}(k\tent-(k+1)\tentau)}e^{-0.008\nu^\f13\tau}\emnt d\tau\\
	&\ls\f{\ueps A'}{(k+1)^2}e^{-10^{-3}\f{\nu^\f13A'}{k+1}}\int_T^t\nr{\sF\rr{\fM_4(\tau,\na)\rho}(\tau,k+1)}e^{\rm\nu^\f13\tau}\nr{\wh{\vf}(k\tent-(k+1)\tentau)}\emnt d\tau.
	\ea\]
	The second inequality uses \eqref{hudlneq}. The third inequality uses $\tau\geq\tentau\geq\f18\f{A'}{k+1}$ on $\nr{A'-(k+1)\tentau}\leq10D_2\rr{\ueps A'}^\f13$. Note that $\f{\ueps A'}{(k+1)^3}e^{-10^{-3}\f{\nu^\f13A'}{k+1}}\ls\f{\ueps A'}{(k+1)^3}\f{(k+1)^3}{\nu{A'}^3}=\f{\ueps}{\nu {A'}^2}\ls\ueps$ since $A'\geq\f{\nu^{-\f12}}{3}$. Thus by Schur's test
	\[\nnr{\fM_4(\tau,\na)\daleth_3(\tau)e^{\rm\nu^\f13\tau}}_{L^2([T,t], L^2_{z})}\ls\ueps\nnr{\fM_4(\tau,\na)\rho(\tau)e^{\rm\nu^\f13\tau}}_{L^2([T,t], L^2_{z})}.\]
	By \eqref{resolvent} we have
	\[\ba
	&\nnr{\fM_4(\tau,\na)\rho(\tau)e^{\rm\nu^\f13\tau}}_{L^2([T,t], L^2_{z})}\ls\nnr{\fM_4(\tau,\na)H(\tau)e^{\rm\nu^\f13\tau}}_{L^2([T,t], L^2_{z})}\\
	&\ls {A'}^4Z(t)+\ueps\nnr{\fM_4(\tau,\na)\rho(\tau)e^{\rm\nu^\f13\tau}}_{L^2([T,t], L^2_{z})}.
	\ea\]
	Note that $\ueps{A'}^6e^{-5D_1\rr{\ueps A'}^\f13}\ls\ueps\nu^{-3}e^{-D_1\rr{\ueps \nu^{-\f12}}^\f13}\ll1$ since $\ueps A'\gs\nu^{-\f16}\gg1$.
	
	\underline{Step 3: Recovering the profile bounds.} We prove \eqref{hrho3}. Arguing as in the proof of \eqref{neqneqh} and \eqref{p0poh} we have
	\[\ba
	&\left\|\fM_{2}(t,\na)\rr{e^{-|\alpha|\nu t}v^\alpha P_{\neq}\th(t)} \right\|_{\lta L^2_{z,v}}e^{\rm\nu^\f13t}
	+\nnr{v^\alpha\<\na_v\>\fM(t,0,\emnt\na_v)P_{0}h(t)}_{\lta L^2_{v}}\ent\\
	\ls&\int^t_T\nnr{\fM_2(\tau,\na)\rr{e^{-\nu|\alpha|\tau}v^\alpha \tQ(\tau)}}_{\lta L^2_{z,v}}e^{\rm\nu^\f13\tau} d\tau+\nnr{\fM_4(\tau,\na)\rho(\tau)e^{\rm\nu^\f13\tau}}_{L^2([T,t], L^2_{z})}\\
	\ls &{A'}^2Z(t)+\nnr{\fM_4(\tau,\na)\rho(\tau)e^{\rm\nu^\f13\tau}}_{L^2([T,t], L^2_{z})}\ls{A'}^4Z(t).
	\ea\]
\end{proof}

\section{Estimates of the finite Taylor expansion}
We now combine the homogeneous and forced linear estimates to control the finite nonlinear expansion. The main task is to retain a factor $B(t)^n$ at order $n$, despite the additional echoes generated by the nonlinear sources. Proposition~\ref{prop:envelope} constructs $B$ from the leading cascade and shows that its powers dominate the relevant forced amplification products. Proposition~\ref{prop:nonlinear-coefficients} then closes the induction by separating each source into a Fourier tail, a localized nonzero-mode part, and a zero-mode part. Since $B(\infty)=A^{-N_1}$, the resulting bounds make the higher-order coefficients smaller than the leading packet and provide the input for the remainder estimates in Section~6.

\subsection{Estimates for $h^{(1,1)}$.}
We begin by choosing the initial amplitude and constructing a continuous envelope for the discrete echo chain. The three parts of the next proposition have distinct roles: part~(1) fixes the final size and quantifies the size reached at the observation time; part~(2) compares the envelope with the leading packet and controls its global norms; part~(3) ensures that additional echoes at order $n$ can be absorbed into $B^n$. The last property is the link between the forced linear analysis and the nonlinear induction.
\begin{prop}[The cascade envelope]\label{prop:envelope}
	There exist $A^{-N_1}e^{-3K}<\eps<A^{-N_1}e^{-1.5K}$ and a nondecreasing function $B(t):[T,\infty)\rightarrow\R_{+}$ satisfying the following properties.
	
	\begin{enumerate}
		\item 
		$B(\infty)=A^{-N_1}$.	For a universal constant $c>0$, one has	
		\beq\label{cnu13} B(c\nu^{-\f13})
		\ge
		A^{-N_1}
		\exp\left(
		-\frac{1.02}{c}K^{\f23}\log K
		\right).\eeq
		
		\item 
		For
		$\llcorner(D_3\nu A^3)^{\f13}\lrcorner+1\le k\le K$,
		we have
		\beq\label{lowerb}
		\nnr{\wh{\th^{(1,1)}}(t_k,k,\xi)\vv{1}_{|\xi-A|\leq20D_2K}}_{L^2_\xi}
		\gs
		K^{-\f12}B(t_k),
		\eeq
		where $t_k$ is defined in \eqref{4.1}.
		Furthermore, for all $t\in [T,\infty)$,
\beq\label{h111}\ba
&\nnr{\fM_4(\tau,\na)\rho^{(1,1)}(\tau)e^{\rm\nu^\f13\tau}}_{L^2([T,t], L^2_{z})}+\left\|\fM_{2}(t,\na)\rr{e^{-|\alpha|\nu t}v^\alpha P_{\neq}\th^{(1,1)}(t)} \right\|_{\lta L^2_{z,v}}e^{\rm\nu^\f13t}\\
&+\nnr{v^\alpha\<\na_v\>\fM(t,0,\emnt\na_v)P_{0}h^{(1,1)}(t)}_{\lta L^2_{v}}e^{\nu t}
\ls\nu^{-\f13} {A}^4KB(t).\ea\eeq
		
		\item 
		Let $n,q,A'$ satisfy
		$
		2\le n\le A^{N'}$, $|q|\le n$,
		$2\mid n-q$, $|A'-qA|\le 40nD_2K$,
		and
		$
		\frac{1-e^{-\nu T}}{\nu}-1
		\le
		A'
		\le
		\frac13\nu^{-\f12}
		$.
		Then, for any $k_1, k_2\in\Z_+$ satisfying
		$
		k_2+1\le k_1
		\le
		\frac{A'+1}{(1-e^{-\nu T})/\nu}
		$,
		the following estimate holds:
		\beq\label{pbqhb}
		\left\{\prod_{i=\max\{k_2+1,\kAt\}}^{\min\{k_1-1,\kAo\}}\max\{
		\frac{\ueps A'}{i^{3}}
		e^{-\frac{\nu A'^3}{3i^3}},1\}\right\}
		B^n(t'_{k_1-1})
		\ls
		K^4B^n(t'_{k_2}),
		\eeq
		where $t'_k$ is defined in \eqref{401}, and $\kAo$ and $\kAt$ are defined in \eqref{kAot}.
	\end{enumerate}
\end{prop}

\begin{proof}[Proof of (1).]
We define $B$ by integrating a continuous majorant of the logarithmic one-step growth and normalize its limiting size. Estimating this integral gives both the admissible initial amplitude and the lower bound at the enhanced-dissipation time scale.

    Define
	\[
	B(t)
	:=
	\varepsilon
	\exp\left\{
	-\frac{\nu A^3}{3K}
	+
	\int_T^t
	\frac{A}{\tau^2}
	\max\left\{
	\log\frac{\ueps \tau^3}{A^2}
	-\frac13\nu\tau^3,
	K^{-2/3}
	\right\}
	d\tau
	\right\}.
	\]	
	Choose $\varepsilon$ so that
	$
	B(\infty)=A^{-N_1}.
	$
	For $x>e$, define $\gamma(x)>e$ by
	\beq\label{gamma}
	\frac{\gamma(x)}{\log \gamma(x)}=x.
	\eeq
	A direct computation yields
	\beq\label{geqv}
	x\log x+x\log\log x
	\le
	\gamma(x)
	\le
	x\log x+x\log\log x+x.
	\eeq
	We set $\check{t}=\left\{\f{A^2}{\ueps}\gamma(\f{3\ueps A}{\nu A^3})\right\}^\f13$.
	For
	$
	T\le t\le \nu^{-\f13},
	$
	we have
	$
	\log\frac{\ueps t^3}{A^2}
	-\frac13\nu t^3
	>
	K^{-2/3}.
	$
	Therefore, for $
	T\le t\le \nu^{-\f13},
	$
	\[\ba
	&\int^{\check{t}}_t\f{A}{\tau^2}\left\{\log\frac{\ueps \tau^3}{A^2}
	-\frac13\nu\tau^3\right\}d\tau\leq\log\f{ B(\infty)}{B(t)}=\int_t^\infty\f{A}{\tau^2}\max\left\{
	\log\frac{\ueps \tau^3}{A^2}
	-\frac13\nu\tau^3,
	K^{-2/3}
	\right\}
	d\tau\\
	\leq&\int^{\check{t}}_t\f{A}{\tau^2}\left\{\log\frac{\ueps \tau^3}{A^2}
	-\frac13\nu\tau^3\right\}d\tau+\int_{\nu^{-\f13}}^\infty\f{AK^{-\f23}}{\tau^2}d\tau=\int^{\check{t}}_t\f{A}{\tau^2}\left\{\log\frac{\ueps \tau^3}{A^2}
	-\frac13\nu\tau^3\right\}d\tau+\f{\nu^\f13A}{K^{\f23}}.	
	\ea\]
	Note that $\nu A^3=1.01^3K^2\rr{1+\mathcal{O}(\f1K)}$. Therefore,
	\[\ba
	&\log\f{ B(\infty)}{B(t)}+\mathcal{O}(1)=\int^{\check{t}}_t\f{A}{\tau^2}\left\{\log\frac{\ueps \tau^3}{A^2}
	-\frac13\nu\tau^3\right\}d\tau=\left.\left\{\f{3A}{\tau}\rr{\log\f{A^{\f23}}{\ueps^\f13\tau}-1 }-\f{\nu A\tau^2}{6}\right\}\right|_{t}^{\check{t}}\\
	&=-3\rr{\ueps A}^\f13\gamma(\f{3\ueps A}{\nu A^3})^{-\f13}\\
	&\times\left[1+\f12\log\gamma(\f{3\ueps A}{\nu A^3})\right]-\f{3A}{t}\rr{\log\f{A^{\f23}}{\ueps^\f13t}-1}+\f{\nu At^2}{6}=-\f{3A}{t}\rr{\log\f{A^{\f23}}{\ueps^\f13t}-1}+\f{\nu At^2}{6}+\mathcal{O}(K^\f23\log^\f23K).
	\ea\]
	Thus we have
	\[\ba
	\log\f{ B(\infty)}{B(T)}=-\log \eps A^{N_1}+\f{\nu A^3}{3K}=\rr{3+\log 1.01}K+\mathcal{O}(K^\f23\log^\f23K).
	\ea\]
	Note that $T=\f{A}{K+\f12}\rr{1+\mathcal{O}(\nu^\f23K^{-\f13})}$. Thus 
	\[\log\eps A^{N_1}=-(3+\log1.01-\f{1.01^3}{3})K+\mathcal{O}(K^\f23\log^\f23K).\]
	Thus we have $A^{-N_1}e^{-3K}<\eps<A^{-N_1}e^{-1.5K}$. Similarly we have
	\[\log\f{B(\infty)}{ B(c\nu^{-\f13})}=\f{1.01}{c}K^\f23\log K+\mathcal{O}(K^\f23\log^\f23K)\leq\f{1.02}{c}K^\f23\log K.\]
	This proves \eqref{cnu13}.
\end{proof}
\begin{proof}[Proof of (2).]
We compare increments of $\log B$ with the logarithms of the discrete transfer factors, then compare $B(T)$ with the actual initial packet. The critical-time lower bound and the global upper bound follow from the homogeneous estimates in Section~3.

	We split the proof into three steps.
	
	\underline{Step 1.} We prove for $\llcorner(D_3\nu A^3)^{\f13}\lrcorner+1\le k\le K$,
	\beq\label{200}\left\{\prod_{i=k+1}^{K}\f{\ueps A}{i^3}e^{-\f{\nu A^3}{3i^3}}\right\}B(T)\sim B(t_k),\eeq
	and for $1\leq k\leq K$,
	\beq\label{300}\left\{\prod_{i=k+1}^{K}\max\{\f{\ueps A}{i^3}e^{-\f{\nu A^3}{3i^3}},1\}\right\}B(T)\ls B(t_k).\eeq
	For $\llcorner(D_3\nu A^3)^{\f13}\lrcorner+1\le i\le K$,
	\[\ba
	&\log B(t_{i-1})-\log B(t_i)=\left.\left\{\f{3A}{\tau}\rr{\log\f{{A}^\f23}{\ueps^\f13\tau}-1}-\f{\nu A\tau^2}{6}\right\}\right|_{t_i}^{t_{i-1}}
	=3(i-\f12)\rr{\log\f{i-\f12}{(\ueps A)^\f13}-1}\\
	&-3(i+\f12)\rr{\log\f{i+\f12}{(\ueps A)^\f13}-1}-\f{\nu A^3}{6(i-\f12)^2}+\f{\nu A^3}{6(i+\f12)^2}+\mathcal{O}(\f{\nu AK^2}{i^3})
	=\log\f{\ueps A}{i^3}+\mathcal{O}(\f1{i^2})-\f{\nu{A}^3}{3i^3}+\mathcal{O}(\f{\nu A^3}{i^5}).
	\ea\]
	Thus for $\llcorner(D_3\nu A^3)^{\f13}\lrcorner+1\le k\le K$
	\[\log B(t_{k})-\log B(T)=\sum_{i=k+1}^{K}\left\{\log B(t_{i-1})-\log B(t_i)\right\}=\sum_{i=k+1}^{K}\left\{\log\f{\ueps A}{i^3}-\f{\nu{A}^3}{3i^3}\right\}+\mathcal{O}(\f1{(\nu A^3)^\f13}).\]
	This proves \eqref{200}. We next prove \eqref{300}. For $\llcorner(\nu A^3)^{\f14}\lrcorner+1\le i\le K$,
	\[
	\log B(t_{i-1})-\log B(t_i)\geq\left.\left\{\f{3A}{\tau}\rr{\log\f{{A}^\f23}{\ueps^\f13\tau}-1}-\f{\nu A\tau^2}{6}\right\}\right|_{t_i}^{t_{i-1}}=\log\f{\ueps A}{i^3}+\mathcal{O}(\f1{i^2})-\f{\nu{A}^3}{3i^3}+\mathcal{O}(\f{\nu A^3}{i^5}).\]
	For $i\leq\llcorner(\nu A^3)^{\f14}\lrcorner$, $\f{\ueps A}{i^3}e^{-\f{\nu A^3}{3i^3}}\leq{\ueps A}e^{-\f{\nu A^3}{3(\nu A^3)^{\f34}}}\ls K^3e^{-\f13K^\f12}\ll1$. Thus
	\[\ba
	&\prod_{i=k+1}^{K}\f{\ueps A}{i^3}e^{-\f{\nu A^3}{3i^3}}\leq\prod_{i=\max\{k,\llcorner(\nu A^3)^{\f14}\lrcorner\}+1}^{K}\f{\ueps A}{i^3}e^{-\f{\nu A^3}{3i^3}}=\exp\left\{\sum_{i=\max\{k,\llcorner(\nu A^3)^{\f14}\lrcorner\}+1}^{K}\log\f{\ueps A}{i^3}-\f{\nu{A}^3}{3i^3}\right\}\\
	&\ls\exp\left\{\log B(t_{\max\{k,\llcorner(\nu A^3)^{\f14}\lrcorner\}})-\log B(T)\right\}\leq\f{B(t_k)}{B(T)}.
	\ea\]
	By the monotonicity of $B(t_k)$ we get \eqref{300}.
	
	\underline{Step 2.} We prove 
	\beq\label{B(T)}\nnr{\wh{\th^{(1,1)}}(T,K,\xi)}_{L^1_\xi}\sim B(T),\quad \nnr{\fM(T,\na)\th^{(1,1)}(T)}_{L^2_{z,v}}e^{\rm\nu^\f13T}\ls B(T).\eeq
	For $|\xi-A|\leq1$, we have
	\[\ba
	&-\log S(T,0,K,\xi)=\nu\int_0^Te^{3\nu s}\nr{\xi-K\f{1-e^{-\nu s}}{\nu}}^2e^{-\nu s}ds\leq\nu e^{3\nu T}\int_0^{\f{A}{K+\f12}}\nr{\xi-Ks}^2ds\\
	&\leq\nu e^{3\nu T}\int_0^{\f{A}{K+\f12}}\nr{A+1-Ks}^2ds=\nu e^{3\nu T}\left\{\f{(A+1)^3}{3K}-\f{(\f{\f12A}{K+\f12}+1)^3}{3K}\right\}\leq\f{\nu A^3}{3K}+\nu^\f13K^\f13.
	\ea\]
	\[\ba
	&-\log S(T,0,K,\xi)\geq\nu\int_0^{\f{A}{K+\f12}}\nr{\xi-Ks}^2ds\geq\nu\int_0^{\f{A}{K+\f12}}\nr{A-1-Ks}^2ds\\
	&=\nu\left\{\f{(A-1)^3}{3K}-\f{(\f{\f12A}{K+\f12}-1)^3}{3K}\right\}\geq\f{\nu A^3}{3K}-\f23K^{-2}.
	\ea\]
	Note that $K^{-2}\gg\nu^\f13K^\f13$, $\log\fM(T,K,\xi)\ls\rm\nu AT\ls\rm\nu^\f13K^\f13$, where we use \eqref{Mxj}. By \eqref{fonemp} we have
	\[\ba
	&\nnr{\wh{\th^{(1,1)}}(T,K,\xi)}_{L^1_\xi}\sim\eps\nnr{\wh{\vf}}_{L^1_\xi}e^{-\f{\nu A^3}{3K}}= B(T),\\
	&\quad\nnr{\fM(T,\na)\th^{(1,1)}(T)}_{L^2_{z,v}}e^{\rm\nu^\f13T}\leq\eps\nnr{\wh{\vf}}_{L^2_\xi}e^{-\f{\nu A^3}{3K}+K^{-\f13}}\ls B(T).
	\ea\]
	
	\underline{Step 3.} We prove \eqref{lowerb} and \eqref{h111}. We first prove \eqref{lowerb}. For
	$\llcorner(D_3\nu A^3)^{\f13}\lrcorner+1\le k\le K$,
	\[\ba
	&\nnr{\wh{\th^{(1,1)}}(t_k,k,\xi)\vv{1}_{|\xi-A|\leq20D_2K}}_{L^2_\xi}\geq\f{1}{\rr{40D_2K}^\f12}\nnr{\wh{\th^{(1,1)}}(t_k,k,\xi)\vv{1}_{|\xi-A|\leq20D_2K}}_{L^1_\xi}\\
	&\geq\f{1}{\rr{40D_2K}^\f12}\left\{\nnr{\wh{\th^{(1,1)}}(t_k,k,\xi)}_{L^1_\xi}-\nnr{\wh{\th^{(1,1)}}(t_k,k,\xi)\vv{1}_{|\xi-A|\geq20D_2K}}_{L^1_\xi}\right\}\\
	&\geq\f{1}{\rr{40D_2K}^\f12}\left\{\nnr{\wh{\th^{(1,1)}}(t_k,k,\xi)}_{L^1_\xi}-e^{-18D_1K}\nnr{e^{\tilde{r}_1\<\na_v-A\>}\fM_2(t_k,\na)P_{\neq}\th^{(1,1)}(t_k)}_{L^2_{z,v}}\right\}\\
	&\geq\f{1}{\rr{40D_2K}^\f12}\left\{\nnr{\wh{\th^{(1,1)}}(t_k,k,\xi)}_{L^1_\xi}-e^{-5D_1K}\right\}\gs\f{1}{K^\f12}B(t_k).
	\ea\]
	Recall $D_2=\f{2D_1}{r_1}$. In the fifth step we use \eqref{hnqenergy}. Note that by \eqref{hklow}, \eqref{200} and \eqref{B(T)} we have for
	$\llcorner(D_3\nu A^3)^{\f13}\lrcorner+1\le k\le K$, $\nnr{\wh{\th^{(1,1)}}(t_k,k,\xi)}_{L^1_\xi}\gs B(t_k)$. We next prove \eqref{h111}. By \eqref{hkupp}, Lemma \ref{lA3}, \eqref{B(T)} and \eqref{300} we have for $1\leq k\leq K$,
	\[\nnr{\sF\rr{\fM(t_{k},\na)\th}(t_{k},k,\xi)}_{L^2_\xi}e^{\rm\nu^\f13t_{k}}\ls\nu^{-\f13} B(t_k).\]
	This, together with \eqref{h11upper}, yields
	\[\ba
	&\nnr{\fM_4(\tau,\na)\rho^{(1,1)}(\tau)e^{\rm\nu^\f13\tau}}_{L^2([T,t], L^2_{z})}+\left\|\fM_{2}(t,\na)\rr{e^{-|\alpha|\nu t}v^\alpha P_{\neq}\th^{(1,1)}(t)} \right\|_{\lta L^2_{z,v}}e^{\rm\nu^\f13t}\\
	&+\nnr{v^\alpha\<\na_v\>\fM(t,0,\emnt\na_v)P_{0}h^{(1,1)}(t)}_{\lta L^2_{v}}e^{\nu t}
	\ls \nu^{-\f13}{A}^4KB(t).
	\ea\]
	Note that by the same proof of \eqref{B(T)} we have $\left\|\fM_{2}(T,\na)\rr{e^{-|\alpha|\nu T}v^\alpha P_{\neq}\th^{(1,1)}(T)} \right\|_{\lta L^2_{z,v}}e^{\rm\nu^\f13T}\ls A^2B(T)$.
\end{proof}

\begin{proof}[Proof of (3).]
We compare the logarithmic growth rate of a possible source cascade with that of $B^n$. A discrete-to-continuous comparison then bounds the product of forced transfer factors by the corresponding envelope ratio. This is the comparison needed when the order-$n$ source generates further echoes.

	We split the proof into three steps.
	
	\underline{Step 1.} For $n, q, A'$ satisfying the desired conditions, we define
	\[\log E(t)=\int^{t}_T\f{A'}{\tau^2}\left\{\log\frac{\ueps \tau^3}{{A'}^2}
	-\frac13\nu\tau^3\right\}d\tau.\]
	We next prove 
	\beq\label{EBn}\f{d}{dt}\log E(t)\leq n\f{d}{dt}\log B(t),\eeq
	namely $\f{A'}{t^2}\log\{\frac{\ueps t^3}{{A'}^2}e^{
	-\frac13\nu t^3}\}\leq\f{nA}{t^2}\max\{\log\{\frac{\ueps t^3}{{A}^2}e^{
	-\frac13\nu t^3}\},K^{-\f23}\}$. Define $\bar{t}>\nu^{-\f13}$ by $\log\frac{\ueps \bar{t}^3}{{A}^2}
	-\frac13\nu \bar{t}^3=K^{-\f23}$. Set $\kappa=\f{40D_2 K}{A}$. Note that $G_t(A'):=\f{A'}{t^2}\log\{\frac{\ueps t^3}{{A'}^2}e^{
		-\frac13\nu t^3}\}$ reaches its maximum $\f{2\sqrt{H(t)}}{et^2}$ at $A'=\f{\sqrt{H(t)}}{e}$, where $H(t)=\ueps t^3e^{
		-\frac13\nu t^3}$. $H(t)$ reaches its maximum $\f{3\ueps}{e\nu}$ at $t=\rr{\f3\nu}^\f13$. Note that $A^2e^{K^{-\f23}}\leq H(t)\leq\f{3\ueps}{e\nu}$ on $T\leq t\leq\bar{t}$, $H(\bar{t})=A^2e^{K^{-\f23}}<H(T)=\ueps T^3e^{-\f{\nu T^3}{3}}$. We split into two cases according to the range of $n$.
	
	\uwave{Case 1: $n\leq\f1\kappa e^{\f12K^{-\f23}-1}$.} We split into the following subcases according to the range of $t$.
	\begin{enumerate}
		\item Subcase 1.1: $T\leq t\leq\bar{t}$, $\f{\sqrt{H(t)}}{e}\leq A(1-n\kappa)$. In this case
		\[\sup_{\substack{|A'-qA|\leq nA\kappa,\\q=0,1,\dots,n}}G_t(A')\le\max\{G_t(nA\kappa),G_t(A(1-n\kappa))\}.\]
		We first prove $G_t(nA\kappa)\leq n\f{d}{dt}\log B(t)$, which is equivalent to $(1-\kappa)\log\frac{\ueps t^3}{{A}^2}e^{
			-\frac13\nu t^3}\geq2\kappa\log\f{1}{n\kappa}$. The latter follows from
		\[2\kappa\log\f{1}{n\kappa}\leq2\kappa\log\f{1}{2\kappa}\leq\f{1-\kappa}{K^\f23}\leq(1-\kappa)\log\frac{\ueps t^3}{{A}^2}e^{
			-\frac13\nu t^3}.\]
		We next prove $G_t(A(1-n\kappa))\leq n\f{d}{dt}\log B(t)$, which is equivalent to $(n-1+n\kappa)\log\f{\ueps t^3}{A^2}e^{-\f13\nu t^3}\geq-2(1-n\kappa)\log(1-n\kappa)$. The latter follows from
		\[(n-1+n\kappa)\log\f{\ueps t^3}{A^2}e^{-\f13\nu t^3}\geq\f{n}{2}K^{-\f23}\geq2n\kappa\geq-2(1-n\kappa)\log(1-n\kappa).\]
		
		\item Subcase 1.2: $T\leq t\leq\bar{t}$, $A(1-n\kappa)\leq \f{\sqrt{H(t)}}{e}\leq nA(1+\kappa)$. In this case
		\[\sup_{\substack{|A'-qA|\leq nA\kappa,\\q=0,1,\dots,n}}G_t(A')\le\f{2\sqrt{H(t)}}{et^2}.\]
		$\f{2\sqrt{H(t)}}{et^2}\leq n\f{d}{dt}\log B(t)$ is equivalent to $\f{1}{ne}\leq\f{A}{\sqrt{H(t)}}\log\f{\sqrt{H(t)}}{A}$. We only need to prove that $\f{1}{ne}\leq\f{1}{e(1-n\kappa)}\log e(1-n\kappa)$ and $\f{1}{ne}\leq\f{1}{en(1+\kappa)}
		\log en(1+\kappa)$, which follow directly.
		
		\item Subcase 1.3: $T\leq t\leq\bar{t}$, $\f{\sqrt{H(t)}}{e}\geq nA(1+\kappa)$. In this case
		\[\sup_{\substack{|A'-qA|\leq nA\kappa,\\q=0,1,\dots,n}}G_t(A')\le G_t(nA(1+\kappa)).\]
		$G_t(nA(1+\kappa))\leq n\f{d}{dt}\log B(t)$ is equivalent to $\kappa\log\f{H(t)}{A^2}\leq2(1+\kappa)\log n(1+\kappa)$. The latter follows from
		\[\kappa\log\f{H(t)}{A^2}\leq\kappa\log\f{3\ueps}{e\nu A^2}<2(1+\kappa)\log 2(1+\kappa)\leq2(1+\kappa)\log n(1+\kappa).\]
		
		\item Subcase 1.4: $t\geq\bar{t}$, $\f{\sqrt{H(t)}}{e}\geq nA\kappa$. In this case 
		\[\sup_{\substack{|A'-qA|\leq nA\kappa,\\q=0,1,\dots,n}}G_t(A')\le\max\{G_t(nA\kappa),G_t(A(1-n\kappa))\}.\]
		$G_t(nA\kappa)\leq\f{nA}{t^2K^\f23}$ is equivalent to $\kappa\rr{\log\f{H(t)}{A^2}-2\log n\kappa}\leq\f{1}{K^\f23}$. The latter follows from
		\[\kappa\log\f{H(t)}{A^2}+2\kappa\log\f1{n\kappa}\leq\f{\kappa}{K^\f23}+2\kappa\log\f{1}{2\kappa}<\f{1}{K^\f23}.\]
		$G_t(A(1-n\kappa))\leq\f{nA}{t^2K^\f23}$ is equivalent to $(1-n\kappa)\rr{\log\f{H(t)}{A^2}-2\log (1-n\kappa)}\leq\f{n}{K^\f23}$. The latter follows from
		\[(1-n\kappa)\left\{\log\f{H(t)}{A^2}-2\log (1-n\kappa)\right\}\leq\f{1-n\kappa}{K^\f23}-2(1-n\kappa)\log (1-n\kappa)\leq\f{n}{K^\f23},\]
		which follows from $2(1-n\kappa)\log \f1{1-n\kappa}\leq 2n\kappa\leq\f{n-1+n\kappa}{K^\f23}$.
		
		\item Subcase 1.5: $t\geq\bar{t}$, $\f{\sqrt{H(t)}}{e}\leq nA\kappa$. In this case 
		\[\sup_{\substack{|A'-qA|\leq nA\kappa,\\q=0,1,\dots,n}}G_t(A')\le\f{2\sqrt{H(t)}}{et^2}.\]
		$\f{2\sqrt{H(t)}}{et^2}\leq\f{nA}{t^2K^\f23}$ is equivalent to $\sqrt{H(t)}\leq\f{eAn}{2K^\f23}$, which follows directly.
	\end{enumerate}
	
	\uwave{Case 2: $n\geq\f1\kappa e^{\f12K^{-\f23}-1}$.} In this case we have
	\[G_t(A')\leq\f{2\sqrt{H(t)}}{et^2}\leq\f{2}{et^2}\sqrt{\f{3\ueps}{e\nu}}<\f{nA}{t^2K^\f23}\leq n\f{d}{dt}\log B(t).\]
	
	\underline{Step 2.} We consider only the nontrivial case $k_2+1<k_1$. We prove 
	\beq\label{100}\prod_{i=\max\{k_2+1,\kAt\}}^{\min\{k_1-1,\kAo\}}
	\frac{\ueps A'}{i^{3}}
	e^{-\frac{\nu A'^3}{3i^3}}
	\ls
	K^4\f{E(t'_{\max\{k_2+1,\kAt\}})}{E(t'_{\min\{k_1-1,\kAo\}})}.\eeq
	For $\max\{k_2+1,\kAt\}+1\leq i\leq\min\{k_1-1,\kAo\}$, 
	\[\ba
	&\log E(t'_{i-1})-\log E(t'_i)=\left.\left\{\f{3A'}{\tau}\rr{\log\f{{A'}^\f23}{\ueps^\f13\tau}-1}-\f{\nu A'\tau^2}{6}\right\}\right|_{t'_i}^{t'_{i-1}}=\log\f{\ueps A'}{i^3}+\mathcal{O}(\f1{i^2})-\f{\nu{A'}^3i}{3(i^2-\f14)^2}\\
	&\geq\log\f{\ueps A'}{(i-1)^3}+\log\f{(i-1)^3}{i^3}-\f{\nu{A'}^3}{3(i-1)^3}+\mathcal{O}(\f1{i^2}).
	\ea\]
	Therefore,
	\[\ba
	&\log E(t'_{\max\{k_2+1,\kAt\}})-\log E(t'_{\min\{k_1-1,\kAo\}})\geq\sum_{i=\max\{k_2+1,\kAt\}}^{\min\{k_1-1,\kAo\}-1}\log\f{\ueps A'}{i^3}-\f{\nu{A'}^3}{3i^3}+3\log\f{\max\{k_2+1,\kAt\}}{\min\{k_1-1,\kAo\}}+\mathcal{O}(1)\\
	&\geq\sum_{i=\max\{k_2+1,\kAt\}}^{\min\{k_1-1,\kAo\}}\log\f{\ueps A'}{i^3}-\f{\nu{A'}^3}{3i^3}-3\log\min\{k_1-1,\kAo\}-\log\f{3\ueps}{\nu {A'}^2}+\mathcal{O}(1),
	\ea\]
	where we use $\log\f{\ueps A'}{i^3}-\f{\nu{A'}^3}{3i^3}\leq\log\f{3\ueps}{\nu{A'}^2}-1$. Note that $\min\{k_1-1,\kAo\}\ls\min\{\f{A'}{T},\rr{\ueps A'}^\f13\}$. This, together with $3\log\min\{\f{A'}{T},\rr{\ueps A'}^\f13\}+\log\f{3\ueps}{\nu {A'}^2}\leq\log K^4+\mathcal{O}(1)$ yields \eqref{100}.
	
	\underline{Step 3.} We prove \eqref{pbqhb}.
	\[\ba
	&\prod_{i=\max\{k_2+1,\kAt\}}^{\min\{k_1-1,\kAo\}}
	\frac{\ueps A'}{i^{3}}
	e^{-\frac{\nu A'^3}{3i^3}}	\ls K^4\exp\left\{\log E(t'_{\max\{k_2+1,\kAt\}})-\log E(t'_{\min\{k_1-1,\kAo\}})\right\}\\
	\leq& K^4\exp\left\{n\log B(t'_{\max\{k_2+1,\kAt\}})-n\log B(t'_{\min\{k_1-1,\kAo\}})\right\}
	\leq K^4\exp\left\{n\log B(t'_{k_2})-n\log B(t'_{k_1-1})\right\}.
	\ea\] 
	The first inequality uses \eqref{100}. The second inequality uses \eqref{EBn}. The third inequality uses the monotonicity of $\log B(t)$. By the monotonicity of $B(t)$ we get \eqref{pbqhb}.
\end{proof}
	
\subsection{Estimates for $h^{(n,q)}$, $n\geq2$.}
The envelope comparison \eqref{pbqhb} allows us to propagate the factor $B^n$ through the forced response. We prove the next proposition by induction on $n$. The rough localization estimate from Section~2 controls the source outside its principal frequency region. Inside that region, a partition into packets permits the use of the three regimes from Section~4, while the zero spatial mode is treated by the homogeneous semigroup. The polynomial powers of $A$ absorb the derivative losses and the number of packets.

\begin{prop}[Bounds for the nonlinear expansion coefficients]\label{prop:nonlinear-coefficients}
	For $1\leq n\leq A^{N'}$, $|q|\leq n$, $2\mid n-q$, $t\geq T$, we have
	\beq\label{hnqu}\ba
	&\nnr{\fM_4(\tau,\na)\rho^{(n,q)}(\tau)e^{\rm\nu^\f13\tau}}_{L^2([T,t], L^2_{z})}+\left\|\fM_{2}(t,\na)\rr{e^{-|\alpha|\nu t}v^\alpha P_{\neq}\th^{(n,q)}(t)} \right\|_{\lta L^2_{z,v}}e^{\rm\nu^\f13t}\\
	&+\nnr{v^\alpha\<\na_v\>\fM(t,0,\emnt\na_v)P_{0}h^{(n,q)}(t)}_{\lta L^2_{v}}e^{\nu t}
	\leq {A}^{10N'(2n-1)}B^n(t).\ea
	\eeq
\end{prop}
\begin{proof}
The induction uses two different bounds for the source. The coarse exponential localization estimate makes the Fourier tail negligible, while the inductive size estimate controls the localized source by $B^n$. The latter must be passed through the forced linear analysis without losing its envelope factor.

	The base case $n=1$ follows from \eqref{h111}. Suppose \eqref{hnqu} holds through order $n-1$. We prove it at order $n$. 
	We set 
	\[\chi(k,\xi)=\vv{1}_{|k-qK|\leq40nD_2K}\int_{|\eta-qA|\leq40nD_2K-\f12}2\wh{\vf}(2(\xi-\eta))d\eta.\]
	Write $h^{(n,q)}=h_1+h_2+h_3$, where the three components have zero initial data at $T$ and satisfy
	\[\ba&\pa_th_1=\cL h_1+Q_1,\quad\tQ_1=\sum_{\substack{n_1+n_2=n,\\q_1+q_2=q}}(1-\chi(\na))P_{\neq}\cN^\dagger(\th^{(n_1,q_1)},\th^{(n_2,q_2)}),\\
	&\pa_th_2=\cL h_2+Q_2,\quad \tQ_2=\sum_{\substack{n_1+n_2=n,\\q_1+q_2=q}}\chi(\na)P_{\neq}\cN^\dagger(\th^{(n_1,q_1)},\th^{(n_2,q_2)}),\\
	&\pa_th_3=\cL h_3+\sum_{\substack{n_1+n_2=n,\\q_1+q_2=q}}P_{0}\cN(h^{(n_1,q_1)},h^{(n_2,q_2)}).
	\ea\]
	
	\underline{Estimate of $h_1$.}	By \eqref{ucm} we have
	\[\ba
	&\int_T^t\left\|e^{\f{3r_1}{4}\<\na_z-qK,\na_v-qA\>}\fA\fM_4(\tau,\na)\rr{e^{-|\alpha|\nu \tau}v^\alpha \sum_{\substack{n_1+n_2=n,\\q_1+q_2=q}}P_{\neq}\cN^\dagger(\th^{(n_1,q_1)},\th^{(n_2,q_2)})(\tau)} \right\|_{\lta L^2_{z,v}}\\
	&\times e^{\rm\nu^\f13\tau+\f12\nu\tau}d\tau\ls n^2A^{4N'+4}e^{10(2n-2)D_1K}\leq A^{6N'+4}e^{10(2n-2)D_1K}.
	\ea\]
	Note that $\mathrm{supp}\chi\subset\{(k,\xi):|k-qK|\leq40nD_2K, |\xi-qA|\leq40nD_2K\}$, $\mathrm{supp}(1-\chi)\subset\{(k,\xi):|k-qK|\geq40nD_2K\ \text{or}\ |\xi-qA|\geq40nD_2K-1\}$. Thus
	\[\ba
	&\int_T^t\left\|\fA\fM_4(\tau,\na)\rr{e^{-|\alpha|\nu \tau}v^\alpha \sum_{\substack{n_1+n_2=n,\\q_1+q_2=q}}\right.\right.\\
	&\quad \left.\left.(1-\chi(\na))P_{\neq}\cN^\dagger(\th^{(n_1,q_1)},\th^{(n_2,q_2)})(\tau)} \right\|_{\lta L^2_{z,v}}e^{\rm\nu^\f13\tau+\f12\nu\tau}d\tau\leq e^{-10nD_1K}.
	\ea\]
	By \eqref{nyb} we have
	\beq\label{weq}\ba&\nnr{\fM_4(\tau,\na)\rho_1(\tau)e^{\rm\nu^\f13\tau}}_{L^2([T,t], L^2_{z})}+\left\|\fM_{2}(t,\na)\rr{e^{-|\alpha|\nu t}v^\alpha P_{\neq}\th_1(t)} \right\|_{\lta L^2_{z,v}}e^{\rm\nu^\f13t}\\
	&+\nnr{v^\alpha\<\na_v\>\fM(t,0,\emnt\na_v)P_{0}h_1(t)}_{\lta L^2_{v}}e^{\nu t}
	\ls e^{-10nD_1K}.\ea\eeq
	
	\underline{Estimate of $h_2$.} We first estimate the source term. By a variant of \eqref{cNk2} we have
	\[\ba
	&\int_T^t\left\|\fM(\tau,\na)\rr{e^{-|\alpha|\nu \tau}v^\alpha \sum_{\substack{n_1+n_2=n,\\q_1+q_2=q}}\chi(\na)P_{\neq}\cN^\dagger(\th^{(n_1,q_1)},\th^{(n_2,q_2)})(\tau)} \right\|_{\lta L^2_{z,v}}e^{\rm\nu^\f13\tau+\f12\nu\tau}d\tau\\
	\ls&\sum_{\substack{n_1+n_2=n,\\q_1+q_2=q}}\nnr{\fM_2(\tau,\na)\rho^{(n_1,q_1)}(\tau)e^{\rm\nu^\f13\tau}}_{L^2([T,t], L^2_{z})}
	\sup_{\tau\in[T,t]}\left\{\nnr{\fM_{1}(\tau,\na)\rr{e^{-|\alpha|\nu\tau}v^\alpha P_{\neq}\th^{(n_2,q_2)} }}_{\lta L^2_{z,v}}e^{2\nu\tau}\right.\\
	&\left.+\nnr{v^\alpha\<\na_v\>\fM(\tau,0,\emntau\na_v)P_{0}h^{(n_2,q_2)} }_{\lta L^2_{v}}e^{\f12\nu\tau}\right\}\ls n^2A^{10N'(2n-2)}B(t)^n\leq A^{N'(20n-18)}B(t)^n.
	\ea\]
	We perform a smooth partition of unity within the set $\{(k,\xi):|k-qK|\leq40nD_2K, |\xi-qA|\leq40nD_2K\}$. Each part of the partition has its support in $\{K'\}\times[A'-1,A'+1]$, where $(K',A')\in\Z^{*}\times\R$ satisfies $|k'-qK|\leq40nD_2K$, $|A'-qA|\leq40nD_2K$. The number of parts is less than $10^4n^2D_2^2K^2$. We reduce the problem to solving the equation
	\[\pa_th=\mathcal{L} h+Q,\quad\mathrm{supp}\wh{\tQ}\subset\{K'\}\times[A'-1,A'+1],\]
	\[\int^t_T\nnr{\fM(\tau,\na)\rr{e^{-\nu|\alpha|\tau}v^\alpha \tQ(\tau)}}_{\lta L^2_{z,v}}e^{\rm\nu^\f13\tau+\f12\nu\tau} d\tau\leq A^{N'(20n-18)}B(t)^n.\]
	By symmetry we only need to consider the $A'\geq0$ case. We apply the three regimes from Section~4.
	
	\uwave{Case 1: $K'<0$ or $K'\f{1-e^{-\nu T}}{\nu}\geq A'+1$.} In this case by \eqref{Case1} we have
	\[\left\|\fM_{2}(t,\na)\rr{e^{-|\alpha|\nu t}v^\alpha P_{\neq}\th(t)} \right\|_{\lta L^2_{z,v}}e^{\rm\nu^\f13t}\ls(nA)^2A^{N'(20n-18)}B(t)^n, \quad\rho(t)=0,\ P_0h(t)=0.\]
	Note that $\<K',A'\>\ls nA$.
	
	\uwave{Case 2: $K'>0$ and $K'\f{1-e^{-\nu T}}{\nu}\leq A'+1$ and $A'\leq\f13\nu^{-\f12}$.} By \eqref{h39}, \eqref{pbqhb} and Lemma \ref{lA3} we have
	\[\ba
	&\nnr{\sF\rr{\fM(t_k',\na)\th}(t'_k,k,\xi)}_{L^2_\xi}e^{\rm\nu^\f13t_k'}\ls\prod_{i=\max\{k+1,\kAt\}}^{\min\{K'-1,\kAo\}}\f{\ueps A'}{i^3}e^{-\f{\nu {A'}^3}{3i^3}+2\rm\nu^\f13(t'_{i-1}-t'_i)}{A'}^2A^{N'(20n-18)}B(t'_{K'-1})^n\\
	&+\ueps^{0.7}{A'}^2e^{-3D_1\rr{\ueps A'}^\f13}A^{N'(20n-18)}B(t'_k)^n\ls\nu^{-\f13} {A'}^2K^4A^{N'(20n-18)}B(t'_k)^n.
	\ea\]
	Thus by \eqref{hrhoh} we have
	\[\ba
	&\nnr{\fM_4(\tau,\na)\rho(\tau)e^{\rm\nu^\f13\tau}}_{L^2([T,t], L^2_{z})}+\left\|\fM_{2}(t,\na)\rr{e^{-|\alpha|\nu t}v^\alpha P_{\neq}\th(t)} \right\|_{\lta L^2_{z,v}}e^{\rm\nu^\f13t}\\
	&+\nnr{v^\alpha\<\na_v\>\fM(t,0,\emnt\na_v)P_{0}h(t)}_{\lta L^2_{v}} e^{\nu t}\ls\nu^{-\f13}{A'}^6K'K^4A^{N'(20n-18)}B(t)^n\ls A^{N'(20n-14)+4}B(t)^n.
	\ea\]
	Note that $\nu^{-\f13}{A'}^6K'K^4\ls\nu^{-\f13}\f{{A'}^7}{T}K^4\ls\nu^{-\f13}\f{\nu^{-\f72}K^4}{T}\ls\nu^{-\f13} A^{\f{19}{2}}K^{-2}\ls A^{4(N'+1)}$.
	
	\uwave{Case 3: $K'>0$ and $K'\f{1-e^{-\nu T}}{\nu}\leq A'+1$ and $A'\geq\f13\nu^{-\f12}$.} By \eqref{hrho3} we have
	\[\ba
	&\nnr{\fM_4(\tau,\na)\rho(\tau)e^{\rm\nu^\f13\tau}}_{L^2([T,t], L^2_{z})}+\left\|\fM_{2}(t,\na)\rr{e^{-|\alpha|\nu t}v^\alpha P_{\neq}\th(t)} \right\|_{\lta L^2_{z,v}}e^{\rm\nu^\f13t}\\
	&+\nnr{v^\alpha\<\na_v\>\fM(t,0,\emnt\na_v)P_{0}h(t)}_{\lta L^2_{v}}e^{\nu t}
	\ls(nA)^4A^{N'(20n-18)}B(t)^n\leq A^{N'(20n-14)+4}B^n(t).
	\ea\]
	
	Combining these cases we obtain 
	\beq\label{Weq}\ba
	&\nnr{\fM_4(\tau,\na)\rho_2(\tau)e^{\rm\nu^\f13\tau}}_{L^2([T,t], L^2_{z})}+\left\|\fM_{2}(t,\na)\rr{e^{-|\alpha|\nu t}v^\alpha P_{\neq}\th_2(t)} \right\|_{\lta L^2_{z,v}}e^{\rm\nu^\f13t}+\nnr{v^\alpha\<\na_v\>\right.\\
		&\left.\fM(t,0,\emnt\na_v)P_{0}h_2(t)}_{\lta L^2_{v}}e^{\nu t}
	\ls n^2D^2_2K^2A^{N'(20n-14)+4}B^n(t)\ls A^{N'(20n-12)+4}K^2B(t)^n.
	\ea\eeq
	
	\underline{Estimate of $h_3$.} By \eqref{cN0} we have
	\[\ba&\int_T^t\nnr{v^\alpha\<\na_v\>\fM(\tau,0,\emntau\na_v)\sum_{\substack{n_1+n_2=n,\\q_1+q_2=q}}P_{0}\cN(h^{(n_1,q_1)},h^{(n_2,q_2)})}_{\lta L^2_{v}}e^{\nu\tau}d\tau\ls\sum_{\substack{n_1+n_2=n,\\q_1+q_2=q}}\nnr{\fM_{3}(\tau,\na)\rho^{(n_1,q_1)}(\tau)\right.\\
	&\left.\times e^{(\sigma+3)\nu\tau}}_{L^2([T,t],L^2_z)}\sup_{\tau\in[T,t]}\nnr{\fM_{2}(\tau,\na)\rr{e^{-|\alpha|\nu\tau}v^\alpha P_{\neq}\th^{(n_2,q_2)} }}_{\lta L^2_{z,v}}\ls A^{N'(20n-18)}B(t)^n.\ea\]
	By \eqref{myb} we have
	\beq\label{WEq}\nnr{v^\alpha\<\na_v\>\fM(t,0,\emnt\na_v)P_{0}h_3}_{\lta L^2_{v}}e^{\nu\tau}\ls A^{N'(20n-18)}B(t)^n,\quad P_{\neq}h_3(t)=0.\eeq
	
	Combining \eqref{weq}, \eqref{Weq}, and \eqref{WEq} proves \eqref{hnqu} and closes the induction.
\end{proof}

\section{Bounds for the error terms}
This section passes from the finite expansion to an exact solution. Proposition~\ref{prop:truncation} bounds the approximate solution in the stronger norms needed for the energy argument and shows that its truncation error is extremely small. We then control the remainder in two stages. Up to the thermal time $\nu^{-1}$, Proposition~\ref{prop:remainder-thermal} uses an energy estimate whose losses are compensated by the high truncation order. After the thermal time, Proposition~\ref{prop6} combines a Volterra estimate for the density, a hypocoercive estimate for the nonzero modes, and a Fokker--Planck estimate for the zero mode. Together, these bounds justify the expansion, preserve the echo lower bound, and provide global relaxation estimates for the exact solution.

Set 
\beq\label{h34} h^\#=\sum_{n=1}^{A^{N'}}h^{(n)}, \quad \wt{h}=h-h^\#.\eeq
Then $\wt{h}$ satisfies
\beq\label{wtheq} \pa_t\wt{h}=\cL\wt{h}+\cN(\wt{h},h^\#)+\cN(h^\#,\wt{h})+\cN(\wt{h},\wt{h})+\wt{Q},\quad \wt{h}(T)=0,\eeq
where
\beq\label{wtQ}\wt{Q}=\sum_{\substack{1\leq n_1,n_2\leq A^{N'}\\n_1+n_2>A^{N'}} }\cN(h^{(n_1)},h^{(n_2)}).\eeq
We first collect bounds for $h^\#$ and $\wt Q$. The size estimates from Section~5 are upgraded to higher regularity by splitting into the localized frequency region and its exponentially small complement. Since every interaction in $\wt Q$ has total order above the truncation order, the bilinear estimates give a much smaller forcing. This smallness will compensate for the energy losses in the remainder equation.
\begin{prop}[Bounds for the approximate solution and truncation error]\label{prop:truncation}
	For $t\in[T,\infty)$, we have
	\beq\label{h3}\ba&\nnr{\fM_{14}(\tau,\na)\rho^{\#}(\tau)e^{\rm\nu^\f13\tau}}_{L^2([T,t], L^2_{z})}+\left\|\fM_{12}(t,\na)\rr{e^{-|\alpha|\nu t}v^\alpha P_{\neq}\th^{\#}(t)} \right\|_{\lta L^2_{z,v}}e^{\rm\nu^\f13t}\\
	&+\nnr{v^\alpha\<\na_v\>\fM_{10}(t,0,\emnt\na_v)P_{0}h^{\#}(t)}_{\lta L^2_{v}}e^{\nu t}\leq A^{-N_1+16},\ea\eeq
	\beq\label{Dwb}\left\|\fM_{N_1-4}(t,\na) P_{\neq}\th^{\#}(t) \right\|_{ L^2_{z,v}}e^{\rm\nu^\f13t}+\nnr{\<\na_v\>\fM_{N_1-5}(t,0,\emnt\na_v)P_{0}h^{\#}(t)}_{ L^2_{v}}e^{\nu t}\ls K^{-1},\eeq
	\beq\label{twQ}\ba&\int^t_T\left\|\fM_{10}(\tau,\na)\rr{e^{-|\alpha|\nu \tau}v^\alpha P_{\neq}\wt{\tQ}(\tau)} \right\|_{\lta L^2_{z,v}}e^{\rm\nu^\f13\tau+\f12\nu\tau}\\
	&+\nnr{v^\alpha\<\na_v\>^{11}\fM(\tau,0,\emntau\na_v)P_{0}\wt{Q}(\tau)}_{\lta L^2_{v}}e^{\nu \tau}d\tau\leq A^{-\f{N_1}{2}A^{N'}}.\ea\eeq
\end{prop}
\begin{proof}
We first upgrade the coefficient bounds to the higher norms in \eqref{h3} by using localization to pay for additional derivatives. Summing over the expansion yields the approximate-solution bounds. We then apply the bilinear estimates to the interactions above the truncation order to obtain \eqref{twQ}.

	We first prove \eqref{h3}. For $n\geq1$, by \eqref{hnqu} and \eqref{hnqenergy} we have
	\[\ba
	&\nnr{\fM_{14}(\tau,\na)\rho^{(n,q)}(\tau)e^{\rm\nu^\f13\tau}}_{L^2([T,t], L^2_{z})}\leq
	\nnr{\vv{1}_{\<\na_z-qK,\na_v-qA\>\leq 40nD_2K}\fM_{14}(\tau,\na)\rho^{(n,q)}(\tau)e^{\rm\nu^\f13\tau}}_{L^2([T,t], L^2_{z})}
	\\&+\nnr{\vv{1}_{\<\na_z-qK,\na_v-qA\>\geq 40nD_2K}\fM_{14}(\tau,\na)\rho^{(n,q)}(\tau)e^{\rm\nu^\f13\tau}}_{L^2([T,t], L^2_{z})}\\
	&\ls(nA)^{10}\nnr{\fM_{4}(\tau,\na)\rho^{(n,q)}(\tau)e^{\rm\nu^\f13\tau}}_{L^2([T,t], L^2_{z})}+e^{-10nD_1K}\ls A^{20N'n+10}B(t)^n\leq A^{-{N_1}n+20N'n+10}.
	\ea\]
	Similarly for $n=1$, by \eqref{h111} and \eqref{hnqenergy} we have
	\[\nnr{\fM_{14}(\tau,\na)\rho^{(1,1)}(\tau)e^{\rm\nu^\f13\tau}}_{L^2([T,t], L^2_{z})}\ls \nu^{-\f13}A^{10}A^4KB(t)\ls \nu^{-\f13}A^{-N_1+14}K.\]
	Thus by \eqref{h34} we have
	\[\ba
	&\nnr{\fM_{14}(\tau,\na)\rho^{\#}(\tau)e^{\rm\nu^\f13\tau}}_{L^2([T,t], L^2_{z})}\leq\sum_{q=\pm1}\nnr{\fM_{14}(\tau,\na)\rho^{(1,q)}(\tau)e^{\rm\nu^\f13\tau}}_{L^2([T,t], L^2_{z})}\\
	&+\sum_{n\geq2}\sum_{|q|\leq n}\nnr{\fM_{14}(\tau,\na)\rho^{(n,q)}(\tau)e^{\rm\nu^\f13\tau}}_{L^2([T,t], L^2_{z})}
	\ls \nu^{-\f13}A^{-N_1+14}K+\sum_{n\geq2}A^{-N'(30n-1)+10}\ls A^{-N_1+15}K.
	\ea\]
	Similarly we can estimate
	\[\ba
	&\left\|\fM_{12}(t,\na)\rr{e^{-|\alpha|\nu t}v^\alpha P_{\neq}\th^{(n,q)}(t)} \right\|_{\lta L^2_{z,v}}e^{\rm\nu^\f13t},\quad \left\|\fM_{12}(t,\na)\rr{e^{-|\alpha|\nu t}v^\alpha P_{\neq}\th^{\#}(t)} \right\|_{\lta L^2_{z,v}}e^{\rm\nu^\f13t},\\
	&\nnr{v^\alpha\<\na_v\>\fM_{10}(\tau,0,\emntau\na_v)P_{0}h^{(n,q)}(\tau)}_{\lta L^2_{v}}e^{\nu\tau}\quad\text{and}\quad \nnr{v^\alpha\<\na_v\>\fM_{10}(\tau,0,\emntau\na_v)P_{0}h^{\#}(\tau)}_{\lta L^2_{v}}e^{\nu \tau}.
	\ea\]
	We obtain \eqref{h3}. \eqref{Dwb} follows by the same argument as \eqref{h3}. We next prove \eqref{twQ}. For fixed $n_1, n_2, q_1, q_2$ satisfying $n_1, n_2\leq A^{N'}$, $n_1+n_2>A^{N'}$, $|q_1|\leq n_1, |q_2|\leq n_2$, by \eqref{cNk2} we have
	\[\ba
	&\int^t_T\left\|\fM_{10}(\tau,\na)\rr{e^{-|\alpha|\nu \tau}v^\alpha P_{\neq}\cN^\dagger(\th^{(n_1,q_1)},\th^{(n_2,q_2)})} \right\|_{\lta L^2_{z,v}}e^{\rm\nu^\f13\tau+\f12\nu\tau}d\tau\\
	&\ls\nnr{\fM_{12}(\tau,\na)\rho^{(n_1,q_1)}(\tau)e^{\rm\nu^\f13\tau}}_{L^2([T,t], L^2_{z})}\\
	&\times\sup_{\tau\in[T,t]}\left\{\left\|\fM_{11}(\tau,\na)\rr{e^{-|\alpha|\nu \tau}v^\alpha P_{\neq}\th^{(n_2,q_2)}(\tau)} \right\|_{\lta L^2_{z,v}}e^{2\nu \tau}\right.\\
	&\left.+\nnr{v^\alpha\<\na_v\>\fM_{10}(\tau,0,\emntau\na_v)P_{0}h^{(n_2,q_2)}(\tau)}_{\lta L^2_{v}}e^{\f12\nu\tau}\right\}\ls A^{-(N_1-20N')(n_1+n_2)+20}.
	\ea\]
	By \eqref{cN0} we have
	\[\ba
	&\int_T^t\nnr{v^\alpha\<\na_v\>^{11}\fM(\tau,0,\emntau\na_v)P_{0}\wt{Q}(\tau)}_{\lta L^2_{v}}e^{\nu \tau}d\tau\ls \nnr{\fM_{13}(\tau,\na)\rho^{(n_1,q_1)}(\tau)e^{\rm\nu^\f13\tau}}_{L^2([T,t], L^2_{z})}\\
	&\times\sup_{\tau\in[T,t]}\left\|\fM_{12}(\tau,\na)\rr{e^{-|\alpha|\nu \tau}v^\alpha P_{\neq}\th^{(n_2,q_2)}(\tau)} \right\|_{\lta L^2_{z,v}}\ls A^{-(N_1-20N')(n_1+n_2)+20}.
	\ea\]
	By \eqref{wtQ} we have
	\[\int^t_T\left\|\fM_{10}(\tau,\na)\rr{e^{-|\alpha|\nu \tau}v^\alpha P_{\neq}\wt{\tQ}(\tau)} \right\|_{\lta L^2_{z,v}}e^{\rm\nu^\f13\tau+\f12\nu\tau}d\tau\ls A^{4N'}A^{-(N_1-20N')A^{N'}+20}\ls A^{-0.55N_1A^{N'}}.\]
	This proves \eqref{twQ}.
\end{proof}

\subsection{Bounds on $[T,\nu^{-1}]$.}
On this finite interval, the shifted variables allow us to control all spatial modes in a single weighted energy estimate. Commutator bounds and velocity diffusion control the derivative terms, while a continuity argument handles the quadratic remainder. Although the resulting bound contains a large factor $e^{\nu^{-2}}$, the truncation error in \eqref{twQ} is small enough to absorb this loss. The estimate below both preserves the transient lower bound and supplies small data at the thermal time.

\begin{prop}[Remainder bounds up to the thermal time]\label{prop:remainder-thermal}
	For $t\in[T,\nu^{-1}]$, we have
	\beq\label{TT1}\ba
	&\nnr{\fM_{10}(\tau,\na)\wt{\rho}(\tau)e^{\rm\nu^\f13\tau}}_{L^2([T,t], L^2_{z})}+\left\|\fM_{10}(t,\na)\rr{e^{-|\alpha|\nu t}v^\alpha P_{\neq}\wt{\th}(t)} \right\|_{\lta L^2_{z,v}}e^{\rm\nu^\f13t}\\
	&+\nnr{v^\alpha\<\na_v\>\fM_{9}(t,0,\emnt\na_v)P_{0}\wt{h}(t)}_{\lta L^2_{v}}e^{\nu t}\leq A^{-\f {N_1}{2}A^{N'}}e^{\nu^{-2}}.
	\ea\eeq
\end{prop}
\begin{proof}
We reduce the desired density and profile bounds to a weighted energy norm of the full remainder in shifted coordinates. The multiplier derivative is controlled by velocity diffusion, and commutator estimates handle the transport terms. The small forcing and zero initial remainder then close a continuity argument up to $\nu^{-1}$.

	Sobolev embedding gives, on $[T,\nu^{-1}]$
	\beq\label{oi}\ba&\nnr{\fM_{10}(\tau,\na)\wt{\rho}(\tau)e^{\rm\nu^\f13\tau}}_{L^2_{z}}\ls\left\|\fM_{10}(\tau,\na)\rr{e^{-|\alpha|\nu \tau}v^\alpha P_{\neq}\wt{\th}(\tau)} \right\|_{\lta L^2_{z,v}}e^{\nu\tau+\rm\nu^\f13\tau}\\
	\ls& e^{\rm\nu^{-\f23}}\left\|\fM_{10}(\tau,\na)\rr{e^{-|\alpha|\nu \tau}v^\alpha P_{\neq}\wt{\th}(\tau)} \right\|_{\lta L^2_{z,v}}.\ea\eeq
	Consequently,
	\beq\label{a1}\ba
	&\nnr{\fM_{10}(\tau,\na)\wt{\rho}(\tau)e^{\rm\nu^\f13\tau}}_{L^2([T,t], L^2_{z})}\leq\nu^{-\f12}\sup_{\tau\in[T,t]}\nnr{\fM_{10}(\tau,\na)\wt{\rho}(\tau)e^{\rm\nu^\f13\tau}}_{L^2_{z}}\\
	\ls&\nu^{-\f12}e^{\rm\nu^{-\f23}}\sup_{\tau\in[T,t]}\left\|\fM_{10}(\tau,\na)\rr{e^{-|\alpha|\nu \tau}v^\alpha P_{\neq}\wt{\th}(\tau)} \right\|_{\lta L^2_{z,v}}.
	\ea\eeq
	Moreover we have
	\beq\label{a2}\left\|\fM_{10}(t,\na)\rr{e^{-|\alpha|\nu t}v^\alpha P_{\neq}\wt{\th}(t)} \right\|_{\lta L^2_{z,v}}e^{\rm\nu^\f13t}\ls e^{\rm\nu^{-\f23}}\left\|\fM_{10}(t,\na)\rr{e^{-|\alpha|\nu t}v^\alpha P_{\neq}\wt{\th}(t)} \right\|_{\lta L^2_{z,v}},\eeq
	\beq\label{a3}\ba&\nnr{v^\alpha\<\na_v\>\fM_{9}(t,0,\emnt\na_v)P_{0}\wt{h}(t)}_{\lta L^2_{v}}e^{\nu t}\ls e^{2\nu t}\left\|\fM_{10}(t,0,\emnt\na)v^\alpha P_{0}\wt{h}(t) \right\|_{\lta L^2_{v}}\\
	&\ls\left\|\fM_{10}(t,\na)\rr{e^{-|\alpha|\nu t}v^\alpha P_{0}\wt{\th}(t)} \right\|_{\lta L^2_{v}}.\ea\eeq
	Thus we only need to estimate $\left\|\fM_{10}(t,\na)\rr{e^{-|\alpha|\nu t}v^\alpha \wt{\th}(t)} \right\|_{\lta L^2_{z,v}}$.
	
	By \eqref{wtheq} we have
	\beq\label{wtheqq}\pa_t\wt{\th}=\cN^\dagger(\wt{\th},\mu(\emnt v)\emnt)+\cN^\dagger(\wt{\th},\uth+\th^\#)+\cN^\dagger(\th^\#,\wt{\th})+\cN^\dagger(\wt{\th},\wt{\th})+\nu {\vect}^2\wt{\th}+\wt{\tQ}.\eeq
	Fix $\alpha\leq2$. By the weighted energy identity we have 
	\beq\label{energy}\ba&\f12\f{d}{dt}\left\|\fM_{10}(t,\na)\rr{e^{-|\alpha|\nu t}v^\alpha \wt{\th}(t)} \right\|_{ L^2_{z,v}}^2+\nu|\alpha|\left\|\fM_{10}(t,\na)\rr{e^{-|\alpha|\nu t}v^\alpha \wt{\th}(t)} \right\|_{L^2_{z,v}}^2\\
	&+\nu\left\|\vect\fM_{10}(t,\na)\rr{e^{-|\alpha|\nu t}v^\alpha \wt{\th}(t)} \right\|_{L^2_{z,v}}^2=\sum_{i=1}^{10}\mathscr{A}_i.\ea\eeq
	The terms $\mathscr{A}_i$ are defined and estimated below.
	\beq\label{A1}\ba&\mathscr{A}_1=\sum_k\int\f{\pa_t\fM}{\fM}\nr{\sF\rr{\fM_{10}(t,\na)\rr{e^{-|\alpha|\nu t}v^\alpha \wt{\th}(t)}}}^2d\xi\\
	\leq&\rm\nu\left\|\fM_{10}(t,\na)\rr{e^{-|\alpha|\nu t}v^\alpha \wt{\th}(t)} \right\|_{ L^2_{z,v}}^2+\rm\nu\left\|\vect\fM_{10}(t,\na)\rr{e^{-|\alpha|\nu t}v^\alpha \wt{\th}(t)} \right\|_{ L^2_{z,v}}^2,\ea\eeq
	Note that the second term can be absorbed by the third term on the left-hand side of \eqref{energy}. We use
	\[\f{\pa_t\fM}{\fM}=\rm\nu\rr{1+|\bar{\eta}(t,k,\xi)|^2}^\f12-\rm\nu\int_0^t\nu e^{-2\nu(t-s)}\rr{e^{-2\nu(t-s)}+|\bar{\eta}(s,k,\xi)|^2}^{-\f12}ds\leq\rm\nu\rr{1+|\bar{\eta}(t,k,\xi)|}\]
	in the estimate.
	\[\ba
	&\mathscr{A}_2=-\sum_k\int\wh{\wt{\rho}}(t,k)\f{1}{k}D^\alpha_\xi\rr{\xi\wh{\mu}(\xi)}(\bar{\eta}(t,k,\xi))\fM_{10}(t,k,\xi)\overline{\sF\rr{\fM_{10}(t,\na)\rr{e^{-|\alpha|\nu t}v^\alpha \wt{\th}(t)}}}d\xi\\
	\leq&\nnr{\fM_{10}(t,\na)\wt{\rho}(t)}_{L^2_z}
	\nnr{\fM_{10}(t,0,\emnt\na_v)v^\alpha\pa_v\mu(v)}_{L^2_v}e^{-\f12\nu t}
	\left\|\fM_{10}(t,\na)\rr{e^{-|\alpha|\nu t}v^\alpha \wt{\th}(t)} \right\|_{ L^2_{z,v}}\\
	\ls&\left\|\fM_{10}(t,\na)\rr{e^{-|\alpha|\nu t}v^\alpha \wt{\th}(t)} \right\|_{\lta L^2_{z,v}}\left\|\fM_{10}(t,\na)\rr{e^{-|\alpha|\nu t}v^\alpha \wt{\th}(t)} \right\|_{ L^2_{z,v}}e^{\f32\nu t}.
	\ea\]
	Note that we use Sobolev embedding to estimate $\nnr{\fM_{10}(t,\na)\wt{\rho}(t)}_{L^2_z}$ just like \eqref{oi}.
	\[\ba
	&\mathscr{A}_3=-\sum_{k,\ell}\int\wh{\wt{\rho}}(t,\ell)\f{\bar{\eta}(t,k,\xi)}{\ell}D^\alpha_\xi\rr{\wh{\uth}+\wh{\th^\#}}(t,k-\ell,\xi-\ell\tent)\fM_{10}(t,k,\xi) \overline{\sF\rr{\fM_{10}(t,\na)\rr{v^\alpha \wt{\th}(t)}}}d\xi\\
	&\times e^{-2|\alpha|\nu t}\ls\nnr{\fM_{10}(t,\na)\wt{\rho}(t)}_{L^2_z}
	\nnr{\fM_{12}(t,\na)\rr{e^{-|\alpha|\nu t}v^\alpha \rr{\uth+\th^\#}(t)}}_{L^2_{z,v}}e^{\nu t}\<t\>\\
	&\quad \left\|\fM_{10}(t,\na)\rr{e^{-|\alpha|\nu t}v^\alpha \wt{\th}(t)} \right\|_{ L^2_{z,v}}\\
	&\ls\ueps e^{3\nu t}\<t\>\left\|\fM_{10}(t,\na)\rr{e^{-|\alpha|\nu t}v^\alpha \wt{\th}(t)} \right\|_{\lta L^2_{z,v}}\left\|\fM_{10}(t,\na)\rr{e^{-|\alpha|\nu t}v^\alpha \wt{\th}(t)} \right\|_{ L^2_{z,v}},
	\ea\]
	where we use $\nr{\bar{\eta}(t,k,\eta)}\leq\ent\<t\>\<k-\ell,\xi-\ell\tent\>$.
	\[\ba
	&\mathscr{A}_4=-\sum_{k,\ell}\int\wh{\rho^\#}(t,\ell)\f{\bar{\eta}(t,k,\xi)}{\ell}D^\alpha_\xi\wh{\th}(t,k-\ell,\xi-\ell\tent)e^{-|\alpha|\nu t}\\
	&\quad \left\{\fM_{10}(t,k,\xi)-\fM_{10}(t,k-\ell,\xi-\ell\tent)\right\}\\
	&\times\overline{\sF\rr{\fM_{10}(t,\na)\rr{e^{-|\alpha|\nu t}v^\alpha \wt{\th}(t)}}}d\xi\ls\nnr{\fM_{12}(t,\na)\rho^\#(t)}_{L^2_z}\ent\left\|\fM_{10}(t,\na)\rr{e^{-|\alpha|\nu t}v^\alpha \wt{\th}(t)} \right\|_{ L^2_{z,v}}^2\\
	&+\rm\nu\nnr{\fM_{3}(t,\na)\rho^\#(t)}_{L^2_z}e^{2\nu t}\left\|\vect\fM_{10}(t,\na)\rr{e^{-|\alpha|\nu t}v^\alpha \wt{\th}(t)} \right\|_{ L^2_{z,v}}\left\|\fM_{10}(t,\na)\rr{e^{-|\alpha|\nu t}v^\alpha \wt{\th}(t)} \right\|_{ L^2_{z,v}}\\
	&\ls\rm^2\nu^2\left\|\vect\fM_{10}(t,\na)\rr{e^{-|\alpha|\nu t}v^\alpha \wt{\th}(t)} \right\|_{ L^2_{z,v}}^2
	+\nnr{\fM_{3}(t,\na)\rho^\#(t)}_{L^2_z}^2e^{4\nu t}\left\|\fM_{10}(t,\na)\rr{e^{-|\alpha|\nu t}v^\alpha \wt{\th}(t)} \right\|_{ L^2_{z,v}}^2\\
	&+\nnr{\fM_{12}(t,\na)\rho^\#(t)}_{L^2_z}e^{\nu t}\left\|\fM_{10}(t,\na)\rr{e^{-|\alpha|\nu t}v^\alpha \wt{\th}(t)} \right\|_{ L^2_{z,v}}^2
	\ea\]
	Note that the first term can be absorbed by the third term on the left-hand side of \eqref{energy}. We use the following commutator estimate
	\[\ba
	&\nr{\fM_{10}(t,k,\xi)-\fM_{10}(t,k-\ell,\xi-\ell\tent)}\leq\fM(t,k,\xi)\nr{\<k,\xi\>^{10}-\<k-\ell,\xi-\ell\tent\>^{10}}\\
&+\<k-\ell,\xi-\ell\tent\>^{10}\nr{\fM(t,k,\xi)-\fM(t,k-\ell,\xi-\ell\tent)}\ls\fM(t,k,\xi)\<\ell,\ell\tent\>\\
	&\times\max\{\<k,\xi\>^9,\<k-\ell,\xi-\ell\tent\>^{9}\}+\<k-\ell,\xi-\ell\tent\>^{10}\rm\nu\int_0^t\nr{\bar{\eta}(s,\ell,\ell\tent)}ds\\
	&\times\fM(t,\ell,\ell\tent)\fM(t,k-\ell,\xi-\ell\tent)
	\ls\fM(t,\ell,\ell\tent)\fM(t,k-\ell,\xi-\ell\tent)\\
	&\times\left\{\<\ell,\ell\tent\>\max\{\<k,\xi\>^9,\<k-\ell,\xi-\ell\tent\>^{9}\}+\rm\nu\<\ell,\ell\tent\>^2e^{2\nu t}\<k-\ell,\xi-\ell\tent\>^{10}\right\}.
	\ea\]
	We next estimate $\mathscr{A}_5$
	\[\ba&\mathscr{A}_5=-\sum_{k,\ell}\int\wh{\wt{\rho}}(t,\ell)\f{\bar{\eta}(t,k,\xi)}{\ell}D^\alpha_\xi\wh{\th}(t,k-\ell,\xi-\ell\tent)e^{-|\alpha|\nu t}\\
	&\quad \left\{\fM_{10}(t,k,\xi)-\fM_{10}(t,k-\ell,\xi-\ell\tent)\right\}\overline{\sF\rr{\fM_{10}(t,\na)\rr{e^{-|\alpha|\nu t}v^\alpha \wt{\th}(t)}}}d\xi.\ea\]
	We use the paraproduct decomposition and split the integral into $\<\ell,\ell\tent\>\leq\<k-\ell,\xi-\ell\tent\>$ part and $\<\ell,\ell\tent\>\geq\<k-\ell,\xi-\ell\tent\>$ part. We denote the two terms as $\mathscr{A}_{5,LH}$ and $\mathscr{A}_{5,HL}$. For $\mathscr{A}_{5,LH}$ we use the same commutator estimate as $\mathscr{A}_4$
	\[\ba
	&\mathscr{A}_{5,LH}\ls\nnr{\fM_{3}(t,\na)\wt{\rho}(t)}_{L^2_z}\ent\left\|\fM_{10}(t,\na)\rr{e^{-|\alpha|\nu t}v^\alpha \wt{\th}(t)} \right\|_{ L^2_{z,v}}^2+\rm\nu\nnr{\fM_{3}(t,\na)\wt{\rho}(t)}_{L^2_z}e^{2\nu t}\\
	&\left\|\vect\fM_{10}(t,\na)\rr{e^{-|\alpha|\nu t}v^\alpha \wt{\th}(t)} \right\|_{ L^2_{z,v}}\left\|\fM_{10}(t,\na)\rr{e^{-|\alpha|\nu t}v^\alpha \wt{\th}(t)} \right\|_{ L^2_{z,v}}\\
	&\ls\rm^2\nu^2\left\|\vect\fM_{10}(t,\na)\rr{e^{-|\alpha|\nu t}v^\alpha \wt{\th}(t)} \right\|_{ L^2_{z,v}}^2\\
	&+\left\|\fM_{10}(t,\na)\rr{e^{-|\alpha|\nu t}v^\alpha \wt{\th}(t)} \right\|_{\lta L^2_{z,v}}^2e^{6\nu t}\left\|\fM_{10}(t,\na)\rr{e^{-|\alpha|\nu t}v^\alpha \wt{\th}(t)} \right\|_{ L^2_{z,v}}^2\\
	&+\left\|\fM_{10}(t,\na)\rr{e^{-|\alpha|\nu t}v^\alpha \wt{\th}(t)} \right\|_{\lta L^2_{z,v}}e^{3\nu t}\left\|\fM_{10}(t,\na)\rr{e^{-|\alpha|\nu t}v^\alpha \wt{\th}(t)} \right\|_{ L^2_{z,v}}^2.
	\ea\]
	For $\mathscr{A}_{5,HL}$ we use a similar estimate to $\mathscr{A}_3$
	\[\ba
	&\mathscr{A}_{5,HL}\ls\nnr{\fM_{10}(t,\na)\wt{\rho}(t)}_{L^2_z}
	\nnr{\fM_{3}(t,\na)\rr{e^{-|\alpha|\nu t}v^\alpha \wt{\th}(t)}}_{L^2_{z,v}}e^{\nu t}\<t\>
	\left\|\fM_{10}(t,\na)\rr{e^{-|\alpha|\nu t}v^\alpha \wt{\th}(t)} \right\|_{ L^2_{z,v}}\\
	&\ls e^{3\nu t}\<t\>\left\|\fM_{10}(t,\na)\rr{e^{-|\alpha|\nu t}v^\alpha \wt{\th}(t)} \right\|_{\lta L^2_{z,v}}\left\|\fM_{10}(t,\na)\rr{e^{-|\alpha|\nu t}v^\alpha \wt{\th}(t)} \right\|_{ L^2_{z,v}}^2.
	\ea\]
	We next estimate $\mathscr{A}_{6}$.
	\[\ba
	&\mathscr{A}_{6}=\int\fM_{10}(t,\na)\rr{e^{-|\alpha|\nu t}v^\alpha \wt{\tQ}(t)}\overline{\fM_{10}(t,\na)\rr{e^{-|\alpha|\nu t}v^\alpha \wt{\th}(t)}}dzdv\\
	\leq&\left\|\fM_{10}(t,\na)\rr{e^{-|\alpha|\nu t}v^\alpha \wt{\tQ}(t)} \right\|_{ L^2_{z,v}}\left\|\fM_{10}(t,\na)\rr{e^{-|\alpha|\nu t}v^\alpha \wt{\th}(t)} \right\|_{ L^2_{z,v}}.
	\ea\]
	The terms $\mathscr{A}_i$, $i=7,8,9,10$, arise only for $\alpha\geq1$; recall that $d=1$. Arguing as in the estimate of $\mathscr{A}_i\ i=3,4,5$ we have
	\[\ba&\mathscr{A}_7=-\sum_{k,\ell}\int\wh{\wt{\rho}}(t,\ell)\f{1}{\ell}D^{\alpha-1}_\xi\rr{\wh{\uth}+\wh{\th^\#}}(t,k-\ell,\xi-\ell\tent)\fM_{10}(t,k,\xi)\overline{\sF\rr{\fM_{10}(t,\na)\rr{v^\alpha \wt{\th}(t)}}}d\xi\\
	&\times e^{-(2|\alpha|+1)\nu t}\ls\nnr{\fM_{10}(t,\na)\wt{\rho}(t)}_{L^2_z}
	\nnr{\fM_{11}(t,\na)\rr{e^{-(\alpha-1)\nu t}v^{\alpha-1} \rr{\uth+\th^\#}(t)}}_{L^2_{z,v}}\\
	&\quad \left\|\fM_{10}(t,\na)\rr{e^{-|\alpha|\nu t}v^\alpha \wt{\th}(t)} \right\|_{ L^2_{z,v}}\\
	&\ls\ueps e^{2\nu t}\left\|\fM_{10}(t,\na)\rr{e^{-|\alpha|\nu t}v^\alpha \wt{\th}(t)} \right\|_{\lta L^2_{z,v}}\left\|\fM_{10}(t,\na)\rr{e^{-|\alpha|\nu t}v^\alpha \wt{\th}(t)} \right\|_{ L^2_{z,v}},
	\ea\]
	\[\ba&\mathscr{A}_8=-\sum_{k,\ell}\int\wh{\rho^\#}(t,\ell)\f{1}{\ell}D^{\alpha-1}_\xi\wh{\wt{\th}}(t,k-\ell,\xi-\ell\tent)e^{-(\alpha-1)\nu t}\fM_{10}(t,k,\xi)\overline{\sF\rr{\fM_{10}(t,\na)\rr{e^{-|\alpha|\nu t}v^\alpha \wt{\th}(t)}}}d\xi\\
	&\ls\nnr{\fM_{10}(t,\na)\rho^\#(t)}_{L^2_z}
	\nnr{\fM_{10}(t,\na)\rr{e^{-(\alpha-1)\nu t}v^{\alpha-1} \wt{\th}(t)}}_{L^2_{z,v}}
	\left\|\fM_{10}(t,\na)\rr{e^{-|\alpha|\nu t}v^\alpha \wt{\th}(t)} \right\|_{ L^2_{z,v}},
	\ea\]
	\[\ba&\mathscr{A}_9=-\sum_{k,\ell}\int\wh{\wt{\rho}}(t,\ell)\f{1}{\ell}D^{\alpha-1}_\xi\wt{\th}(t,k-\ell,\xi-\ell\tent)e^{-(\alpha-1)\nu t}\fM_{10}(t,k,\xi)\overline{\sF\rr{\fM_{10}(t,\na)\rr{e^{-|\alpha|\nu t}v^\alpha \wt{\th}(t)}}}d\xi\\
	&\ls\nnr{\fM_{10}(t,\na)\wt{\rho}(t)}_{L^2_z}
	\nnr{\fM_{10}(t,\na)\rr{e^{-(\alpha-1)\nu t}v^{\alpha-1} \wt{\th}(t)}}_{L^2_{z,v}}
	\left\|\fM_{10}(t,\na)\rr{e^{-|\alpha|\nu t}v^\alpha \wt{\th}(t)} \right\|_{ L^2_{z,v}}.
	\ea\]
	We next estimate $\mathscr{A}_{10}$. For $\alpha=1$
	\[\ba
	&\mathscr{A}_{10}=-2\nu\int\vect\fM_{10}(t,\na)\wt{\th}(t)\overline{\fM_{10}(t,\na)\rr{e^{-|\alpha|\nu t}v^\alpha \wt{\th}(t)}}dzdv\leq2\nu\nnr{\vect\fM_{10}(t,\na)\rr{e^{-|\alpha|\nu t}v^\alpha \wt{\th}(t)}}_{L^2_{z,v}}\\
    &\times\nnr{\fM_{10}(t,\na)\wt{\th}(t)}_{L^2_{z,v}}\leq\f\nu8\nnr{\vect\fM_{10}(t,\na)\rr{e^{-|\alpha|\nu t}v^\alpha \wt{\th}(t)}}_{L^2_{z,v}}^2+8\nu\nnr{\fM_{10}(t,\na)\wt{\th}(t)}_{L^2_{z,v}}^2.
	\ea\]
	For $\alpha=2$
	\[\ba
	&\mathscr{A}_{10}=2\nu\int\fM_{10}(t,\na)\wt{\th}(t)\overline{\fM_{10}(t,\na)\rr{e^{-|\alpha|\nu t}v^\alpha \wt{\th}(t)}}dzdv-4\nu\int\vect\fM_{10}(t,\na)\rr{\emnt v\wt{\th}(t)}\\
	&\times\overline{\fM_{10}(t,\na)\rr{e^{-|\alpha|\nu t}v^\alpha \wt{\th}(t)}}dzdv\leq\f18\nu\nnr{\vect\fM_{10}(t,\na)\rr{e^{-|\alpha|\nu t}v^\alpha \wt{\th}(t)}}_{L^2_{z,v}}^2\\
	&+32\nu\nnr{\fM_{10}(t,\na)\rr{\emnt v\wt{\th}(t)}}_{L^2_{z,v}}^2\\
	&+\nu\nnr{\fM_{10}(t,\na)\wt{\th}(t)}_{L^2_{z,v}}^2+\nu\nnr{\fM_{10}(t,\na)\rr{e^{-|\alpha|\nu t}v^\alpha \wt{\th}(t)}}_{L^2_{z,v}}^2
	\ea\]
	Note that the first term can be absorbed by the third term on the left-hand side of \eqref{energy}.
	
	We use a standard continuity argument to prove that 
	\beq\label{Pop}\nnr{\fM_{10}(t,\na)\rr{e^{-|\alpha|\nu t}v^\alpha \wt{\th}(t)}}_{L^2_{z,v}}\leq A^{-\f {N_1}{2} A^{N'}}e^{0.8\nu^{-2}}\eeq
	for $t\in[T,\nu^{-1}]$.	Note that $\wt{h}(T)=0$. We suppose that \eqref{Pop} holds on $[T,T^*]$, $T<T^*\leq\nu^{-1}$. Then combining the previous estimates we obtain 
	\[\ba
	&\f12\f{d}{dt}\nnr{\fM_{10}(t,\na)\rr{e^{-|\alpha|\nu t}v^\alpha \wt{\th}(t)}}_{\lta L^2_{z,v}}^2\leq C\left\{\ueps t+\nnr{\fM_{12}(t,\na)\rho^\#(t)}_{L^2_z}^2+\nnr{\fM_{12}(t,\na)\rho^\#(t)}_{L^2_z}\right\}\\
	&\times\nnr{\fM_{10}(t,\na)\rr{e^{-|\alpha|\nu t}v^\alpha \wt{\th}(t)}}_{\lta L^2_{z,v}}^2\\
	&+\nnr{\fM_{10}(t,\na)\rr{e^{-|\alpha|\nu t}v^\alpha \wt{\th}(t)}}_{\lta L^2_{z,v}}\nnr{\fM_{10}(t,\na)\rr{e^{-|\alpha|\nu t}v^\alpha \wt{\tQ}(t)}}_{\lta L^2_{z,v}},
	\ea\]
	on $[T,T^*]$, which implies that
	\[\ba
	&\nnr{\fM_{10}(t,\na)\rr{e^{-|\alpha|\nu t}v^\alpha \wt{\th}(t)}}_{\lta L^2_{z,v}}\leq\int_T^t\nnr{\fM_{10}(\tau,\na)\rr{e^{-|\alpha|\nu \tau}v^\alpha \wt{\tQ}(\tau)}}_{\lta L^2_{z,v}}d\tau\ \exp\left\{\int_T^tC\left\{\ueps \tau\right.\right.\\
	&\left.\left.+\nnr{\fM_{12}(\tau,\na)\rho^\#(\tau)}_{L^2_z}^2+\nnr{\fM_{12}(\tau,\na)\rho^\#(\tau)}_{L^2_z}\right\}d\tau\right\}\\
	&\leq A^{-\f {N_1}{2}A^{N'}}\exp\left\{C(\ueps\nu^{-2}+A^{-2N_1+30}+\nu^{-\f12}A^{-N_1+15})\right\},
	\ea\]
	where we use \eqref{twQ} and \eqref{h3}. This implies that \eqref{Pop} holds on a longer time interval. A standard continuity argument shows that \eqref{Pop} holds on $[T,\nu^{-1}]$. This, together with \eqref{a1}, \eqref{a2}, and \eqref{a3}, yields \eqref{TT1}.
\end{proof}

\subsection{Bounds on $[\nu^{-1},\infty)$.}
The finite-time energy estimate is no longer sufficient for global control. We therefore use separate decay mechanisms for the density, the nonzero spatial modes, and the zero mode. The next proposition propagates the small remainder beyond the thermal time. Its proof closes three coupled bootstrap bounds: the density estimate supplies the electric-field control needed for hypocoercivity, the hypocoercive estimate gives enhanced dissipation, and the homogeneous semigroup gives thermal relaxation of the zero mode.

\begin{prop}[Global remainder bounds after the thermal time]\label{prop6}
	There exists a universal constant $0<\rm_1\ll\rm$ such that the following holds on $[\nu^{-1},\infty)$.
	\beq\label{T1infty}\ba
	&\nnr{\fM_{8.5}(\tau,\na)\wt{\rho}(\tau)e^{\rm_1\nu^\f13\tau}}_{L^2([\nu^{-1},t], L^2_{z})}+\left\|\fM_{7.5}(t,\na)\rr{e^{-|\alpha|\nu t}v^\alpha P_{\neq}\wt{\th}(t)} \right\|_{\lta L^2_{z,v}}e^{\rm_1\nu^\f13t}\\
	&+\nnr{v^\alpha\<\na_v\>\fM_{7.5}(t,0,\emnt\na_v)P_{0}\wt{h}(t)}_{\lta L^2_{v}}e^{\nu t}\leq A^{-\f {N_1}{2}A^{N'}}e^{2\nu^{-2}}.
	\ea\eeq
\end{prop}

We prove Proposition~\ref{prop6}. The strategy follows \cite{LY}; the argument is simpler because of the smallness of $\wt{h}$  and because the initial time is already at the thermal scale $\nu^{-1}$. For convenience, set 
\[\eps_1=A^{-\f {N_1}{2}A^{N'}}e^{\nu^{-2}}.\]
We introduce the following bootstrap hypotheses
\beq\label{wtrho}\nnr{\fM_{8.5}(\tau,\na)\wt{\rho}(\tau)e^{\rm\nu^\f13\tau}}_{L^2([\nu^{-1},t], L^2_{z})}\leq K_1\eps_1,\eeq
\beq\label{wthed}\left\|\pa_x^2\fM_{7.5}(t,\na)\rr{e^{-|\alpha|\nu t}v^\alpha \wt{\th}(t)} \right\|_{\lta L^2_{z,v}}+\nu^\f13\left\|\pa_x\vect\fM_{7.5}(t,\na)\rr{e^{-|\alpha|\nu t}v^\alpha \wt{\th}(t)} \right\|_{\lta L^2_{z,v}}\leq K_2\eps_1e^{-\rm_1\nu^\f13t},\eeq
\beq\label{wth0}\nnr{v^\alpha\<\na_v\>\fM_{7.5}(t,0,\emnt\na_v)P_{0}\wt{h}(t)}_{\lta L^2_{v}}\leq K_3\eps_1e^{-\nu t}.\eeq
$K_1, K_2, K_3$ are uniform constants. Estimate~\eqref{T1infty} follows once we show that (\ref{wtrho}-\ref{wth0}) hold for any $t\in[\nu^{-1},\infty)$. We will improve these estimates under the bootstrap hypotheses and close the argument by continuity.

\subsubsection{Improvement of \eqref{wtrho}.}
We first estimate the density through its Volterra resolvent. The data at $\nu^{-1}$ and the truncation source contribute at the scale $\eps_1$. For the nonlinear terms, we separate zero and nonzero profile modes and use a paraproduct decomposition to place the higher frequency weight on the appropriate factor. The small remainder and the decay of the approximate solution then yield an improvement of the density bootstrap.
 By \eqref{wtheqq} and Duhamel's formula we have
\[\ba
&\wh{\wt{\rho}}(t,k)=\wh{\wt{H}}(t,k)-\int_{\nu^{-1}}^t\wh{\wt{H}}(\tau,k)\wh{R}(t-\tau,k)d\tau,\\
	&\quad 
\wh{\wt{H}}(t,k)=\wh{\wt{\th}}(\nu^{-1},k,k\tent)S(t-\nu^{-1},k)+\int_{\nu^{-1}}^t\wh{\wt{\tQ}}(\tau,k,k\tent)\\
&\times S(t-\tau,k)d\tau-\sum_{\ell}\int_{\nu^{-1}}^t\wh{\wt{\rho}}(\tau,\ell)\f k\ell\f{1-e^{-\nu(t-\tau)}}{\nu}\sF\rr{\uth+\th^\#}(\tau,k-\ell,k\tent-\ell\tentau)S(t-\tau,k)d\tau\\
&-\sum_{\ell}\int_{\nu^{-1}}^t\wh{\rho^\#}(\tau,\ell)\f k\ell\f{1-e^{-\nu(t-\tau)}}{\nu}\wh{\wt{\th}}(\tau,k-\ell,k\tent-\ell\tentau)S(t-\tau,k)d\tau\\
&-\sum_{\ell}\int_{\nu^{-1}}^t\wh{\wt{\rho}}(\tau,\ell)\f k\ell\f{1-e^{-\nu(t-\tau)}}{\nu}\wh{\wt{\th}}(\tau,k-\ell,k\tent-\ell\tentau)S(t-\tau,k)d\tau,
\ea\]
where $\wh{R}(t,k)$ is defined in \eqref{definereso}. We define the five terms in $\wh{\wt{H}}(t,k)$ as $\mathscr{B}_i$, $i=1,2,3,4,5$. By \eqref{resolvent} (similar to \eqref{egz}) we only need to estimate $\nnr{\fM_{8.5}(\tau,\na)\wt{H}(\tau)e^{\rm\nu^\f13\tau}}_{L^2([\nu^{-1},t], L^2_{z})}$.
As in \eqref{egx} and \eqref{egc} we have
\[\nnr{\fM_{8.5}(\tau,\na)\mathscr{B}_1e^{\rm\nu^\f13\tau}}_{L^2([\nu^{-1},t], L^2_{z})}\ls\nnr{\fM_{8.5}(\nu^{-1},\na)P_{\neq}\wt{\th}(\nu^{-1})}_{L^2_{z,v}}e^{\rm\nu^{-\f23}}\leq\eps_1,\]
\[\nnr{\fM_{8.5}(\tau,\na)\mathscr{B}_2e^{\rm\nu^\f13\tau}}_{L^2([\nu^{-1},t], L^2_{z})}\ls\int^t_{\nu^{-1}}\nnr{\fM_{8.5}(\tau,\na)P_{\neq}\wt{\tQ}(\tau)}_{L^2_{z,v}}e^{\rm\nu^\f13\tau+\f12\nu\tau}d\tau\leq\eps_1,\]
where we use \eqref{TT1} and \eqref{twQ}. 

We focus on the more technical nonlinear term $\mathscr{B}_5$. We split $\mathscr{B}_5$ into $\mathscr{B}_{5,\neq}$ and $\mathscr{B}_{5,0}$, corresponding to $\ell\neq k$ and $\ell=k$, respectively. We use the standard paraproduct decomposition and split the integral into the two regions $\<\ell,\ell\tentau\>\leq\<k-\ell,k\tent-\ell\tentau\>$ and $\<\ell,\ell\tentau\>\geq\<k-\ell,k\tent-\ell\tentau\>$ and denote the corresponding contributions by LH and HL.

We first consider $\mathscr{B}_{5,\neq,HL}$.
\[\ba
&\nr{\sF\rr{\fM_{8.5}(t,\na)\mathscr{B}_{5,\neq,HL}}(t,k)}e^{\rm\nu^\f13t}
\ls\sum_{\ell}\int_{\nu^{-1}}^t\nr{\sF\rr{\fM_{8.5}(\tau,\na)\wt{\rho}}(\tau,\ell)}e^{\rm\nu^\f13\tau}\f{|k|}{|\ell|}\f{1-e^{-\nu(t-\tau)}}{\nu}\\
	&\quad \nr{\sF\rr{\fM(\tau,\na)P_{\neq}\wt{\th}}(\tau,k-\ell,k\tent-\ell\tentau)}S^{0.9}(t-\tau,k)d\tau\\
	&\ls\nu^{-\f13}\sum_{\ell}\int_{\nu^{-1}}^t\nr{\sF\rr{\fM_{8.5}(\tau,\na)\wt{\rho}}(\tau,\ell)}e^{\rm\nu^\f13\tau}\\
&\times\nr{\sF\rr{\na_x\fM(\tau,\na)P_{\neq}\wt{\th}}(\tau,k-\ell,k\tent-\ell\tentau)}S^{0.8}(t-\tau,k)d\tau\\
	&\ls\nu^{-\f13}\sum_{\ell}\int_{\nu^{-1}}^t\nr{\sF\rr{\fM_{8.5}(\tau,\na)\wt{\rho}}(\tau,\ell)}\\
&\times e^{\rm\nu^\f13\tau}\<k-\ell,k\tent-\ell\tentau\>^{-4}\emnt d\tau\ \sup_{\tau\in[\nu^{-1},t]}\nnr{\na_x\fM_4(\tau,\na)\rr{e^{-|\alpha|\nu\tau}v^\alpha P_{\neq}\wt{\th}} }_{\lta L^2_{z,v}}e^{3\nu\tau}.
\ea\]
The second inequality uses $\f{1-e^{-\nu(t-\tau)}}{\nu}S^{0.1}(t-\tau,k)\ls\nu^{-\f13}$ and $|k|\ls|k-\ell||\ell|$. (Note that $\ell\neq0$, $k-\ell\neq0$.) The third inequality uses Sobolev embedding. By the bootstrap bound \eqref{wthed} and Schur's test we have
\[\nnr{\fM_{8.5}(\tau,\na)\mathscr{B}_{5,\neq,HL}e^{\rm\nu^\f13\tau}}_{L^2([\nu^{-1},t], L^2_{z})}\ls\eps_1\nu^{-\f13}\nnr{\fM_{8.5}(\tau,\na)\wt{\rho}e^{\rm\nu^\f13\tau}}_{L^2([\nu^{-1},t], L^2_{z})}\ls\eps_1^2\nu^{-\f13}. \]

We next consider $\mathscr{B}_{5,\neq,LH}$. As in the estimate of $\mathscr{B}_{5,\neq,HL}$ we have
\[\ba
&\nr{\sF\rr{\fM_{8.5}(t,\na)\mathscr{B}_{5,\neq,LH}}(t,k)}e^{\rm\nu^\f13t}\ls
\nu^{-\f13}\sum_{\ell}\int_{\nu^{-1}}^t\nr{\sF\rr{\fM(\tau,\na)\wt{\rho}}(\tau,\ell)}e^{\rm\nu^\f13\tau}\\
&\times e^{\f12\nu\tau}\nr{\sF\rr{\na_x\fM_{8.5}(\tau,\na)P_{\neq}\wt{\th}}(\tau,k-\ell,k\tent-\ell\tentau)}S^{0.8}(t-\tau,k)e^{-\f12\nu t}d\tau
\ea\]
By Minkowski's inequality we have
\[\ba&\nnr{\fM_{8.5}(\tau,\na)\mathscr{B}_{5,\neq,LH}e^{\rm\nu^\f13\tau}}_{L^2([\nu^{-1},t], L^2_{z})}\ls\nu^{-\f13}\nnr{\fM_1(\tau,\na)\wt{\rho}e^{\rm\nu^\f13\tau-\f12\nu\tau}}_{L^1([\nu^{-1},t], L^2_{z})}\\
	&\quad \sup_{\tau\in[\nu^{-1},t]}\nnr{\na_x\fM_{8.5}(\tau,\na) P_{\neq}\wt{\th} }_{ L^2_{z,v}} e^{\nu\tau}\\
&\ls\nu^{-\f13}\eps_1\nnr{\fM_2(\tau,\na)\wt{\rho}e^{\rm\nu^\f13\tau}}_{L^2([\nu^{-1},t], L^2_{z})}\nnr{\<\tentau\>^{-1}e^{-\f12\nu\tau}}_{L^2([\nu^{-1},t])}\ls\nu^{-\f13}\eps_1^2.
\ea\]
For the estimate of $\nnr{\na_x\fM_{8.5}(\tau,\na)\rr{e^{-|\alpha|\nu\tau}v^\alpha P_{\neq}\wt{\th}} }_{\lta L^2_{z,v}}$ we use
\beq\label{519}\ba&\nnr{\na_x\fM_{8.5}(\tau,\na)\rr{e^{-|\alpha|\nu\tau}v^\alpha P_{\neq}\wt{\th}} }_{\lta L^2_{z,v}}\ls\nnr{\na_x^2\fM_{7.5}(\tau,\na)\rr{e^{-|\alpha|\nu\tau}v^\alpha P_{\neq}\wt{\th}} }_{\lta L^2_{z,v}}\f{\entau-1}{\nu}\\
&+\nnr{\na_x\vect\fM_{7.5}(\tau,\na)\rr{e^{-|\alpha|\nu\tau}v^\alpha P_{\neq}\wt{\th}} }_{\lta L^2_{z,v}}
\ls\nu^{-\f13}\f{\entau-1}{\nu}e^{-\rm_1\nu^\f13\tau}\eps_1\ls\eps_1 e^{-0.5\rm_1\nu^\f13\tau}
\ea\eeq
on $[\nu^{-1},\infty)$. Note that we use \eqref{wthed} and $|\xi|\leq|\bar{\eta}|+|k|\f{\ent-1}{\nu}$.

We next consider $\mathscr{B}_{5,0,HL}$. By $\f{1-e^{-\nu(t-\tau)}}{\nu}S^{0.1}(t-\tau,k)\ls\nu^{-\f13}$ and Sobolev embedding we have
\[\ba
&\nr{\sF\rr{\fM_{8.5}(t,\na)\mathscr{B}_{5,0,HL}}(t,k)}e^{\rm\nu^\f13t}\ls
\nu^{-\f13}\int_{\nu^{-1}}^t\nr{\sF\rr{\fM_{8.5}(\tau,\na)\wt{\rho}}(\tau,k)}e^{\rm\nu^\f13\tau}S^{0.8}(t-\tau,k)d\tau\\ &\times\sup_{\tau\in[\nu^{-1},t]}\nnr{v^\alpha\fM(\tau,0,\emntau\na_v)P_0\wt{h}}_{\lta L^2_v}.
\ea\]
Thus by Schur's test and \eqref{wth0} we have
\[\nnr{\fM_{8.5}(\tau,\na)\mathscr{B}_{5,0,HL}e^{\rm\nu^\f13\tau}}_{L^2([\nu^{-1},t], L^2_{z})}\ls\eps_1\nu^{-\f23}\nnr{\fM_{8.5}(\tau,\na)\wt{\rho}e^{\rm\nu^\f13\tau}}_{L^2([\nu^{-1},t], L^2_{z})}\ls\eps_1^2\nu^{-\f23}. \]
Note that we use $\int_0^tS^{0.5}(t-\tau,k)d\tau\ls\nu^{-\f13}$.

We next consider $\mathscr{B}_{5,0,LH}$. 
\[\ba
&\nr{\sF\rr{\fM_{8.5}(t,\na)\mathscr{B}_{5,0,LH}}(t,k)}e^{\rm\nu^\f13t}\ls
\nu^{-\f13}\int_{\nu^{-1}}^t\nr{\sF\rr{\fM(\tau,\na)\wt{\rho}}(\tau,k)}e^{\rm\nu^\f13\tau}\\
&\times \nr{\sF\rr{\fM_{8.5}(\tau,0,\emntau\na)P_{0}\wt{h}}(\tau,0,k\f{1-e^{-\nu(t-\tau)}}{\nu} )}e^{-\f12\nu (t-\tau)}d\tau
\ea\]
Thus by Minkowski's inequality and Cauchy--Schwarz inequality we have
\[\ba
&\nnr{\fM_{8.5}(\tau,\na)\mathscr{B}_{5,0,LH}e^{\rm\nu^\f13\tau}}_{L^2([\nu^{-1},t], L^2_{z})}\ls\nu^{-\f13}\nnr{\fM_1(\tau,\na)\wt{\rho}e^{\rm\nu^\f13\tau-\f12\nu\tau}}_{L^1([\nu^{-1},t], L^2_{z})}\\
&\times\sup_{\tau\in[\nu^{-1},t]}\nnr{\fM_{8.5}(\tau,0,\emntau\na)P_{0}\wt{h}}_{L^2_v}e^{\f12\nu\tau}
\ls\nu^{-\f13}\eps_1\nnr{\fM_2(\tau,\na)\wt{\rho}e^{\rm\nu^\f13\tau}}_{L^2([\nu^{-1},t], L^2_{z})}\ls\nu^{-\f13}\eps_1^2.
\ea\]
This completes the estimate of $\mathscr{B}_5$. $\mathscr{B}_3$ and $\mathscr{B}_4$ can be bounded by the same argument. Note that by \eqref{h3}, $h^\#$ satisfies the estimates (\ref{wtrho}-\ref{wth0}) with $\eps_1=A^{-N_1+16}\ll\nu^{3}$. $P_{\neq}\uh$ satisfies
\beq\label{Hg}\left\|\fM_{10}(t,\na)\rr{e^{-|\alpha|\nu t}v^\alpha P_{\neq}\uth(t)} \right\|_{\lta L^2_{z,v}}\ls\ueps e^{-\rm\nu^\f13t}\ls \ueps e^{-0.5\rm\nu^{-\f23}} e^{-\rm_1\nu^\f13t}.\eeq
$P_0\udl{h}$ has sufficiently high regularity so the term containing $P_0\udl{h}$ can be bounded by the same argument as \eqref{egv}. There exist a uniform constant $K'_1$ and $K'_{11}(K_1,K_2,K_3)$ depending on $K_1,K_2,K_3$ such that if (\ref{wtrho}-\ref{wth0}) hold then
\beq\label{imp1}\nnr{\fM_{8.5}(\tau,\na)\wt{\rho}e^{\rm\nu^\f13\tau}}_{L^2([\nu^{-1},t], L^2_{z})}\leq K'_1\eps_1+K'_2(K_1,K_2,K_3)\nu^{-2}\rr{\eps_1^2+A^{-N_1+16}\eps_1}.\eeq

\subsubsection{Improvement of \eqref{wthed}.}
We next control the nonzero modes by a hypocoercive energy. The mixed term couples spatial and velocity derivatives and produces dissipation at the enhanced-dissipation scale. Terms arising from the multiplier and the differentiated characteristics are absorbed into the dissipation, while the improved density bound controls the electric-field terms. We detail the quadratic remainder interaction; the terms involving the approximate solution are treated with its established bounds.
 We define the following hypocoercive energy functional
\[\ba
\sE_\alpha(t&)=\left\|\pa_x^2\fM_{7.5}(t,\na)\rr{e^{-|\alpha|\nu t}v^\alpha \wt{\th}(t)} \right\|_{ L^2_{z,v}}^2e^{2\rm_1\nu^\f13t}+\nu^\f23\left\|\pa_x\vect\fM_{7.5}(t,\na)\rr{e^{-|\alpha|\nu t}v^\alpha \wt{\th}(t)} \right\|_{L^2_{z,v}}^2e^{2\rm_1\nu^\f13t}\\
&+\nu^\f13\int\pa_x^2\fM_{7.5}(t,\na)\rr{e^{-|\alpha|\nu t}v^\alpha \wt{\th}(t)}\cdot\pa_x\vect\fM_{7.5}(t,\na)\rr{e^{-|\alpha|\nu t}v^\alpha \wt{\th}(t)}dzdve^{2\rm_1\nu^\f13t}
\ea\]
for $|\alpha|\leq2$, which has the same structure as that of \cite{LY}. We set
\[\sE(t)=\sum_{|\alpha|\leq2}\sE_\alpha(t).\]
Let
\[\ba
\sD_\alpha(t)=&\nu\left\|\vect\pa_x^2\fM_{7.5}(t,\na)\rr{e^{-|\alpha|\nu t}v^\alpha \wt{\th}(t)} \right\|_{ L^2_{z,v}}^2e^{2\rm_1\nu^\f13t}+\nu^\f53\left\|\pa_x{\vect}^2\fM_{7.5}(t,\na)\rr{e^{-|\alpha|\nu t}v^\alpha \wt{\th}(t)} \right\|_{L^2_{z,v}}^2e^{2\rm_1\nu^\f13t}\\
&+\f{\nu^\f13}2\left\|\pa_x^2\fM_{7.5}(t,\na)\rr{e^{-|\alpha|\nu t}v^\alpha \wt{\th}(t)} \right\|_{ L^2_{z,v}}^2e^{2\rm_1\nu^\f13t}.
\ea\]
Direct comparison of the quadratic forms gives
\[\ba
	&\sE_\alpha(t)\sim\left\{\left\|\pa_x^2\fM_{7.5}(t,\na)\rr{e^{-|\alpha|\nu t}v^\alpha \wt{\th}(t)} \right\|_{ L^2_{z,v}}^2\right.\\
	&\left.+\nu^\f23\left\|\pa_x\vect\fM_{7.5}(t,\na)\rr{e^{-|\alpha|\nu t}v^\alpha \wt{\th}(t)} \right\|_{L^2_{z,v}}^2\right\}e^{2\rm_1\nu^\f13t}\ls\nu^{-\f13}\sD_\alpha(t).
	\ea\]
Fix $|\alpha|\leq2$. By \eqref{wtheqq} we have
\[\ba
&\f12\f{d}{dt}\sE_\alpha(t)+\sD_\alpha(t)+|\alpha|\nu\left\|\pa_x^2\fM_{7.5}(t,\na)\rr{e^{-|\alpha|\nu t}v^\alpha \wt{\th}(t)} \right\|_{ L^2_{z,v}}^2e^{2\rm_1\nu^\f13t}\\
&+|\alpha|\nu^\f53\left\|\pa_x\vect\fM_{7.5}(t,\na)\rr{e^{-|\alpha|\nu t}v^\alpha \wt{\th}(t)} \right\|_{L^2_{z,v}}^2e^{2\rm_1\nu^\f13t}=\sum_{i=1}^{10}\sC_i.
\ea\]
The terms $\sC_i$ are defined and estimated below.
\[
\sC_1=\rm_1\nu^\f13\sE_\alpha(t)\ls\rm_1\sD_\alpha(t).
\]
For $\rm_1\ll1$ this term can be absorbed by $\sD_\alpha(t)$.
\[\ba&\sC_2=\sum_{k}\int\f{\pa_t\fM}{\fM}\left|\sF\rr{\pa_x\fM_{7.5}(t,\na)\rr{e^{-|\alpha|\nu t}v^\alpha \wt{\th}(t)}} \right|^2\rr{k^2+\nu^\f23\bar{\eta}^2+\nu^\f13k\bar{\eta}}d\xi e^{2\rm_1\nu^\f13t}\ls\rm\sD_\alpha(t)+\rm\nu\sE_\alpha(t).
\ea\]
We use the multiplier-derivative estimate from \eqref{A1}. For $\rm\ll1$ this term can be absorbed by $\sD_\alpha(t)$.
\[\ba
&\sC_3=\sum_k\int\left|\sF\rr{\pa_x\fM_{7.5}(t,\na)\rr{e^{-|\alpha|\nu t}v^\alpha \wt{\th}(t)}} \right|^2\rr{\nu^\f23\left\{\nu\bar{\eta}-k\right\}\bar{\eta}+\f{\nu^\f43}{2}\bar{\eta}k}d\xi e^{2\rm_1\nu^\f13t}\leq\rr{\f{2\nu}{3}+\nu^\f43}\\
&\times\left\|\pa_x\vect\fM_{7.5}(t,\na)\rr{e^{-|\alpha|\nu t}v^\alpha \wt{\th}(t)} \right\|_{ L^2_{z,v}}^2e^{2\rm_1\nu^\f13t}+\rr{\f{3\nu^\f13}{8}+\nu^\f43}\left\|\pa_x^2\fM_{7.5}(t,\na)\rr{e^{-|\alpha|\nu t}v^\alpha \wt{\th}(t)} \right\|_{ L^2_{z,v}}^2e^{2\rm_1\nu^\f13t}
\ea\]
This term can be absorbed by $\sD_\alpha(t)$.
\[\ba
\sC_4&=-\sum_k\int\wh{\wt{\rho}}(t,k)D^\alpha_\xi\rr{\xi\wh{\mu}(\xi)}(\bar{\eta}(t,k,\xi))\rr{k^2+\nu^\f23\bar{\eta}^2+\nu^\f13k\bar{\eta}}\overline{\sF\rr{\pa_x\fM_{7.5}(t,\na)\rr{e^{-|\alpha|\nu t}v^\alpha \wt{\th}(t)}}(t,k,\xi)}\\
&\times\fM_{7.5}(t,k,\xi)d\xi e^{2\rm_1\nu^\f13t}
\ls\nnr{\fM_{8.5}(t,\na)\wh{\wt{\rho}}e^{\rm\nu^\f13t}}_{L^2_z}\nnr{\pa_x^2\fM_{7.5}(t,\na)\rr{e^{-|\alpha|\nu t}v^\alpha \wt{\th}(t)}}_{L^2_{z,v}}e^{\rm_1\nu^\f13t}e^{-(\rm-\rm_1)\nu^\f13t}.
\ea\]
Note that under the bootstrap hypotheses \eqref{wtrho} and \eqref{wthed}
\[\ba
\int^{t}_{\nu^{-1}}\sC_4d\tau&\ls\nnr{\fM_{8.5}(\tau,\na)\wh{\wt{\rho}}e^{\rm\nu^\f13\tau}}_{L^2([\nu^{-1},t],L^2_z)}\sup_{\tau\in[\nu^{-1},t]}\nnr{\pa_x^2\fM_{7.5}(\tau,\na)\rr{e^{-|\alpha|\nu \tau}v^\alpha \wt{\th}(\tau)}}_{L^2_{z,v}}e^{\rm_1\nu^\f13\tau}\\
&\times\nnr{e^{-(\rm-\rm_1)\nu^\f13\tau}}_{L^2([\nu^{-1},t])}\ls\eps_1^2.
\ea\]
\[\ba
\sC_5=&\int\sF\rr{\pa_x\fM_{7.5}(t,\na)e^{-|\alpha|\nu t}v^\alpha\cN^\dagger(\wt{\th},\uth+\th^\#)}(t,k,\xi)\overline{\sF\rr{\pa_x\fM_{7.5}(t,\na)\rr{e^{-|\alpha|\nu t}v^\alpha \wt{\th}(t)}}(t,k,\xi)}\\
&\times\rr{k^2+\nu^\f23\bar{\eta}^2+\nu^\f13k\bar{\eta}}d\xi e^{2\rm_1\nu^\f13t}.
\ea\]
\[\ba
\sC_6=&\int\sF\rr{\pa_x\fM_{7.5}(t,\na)e^{-|\alpha|\nu t}v^\alpha\cN^\dagger(\th^\#,\wt{\th})}(t,k,\xi)\overline{\sF\rr{\pa_x\fM_{7.5}(t,\na)\rr{e^{-|\alpha|\nu t}v^\alpha \wt{\th}(t)}}(t,k,\xi)}\\
&\times\rr{k^2+\nu^\f23\bar{\eta}^2+\nu^\f13k\bar{\eta}}d\xi e^{2\rm_1\nu^\f13t}.
\ea\]
\[\ba
\sC_7=&\int\sF\rr{\pa_x\fM_{7.5}(t,\na)e^{-|\alpha|\nu t}v^\alpha\cN^\dagger(\wt{\th},\wt{\th})}(t,k,\xi)\overline{\sF\rr{\pa_x\fM_{7.5}(t,\na)\rr{e^{-|\alpha|\nu t}v^\alpha \wt{\th}(t)}}(t,k,\xi)}\\
&\times\rr{k^2+\nu^\f23\bar{\eta}^2+\nu^\f13k\bar{\eta}}d\xi e^{2\rm_1\nu^\f13t}.
\ea\]
We focus on the nonlinear term $\sC_7$. By
	\[
	e^{-|\alpha|\nu t}v^\alpha\cN^\dagger(\wt{\th},\wt{\th})=\cN^\dagger(\wt{\th},e^{-|\alpha|\nu t}v^\alpha\wt{\th})+\wt{E}(t,z+\tent v)e^{-(\alpha-1)\nu t}v^{\alpha-1}\wt{\th}
	\]
	we split $\sC_{7}$ into the top order part and the remaining part. We use $\sC_{7,top,x}$, $\sC_{7,top,v}$, $\sC_{7,top,hc}$  to denote the term containing $k^2$, $\nu^\f23\bar{\eta}^2$, $\nu^\f13k\bar{\eta}$, respectively. By $\cN^\dagger(\wt{\th},\wt{\th})=\cN^\dagger(\wt{\th},P_{\neq}\wt{\th})+\cN^\dagger(\wt{\th},P_0\wt{\th})$ we split $\sC_{7,top,x}$ into $\sC_{7,top,x,\neq}$ and $\sC_{7,top,x,0}$. We further apply the paraproduct decomposition and split $\sC_{7,top,x,\neq}$ into $\sC_{7,top,x,\neq,LH}$ and $\sC_{7,top,x,\neq,HL}$.
\[\ba
&\sC_{7,top,x,\neq,HL}\ls\sum_{k,\ell}\int\nr{\sF\rr{\fM_{8.5}(t,\na)\wt{\rho}}(t,\ell)}\nr{\sF\rr{\pa^2_x\fM(t,\na)e^{-|\alpha|\nu t}v^\alpha\wt{\th}}(t,k-\ell,\xi-\ell\tent) }\\
&\times\nr{\sF\rr{\vect\pa_x^2\fM_{7.5}(t,\na)\rr{e^{-|\alpha|\nu t}v^\alpha \wt{\th}(t)}}(t,k,\xi)}d\xi e^{2\rm_1\nu^\f13t}
\ls\nnr{\fM_{8.5}(t,\na)\wt{\rho}}_{L^2_z}\nnr{\pa^2_x\fM_1(t,\na)e^{-|\alpha|\nu t}v^\alpha\wt{\th}}_{L^2_{z,v}}\\
&\times\nnr{\vect\pa_x^2\fM_{7.5}(t,\na)\rr{e^{-|\alpha|\nu t}v^\alpha \wt{\th}(t)}}_{L^2_{z,v}}e^{2\rm_1\nu^\f13t}
\ls\nu^{1.01}\nnr{\vect\pa_x^2\fM_{7.5}(t,\na)\rr{e^{-|\alpha|\nu t}v^\alpha \wt{\th}(t)}}_{L^2_{z,v}}^2e^{2\rm_1\nu^\f13t}\\
&+\nu^{-1.01}\nnr{\fM_{8.5}(t,\na)\wt{\rho}}_{L^2_z}^2\nnr{\pa^2_x\fM_1(t,\na)e^{-|\alpha|\nu t}v^\alpha\wt{\th}}_{L^2_{z,v}}^2e^{2\rm_1\nu^\f13t}.
\ea\]
The first term can be absorbed by $\sD_\alpha(t)$. The second term satisfies
\[\int_{\nu^{-1}}^t\nu^{-1.01}\nnr{\fM_{8.5}(\tau,\na)\wt{\rho}}_{L^2_z}^2\nnr{\pa^2_x\fM_1(\tau,\na)e^{-|\alpha|\nu \tau}v^\alpha\wt{\th}}_{L^2_{z,v}}^2e^{2\rm_1\nu^\f13\tau}d\tau
\ls\nu^{-1.01}\eps_1^4
\]
under the hypotheses \eqref{wtrho} and \eqref{wthed}. Similarly for $\sC_{7,top,x,\neq,LH}$,
\[\ba
&\sC_{7,top,x,\neq,LH}\ls\sum_{k,\ell}\int\nr{\sF\rr{\fM_{1}(t,\na)\wt{\rho}}(t,\ell)}\nr{\sF\rr{\pa^2_x\fM_{7.5}(t,\na)e^{-|\alpha|\nu t}v^\alpha\wt{\th}}(t,k-\ell,\xi-\ell\tent) }\\
&\times\nr{\sF\rr{\vect\pa_x^2\fM_{7.5}(t,\na)\rr{e^{-|\alpha|\nu t}v^\alpha \wt{\th}(t)}}(t,k,\xi)}d\xi e^{2\rm_1\nu^\f13t}\ls\nu^{1.01}\nnr{\vect\pa_x^2\fM_{7.5}(t,\na)\rr{e^{-|\alpha|\nu t}v^\alpha \wt{\th}(t)}}_{L^2_{z,v}}^2\\
&\times e^{2\rm_1\nu^\f13t}+\nu^{-1.01}\nnr{\fM_{2}(t,\na)\wt{\rho}}_{L^2_z}^2\nnr{\pa^2_x\fM_{7.5}(t,\na)e^{-|\alpha|\nu t}v^\alpha\wt{\th}}_{L^2_{z,v}}^2e^{2\rm_1\nu^\f13t}.
\ea\]
For $\sC_{7,top,x,0}$,
\[\ba
&\sC_{7,top,x,0}\ls\sum_k\int\nr{\sF\rr{\fM_{8.5}(t,\na)\wt{\rho}}(t,k)}\nr{\sF\rr{v^\alpha\fM_{7.5}(t,0,\emnt\na)P_0\wt{h}}(t,0,\bar{\eta}(t,k,\xi)) }\\
&\times\nr{\sF\rr{\vect\pa_x^2\fM_{7.5}(t,\na)\rr{e^{-|\alpha|\nu t}v^\alpha \wt{\th}(t)}}(t,k,\xi)}d\xi e^{2\rm_1\nu^\f13t}\ls\nnr{\fM_{8.5}(t,\na)\wt{\rho}}_{L^2_z}\nnr{v^\alpha\fM_{7.5}(t,0,\emnt\na)P_0\wt{h}}_{L^2_{v}}\\
&\times\nnr{\vect\pa_x^2\fM_{7.5}(t,\na)\rr{e^{-|\alpha|\nu t}v^\alpha \wt{\th}(t)}}_{L^2_{z,v}}e^{2\rm_1\nu^\f13t}\ls\nu^{1.01}\nnr{\vect\pa_x^2\fM_{7.5}(t,\na)\rr{e^{-|\alpha|\nu t}v^\alpha \wt{\th}(t)}}_{L^2_{z,v}}^2e^{2\rm_1\nu^\f13t}\\
&+\nu^{-1.01}\nnr{\fM_{8.5}(t,\na)\wt{\rho}}_{L^2_z}^2\nnr{v^\alpha\fM_{7.5}(t,0,\emnt\na)P_0\wt{h}}_{L^2_{v}}^2e^{2\rm_1\nu^\f13t}.
\ea\]
The second term satisfies
\[\int_{\nu^{-1}}^{t}\nu^{-1.01}\nnr{\fM_{8.5}(\tau,\na)\wt{\rho}e^{\rm\nu^\f13\tau}}_{L^2_z}^2\nnr{v^\alpha\fM_{7.5}(\tau,0,\emntau\na)P_0\wt{h}}_{L^2_{v}}^2e^{2(\rm_1-\rm)\nu^\f13\tau}d\tau\ls\nu^{-1.01}\eps_1^4.\]
We next consider $\sC_{7,top,v}$.
\[\ba
&\sC_{7,top,v,\neq}\ls\nu^\f23\sum_{k,\ell}\int\nr{\sF\rr{\fM_{8.5}(t,\na)\wt{\rho}}(t,\ell)}\f{1}{|\ell|}\nr{\sF\rr{\pa_x\vect\fM_{7.5}(t,\na)e^{-|\alpha|\nu t}v^\alpha\wt{\th}}(t,k-\ell,\xi-\ell\tent) }\\
&\times\nr{\sF\rr{{\vect}^2\pa_x\fM_{7.5}(t,\na)\rr{e^{-|\alpha|\nu t}v^\alpha \wt{\th}(t)}}(t,k,\xi)}d\xi e^{2\rm_1\nu^\f13t}\ls\nu^{\f53+0.01}\nnr{{\vect}^2\pa_x\fM_{7.5}(t,\na)\rr{e^{-|\alpha|\nu t}v^\alpha \wt{\th}(t)}}_{L^2_{z,v}}^2\\
&\times e^{2\rm_1\nu^\f13t}+\nu^{-1.01}\nnr{\fM_{8.5}(t,\na)\wt{\rho}}_{L^2_z}^2\nu^{\f23}\nnr{\pa_x\vect\fM_{7.5}(t,\na)e^{-|\alpha|\nu t}v^\alpha\wt{\th}}_{L^2_{z,v}}^2e^{2\rm_1\nu^\f13t}.
\ea\]
Under the bootstrap hypotheses, the time integral of the second term is bounded by $\nu^{-1.01}\eps_1^4$.
\[\ba
&\sC_{7,top,v,0}\ls\nu^\f23\sum_k\int\nr{\sF\rr{\fM_{8.5}(t,\na)\wt{\rho}}(t,k)}\f{1}{|k|}\nr{\sF\rr{v^\alpha\na_v\fM_{7.5}(t,0,\emnt\na)P_0\wt{h}}(t,0,\bar{\eta}(t,k,\xi)) }\\
&\times\nr{\sF\rr{{\vect}^2\pa_x\fM_{7.5}(t,\na)\rr{e^{-|\alpha|\nu t}v^\alpha \wt{\th}(t)}}(t,k,\xi)}d\xi e^{2\rm_1\nu^\f13t}\ls\nu^{\f53+0.01}\nnr{{\vect}^2\pa_x\fM_{7.5}(t,\na)\rr{e^{-|\alpha|\nu t}v^\alpha \wt{\th}(t)}}_{L^2_{z,v}}^2\\
&\times e^{2\rm_1\nu^\f13t}+\nu^{-1.01}\nnr{\fM_{8.5}(t,\na)\wt{\rho}}_{L^2_z}^2\nu^\f23\nnr{v^\alpha\na_v\fM_{7.5}(t,0,\emnt\na)P_0\wt{h}}_{L^2_{v}}^2e^{2\rm_1\nu^\f13t}.
\ea\]
Similarly for $\sC_{7,top,hc}$ we have
\[\ba
\sC_{7,top,hc}\ls\nu^{1.01}\nnr{\vect\pa_x^2\fM_{7.5}(t,\na)\rr{e^{-|\alpha|\nu t}v^\alpha \wt{\th}(t)}}_{L^2_{z,v}}^2e^{2\rm_1\nu^\f13t}+\nu^{-1.01}\nnr{\fM_{8.5}(t,\na)\wt{\rho}}_{L^2_z}^2\nu^\f23\\\times\left\{\nnr{\pa_x\vect\fM_{7.5}(t,\na)e^{-|\alpha|\nu t}v^\alpha\wt{\th}}_{L^2_{z,v}}^2+\nnr{v^\alpha\na_v\fM_{7.5}(t,0,\emnt\na)P_0\wt{h}}_{L^2_{v}}^2\right\}e^{2\rm_1\nu^\f13t}.
\ea\]
This completes the estimate of $\sC_{7,top}$. $\sC_{7,re}$ can be bounded by the same argument. $\sC_5$ and $\sC_6$ can be estimated similarly thanks to \eqref{h3} and \eqref{Hg}. Note that $P_0\uh$ has sufficiently high regularity so the term containing $P_0\uh$ can be bounded by the same argument as $\sC_4$.
\[\ba
&\sC_8=-\sum_k\int\nu^\f43k\bar{\eta}^3\left|\sF\rr{\pa_x\fM_{7.5}(t,\na)\rr{e^{-|\alpha|\nu t}v^\alpha \wt{\th}(t)}} \right|^2d\xi e^{2\rm_1\nu^\f13t}\\
&\leq0.3\nu\nnr{\vect\pa_x^2\fM_{7.5}(t,\na)\rr{e^{-|\alpha|\nu t}v^\alpha \wt{\th}(t)}}_{L^2_{z,v}}^2e^{2\rm_1\nu^\f13t}+\f56\nu^\f53\nnr{{\vect}^2\pa_x\fM_{7.5}(t,\na)\rr{e^{-|\alpha|\nu t}v^\alpha \wt{\th}(t)}}_{L^2_{z,v}}^2e^{2\rm_1\nu^\f13t},
\ea\]
which can be absorbed by $\sD_\alpha(t)$.
\[\ba
\sC_{9}=&-\nu\sum_k\int\sF\rr{\pa_x\fM_{7.5}(t,\na)e^{-|\alpha|\nu t}\rr{v^\alpha{\vect}^2-{\vect}^2v^\alpha} \wt{\th}(t)}\overline{\sF\rr{\pa_x\fM_{7.5}(t,\na)\rr{e^{-|\alpha|\nu t}v^\alpha \wt{\th}(t)}}(t,k,\xi)}\\
&\times\rr{k^2+\nu^\f23\bar{\eta}^2+\nu^\f13k\bar{\eta}}d\xi e^{2\rm_1\nu^\f13t}\leq0.01\sD_{\alpha}(t)+1000\nu\sE(t).
\ea\]
We argue as in the estimate of $\mathscr{A}_{10}$ in \eqref{energy}. The second term can be absorbed by $\sD_\beta, \beta\leq\alpha$.
\[\ba
\sC_{10}=&\int\sF\rr{\pa_x\fM_{7.5}(t,\na)e^{-|\alpha|\nu t}v^\alpha\wt\tQ}(t,k,\xi)\overline{\sF\rr{\pa_x\fM_{7.5}(t,\na)\rr{e^{-|\alpha|\nu t}v^\alpha \wt{\th}(t)}}(t,k,\xi)}\\
&\times\rr{k^2+\nu^\f23\bar{\eta}^2+\nu^\f13k\bar{\eta}}d\xi e^{2\rm_1\nu^\f13t}\ls\sE_\alpha(t)^\f12\left\|\fM_{10}(t,\na)\rr{e^{-|\alpha|\nu t}v^\alpha P_{\neq}\wt{\tQ}(t)} \right\|_{\lta L^2_{z,v}}e^{\rm\nu^\f13t+\f12\nu t},
\ea\]

Combining these estimates and summing over $|\alpha|\leq2$, we find a uniform constant $K_2'$ and $K_{21}'(K_1,K_2,K_3)$ depending on $K_1,K_2,K_3$ such that if (\ref{wtrho}-\ref{wth0}) hold then
\[\sE(t)\leq\sE(\nu^{-1})+K_2'\eps_1^2+K_{21}'(K_1,K_2,K_3)\eps_1^2\nu^{-1.01}\rr{A^{-2N_1+32}+\eps_1^2}.\]
By adjusting the constants we obtain
\beq\label{imp2}\ba&\left\{\left\|\pa_x^2\fM_{7.5}(t,\na)\rr{e^{-|\alpha|\nu t}v^\alpha \wt{\th}(t)} \right\|_{\lta L^2_{z,v}}+\nu^\f13\left\|\pa_x\vect\fM_{7.5}(t,\na)\rr{e^{-|\alpha|\nu t}v^\alpha \wt{\th}(t)} \right\|_{\lta L^2_{z,v}}\right\}e^{\rm_1\nu^\f13t}\\
&\leq K_2'\eps_1+K_{21}'(K_1,K_2,K_3)\eps_1\nu^{-0.6}\rr{A^{-N_1+16}+\eps_1}.\ea\eeq
Note that the first term on the right-hand side contains the contribution of $\sE(\nu^{-1})$ and the linear electric field term.

\subsubsection{Improvement of \eqref{wth0}.}
Finally, the zero spatial mode is estimated by the homogeneous Fokker--Planck semigroup. The initial contribution has zero mass and therefore decays at the thermal rate. For the nonlinear velocity-divergence sources, part~(2) of Proposition~\ref{etl} supplies an integrable regularization kernel, and the decay of the nonzero modes makes the source time integral small.
 By \eqref{wtheq} we have
\beq\label{Gk}\ba
&P_0\wt{h}(t)=e^{\nu(t-\nu^{-1})L}P_0\wt{h}(\nu^{-1})-\int_{\nu^{-1}}^te^{\nu(t-\tau)L}\pa_vP_0\rr{\wt{E}\rr{P_{\neq}\uh+h^\#}}d\tau\\
&-\int_{\nu^{-1}}^te^{\nu(t-\tau)L}\pa_vP_0\rr{E^\#\wt{h}}d\tau-\int_{\nu^{-1}}^te^{\nu(t-\tau)L}\pa_vP_0\rr{\wt{E}\wt{h}}d\tau+\int_{\nu^{-1}}^te^{\nu(t-\tau)L}P_0\wt{Q}(\tau)d\tau. \ea\eeq
By \eqref{TT1} and Proposition~\ref{etl} we have
\[\nnr{v^\alpha\<\na_v\>\fM_{7.5}(t,0,\emnt\na_v)e^{\nu(t-\nu^{-1})L}P_0\wt{h}(\nu^{-1})}_{\lta L^2_v}\ent\ls\nnr{v^\alpha\<\na_v\>\fM_{7.5}(\nu^{-1},0,e^{-1}\na_v)P_0\wt{h}(\nu^{-1})}_{\lta L^2_v}\ls\eps_1.\]
By \eqref{twQ} and Proposition~\ref{etl} we have
\[\ba\nnr{v^\alpha\<\na_v\>\fM_{7.5}(t,0,\emnt\na_v)\int_{\nu^{-1}}^te^{\nu(t-\tau)L}P_0\wt{Q}(\tau)d\tau}_{\lta L^2_v}\ent\\
\ls\int^t_{\nu^{-1}}\nnr{v^\alpha\<\na_v\>\fM_{7.5}(\tau,0,e^{-\nu\tau}\na_v)P_0\wt{Q}(\tau)}_{\lta L^2_v}\entau d\tau\ls\eps_1.\ea\]
We focus on the most technical nonlinear term. We first give an estimate of $P_0\rr{\wt{E}\wt{h}}$. Fix $|\alpha|\leq2$. By $\<\xi\>\ls\ent\rr{\<-\ell,\emnt\xi-\ell\tent\>+\<\ell,\ell\tent\>}$ we have
\[\ba&\nr{\sF\left\{v^\alpha\<\na_v\>\fM_{7.5}(t,0,\emnt\na_v)P_0\rr{\wt{E}\wt{h}}\right\} (t,\xi)}\ent\ls\sum_{\ell}\nr{\sF\rr{\fM_{8.5}(t,\na)\wt{\rho}} (t,\ell)}\f{e^{2\nu t}}{|\ell|}\\
&\times\nr{\sF\rr{\fM_{8.5}(t,\na)\rr{e^{-|\alpha|\nu t}v^\alpha\wt{\th}}}(t,-\ell,\emnt\xi-\ell\tent)}.\ea\]
Thus by Sobolev embedding and \eqref{519} we have
\[\ba
&\nnr{v^\alpha\<\na_v\>\fM_{7.5}(t,0,\emnt\na_v)P_0\rr{\wt{E}\wt{h}}}_{\lta L^2_v}\ent\ls e^{\f52\nu t}\nnr{\fM_{8.5}(t,\na)\wt{\rho}}_{L^2_z}\nnr{\fM_{8.5}(\tau,\na)\rr{e^{-|\alpha|\nu\tau}v^\alpha P_{\neq}\wt{\th}} }_{\lta L^2_{z,v}}\\
&\ls e^{\f92\nu t}\nnr{\fM_{8.5}(\tau,\na)\rr{e^{-|\alpha|\nu\tau}v^\alpha P_{\neq}\wt{\th}} }_{\lta L^2_{z,v}}^2\ls e^{\f92\nu t}\eps_1^2 e^{-\rm_1\nu^{\f13}t}.
\ea\]
By Proposition~\ref{etl} (2) we have
\[\ba
&\nnr{v^\alpha\<\na_v\>\fM_{7.5}(t,0,\emnt\na_v)\int_{\nu^{-1}}^te^{\nu(t-\tau)L}\pa_vP_0\rr{\wt{E}\wt{h}}d\tau}_{\lta L^2_v}\ent\leq \int_{\nu^{-1}}^t\nnr{v^\alpha\<\na_v\>\fM_{7.5}(t,0,\emnt\na_v)\right.\\
	&\left.e^{\nu(t-\tau)L}\pa_vP_0\rr{\wt{E}\wt{h}}}_{\lta L^2_v}\ent d\tau
\ls\int_{\nu^{-1}}^t\nnr{v^\alpha\<\na_v\>\fM_{7.5}(\tau,0,\emntau\na_v)P_0\rr{\wt{E}\wt{h}}}_{\lta L^2_v}\entau \f{d\tau}{\sqrt{1-e^{-2\nu(t-\tau)}}}\\
&\ls\eps_1^2\int_{\nu^{-1}}^te^{-0.5\rm_1\nu^{\f13}\tau}\f{d\tau}{\sqrt{1-e^{-2\nu(t-\tau)}}}\ls\eps_1^2.
\ea\]
The second and third terms on the right-hand side of \eqref{Gk} can be bounded by the same argument thanks to \eqref{h3} and \eqref{Hg}. There exist a uniform constant $K'_3$ and $K'_{31}(K_2)$ depending on $K_2$ such that if (\ref{wtrho}-\ref{wth0}) hold then
\beq\label{imp3}\nnr{v^\alpha\<\na_v\>\fM_{7.5}(t,0,\emnt\na_v)P_0\wt{h}(t)}_{\lta L^2_v}\ent\leq K'_3\eps_1+K'_{31}(K_2)\rr{\eps_1^2+A^{-N_1+15}\eps_1}.\eeq

The three improvements now close the bootstrap for Proposition~\ref{prop6}.

\begin{proof}[Proof of Proposition~\ref{prop6}]
	In (\ref{wtrho}-\ref{wth0}) we set $K_1=2K'_1$, $K_2=2K'_2$, $K_3=2K'_3$, where $K'_1, K'_2, K'_3$ are defined in \eqref{imp1}, \eqref{imp2}, and \eqref{imp3}. By standard local well-posedness and energy estimates, together with \eqref{TT1}, (\ref{wtrho}-\ref{wth0}) hold on $[\nu^{-1}, T^*]$ for some $T^*>\nu^{-1}$. Suppose  (\ref{wtrho}-\ref{wth0}) hold on $[\nu^{-1}, T^*]$, $T^*>\nu^{-1}$. By \eqref{imp1}, \eqref{imp2}, and \eqref{imp3} and the fact that $\eps_1\ll\nu^3$, $A^{-N_1+15}\ll\nu^{3}$ the estimates (\ref{wtrho}-\ref{wth0}) hold with $K_i$ replaced by $0.9K_i$, $i=1,2,3$ on $[\nu^{-1}, T^*]$. Thus (\ref{wtrho}-\ref{wth0}) hold on a longer time interval. A standard continuity argument shows that (\ref{wtrho}-\ref{wth0}) hold on $[\nu^{-1},\infty)$.
\end{proof}

\section{Proof of the instability and relaxation results}

We now assemble the estimates to prove Proposition~\ref{P1} and Theorem~\ref{T1}. For Proposition~\ref{P1}, the leading packet gives the lower bound at a critical time, while Fourier localization, the higher-order estimates, and the remainder bounds prevent cancellation. The smoothing multiplier then converts the global bounds into the stated relaxation estimates in the original variables. To obtain Theorem~\ref{T1}, we add the explicit background, shift the initial time, and verify positivity. Finally, separate choices of $\ueps$ balance the smallness of the initial data against the cascade amplification in the Sobolev and Gevrey cases.

\begin{proof}[Proof of Proposition~\ref{P1}]
	\underline{Proof of (1).} We set $\check{T}=t_{\llcorner\rr{D_3\nu A^3}^\f13\lrcorner+1}$, where $t_{\llcorner\rr{D_3\nu A^3}^\f13\lrcorner+1}$ is defined in \eqref{4.1}. We set $c=\f{D_3^{-\f13}}{2}$. By \eqref{lowerb}, \eqref{hnqenergy}, \eqref{hnqu}, and \eqref{TT1} we have
	\[\ba
	&\rr{\sum_k\int_{|\xi-A|\leq20D_2K} \nr{\wh{P_{\neq}\th}(\check{T},k,\xi) }^2d\xi}^\f12\geq\rr{\sum_k\int_{|\xi-A|\leq20D_2K} \nr{\wh{P_{\neq}\th^{(1,1)}}(\check{T},k,\xi) }^2d\xi}^\f12\\
	&-\rr{\sum_k\int_{|\xi-A|\leq20D_2K} \nr{\wh{P_{\neq}\th^{(1,-1)}}(\check{T},k,\xi) }^2d\xi}^\f12-\sum_{n=2}^{A^{N'}}\sum_{|q|\leq n}\nnr{P_{\neq}\th^{(n,q)}(\check{T})}_{L^2_{z,v}}-\nnr{P_{\neq}\wt{\th}(\check{T})}_{L^2_{z,v}}\\
	&\geq K^{-0.6}B(\check{T})\\
	&-e^{-10D_1K}-\sum_{n=2}^{A^{N'}}(2n+1)A^{10N'(2n-1)}B^n(\check{T})-A^{-\f{N_1}{2}A^{N'}}e^{\nu^{-2}}
	\geq K^{-1}B(\check{T})\geq K^{-1}\exp\rr{-1.02c^{-1}K^\f23\log K}A^{-N_1}.
	\ea\]
	Note that we use \eqref{hnqenergy} to control the second term. The last inequality uses \eqref{cnu13}.
	
	\underline{Proof of (2).}  We focus on the more complicated case $t\geq T$. We first estimate $P_{\neq}h$. We choose $\mathsf{r}\ll\rm_1$. By \eqref{Mxj} and the coordinate shift \eqref{coordinateshift} we have
	\[\ba
	&\nnr{\<\na_x,\emnt\na_v-\tent\na_x\>^{N_1-4}e^{\mathsf{r}\min\{1,\nu t\}\rr{|\na_v|+|\na_x|t}} P_{\neq}h(t)}_{L^2_{x,v}}e^{\mathsf{r}\nu^\f13 t}
	\leq\nnr{\<\na\>^{N_1-4}\fM^\f12(t,\na)P_{\neq}\th(t)}_{L^2_{z,v}}\\
	&\times e^{\mathsf{r}\nu^\f13t+\f{\nu}{2}t}
	\leq\nnr{\fM_{N_1-4}(t,\na)P_{\neq} \th^{\#}}_{L^2_{z,v}}e^{\rm\nu^\f13t+\f{\nu}{2}t}+\nnr{\<\na\>^{N_1-4}\fM^\f12(t,\na)P_{\neq}\wt{ \th}}_{L^2_{z,v}}e^{\mathsf{r}\nu^\f13t+\f{\nu}{2}t}.
	\ea\]
	The first term is controlled by \eqref{Dwb}. It remains to estimate the remainder term. By \eqref{Mxj} we have
	\[\log\fM^\f12(t,k,\xi)\gs\nu T\rr{|kt|+|\bar{\eta}(t,k,\xi)|}\geq\f{\nu T}{2}|k|+ \f{\nu T}{2}\rr{\ent\nr{\xi-k\tent}+\nr{k\tent}}\geq\f{\nu T}{2}\rr{|k|+|\xi|},\]
	which implies that $\<k,\xi\>^{N_1-4}\ls \rr{\nu T}^{-N_1+4}\fM^\f12(t,k,\xi)$. Thus
	\[\nnr{\<\na\>^{N_1-4}\fM^\f12(t,\na)P_{\neq}\wt{ \th}}_{L^2_{z,v}}e^{\rm_1\nu^\f13t}\ls\rr{\nu T}^{-N_1+4}\nnr{\fM(t,\na)P_{\neq}\wt{ \th}}_{L^2_{z,v}}e^{\rm_1\nu^\f13t}\ls\rr{\nu^{-\f23}K^{\f13}}^{N_1-4}A^{-\f{N_1}{2}A^{N'}}e^{2\nu^{-2}}.\]
	This, together with \eqref{Dwb}, yields \eqref{zx}.
	
	We next estimate $P_0h$. We first consider the time interval $[T,\nu^{-1}]$. By \eqref{M0} we have
	\[\ba
	&\nnr{\<\na_v\>^{N_1-4}e^{\mathsf{r}\min\{1,\nu t\}|\na_v|} P_{0}h(t)}_{L^2_{v}}\ent\\
	&\ls\nnr{e^{(N_1-4)\nu t}\<\emnt\na_v\>^{N_1-4}\fM^\f12(t,0,\emnt\na_v) P_{0}h(t)}_{L^2_{v}}\ent
	\ls\nnr{\<\emnt\na_v\>^{N_1-4}\right.\\
	&\left.\fM^\f12(t,0,\emnt\na_v) P_{0}h(t)}_{L^2_{v}}\leq\nnr{\fM_{N_1-4}(t,0,\emnt\na_v) P_{0}h^\#(t)}_{L^2_{v}}+\nnr{\<\na_v\>^{N_1-4}\fM^\f12(t,0,\emnt\na_v) P_{0}\wt{h}(t)}_{L^2_{v}}.
	\ea\]
	The first term is controlled by \eqref{Dwb}. For the remainder term,
	\[\ba&
	\nnr{\<\na_v\>^{N_1-4}\fM^\f12(t,0,\emnt\na_v) P_{0}\wt{h}(t)}_{L^2_{v}}\ls(\nu T)^{-N_1+4}\nnr{\fM(t,0,\emnt\na_v) P_{0}\wt{h}(t)}_{L^2_{v}}\\
	&\ls\rr{\nu^{-\f23}K^{\f13}}^{N_1-4}A^{-\f{N_1}{2}A^{N'}}e^{2\nu^{-2}},
	\ea\]
    where we use $\<\xi\>^{N_1-4}\ls(\nu T)^{-N_1+4}\fM^\f12(t,0,\emnt\xi)$ for $t\geq T$, which follows from \eqref{M0}. Thus \eqref{zX} holds on $[T,\nu^{-1}]$. We next consider the time interval $[\nu^{-1},\infty)$.
    \[\ba
    &\nnr{\<\na_v\>^{N_1-4}e^{\mathsf{r}\min\{1,\nu t\}|\na_v|} P_{0}h(t)}_{L^2_{v}}\ent\ls \nnr{e^{2\mathsf{r}|\na_v|}P_{0}h(t) }_{L^2_v}\ent\leq \nnr{\fM(t,0,\emnt\na_v)P_0h(t)}_{L^2_v}\ent\\
    &\ls\nnr{\fM(t,0,\emnt\na_v)P_0h^\#(t)}_{L^2_v}\ent+\nnr{\fM(t,0,\emnt\na_v)P_0\wt{h}(t)}_{L^2_v}\ent\ls K^{-1},
    \ea\]
    where we use $\<\xi\>^{N_1-4}\ls e^{\mathsf{r}|\na_v|}$ and \eqref{Dwb} and \eqref{T1infty}. On $[1,T]$, the explicit formula \eqref{Vb} gives \eqref{zx} and \eqref{zX} directly.
\end{proof}

Proposition~\ref{P1} gives both transient amplification and eventual relaxation for the perturbation. We now add the background to obtain a solution of the original system. The global bounds and positivity are verified first, before choosing the parameters that yield the two instability statements.

\begin{proof}[Proof of Theorem \ref{T1}]
	We set 
	\beq\label{deff}f(t)=h(t+1)+\uh(t+1)\eeq
    for $t\geq0$. We verify that $f(t)$ satisfies the required properties.
    
    \underline{Global relaxation estimates.} We verify that $f(t)$ satisfies the estimates \eqref{edted} and \eqref{edt0}. We first establish two norm comparisons. 
    \beq\label{eqv1} \f12\<k,\emnt\xi-k\tent\>\leq\<k,e^{-\nu(t-1)}\xi-k\f{1-e^{-\nu(t-1)}}{\nu}\>\leq 2\<k,\emnt\xi-k\tent\>.\eeq
    For $0\leq t\leq \nu^{-\f13}$ we have
    \beq\label{eqv}\f{1}{1.1}\<k,\xi-kt\>\leq\<k,\emnt\xi-k\tent\>\leq1.1\<k,\xi-kt\>.\eeq
     Note that
    \[\nr{e^{-\nu(t-1)}\xi-k\f{1-e^{-\nu(t-1)}}{\nu} }\leq e^\nu \nr{\emnt\xi-k\tent}+|k|\f{e^\nu-1}{\nu}\leq e^\nu\rr{\nr{\emnt\xi-k\tent}+|k|},\]
    which implies that $\<k,e^{-\nu(t-1)}\xi-k\f{1-e^{-\nu(t-1)}}{\nu}\>\leq 2\<k,\emnt\xi-k\tent\>$. $\<k,e^{-\nu(t-1)}\xi-k\f{1-e^{-\nu(t-1)}}{\nu}\>\geq \f1 2\<k,\emnt\xi-k\tent\>$ follows by the same argument. Note that for $0\leq t\leq \nu^{-\f13}$,
    \[\nr{\emnt\xi-k\tent}\leq\emnt\nr{\xi-kt}+|k|\nr{t\emnt-\tent}\leq|\xi-kt|+\nu^\f13|k|,\]
    which implies that $\<k,\emnt\xi-k\tent\>\leq1.1\<k,\xi-kt\>$. $\<k,\emnt\xi-k\tent\>\geq\f1{1.1}\<k,\xi-kt\>$ follows by the same argument.
    
    We first prove \eqref{edted}. By \eqref{deff}, it suffices to estimate $h(t+1)$ and $\uh(t+1)$. By \eqref{eqv} and \eqref{zx} we have
    \[\ba
    &\sup_{t\ge0}\nnr{\<\na_x,\emnt\na_v-\tent\na_x\>^{N_1-4}e^{\mathsf{r}\min\{1,\nu t\}\rr{|\na_v|+|\na_x|t}} P_{\neq}h(t+1)}_{L^2_{x,v}}e^{\mathsf{r}\nu^\f13 t}\\
    \ls&\sup_{t\ge1}\nnr{\<\na_x,\emnt\na_v-\tent\na_x\>^{N_1-4}e^{\mathsf{r}\min\{1,\nu t\}\rr{|\na_v|+|\na_x|t}} P_{\neq}h(t)}_{L^2_{x,v}}e^{\mathsf{r}\nu^\f13 t}\ls K^{-1}.
    \ea\]
    By \eqref{963} and \eqref{Mxj} we have
    \[\ba
    &\sup_{t\ge0}\nnr{\<\na_x,\emnt\na_v-\tent\na_x\>^{N_1-4}e^{\mathsf{r}\min\{1,\nu t\}\rr{|\na_v|+|\na_x|t}} P_{\neq}\uh(t+1)}_{L^2_{x,v}}e^{\mathsf{r}\nu^\f13 t}\\
	&\ls\sup_{t\ge1}\nnr{\<\na_x,\emnt\na_v-\tent\na_x\>^{N_1-4}\right.\\
    	&\left.e^{\mathsf{r}\min\{1,\nu t\}\rr{|\na_v|+|\na_x|t}} P_{\neq}\uh(t)}_{L^2_{x,v}}e^{\mathsf{r}\nu^\f13 t}
    \ls\sup_{t\ge1}\nnr{\fM_{N_1-4}(t,\na)P_{\neq}\uth(t)}_{L^2_{z,v}}e^{\mathsf{r}\nu^\f13 t+\f12\nu t}\ls\ueps.
    \ea\]    
    This proves \eqref{edted}. We next prove \eqref{edt0}. By \eqref{zX} we have
    \[\sup_{t\ge0}\nnr{\<\na_v\>^{N_1-4}e^{\mathsf{r}\min\{1,\nu t\}|\na_v|} P_{0}h(t+1)}_{L^2_{v}}\ent\leq\sup_{t\ge1}\nnr{\<\na_v\>^{N_1-4}e^{\mathsf{r}\min\{1,\nu t\}|\na_v|} P_{0}h(t)}_{L^2_{v}}\ent\ls K^{-1}.\]
    By \eqref{957} and \eqref{M0} and Proposition~\ref{etl} we have
    \[\nnr{\<\na_v\>^{N_1-4}e^{\mathsf{r}\min\{1,\nu t\}|\na_v|} P_{0}\uh(t+1)}_{L^2_{v}}\ent\ls \ueps.\]
    This proves \eqref{edt0}.
    
    \underline{Positivity.} Since \eqref{VPF} preserves positivity, it suffices to prove $f(0)+\sM(v)>0$. By \eqref{udlh}, \eqref{h1T} and \eqref{vfdecay} we have
    \[\ba f(0)+\sM(v)\geq(1-\udl{\eps})\sM(v)+\udl{\eps}\sP(v)-3\udl{\eps}C\int e^{-|v-u|^{0.99}}\f{1}{(2\pi)^\f12(1-e^{-2\nu})^\f12}\exp\rr{-\f{u^2}{2(1-e^{-2\nu})}}du\\
    -3\eps C\int e^{-|v-u|^{0.99}}\f{1}{(2\pi)^\f12(1-e^{-2\nu})^\f12}\exp\rr{-\f{u^2}{2(1-e^{-2\nu})}}du
    \geq(1-\udl{\eps})\sM(v)+\udl{\eps}\sP(v)-6\udl{\eps}Ce^{-|v|^{0.99}}.
    \ea\]
    On the region $\{|v|\leq\rr{2\log\f{1}{10^5C\udl{\eps}}}^\f12\}$ it holds that $6\udl{\eps}Ce^{-|v|^{0.99}}\leq\f12\sM(v)$.\linebreak On the region $\{|v|\geq\rr{2\log\f{1}{10^5C\udl{\eps}}}^\f12\}$ it holds that $6\udl{\eps}Ce^{-|v|^{0.99}}\leq\udl{\eps}\sP(v)$.	
    
    It remains to choose $\ueps$ so that the initial perturbation is small in the prescribed topology while the echo lower bound gives the desired norm inflation.
    
    \underline{Proof of instability for Sobolev initial data.}For $N\geq10^4$, define
    \[
    K_*=\frac{\gamma\bigl(\frac{2N}{3}1.01^{3/2}\nu^{-1/2}\bigr)}
    {1.01^{3/2}\nu^{-1/2}},
    \qquad
    \ueps=\nu^{1/3}\lfloor K_*\rfloor^{7/3},
    \]
    where $\gamma$ is defined in \eqref{gamma}. By \eqref{KAT},
    $K=\lfloor K_*\rfloor$ and $A=1.01\nu^{-1/3}K^{2/3}$.
    By \eqref{gamma},
    \[
    K_*=\frac{2N}{3}\log\bigl(1.01^{3/2}\nu^{-1/2}K_*\bigr)
    =N\log\bigl(1.01\nu^{-1/3}K_*^{2/3}\bigr).
    \]
    Choose $\nu(N)>0$ sufficiently small that $K\geq N$ whenever
    $0<\nu<\nu(N)$. On $[K,K_*]$, the function
    $F(x)=x-N\log(1.01\nu^{-1/3}x^{2/3})$ satisfies
    $0<F'(x)=1-\frac{2N}{3x}<1$. Since $F(K_*)=0$,
    \[
    0\leq N\log A-K=\int_K^{K_*}F'(x)\,dx\leq K_*-K<1.
    \]
    
    We first prove \eqref{initial}.
    \[\ba
    &\nnr{\<\na\>^Nf(0)}_{L^2_{x,v}}\leq\nnr{\<\na\>^NP_0\uh(1)}_{L^2_{x,v}}+\nnr{\<\na\>^NP_{\neq}\uh(1)}_{L^2_{x,v}}+\nnr{\<\na\>^NP_{\neq} h(1)}_{L^2_{x,v}}\\
    &\ls_N\ueps+\nnr{\<\na\>^NP_{\neq} \th(1)}_{L^2_{z,v}}+\nnr{\<\na\>^NP_{\neq} \uth(1)}_{L^2_{z,v}}\ls_N \ueps+\eps A^N\ls_N \nu^\f13\rr{\log\f1\nu}^\f73.\ea\]
    The second inequality uses \eqref{eqv1} and the coordinate shift \eqref{coordinateshift}. More explicitly, \[\nnr{\<\na\>^NP_{\neq}\uh(1)}_{L^2_{x,v}}\ls\nnr{\<\na_x, e^{-\nu}\na_v-\f{1-e^{-\nu}}{\nu}\na_x\>^NP_{\neq}\uh(1)}_{L^2_{x,v}}\leq\nnr{\<\na\>^NP_{\neq} \th(1)}_{L^2_{z,v}}e^{\f\nu2}. \]
    The third inequality uses \eqref{963} and \eqref{Vb}. Note that by \eqref{geqv} we have $\ueps\sim \nu^\f13\rr{\f{N}{3}\log\f{1.01^3\times4N^2}{9\nu}}^\f73$. By $|K-N\log A|\leq1$ and $A^{-N_1}e^{-3K}<\eps<A^{-N_1}e^{-1.5K}$ we have $\eps A^N\ls \eps e^K\leq A^{-N_1}\ll \nu$.

    We next prove \eqref{checkT}. For $\check{T}$ defined in Proposition~\ref{P1}(1), we have
    \[\ba
    &\nnr{\<\na_x,\na_v-(\check{T}-1)\na_x\>^{N_1+4}P_{\neq}f(\check{T}-1) }_{L^2_{x,v}}\gs \nnr{\<\na_x,e^{-\nu\check{T}}\na_v-\f{1-e^{-\nu\check{T}} }{\nu}\na_x\>^{N_1+4}P_{\neq}f(\check{T}-1) }_{L^2_{x,v}}\\
    &\gs\nnr{\<\na\>^{N_1+4}P_{\neq}\th(\check{T})}_{L^2_{z,v}}-\nnr{\<\na\>^{N_1+4}P_{\neq}\uth(\check{T})}_{L^2_{z,v}}\gs A^4K^{-1}\exp\rr{-1.02c^{-1}K^\f23\log K}\gs\nu^{-1}K.
    \ea\]
    The first inequality uses \eqref{eqv} and \eqref{eqv1}. The third inequality uses \eqref{Zx}. Note that by $|K-N\log A|\leq1$ we have $\nr{K^\f23\log K-\rr{N\log A}^\f23\log \rr{N\log A}}\ls K^{-\f13}\log K$, which implies that $K^\f23\log K\ll\log A$. Replace $\check{T}-1$ by $\check{T}$ and we obtain \eqref{checkT}.
    
    \underline{Proof of instability for Gevrey initial data.} For $0<s<\frac13$, define
    \[
    K_*=\bigl(3^{1/s}\times1.01\nu^{-1/3}\bigr)^{\frac{3s}{3-2s}},
    \qquad
    \ueps=\nu^{1/3}\lfloor K_*\rfloor^{7/3}.
    \]
    By \eqref{KAT}, $K=\lfloor K_*\rfloor$ and $A=1.01\nu^{-1/3}K^{2/3}$.
    The definition of $K_*$ gives
    \[
    \frac{K_*}{3}=\bigl(1.01\nu^{-1/3}K_*^{2/3}\bigr)^s,
    \qquad
    A^s=\frac{K_*}{3}\left(\frac{K}{K_*}\right)^{2s/3}.
    \]
    Since $0\leq K/K_*\leq1$ and $0<2s/3<1$, we have
    $K/K_*\leq(K/K_*)^{2s/3}\leq1$. It follows that
    \[
    0\leq A^s-\frac{K}{3}\leq\frac{K_*-K}{3}<\frac13.
    \]
    
    We first prove \eqref{checkT2}. By \eqref{eqv1}, \eqref{963}, and \eqref{Vb} we have
    \[\ba
    &\nnr{e^{\<\na\>^s}f(0)}_{L^2_{x,v}}\leq\nnr{e^{\<\na\>^s}P_0\uh(1)}_{L^2_{x,v}}+\nnr{e^{\<\na\>^s}P_{\neq}\uh(1)}_{L^2_{x,v}}+\nnr{e^{\<\na\>^s}P_{\neq} h(1)}_{L^2_{x,v}}\\
    &\ls\ueps+\nnr{e^{1.3\<\na\>^s}P_{\neq} \th(1)}_{L^2_{z,v}}+\nnr{e^{1.3\<\na\>^s}P_{\neq} \uth(1)}_{L^2_{z,v}}\ls \ueps+\eps e^{1.4A^s}\ls \nu^{\f{1-3s}{3-2s}}.\ea\]
    Note that $\eps e^{1.4A^s}\ls A^{-N_1}e^{-K}\ll\nu$.
    
    We next prove \eqref{checkT2}. By \eqref{eqv}, \eqref{eqv1}, and \eqref{Zx} we have
    \[\ba
    &\nnr{e^{10c^{-1}\<\na_x,\na_v-(\check{T}-1)\na_x\>^{\f{2s}{3}}\log\<\na_x,\na_v-(\check{T}-1)\na_x\>}P_{\neq}f(\check{T}-1) }_{L^2_{x,v}}\gs \nnr{\exp\left\{5c^{-1}\<\na_x,e^{-\nu\check{T}}\na_v-\f{1-e^{-\nu\check{T}} }{\nu}\na_x\>^{\f{2s}{3}}\right.\right.\\
    &\left.\left.\times\log\<\na_x,e^{-\nu\check{T}}\na_v-\f{1-e^{-\nu\check{T}} }{\nu}\na_x\>\right\}P_{\neq}f(\check{T}-1) }_{L^2_{x,v}}
    \gs\nnr{e^{5c^{-1}\<\na\>^{\f{2s}{3}}\log\<\na\>}P_{\neq}\th(\check{T})}_{L^2_{x,v}}\\
    &- \nnr{e^{5c^{-1}\<\na\>^{\f{2s}{3}}\log\<\na\>}P_{\neq}\uth(\check{T})}_{L^2_{x,v}}
    \gs e^{4c^{-1}A^{\f{2s}{3}}\log A-1.02c^{-1}K^\f23\log K}e^{-\log K-N_1\log A}\gs\f1{\nu^2}.
    \ea\]
    Note that $|\f13K-A^s|\leq1$ implies $\nr{N_1\log A-N_1\log\rr{\f13K}^\f1s}\leq1$, which implies that $N_1\log A+\log K\ll K^\f23$. Replace $\check{T}-1$ by $\check{T}$ and we obtain \eqref{checkT2}.
\end{proof}

	\appendix
\section{Appendix A: Auxiliary estimates for resonance and smoothing}
This appendix contains two elementary estimates used in the weighted analysis. Lemma~\ref{resonance} shows that the time variation of $\fA$ compensates for the resonant interaction coefficient; this is the small kernel bound used to close the density estimate in Proposition~\ref{3.2}. Lemma~\ref{lm7} gives a lower bound for an integral along a characteristic and is used in Appendix~B to identify the regularity measured by the smoothing multiplier at and beyond the thermal time.

We first estimate the multiplier ratio on the resonant support. The proof separates the region where the interaction coefficient is already small from the region where the decrease of $\fA$ supplies the required smallness.

\begin{lem}\label{resonance}
	Recall $\fA(t,k,\xi)$ defined in \eqref{deffA}.	For any $k\geq1$, $k\in\Z$, $t\geq T$, $D_1>10^3$, we have
	\[\ueps\f{1-\emnt}{\nu k^2}\vv{1}_{|k\tent-(k+1)\tentau|\leq1}\f{\fA(t,k,k\tent)}{\fA(\tau,k,k\tent)}\ls D_1^{-\f14}.\]
\end{lem}
\begin{proof}
	To simplify the notation we set $\mathsf{t}=\tent$, $\uptau=\tentau$. We only need to prove 
	\beq\label{734}\ba
	&\f{\ueps \mathsf{t}}{k^2}\vv{1}_{|k \mathsf{t}-(k+1)\uptau|\leq1}\exp \left\{D_1\left[\rr{1+\ueps \mathsf{t}\<\ueps|k,k\mathsf{t}|\>^{-\f23}}^{-0.01}\right.\right.\\
	&\quad \left.\left.-\rr{1+\ueps \uptau\<\ueps|k,k\mathsf{t}|\>^{-\f23}}^{-0.01}\right]\<\ueps|k,k\mathsf{t}|\>^{\f13}\right\}\ls D_1^{-\f14}.
	\ea\eeq
	We split into two cases.
	
	\underline{Case 1: $\ueps \mathsf{t}\<\ueps|k,k\mathsf{t}|\>^{-\f23}\leq D_1^{-\f14}$.} In this case either $\ueps \mathsf{t}\ls D_1^{-\f14}$ or $\ueps \mathsf{t}\ls D_1^{-\f14}\rr{\ueps|k,k\mathsf{t}|}^{\f23}$ holds. In the first case, \eqref{734} holds immediately. The second case implies that $\f{\ueps^\f13\mathsf{t}^\f13}{k^\f23}\ls D_1^{-\f14}$, which implies that $\f{\ueps \mathsf{t}}{k^2}\ls D_1^{-\f34}$. 
	
	\underline{Case 2: $\ueps \mathsf{t}\<\ueps|k,k\mathsf{t}|\>^{-\f23}\geq D_1^{-\f14}$.} On $|k\mathsf{t}-(k+1)\uptau|\leq1$ we have
	\[\ba
	&D_1\left[\rr{1+\ueps \uptau\<\ueps|k,k\mathsf{t}|\>^{-\f23}}^{-0.01}-\rr{1+\ueps \mathsf{t}\<\ueps|k,k\mathsf{t}|\>^{-\f23}}^{-0.01}\right]\<\ueps|k,k\mathsf{t}|\>^{\f13}\\
	&\geq \f{0.01D_1(\mathsf{t}-\uptau)\ueps\<\ueps|k,k\mathsf{t}|\>^{-\f13}}{\rr{1+\ueps \mathsf{t}\<\ueps|k,k\mathsf{t}|\>^{-\f23}}^{1.01}}\\
	\geq& \f{10^{-3}D_1\f{\ueps \mathsf{t}^2}{|k,k\mathsf{t}|}\<\ueps|k,k\mathsf{t}|\>^{-\f13}}{\rr{1+\ueps \mathsf{t}\<\ueps|k,k\mathsf{t}|\>^{-\f23}}^{1.01}}\\
	&=10^{-3}D_1\f{\<\ueps|k,k\mathsf{t}|\>}{\ueps|k,k\mathsf{t}|}\f{\rr{\ueps \mathsf{t}\<\ueps|k,k\mathsf{t}|\>^{-\f23}}^2}{\rr{1+\ueps \mathsf{t}\<\ueps|k,k\mathsf{t}|\>^{-\f23}}^{1.01}}\\
	&\geq10^{-3}D_1\f{\rr{\ueps \mathsf{t}\<\ueps|k,k\mathsf{t}|\>^{-\f23}}^2}{\rr{1+\ueps \mathsf{t}\<\ueps|k,k\mathsf{t}|\>^{-\f23}}^{1.01}}\\
	\geq&10^{-4}D_1^\f23\rr{\ueps \mathsf{t}\<\ueps|k,k\mathsf{t}|\>^{-\f23}}^{0.99}.
	\ea\]
	The first inequality uses the mean value theorem. The last inequality uses $\ueps \mathsf{t}\<\ueps|k,k\mathsf{t}|\>^{-\f23}\geq D_1^{-\f14}$. We consider the following subcases. \uwave{Subcase 2.1: $\ueps|k,k\mathsf{t}|\leq1$.} In this subcase we have
	\[\ueps\f{1-\emnt}{\nu k^2}\vv{1}_{|k\tent-(k+1)\tentau|\leq1}\f{\fA(t,k,k\tent)}{\fA(\tau,k,k\tent)}\leq\f{\ueps \mathsf{t}}{k^2}e^{-10^{-5}D_1^\f23\rr{\ueps \mathsf{t}}^{0.99}}\ls D_1^{-\f23}.\]
	\uwave{Subcase 2.2: $\ueps|k,k\mathsf{t}|\geq1$.} In this subcase we have
	\[\ueps\f{1-\emnt}{\nu k^2}\vv{1}_{|k\tent-(k+1)\tentau|\leq1}\f{\fA(t,k,k\tent)}{\fA(\tau,k,k\tent)}\leq\f{\ueps \mathsf{t}}{k^2}e^{-10^{-5}D_1^\f23\rr{\f{\ueps \mathsf{t}}{k^2}}^{0.33}}\ls D_1^{-2}.\]
\end{proof}

The next estimate addresses a different issue: cancellation along a characteristic must not destroy the lower bound for the accumulated smoothing weight. After a change of variables, the integral in Lemma~\ref{lm7} is the scalar expression that appears in the large-time analysis of $\fM$. Its proof separates the cases in which the integrand keeps or changes sign.

\begin{lem}\label{lm7}
	It holds that
	\[\int_{e^{-u}}^1\nr{x+\f{y}{s}}ds\gs|x|+|y|u\]
	for any $x,y\in\R$, $u\geq1$.
\end{lem}
\begin{proof}
	The case $y=0$ is immediate. By symmetry, we may therefore assume $y>0$. We split into two cases.
	
	\underline{Case 1: $x+y\geq0$.} 
	\[\int_{e^{-u}}^1\nr{x+\f{y}{s}}ds=(1-e^{-u})x+uy\gs|x|+uy.\]
	
	\underline{Case 2: $x+y<0$.} 
	\[\int_{e^{-u}}^1\nr{x+\f{y}{s}}ds=\left\{\ba&(1-e^{-u})(-x)-uy,\quad u\leq\log\f{-x}{y}\\&(1+e^{-u})(-x)+\rr{u-2-2\log\f{-x}{y}}y,\quad u\geq\log\f{-x}{y}  \ea\right.\]
	We set $v=\f{-x}{y}\geq1$. We only need to prove for any $u,v\geq1$,
	\[u+v\ls\left\{\ba&(1-e^{-u})v-u,\quad 1\leq u\leq\log v\\&(1+e^{-u})v+u-2-2\log v,\quad u\geq\max\{1,\log v\}.\ea\right.\]
	We split into two subcases.
	
	\uwave{Subcase 2.1: $1\leq v\leq e$.} We only need to prove $(1+e^{-u})v+u-2-2\log v\gs u$ for $u\geq1$ since $v\ls u$. We set
	\[G_1(u)=(1+e^{-u})\f vu+\f{u-2-2\log v}{u}.\]
	We next show that $G_1(u)>0.001$ for $u\geq1$, $1\leq v\leq e$. Direct differentiation gives
	\[G'_1(u)=\f{v}{u^2}\rr{-e^{-u}(1+u)+2\f{\log v+1}{v}-1}.\]
	We split into two subcases.
	\begin{itemize}
		\item Subcase 2.1.1: $1+\f2e\leq2\f{\log v+1}{v}\leq2$. In this case $G'_1(u)\geq0$, which implies that $G_1(u)\geq G_1(1)=\rr{1+\f1e}v-2\log v-1>0.01$.
		\item Subcase 2.1.2: $\f4e\leq2\f{\log v+1}{v}<1+\f2e$. In this case there exists a unique $u_0>1$ such that $G'_1(u_0)=0$, which is the minimum point of $G_1(u)$ on $[1,\infty)$. We have
		\[G_1(u)\geq G_1(u_0)=1-\f{2\rr{\log v+1}-v}{u_0+1}>1-\f{2\rr{\log v+1}-v}{2}\geq1-\log2.\]
		The second inequality uses $e^{-u_0}=\f{1}{u_0+1}\rr{2\f{\log v+1}{v}-1}$. The third inequality uses $u_0>1$. The fourth step follows from basic calculus.
	\end{itemize}
	
	\uwave{Subcase 2.2: $v>e$.} We first prove $(1-e^{-u})v-u\gs v+u$ on $1\leq u\leq\log v$. We only need to prove $(1-e^{-u})v-u\gs u$ since $u\leq\log v\leq v$. We set
	\[G_2(u)=(1-e^{-u})v-u.\]
	Then we have $G_2'(u)=e^{-u}v-1\geq0$ on $1\leq u\leq\log v$. Thus $G_2(u)\geq G_2(1)=(1-e^{-1})v-1\geq0.01v$.
	
	We next prove $(1+e^{-u})v+u-2-2\log v\gs v+u$ on $u\geq\log v$. First we prove $(1+e^{-u})v+u-2-2\log v\gs v$. We set
	\[G_3(u)=(1+e^{-u})v+u-2-2\log v.\]
	Then we have $G'_3(u)=1-ve^{-u}\geq0$ on $u\geq\log v$. Thus $G_3(u)\geq G_3(\log v)=v-\log v-1\geq0.01v$. We next prove $(1+e^{-u})v+u-2-2\log v\gs u$. Again we set
	\[G_1(u)=(1+e^{-u})\f vu+\f{u-2-2\log v}{u}.\quad \text{Then}\quad G'_1(u)=\f{v}{u^2}\rr{-e^{-u}(1+u)+2\f{\log v+1}{v}-1}.\]
	We split into two subcases.
	\begin{itemize}
		\item Subcase 2.2.1: $1<2\f{\log v+1}{v}\leq\f4e$. In this case there exists a unique $u_0>\log v$ such that $G'_1(u_0)=0$, which is the minimum point of $G_1(u)$ on $[\log v,\infty)$. We have
		\[G_1(u)\geq G_1(u_0)=1-\f{2\rr{\log v+1}-v}{u_0+1}>\f{v}{\log v+1}-1>0.01.\]
		\item Subcase 2.2.2: $2\f{\log v+1}{v}\leq1$. In this case $G'_1(u)\leq0$. Thus $G_1(u)\geq G_1(+\infty)=1$.
	\end{itemize}
\end{proof}
\begin{lem}\label{lA3}
	Suppose $\nu$ and $\ueps$ satisfy the conditions in Proposition \ref{P1}.  Let $0<A'\le\f13\nu^{-\f13}$, $0<\rm<10^{-5}$. Recall the definition of  $t_i'$ for $i\in\Z_{\ge0}$ in \eqref{401}. For $k_1, k_2\in\Z_+$, $k_1\le k_2$, it holds that
	\beq\label{653}\prod_{i=k_1}^{k_2}\f{\ueps A'}{i^3}e^{-\f{\nu {A'}^3}{3i^3}+2\rm\nu^\f13(t'_{i-1}-t'_i)}\le\nu^{-\f13}\prod_{i=k_1}^{k_2}\max\{\f{\ueps A'}{i^3}e^{-\f{\nu {A'}^3}{3i^3}},1\}.\eeq
\end{lem}
\begin{proof} We note that $2\nu^\f13(t'_{i-1}-t'_i)\le5\f{\nu^\f13A'}{i^2}$.
	
	For the case $A'\le\nu^{-\f13}$, we have
	\[\prod_{i=k_1}^{k_2}e^{2\rm\nu^\f13(t'_{i-1}-t'_i)}\le e^{2\rm\nu^\f13 t'_0}\le e^{8\rm}.\]
	
	Then we consider the case $\nu^{-\f13}\le A'<\sqrt{\f{6\ueps}{e\nu}}$. Recall the definition of $\gamma(x)$ in \eqref{gamma}. We set 
	\[i_1=\llcorner\rr{\f{\ueps A'}{\gamma(\f{6\ueps}{\nu {A'}^2})}}^\f13\lrcorner.\]
	By \eqref{geqv} we have
	\[\prod_{i\ge i_1}e^{5\rm\nu^\f13\f{A'}{i^2}}\le e^{10^2\rm\nu^\f13\f{A'}{i_1}}\le e^{10^3\rm\rr{\log\f{6\ueps}{\nu {A'}^2}}^\f13 }\le e^{10^3\rm\rr{\log\f{6\ueps}{\nu^\f13}}^\f13 }\le\nu^{-\f13}.\]
	For $i<i_1$, since $\f{\ueps A'}{i^3}e^{-\f{\nu {A'}^3}{6i^3}}\le1$, $-\f{\nu {A'}^3}{6i^3}+5\rm\nu^\f13\f{A'}{i^2}\le0$ on this region, we have
	\[\f{\ueps A'}{i^3}e^{-\f{\nu {A'}^3}{3i^3}+2\rm\nu^\f13(t'_{i-1}-t'_i)}\le1.\]
	Thus \eqref{653} holds.
	
	Then we consider the case $A'\ge\sqrt{\f{6\ueps}{e\nu}}$. It holds that
	\[\prod_{i>\f{(\nu^\f13A')^2}{5\rm}}e^{5\rm\nu^\f13\f{A'}{i^2}}\le e^{100\rm\nu^\f13A'\f{5\rm}{(\nu^\f13A')^2}}=e^{\f{500\rm^2}{\nu^\f13A'}}\le2.\]
	For $i\le\f{(\nu^\f13A')^2}{5\rm}$, since $\f{\ueps A'}{i^3}e^{-\f{\nu {A'}^3}{6i^3}}\le1$, $-\f{\nu {A'}^3}{6i^3}+5\rm\nu^\f13\f{A'}{i^2}\le0$ on this region, we have
	\[\f{\ueps A'}{i^3}e^{-\f{\nu {A'}^3}{3i^3}+2\rm\nu^\f13(t'_{i-1}-t'_i)}\le1.\]
	Thus \eqref{653} holds.
\end{proof}

\section{Appendix B: Properties of the multipliers and semigroups}
We collect the multiplier and semigroup bounds used in the preceding sections, stating them in dimension $d$. Proposition~\ref{mproperty} gives the product, derivative, commutator, and propagation bounds for $\fM$, and identifies the exponential Fourier weight it supplies in the original variables. Proposition~\ref{prop:semigroup-symbol} records the explicit formulas, derivative bounds, and decay estimate for the Fourier semigroup. Proposition~\ref{etl} treats the zero spatial mode, providing thermal decay for mean-zero data and regularization for differentiated sources. These results justify the weighted Duhamel estimates and the global remainder argument.

We begin with $\fM$. Its submultiplicativity permits the weight to be distributed across bilinear interactions, its derivative and difference bounds control velocity moments and commutators, and its compatibility with the semigroup allows the smoothing weight to be propagated in time.

\begin{prop}[Properties of the smoothing multiplier]\label{mproperty}
	For 
	\[\fM(t,k,\xi)=\exp\rr{\mathrm{r}\nu\int_0^t\rr{e^{-2\nu(t-s)}+|\bar{\eta}(s,k,\xi)|^2}^\f12ds},\]
	the following properties hold:
	\begin{itemize}
		\item[(1)] It holds that
	    \beq\label{M1}\fM(t,k+\ell,\xi+\eta)\leq\fM(t,k,\xi)\fM(t,\ell,\eta).\eeq
	    \item[(2)] For every $\alpha\in\N^d$, we have
	    \beq\label{Mds}\nr{\pa_\xi^\alpha\fM(t,k,\xi)}\ls_{\alpha}e^{|\alpha|\nu t}\fM(t,k,\xi).\eeq
	    \item[(3)] It holds that
	    \beq\label{Mjhz} \nr{\fM(t,k,\xi)-\fM(t,\ell,\eta)}\leq\rm\nu\int_0^t\nr{\bar{\eta}(s,k-\ell,\xi-\eta)}ds\ \max\{\fM(t,k,\xi),\fM(t,\ell,\eta)\}.\eeq
	    \item[(4)] For $k\neq0$, we have
	    \beq\label{MSk}\fM(t,k,\xi)S^\rm(t,\tau,k,\xi)\leq\fM(\tau,k,\xi)\exp\rr{\rm(1-\emnt+4\nu)},\eeq
	    \beq\label{MSkt}\fM(t,k,k\tent)\leq\fM(\tau,k,k\tent)\exp\rr{\rm\min\{\nu|k|(t-\tau)^2,|k|(t-\tau)\}+1}.\eeq
	    For $k=0$, we have
	    \beq\label{MS0}\fM(t,0,\xi)S^\rm(t,\tau,0,\xi)\ls\fM(\tau,0,\xi),\eeq
	    \beq\label{M0}\fM(t,0,\emnt\xi)=\exp\rr{\rm\rr{1-\emnt}\<\xi\>}.\eeq
	    \item[(5)] It holds that
	    \beq\label{Mxj}\log\fM(t,k,\bar{\eta}(-t,k,\xi))\sim \rm\min\{1,\nu t\}\rr{|\xi|+|kt|+1}.\eeq
	\end{itemize}
\end{prop}
\begin{proof}
We prove the properties in the order stated. The product and difference estimates follow from the structure of the exponential weight, and differentiation is handled by the Fa\`a di Bruno formula. Comparison with the dissipative integral gives the propagation bounds. The final lower bound is proved separately before and after the thermal time, using Lemma~\ref{lm7} in the latter regime.

	We begin with part~(1).	Note that $|\bar{\eta}(t,k+\ell,\xi+\eta)|\leq|\bar{\eta}(t,k,\xi)|+|\bar{\eta}(t,\ell,\eta)|$. Thus \eqref{M1} holds. For part~(2), the Fa\`a di Bruno formula gives
	\[\ba
	&\nr{\pa_\xi^\alpha\fM(t,k,\xi)}\ls_\alpha \sum_{\substack{\sum m_\beta\beta=\alpha,\\m_\beta\in\N,\beta\in\N^d}}\fM(t,k,\xi)\prod_\beta\rr{\nu\int_0^te^{-\nu(t-s)}\nr{\pa_\xi^\beta\<e^{\nu(t-s)}\bar{\eta}(s,k,\xi)\>}ds }^{m_\beta}\\
	\ls_\alpha& \sum_{\substack{\sum m_\beta\beta=\alpha,\\m_\beta\in\N,\beta\in\N^d}}\fM(t,k,\xi)\prod_\beta\rr{\nu\int_0^te^{-\nu(t-s)}e^{|\beta|\nu t}ds }^{m_\beta}\ls_\alpha \fM(t,k,\xi)e^{|\alpha|\nu t}.
	\ea\]
	For part~(3), the mean value theorem gives
	\[\ba&\nr{\fM(t,k,\xi)-\fM(t,\ell,\eta)}\leq\rm\nu\int_0^te^{-\nu(t-s)}\nr{\<e^{\nu(t-s)}\bar{\eta}(s,k,\xi)\>-\<e^{\nu(t-s)}\bar{\eta}(s,\ell,\eta)\>}ds\\&\times \max\{\fM(t,k,\xi),\fM(t,\ell,\eta)\}
	\leq \rm\nu\int_0^t\nr{\bar{\eta}(s,k-\ell,\xi-\eta)}ds\ \max\{\fM(t,k,\xi),\fM(t,\ell,\eta)\}.
	\ea\]
	For part~(4), note that
	\[\ba
	&\nu\int_0^t\rr{e^{-2\nu(t-s)}+|\bar{\eta}(s,k,\xi)|^2}^\f12ds\leq\nu\int_0^t e^{-\nu(t-s)}+|\bar{\eta}(s,k,\xi)| ds=1-e^{-\nu t}+\nu\int^\tau_0|\bar{\eta}(s,k,\xi)| ds\\
	&+\nu\int_\tau^t|\bar{\eta}(s,k,\xi)| ds
	=1-e^{-\nu t}+\nu\int^\tau_0|\bar{\eta}(s,k,\xi)| ds+\nu\int^t_\tau|\bar{\eta}(s,k,\xi)|\rr{\vv{1}_{|\bar{\eta}(s,k,\xi)|\geq1}+\vv{1}_{|\bar{\eta}(s,k,\xi)|\leq1}} ds\\
	\leq&1-e^{-\nu t}+\nu\int^\tau_0|\bar{\eta}(s,k,\xi)| ds+\nu\int^t_\tau|\bar{\eta}(s,k,\xi)|^2 ds+\nu\int^t_\tau|\bar{\eta}(s,k,\xi)|\vv{1}_{|\bar{\eta}(s,k,\xi)|\leq1} ds
	\ea\]
	Now we estimate the last term. For $k\neq0$, $|\bar{\eta}(s,k,\xi)|=|e^{\nu s}\rr{\xi-\f{k}{\nu}}+\f{k}{\nu}|\leq1$ implies that $|\dot{\bar{\eta}}(s,k,\xi)|=|e^{\nu s}\rr{\nu\xi-k}|\geq|k|-\nu\geq\f12$. Thus
	\[\nu\int^t_\tau|\bar{\eta}(s,k,\xi)|\vv{1}_{|\bar{\eta}(s,k,\xi)|\leq1} ds\leq4\nu.\]
	Thus
	\[\log\f{\fM(t,k,\xi)}{\fM(\tau,k,\xi)}\leq \rm(1-e^{-\nu t}+4\nu)-\rm\log S(t,\tau,k,\xi).
	\]
	Thus \eqref{MSk} holds.	For $k=0$,
	\[\ba
	\nu\int_\tau^te^{\nu s}|\xi|\vv{1}_{e^{\nu s}|\xi|\leq1}ds\leq\nu\int_0^{\f1\nu\log\f1{|\xi|} }e^{\nu s}|\xi|ds\leq1-|\xi|\leq1.
	\ea\] 
	Thus \eqref{MS0} holds. We next prove \eqref{MSkt}.
	\[\ba
	&\log\f{\fM(t,k,k\tent)}{\fM(\tau,k,k\tent)}\le\rm\nu\int_\tau^t\rr{e^{-2\nu(t-s)}+|\bar{\eta}(s,k,k\tent)|^2}^\f12ds\\
	=&\rm\nu\int_0^{t-\tau}\rr{e^{-2\nu s}+|k\f{1-e^{-\nu s}}{\nu} |^2}^\f12ds\leq \rm\min\{\nu|k|(t-\tau)^2,|k|(t-\tau)\}+1.
	\ea\]
	\eqref{M0} holds by direct computation. Finally, we prove part~(5). We only prove the lower bound. The upper bound follows from the triangle inequality.
	\[\ba&\f1\rm\log\fM(t,k,\bar{\eta}(-t,k,\xi))\geq\f{1}{\sqrt{2}}\nu\int_0^t\nr{\bar{\eta}(s,k,\bar{\eta}(-t,k,\xi))}ds+\f{1}{\sqrt{2}}\nu\int_0^te^{-\nu(t-s)}ds\\
	&=\f{1}{\sqrt{2}}\nu\int_0^t\nr{\bar{\eta}(-s,k,\xi)}ds+\f{1}{\sqrt{2}}(1-\emnt).\ea\]
	We only have to estimate the first term. We split into two cases. We denote $\Pi(k)=I-\f{k}{|k|}\otimes\f{k}{|k|}$. The inequality trivially holds for $k=0$ so we restrict ourselves to the $k\neq 0$ case.
	
    \underline{Case 1: $\nu t\leq1$.}
	\[\ba
	&\nu\int_0^t\nr{e^{-\nu s}\xi+k\f{1-e^{-\nu s}}{\nu}}ds=\nu\int_0^{\f{\ent-1}{\nu}}\f{|\xi+ks|}{(1+\nu s)^2}ds	\geq\f{\nu}{e^2}\int_0^{\f{\ent-1}{\nu}}|\xi+ks|ds
	\geq\f{\nu}{\sqrt{2}e^2}\int_0^{\f{\ent-1}{\nu}}|\Pi(k)\xi|\\
	&+|\rr{I-\Pi(k)}\xi+ks|ds\gs\nu\rr{|k|t^2+|\Pi(k)\xi|t+\nr{(I-\Pi(k))\xi}t}\gs\nu t\rr{|k|t+|\xi|}.
	\ea\]
	For the first equality we use the change of variables $s'=\f{e^{\nu s}-1}{\nu}$. The fourth inequality uses $\f{\ent-1}{\nu}\sim t$ and the fact that $\int_0^t\nr{x+ys}ds\sim|x|t+|y|t^2$, for any $x,y\in\R$, $t>0$, which follows by direct integration.
	
	\underline{Case 2: $\nu t\geq1$. }
	\[\ba&\nu\int_0^t\nr{e^{-\nu s}\rr{\xi-\f{k}{\nu}}+\f{k}{\nu}}ds=\int_{\emnt}^1\nr{\xi-\f{k}{\nu}+\f{k}{\nu s}}ds\geq\f{1}{\sqrt{2}}\int_{\emnt}^1\nr{\rr{I-\Pi(k)}\xi-\f{k}{\nu}+\f{k}{\nu s}}\\
	&+\nr{\Pi(k)\xi}ds
	\gs\nr{\Pi(k)\xi}+\nr{\rr{I-\Pi(k)}\xi-\f{k}{\nu}}+|kt|\gs|\xi|+|kt|.\ea\]
	For the first equality we use the change of variables $s'=e^{-\nu s}$. The third inequality uses Lemma~\ref{lm7}.

\end{proof}

We next record the complementary properties of the Fourier semigroup. The explicit formulas quantify collisional damping, the derivative bound permits velocity moments to be propagated, and the final estimate is the decay bound used with the enhanced-dissipation time weights in the main argument.

\begin{prop}[Properties of the Fourier semigroup]\label{prop:semigroup-symbol}
	For $S(t,\tau,k,\xi)$ and $S(t,k)$ defined in \eqref{deS} and \eqref{deSt}, the following properties hold:
	\begin{itemize}
		\item[(1)] It holds that
		\beq\label{S1}\log S(t,0,k,\xi)=-\f{e^{2\nu t}-1}{2}\nr{\xi-\f{\ent-1}{\nu(\ent+1)}k}^2+\rr{\f{2(\ent-1)}{\nu^2(\ent+1)}-\f{t}{\nu}}|k|^2,\eeq
		\beq\label{St}-\log S(t,k)=\f{|k|^2}{\nu}\rr{t-\f2\nu\rr{1-\emnt}+\f{1}{2\nu}\rr{1-e^{-2\nu t}}}\sim\min\{\nu|k|^2t^3,\f{|k|^2}{\nu}t\}.\eeq
		\item[(2)] For every $\alpha\in\N^d$, we have
		\beq\label{Sds}\nr{\pa_\xi^\alpha S(t,\tau,k,\xi)}\ls_{\alpha}e^{|\alpha|\nu t}S^{0.9}(t,\tau,k,\xi).\eeq
		\item[(3)] For $|k|\neq0$, for any $0\leq\tau\leq t$, we have
		\beq\label{Snu13}S(t,\tau,k,\xi)\le e^{-0.01\nu (t-\tau)^3}.\eeq
	\end{itemize}
\end{prop}
\begin{proof}
The explicit formulas follow by integration along the shifted characteristics. We use the Fa\`a di Bruno formula for the derivative estimate and absorb the resulting polynomial factors into a small part of the semigroup damping. The decay estimate is obtained by completing the square in the characteristic integral.

   	Part~(1) follows by direct computation. For part~(2), the Fa\`a di Bruno formula gives
   	\[\ba
   	&\nr{\pa_\xi^\alpha S(t,\tau,k,\xi)}\ls_\alpha \sum_{\substack{\sum m_\beta\beta=\alpha,\\m_\beta\in\N,\beta\in\N^d}}S(t,\tau,k,\xi)\prod_\beta\nr{\nu\int_\tau^t\pa_\xi^\beta\nr{\bar{\eta}(s,k,\xi)}^2ds }^{m_\beta}\\
   	\ls_\alpha &\sum_{\substack{\sum m_\beta\beta=\alpha,\\m_\beta\in\N,\beta\in\N^d}}S(t,\tau,k,\xi)\prod_\beta\left\{e^{|\beta|\nu t}\rr{\nu\int_\tau^t|\bar{\eta}(s,k,\xi)|ds+1}\right\}^{m_\beta}\ls_\alpha S^{0.9}(t,\tau,k,\xi)e^{|\alpha|\nu t}.
   	\ea\]
   	The second inequality uses
	\[\nu\int_\tau^t\pa_\xi^\beta\nr{\bar{\eta}(s,k,\xi)}^2ds\leq e^{|\beta|\nu t}\rr{\nu\int_\tau^t|\bar{\eta}(s,k,\xi)|ds+1},\]
	which follows from
	\[\pa_\xi^\beta\nr{\bar{\eta}(s,k,\xi)}^2\leq e^{\nu|\beta|s}\rr{\nr{\bar{\eta}(s,k,\xi)}+1}.\]
	The last inequality uses $\nu\int_\tau^t|\bar{\eta}(s,k,\xi)|ds\ls\nu\int_\tau^t|\bar{\eta}(s,k,\xi)|^2ds+10$, as shown in the proof of Proposition~\ref{mproperty} (4). For part~(3),
   	\[\ba
   	&-\log S(t,\tau,k,\bar{\eta}(-t,k,\xi))=\nu\int^t_\tau|\bar{\eta}(s-t,k,\xi)|^2ds=\nu\int^{t-\tau}_0|\bar{\eta}(-s,k,\xi)|^2ds=\f{1-e^{-2\nu (t-\tau)}}{2}\left|\xi\right.\\
   	&\left.+\f{1-e^{-\nu (t-\tau)}}{\nu(1+e^{-\nu (t-\tau)}) }k\right|^2-\rr{\f{2(1-e^{-\nu (t-\tau)})}{\nu^2(1+e^{-\nu (t-\tau)})}-\f{t-\tau}{\nu}}|k|^2\geq\f{t-\tau}{\nu}-\f{2(1-e^{-\nu (t-\tau)})}{\nu^2(1+e^{-\nu (t-\tau)})}\geq0.01\nu(t-\tau)^3.
   	\ea\]
   	Thus \eqref{Snu13} holds.
\end{proof}

The zero spatial mode requires a separate estimate in the original velocity variable. Part~(1) below uses the mean-zero condition to obtain thermal decay. Part~(2) provides the regularization needed for sources in velocity-divergence form, as in the zero-mode remainder equation. We first prove both estimates for $\tau=0$ and then transport the multiplier through the semigroup to handle arbitrary starting times.

\begin{prop}[Weighted Fokker--Planck semigroup estimates]\label{etl}
	\begin{enumerate}
	\item	Let $m>\f d2+1$ be an integer, $\sigma\geq0$. For any mean-zero $g(v)$, $0\leq\tau\leq t$, we have
	\[\nnr{v^\alpha\<\na_v\>^\sigma\fM(t,0,\emnt\na_v)e^{\nu(t-\tau)L}g}_{\ell^2(|\alpha|\leq m)L^2_v}\ls_{\sigma,m}\nnr{v^\alpha\<\na_v\>^\sigma\fM(\tau,0,\emntau\na_v)g}_{\ell^2(|\alpha|\leq m)L^2_v}e^{-\nu(t-\tau)}.\]
	\item Let $m>\f d2$ be an integer, $\sigma\geq0$, $\beta\in\Z_{\geq0}^d$. For any $g(v)$, $0\leq\tau\leq t$, we have
	\[\ba&\nnr{v^\alpha\<\na_v\>^\sigma\fM(t,0,\emnt\na_v)e^{\nu(t-\tau)L}D^\beta_vg}_{\ell^2(|\alpha|\leq m)L^2_v}\\
	\ls_{\sigma,m,\beta}&\nnr{v^\alpha\<\na_v\>^\sigma\fM(\tau,0,\emntau\na_v)g}_{\ell^2(|\alpha|\leq m)L^2_v}\f{1}{\rr{e^{2\nu(t-\tau)}-1}^{\f{|\beta|}{2}} }.\ea\]
\end{enumerate}
\end{prop}
\begin{proof}
	\underline{Proof of (1).}
	We first prove the $\tau=0$ case. By \eqref{M0}, \eqref{Sds} and \eqref{MS0} we have
	\beq\label{Poi}\ba
	&\nnr{v^\alpha\<\na_v\>^\sigma\fM(t,0,\emnt\na_v)e^{\nu t L}g}_{\ell^2(|\alpha|\leq m)L^2_v}\ls\nnr{\<\na_v\>^\sigma\fM(t,0,\emnt\na_v)v^\alpha e^{\nu tL}g}_{\ell^2(|\alpha|\leq m)L^2_v}\\
	\ls&\sum_{|\alpha_1|+|\alpha_2|\leq m}\nnr{\<\xi\>^\sigma\fM(t,0,\emnt\xi)D^{\alpha_1}_\xi \wh{g}(\emnt\xi)D^{\alpha_2}_\xi S(t,0,0,\emnt\xi)}_{L^2_\xi}e^{-(|\alpha_1|+|\alpha_2|)\nu t}\\
	\ls&\nnr{\<\xi\>^\sigma\fM(t,0,\emnt\xi)D^{\alpha}_\xi \wh{g}(\emnt\xi) S^{0.9}(t,0,0,\emnt\xi)e^{-|\alpha|\nu t}}_{\ell^2(|\alpha|\leq m)L^2_\xi}\\
	\ls&\nnr{\<\xi\>^\sigma D^{\alpha}_\xi \wh{g}(\emnt\xi) S^{0.8}(t,0,0,\emnt\xi)e^{-|\alpha|\nu t}}_{\ell^2(|\alpha|\leq m)L^2_\xi}.
	\ea\eeq
	We distinguish the following cases.
	
	\uwave{Case 1: $\nu t\leq1$.} In this case we have
	\beq\label{poi}\ba
	&\nnr{\<\xi\>^\sigma D^{\alpha}_\xi \wh{g}(\emnt\xi) S^{0.8}(t,0,0,\emnt\xi)e^{-|\alpha|\nu t}}_{\ell^2(|\alpha|\leq m)L^2_\xi}\leq\nnr{\<\xi\>^\sigma D^{\alpha}_\xi \wh{g}(\emnt\xi) }_{\ell^2(|\alpha|\leq m)L^2_\xi}\\
	=&\nnr{\<\ent\xi\>^\sigma D^{\alpha}_\xi \wh{g}(\xi) }_{\ell^2(|\alpha|\leq m)L^2_\xi}e^{\f{d\nu t}{2}}
	\ls\nnr{\<\xi\>^\sigma D^{\alpha}_\xi \wh{g}(\xi) }_{\ell^2(|\alpha|\leq m)L^2_\xi}\ls\nnr{v^\alpha\<\na_v\>^\sigma g}_{\ell^2(|\alpha|\leq m)L^2_v}\emnt.
	\ea\eeq
	
	\uwave{Case 2: $\nu t\geq1$.} We first consider $\alpha=0$.
	\[\ba
	&\nnr{\<\xi\>^\sigma \wh{g}(\emnt\xi) S^{0.8}(t,0,0,\emnt\xi)}_{L^2_\xi}\ls\nnr{\<\xi\>^\sigma \f{\wh{g}(\emnt\xi)}{|\emnt\xi|}|\emnt\xi| e^{-0.1|\xi|^2}}_{L^2_\xi}\\
	\ls&\nnr{v^\alpha g}_{\ell^2(|\alpha|\leq m)L^2_v}\nnr{\<\xi\>^{\sigma+1} e^{-0.1|\xi|^2}}_{L^2_\xi}\emnt\ls\nnr{v^\alpha g}_{\ell^2(|\alpha|\leq m)L^2_v}\emnt.
	\ea\]
	The second inequality uses $\f{|\wh{g}(\xi)|}{|\xi|}\ls\nnr{v^\alpha g}_{\ell^2(|\alpha|\leq m)L^2_v}$, which follows from $\int g(v)dv=0$ and Sobolev embedding. We next consider $|\alpha|=1$.
	\[
	\ba
	&\nnr{\<\xi\>^\sigma D \wh{g}(\emnt\xi) S^{0.8}(t,0,0,\emnt\xi)e^{-\nu t}}_{L^2_\xi}\\
	&\ls\nnr{v^\alpha g}_{\ell^2(|\alpha|\leq m)L^2_v}\nnr{\<\xi\>^{\sigma+1} e^{-0.1|\xi|^2}}_{L^2_\xi}\emnt\ls\nnr{v^\alpha g}_{\ell^2(|\alpha|\leq m)L^2_v}\emnt.
	\ea\]
	The second inequality uses $|D\wh{g}(\xi)|\ls\nnr{v^\alpha g}_{\ell^2(|\alpha|\leq m)L^2_v}$, which follows from Sobolev embedding. We next consider the case $2\leq |\alpha|\leq m$.
	\[\ba
	&\nnr{\<\xi\>^\sigma D^{\alpha}_\xi \wh{g}(\emnt\xi) S^{0.8}(t,0,0,\emnt\xi)e^{-|\alpha|\nu t}}_{L^2_\xi}\ls\nnr{\<\xi\>^\sigma D^{\alpha}_\xi \wh{g}(\emnt\xi)e^{-0.1|\xi|^2} }_{L^2_\xi}e^{-|\alpha|\nu t}\\
	\ls&\nnr{D^{\alpha}_\xi \wh{g}(\emnt\xi)}_{L^{\max\{\f{d}{|\alpha|-1},2\}}_\xi}e^{-|\alpha|\nu t}\\
	\ls&\nnr{D^{\alpha}_\xi \wh{g}(\xi)}_{L^{\max\{\f{d}{|\alpha|-1},2\}}_\xi}\emnt\ls\nnr{v^\alpha g}_{\ell^2(|\alpha|\leq m)L^2_v}\emnt.
	\ea\]
	The last inequality uses Sobolev embedding.
	
	We next consider general $0\leq\tau\leq t$. We reduce it to the $\tau=0$ case.
	\beq\label{Oi}\ba
	&\nnr{v^\alpha\<\na_v\>^\sigma\fM(t,0,\emnt\na_v)e^{\nu(t-\tau)L}g}_{\ell^2(|\alpha|\leq m)L^2_v}\ls\nnr{\<\na_v\>^\sigma\fM(t,0,\emnt\na_v)v^\alpha e^{\nu(t-\tau)L}g}_{\ell^2(|\alpha|\leq m)L^2_v}\\
	\ls&\nnr{\<\na_v\>^\sigma\fM(t-\tau,0,e^{-\nu(t-\tau)}\na_v)\fM(\tau,0,\emnt\na_v)v^\alpha e^{\nu(t-\tau)L}g}_{\ell^2(|\alpha|\leq m)L^2_v}\\
	\ls&\nnr{\<\na_v\>^\sigma\fM(t-\tau,0,e^{-\nu(t-\tau)}\na_v)v^\alpha\fM(\tau,0,\emnt\na_v) e^{\nu(t-\tau)L}g}_{\ell^2(|\alpha|\leq m)L^2_v}\\
	\ls&\nnr{\<\na_v\>^\sigma\fM(t-\tau,0,e^{-\nu(t-\tau)}\na_v)v^\alpha e^{\nu(t-\tau)L}\fM(\tau,0,\emntau\na_v)g}_{\ell^2(|\alpha|\leq m)L^2_v}\\
	\ls&\nnr{v^\alpha\<\na_v\>^\sigma\fM(\tau,0,\emntau\na_v)g}_{\ell^2(|\alpha|\leq m)L^2_v}e^{-\nu(t-\tau)}.
	\ea\eeq
	The second inequality uses the following inequality
	\[\ba
	&\f{1}{\rm}\log \fM(t,0,\emnt\xi)\leq1-\emnt+\nu\int^\tau_0|\bar{\eta}(s,0,\emnt\xi)| ds+\nu\int_\tau^t|\bar{\eta}(s,0,\emnt\xi)| ds\\
	\leq&1+\f1\rm\log \fM(\tau,0,\emnt\xi)+\rr{1-e^{-\nu(t-\tau)}}|\xi|\leq1+\f1\rm\log \fM(\tau,0,\emnt\xi)+\f1\rm\log\fM(t-\tau,0,e^{-\nu(t-\tau)}\xi),
	\ea\]
	which implies that $\fM(t,0,\emnt\xi)\ls\fM(t-\tau,0,e^{-\nu(t-\tau)}\xi)\fM(\tau,0,\emnt\xi)$. The third inequality uses \eqref{Mds}. The last inequality uses the inequality for $\tau=0$.
	
	\underline{Proof of (2).} We first prove the $\tau=0$ case. As in \eqref{Poi} we have
	\[\ba
	&\nnr{v^\alpha\<\na_v\>^\sigma\fM(t,0,\emnt\na_v)e^{\nu t L}D^\beta_vg}_{\ell^2(|\alpha|\leq m)L^2_v}\ls\nnr{\<\xi\>^\sigma\fM(t,0,\emnt\xi)D_\xi^\alpha \rr{e^{-|\beta|\nu t}\xi^\beta\wh{g}(\emnt\xi)\right.\right.\\
			&\left.\left.\times S(t,0,0,\emnt\xi)}}_{\ell^2(|\alpha|\leq m)L^2_\xi}
	\ls\nnr{\<\xi\>^{\sigma+|\beta|}e^{-|\alpha|\nu t} D_\xi^\alpha \wh{g}(\emnt\xi)S^{0.8}(t,0,0,\emnt\xi)}_{\ell^2(|\alpha|\leq m)L^2_\xi}e^{-|\beta|\nu t}.
	\ea\]
	We distinguish the following cases.
	
	\uwave{Case 1: $\nu t\leq1$.} As in \eqref{poi} we have
	\[\ba
	&\nnr{\<\xi\>^{\sigma+|\beta|}e^{-|\alpha|\nu t} D_\xi^\alpha \wh{g}(\emnt\xi)S^{0.8}(t,0,0,\emnt\xi)}_{\ell^2(|\alpha|\leq m)L^2_\xi}e^{-|\beta|\nu t}\\
	\ls&\f{1}{\rr{1-e^{-\nu t}}^{\f{|\beta|}{2}} }\nnr{\<\xi\>^{\sigma} D_\xi^\alpha \wh{g}(\emnt\xi)}_{\ell^2(|\alpha|\leq m)L^2_\xi}
	\ls\nnr{\<\na_v\>^{\sigma} v^\alpha g}_{\ell^2(|\alpha|\leq m)L^2_v}\f{1}{\rr{e^{2\nu t}-1}^{\f{|\beta|}{2}} }.
	\ea\]
	The first inequality uses $\<\xi\>^{|\beta|} e^{-0.4\rr{1-e^{-2\nu t}}|\xi|^2 }\ls_\beta\f{1}{\rr{1-e^{-2\nu t}}^{\f{|\beta|}{2}} }$.
	
	\uwave{Case 2: $\nu t\geq1$.} We first consider $\alpha=0$.
	\[\ba&
	\nnr{\<\xi\>^{\sigma+|\beta|} \wh{g}(\emnt\xi)S^{0.8}(t,0,0,\emnt\xi)}_{L^2_\xi}e^{-|\beta|\nu t}
	\ls\nnr{\wh{g}}_{L^\infty_\xi}\nnr{\<\xi\>^{\sigma+|\beta|}e^{-0.1|\xi|^2} }_{L^2_\xi}e^{-|\beta|\nu t}\\
	&\ls\nnr{v^\alpha g}_{\ell^2(|\alpha|\leq m)L^2_v}e^{-|\beta|\nu t}.
	\ea\]
	We next consider $1\leq|\alpha|\leq m$.
	\[\ba
	&\nnr{\<\xi\>^{\sigma+|\beta|}e^{-|\alpha|\nu t} D_\xi^\alpha \wh{g}(\emnt\xi)S^{0.8}(t,0,0,\emnt\xi)}_{L^2_\xi}e^{-|\beta|\nu t}\ls\nnr{D^{\alpha}_\xi \wh{g}(\emnt\xi)}_{L^{\max\{\f{d}{|\alpha|},2\}}_\xi}e^{-\rr{|\alpha|+|\beta|}\nu t}\\
	&\ls\nnr{D^{\alpha}_\xi \wh{g}(\xi)}_{L^{\max\{\f{d}{|\alpha|},2\}}_\xi}e^{-|\beta|\nu t}\ls\nnr{v^\alpha g}_{\ell^2(|\alpha|\leq m)L^2_v}\f{1}{\rr{e^{2\nu t}-1}^{\f{|\beta|}{2}} }.
	\ea\]
	
	We next consider the general $0\leq\tau\leq t$ case. As in \eqref{Oi} we have
	\[\ba
	&\nnr{v^\alpha\<\na_v\>^\sigma\fM(t,0,\emnt\na_v)e^{\nu(t-\tau)L}D^\beta_vg}_{\ell^2(|\alpha|\leq m)L^2_v}\ls\nnr{\<\na_v\>^\sigma\fM(t-\tau,0,e^{-\nu(t-\tau)}\na_v)v^\alpha e^{\nu(t-\tau)L}\right.\\
		&\left.\times D^\beta_v\fM(\tau,0,\emntau\na_v)g}_{\ell^2(|\alpha|\leq m)L^2_v}
	\ls\nnr{v^\alpha\<\na_v\>^\sigma\fM(\tau,0,\emntau\na_v)g}_{\ell^2(|\alpha|\leq m)L^2_v}\f{1}{\rr{e^{2\nu(t-\tau)}-1}^{\f{|\beta|}{2}} }.
	\ea\]
\end{proof}
	
\textbf{AI declaration} This work was completed without the assistance of artificial intelligence (AI) tools. The author subsequently used ChatGPT solely for language editing and polishing of the manuscript. The author takes full responsibility for the correctness and validity of all mathematical results presented in this paper.

\end{document}